\documentclass[reqno]{amsart}

\usepackage{geometry}
\usepackage{amsmath,amssymb,amsthm,amsfonts}
\usepackage{physics}
\allowdisplaybreaks[4]
\usepackage{hyperref}
\usepackage{mathrsfs}
\usepackage{appendix}
\usepackage{graphicx}
\usepackage{setspace}
\usepackage{mathtools}
\usepackage{latexsym}
\usepackage{mleftright}
\usepackage{cite}
\usepackage{color}

\newtheorem{lemma}{Lemma}[section]
\newtheorem{theorem}{Theorem}[section]

\newtheorem{remark}{Remark}[section]

\numberwithin{equation}{section}
\begin{document}
\title[Incompressible Limit for Diffuse Interface Models]{Incompressible Limit of Strong Solutions to the \\ Diffuse Interface Model for Two-phase Flows}
\thanks{$^*$Corresponding author}
\author[Yinghua Li]{Yinghua Li}
\address[Y. Li]{School of Mathematical Sciences, South China Normal University,
Guangzhou, 510631, China}
\email{yinghua@scnu.edu.cn}
\author[Manrou Xie]{Manrou Xie$^*$}
\address[M. Xie]{School of Mathematical Sciences, South China Normal University,
Guangzhou, 510631, China}
\email{2550743137@qq.com}

\date{\today}

\begin{abstract}
This paper is concerned with the incompressible limit problem for strong solutions of compressible two-phase flow models under periodic boundary conditions, where the Navier-Stokes equations are nonlinearly coupled with either Cahn-Hilliard equations or Allen-Cahn equations. The viscosity coefficients are allowed to depend both on the density and the phase field variable.
We establish rigorous convergence of both local and global strong solutions of compressible systems to their incompressible systems as the Mach number tends to zero.
This theoretical framework establishes an essential linkage between compressible and incompressible phase field models, demonstrating that both formulations exhibit consistent physical fidelity in capturing interfacial flow dynamics.
Furthermore, we provide some convergence rate estimates of the solutions.
\end{abstract}

\maketitle
\tableofcontents
\vspace{-5mm}

\section{Introduction}
\subsection{The problem}
The diffuse interface model for compressible viscous fluid mixtures
can be used to describe the contact interface of immiscible two-phase flows with different densities.
The model takes into account the capillarity effect and the partial mixing of the fluids in the interfacial region.
The interface is viewed as a thin diffusive layer.
This model turns out to be effective in capturing the complexities of immiscible two-phase flow interface.
The dynamic of the fluids is governed by Navier-Stokes equations, and the diffusive effect can be described by
Cahn-Hilliard equations or Allen-Cahn equations.
Then it leads to Navier-Stokes/Cahn-Hilliard (NSCH) system and Navier-Stokes/Allen-Cahn (NSAC) system.
Let $\mathbb{T}^N$ be a torus in $\mathbb{R}^N$ $(N=2, 3)$, and given a time $T>0$, for $(x, t)\in\mathbb T\times(0, T)$,
$\rho=\rho(x, t)$ is the total density, $u=u(x, t)$ represents the mean velocity of the fluid mixture,
$\phi=\phi(x, t)$ is the concentration difference of the two partly miscible compressible fluids, $\mu=\mu(x, t)$ denotes the chemical potential.
Then the dynamics of the mixture of fluids can be simulated by the following system
\begin{align}\label{P}
\begin{cases}
\partial_t\rho+{\rm div}(\rho u)=0,
\\
\partial_t(\rho u)+{\rm div}(\rho u\otimes u)-{\rm div}{\mathbb S}+\dfrac1M\nabla P=
-{\rm div}\left(\nabla\phi\otimes\nabla\phi-\frac{|\nabla\phi|^2}{2}{\mathbb I}_N\right),
\\
\partial_t(\rho\phi)+{\rm div}(\rho u\phi)=A_K\mu,
\\
\rho\mu=-\Delta\phi+\rho\dfrac1M\dfrac{\partial f}{\partial\phi},
\end{cases}
\end{align}
where $P=\rho^2\frac{\partial f}{\partial\rho}$ and
$$
{\mathbb S}=\nu {\mathbb D}u+\tilde\eta{\rm div}u{\mathbb I}, \qquad
{\mathbb D}u=\nabla u+\nabla u^T-\frac2N{\rm div}u\mathbb I.
$$
Here, $f=f(\rho, \phi)$ is a homogeneous free energy density of the mixture and $M>0$ represents an analogue of a Mach number.
Moreover, the operator $A_K$ has a twofold definition according to the case $K = CH$ (Cahn-Hilliard) or $K = AC$ (Allen-Cahn), namely,
$$
A_{CH}\mu=\Delta\mu, \qquad A_{AC}\mu=-\mu.
$$
And the chemical potential $\mu$ is related with the variational derivative of the following free energy functional
$$
E_{free}(\rho, \phi)=\int_{\mathbb T^N}\left(\frac1M\rho f(\rho, \phi)+\frac12|\nabla\phi|^2\right){\rm d}x.
$$

As in \cite{A-L-N}, we assume that $f$ is given in the form
$$
f(\rho, \phi)=f_e(\rho)+MG(\phi).
$$
This choice coincides with the assumption in \cite{A-F} with $H\equiv0$ and by adding the factor $M$ in front for $G$ therein for NSCH system.
It is also coincides with the hypotheses in \cite{F-P-R-S} with $g_i\equiv0(i=1,2)$ and inserting $M$ in front of $W$ therein for NSAC system.
Then
$$
P=\rho^2\frac{\partial f}{\partial\rho}=p_e(\rho), \qquad f_e(\rho)=\int_1^\rho\frac{p_e(z)}{z^2}{\rm d}z.
$$
We assume $p_e(\cdot)$ is smooth and satisfy $p_e^{\prime}(\rho)>0$ for $\rho>0$. Choose $G(\phi)=\frac14(\phi^2-1)^2$, set $\eta=(1-\frac2N)\nu+\tilde\eta$,
and denote Mach number by $M=\varepsilon^2$.
Then, the equations in (\ref{P}) reduce to
\begin{align}\label{q.1.1}
\begin{cases}
\partial_t\rho^\varepsilon+{\rm div}\left(\rho^\varepsilon u^\varepsilon\right)=0,
\\
\partial_t\left(\rho^\varepsilon u^\varepsilon\right)
+{\rm div}\left(\rho^\varepsilon u^\varepsilon\otimes u^\varepsilon\right)
+\dfrac1{\varepsilon^2}\nabla P\left(\rho^\varepsilon\right)
\\
\quad~\quad
=\nu(\rho^\varepsilon, \phi^\varepsilon)\Delta u^\varepsilon
+\eta(\rho^\varepsilon, \phi^\varepsilon)\nabla\left({\rm {\rm div}} u^\varepsilon\right)
-{\rm div}\left(\nabla\phi^\varepsilon\otimes\nabla\phi^\varepsilon-\dfrac{|\nabla\phi^\varepsilon|^2}{2}\mathbb I_N\right),
\\
\partial_t\left(\rho^\varepsilon\phi^\varepsilon\right)+{\rm {\rm div}}\left(\rho^\varepsilon u^\varepsilon\phi^\varepsilon\right)
={A_K}\mu^\varepsilon,
\\
\rho^\varepsilon\mu^\varepsilon
=-\Delta\phi^\varepsilon+\rho^\varepsilon\left(\left(\phi^\varepsilon\right)^3-\phi^\varepsilon\right).\\	
\end{cases}	
\end{align}

From a mathematical perspective, it is reasonable to expect that as $\rho^\varepsilon\rightarrow1$, the first equation in (\ref{q.1.1})
reduces to the incompressible condition $\nabla\cdot u=0$.
Assuming that the limits $u^\varepsilon\rightarrow u$ and $\raisebox{0.5ex}{$\phi^\varepsilon$}\rightarrow \raisebox{0.5ex}{$\phi$}$ exist
as $\varepsilon\rightarrow 0$, and that the viscosity coefficients $\nu(\cdot, \cdot)$ and $\eta(\cdot, \cdot)$ are sufficiently smooth
and satisfy the uniform bounds $\nu_*\le\nu\le \nu^*$ and $\eta_*\le\eta\le \eta^*$ for positive constants $\nu_*$, $\nu^*$, $\eta_*$ and $\eta^*$,
passing to the limit $\varepsilon\rightarrow 0$ in (\ref{q.1.1}) yield the following incompressible model:
	\begin{align}\label{q.1.2}
		\begin{cases}
			\nabla\cdot u=0,\\
			u_t+\left(u\cdot\nabla\right)u +\nabla p=\nu(\phi)\Delta u-\nabla\cdot\left(\nabla\phi\otimes\nabla\phi\right),\\
			\phi_t+\left(u\cdot\nabla\right)\phi=A_K\mu,
\\
\mu=-\Delta\phi+\left(\phi^3-\phi\right),
		\end{cases}
	\end{align}
	where $\nu(\phi)=\nu(1, \phi)$ and $\nabla p$ is the limit of $\frac{1}{\varepsilon^2\rho^\varepsilon}\nabla\left(P(\rho^\varepsilon)\right)-\frac{1}{\rho^\varepsilon}\nabla\left(\frac{|\nabla\phi^\varepsilon|^2}{2}\right)$.

\subsection{Literatures review}
The compressible models in (\ref{q.1.1}) originate from the work of Abels-Feireisl \cite{A-F} (NSCH system) and Blesgen \cite{B} (NSAC system).
There have been lots of studies on one-dimensional case, we refer to \cite{D-L, D-L-L, C-G, C-H-M-S, C-Z} and the references therein.
In contrast, for multidimensional settings, there are only a few analysis results available.
Abels-Feireisl \cite{A-F} established the existence of weak solutions for the NSCH system,
while Feireisl et al. \cite{F-P-R-S} addressed the NSAC system.
Kotschote \cite{Kotschote, Kotschote1} demonstrated the local existence of strong solutions for the Lowengrub-Truskinovsky NSCH model \cite{L-T} and the non-isentropic NSAC system, respectively.
Chen-Hong-Shi \cite{C-H-S} proved the global existence and stability of strong solutions for the NSCH system with periodic boundary conditions near the phase separation state, without the detail proof of the existence of local solutions.

Regarding equation (\ref{q.1.2}), the incompressible NSCH system is the so called ``Model H" originally introduced by Hohenberg-Halperin \cite{H-H-77}.
Initial progress on the existence of strong solutions in $\mathbb R^2$ was established by Starovo\u itov \cite{SVN} under the assumption of sufficiently smooth double-well potentials.
More comprehensive results were subsequently demonstrated by Boyer \cite{B-99} in the context of periodic channel domains with smooth double-well potentials. Specifically, the author established the existence of global weak solutions while proving their local in time strong uniqueness in both two-dimensional and three-dimensional settings.
Moreover, the case of the physical relevant logarithmic potential is also considered in connection with a degenerate mobility.
In this case the existence of weak solutions is shown.
The investigations on the asymptotic behavior of solutions have been reported in \cite{G-G-10, G-G-11}.
The studies addressing NSCH systems with singular free energy densities can be found in \cite{A-09, G-M-T}.
The incompressible NSAC system can be seen as the simplified version of liquid crystal system.
For the well-posedness results, one can see \cite{L-L-95, L-L-96} for example.
While in liquid crystal system, the viscosity generally do not depend on the crystal material.

Since the pioneering work of Klainerman-Majda \cite{K-M,K-M-2} on the incompressible limit of quasilinear hyperbolic systems and viscous equations, much progress has been achieved in studying incompressible limits for various compressible flow models. Examples include:
Magnetohydrodynamic equations \cite{H-W, J-J-L ,J-J-L2,J-J-L-X},
Oldroyd/Oldroyd-B models \cite{Lei, L-Z},
liquid-gas two-phase flow models \cite{H-L, Y-Z-Z},
Navier-Stokes-Smoluchowski system \cite{H-H-W},
liquid crystal system \cite{D-H-W-Z},
Navier-Stokes-Fourier system \cite{F-N,F-J-A,F-N-P,Feireisl}, among others.
Recently, Abels-Liu-Ne$\rm\breve{c}$asov$\rm\acute{a}$ \cite{A-L-N} established the low Mach number limit for weak solutions of the compressible NSCH system through a relative entropy method, conditional on the existence of smooth solutions to the corresponding incompressible system.
In parallel, Feireisl-Petcu-Pra\v z\'ak \cite{F-P-P} investigated a simplified compressible NSAC model, demonstrating the convergence of weak solutions under the low Mach number limit, similarly relying on the assumption of smooth incompressible solutions.

This paper rigorously justifies the incompressible limit (i.e., the low Mach number limit) for both local in time
and global strong solutions of the compressible isentropic NSCH/NSAC models,
where the viscosity coefficients depend on the density and phase field variable.
Specifically, under some smallness assumptions on initial data, we
establish the existence of global strong solutions for both compressible equations (with small Mach number) and incompressible systems,
and derive some convergence rate estimates for strong solutions in some Sobolev norms.

\subsection{Statements on the proof}

The specific form of the linearized problem (\ref{q.3.2}) is motivated by the anticipated cancellation mechanism for the singular term $\frac1{\varepsilon^2}\frac{P^{\prime}(\xi^\varepsilon)}{\xi^\varepsilon}\nabla\rho^\varepsilon$, where the chosen structure enables essential term elimination through energy estimates (\ref{odd}).
Similar cancellation mechanism can also be found in the proof of Lemma \ref{L3.3},
see the terms in (\ref{C-1}) and (\ref{C-11}), (\ref{C-2}) and (\ref{qq.3.630}).
These energy cancellation techniques align with the methods in \cite{Lei, D-H-W-Z}.

\vskip1mm
Beyond distinctions in derivative order, the fundamental difference between compressible NSCH and NSAC systems lies in their density coupling mechanisms. The Allen-Cahn system directly couples the density field itself into its evolution equation, whereas the Cahn-Hilliard system features coupling by the density itself and its second-order spatial derivative.
Hence for the estimate of $\phi$, the fourth term
$-\int_{\mathbb{T}^N}\nabla^{\alpha_1}[\xi\Delta(\frac1{\xi})\Delta\psi]\nabla^{\alpha_1}\phi_t$
on the right hand side of (\ref{q.3.5}), unlike its counterpart in (\ref{qqq.3.5}),
requires the density's derivative order to be elevated to $s+2$.
Furthermore, we observe that in the proof of local existence for NSAC equations, we do not need to use the smallness of the Mach number,
except for the restriction on the initial data, see (\ref{qqq.3.55}) and (\ref{qqq.3.57}).
In contrast, the proof of local existence for NSCH system deeply rely on the smallness of the Mach number,
see (\ref{q.3.6}), (\ref{q.3.7}), (\ref{q.3.80000}) for example.

\vskip1mm
The crucial point in the proof of Theorem \ref{th2} is rewriting the equations $(\ref{q.1.1})_{3,4}$ into the following form
$$
\phi_t+\frac1{\rho}\Delta\left(\frac{\Delta\phi}{\rho}\right)-\frac2{\rho}\Delta\phi
=\frac3{\rho}\left(\phi^2-1\right)\Delta\phi-u\cdot\nabla\phi+\frac{6\phi}{\rho}|\nabla\phi|^2.
$$
With the help of the term $-\frac2{\rho}\Delta\phi$ in above equality, we can deduce another dissipative term $\int_{\mathbb{T}^N}\frac2{\rho}|\nabla^\alpha\Delta\phi|^2$
in (\ref{qq.4.8}) which is powerful to absorb the estimates from the right hand side terms
and permits the establishment of the key inequality (\ref{qq.4.16}).

\vskip1mm
In the proof of Theorem \ref{th3}, for showing the convergent rate estimates,
we first need to treat the density and phase variable dependent viscosities carefully, such as $Y_4$ in (\ref{q.5.10}).
When handling the phase field part, we find the obstacle comes from the low-order terms
$\int_{\mathbb{T}^N}\frac{1}{4}\rho^\varepsilon|(\phi^\varepsilon)^2-1)|^2$ and $\int_{\mathbb{T}^N}\frac{1}{4}|\phi^2-1|^2$
in (\ref{q.5.10}). It's hard to represent them in the form of subtraction.
But after resetting the cross term $2\nabla\mu\nabla\mu^\varepsilon$,
we demonstrate the precise cancellation mechanism aligning with the terms in $Y_8$ (\ref{Y8s}),
thereby resolving this problem.
Moreover, the order of the decay rate $\varepsilon$ is limited by the decay rate of the term $\|\rho^\varepsilon-1\|$
in (\ref{q.5.11})--(\ref{q.5.17}).

\subsection{Preliminaries}
We first give the notations to be used throughout this paper.
Let $\int_{\mathbb{T}^N}f\left(x\right)=\int_{\mathbb{T}^N}f\left(x\right)\rm dx$ and
$\int_{0}^{t}g\left(s\right)=\int_{0}^{t}g\left(s\right)\rm ds$.
Denote $\left\|\cdot\right\|,\ \left\|\cdot\right\|_{s}$ and $\left\|\cdot\right\|_\infty$
as the norms in $L^2\left(\mathbb{T}^N\right),\ H^s\left(\mathbb{T}^N\right)$ and $L^\infty\left(\mathbb{T}^N\right)$, respectively.

The following lemma will be frequently used in our calculations.
\begin{lemma} (see \cite{H-L-95, K-M, L-Z}) \label{L3.1}
Let $s>\frac{N}{2}$. For any functions $f, g$ (possibly vector-valued in $\mathbb{R}^n$) in Sobolev space $H^s(\mathbb{T}^N, \mathbb{R}^n)$, $\Psi \in C^s(\mathbb{R}^n) $, $\|\nabla^j\Psi\|_\infty<\infty~(j=1,2,...,s)$, and multi-index $\alpha$ satisfying $|\alpha|\le s$, we have the Sobolev inequality:
\begin{equation*}
  \|f\|_\infty\le \|f\|_s,
\end{equation*}
the estimates based on the chain rule:
\begin{equation*}
  \|\nabla^r(\Psi\circ f)\|\le C(1+ \|f\|_\infty^{s-1})\|\nabla f\|_{s-1},\qquad for\ 1\le r \le s,
\end{equation*}
the estimate based on Leibniz's rule:
\begin{equation*}
  \|\nabla^\alpha( fg)\|\le C(\|f\|_\infty\|\nabla^\alpha g\|+\|g\|_\infty\|\nabla^\alpha f\|),
\end{equation*}
\begin{equation*}
  \|\nabla^\alpha( fg)-f\nabla^\alpha g\|\le C(\|\nabla f\|_\infty\|g\|_{s-1}+\|g\|_\infty\|\nabla f\|_{s-1}),
\end{equation*}
where the constant $C$ are independent of $f, g$, but may depend on $|\alpha|$ and $\Psi$.
\end{lemma}

We also need the following embedding lemma.
\begin{lemma}\label{L3.5}
 Assume that $X\subset E\subset Y$ are Banach spaces and $X\hookrightarrow\hookrightarrow E$. Then the following embedding are compact:

${\rm(i)}~\left \{ \varphi: \varphi\in L^q\left(0,T; X\right), \varphi_t\in L^1\left(0, T; Y\right)\right \}\hookrightarrow \hookrightarrow L^q\left(0, T; E\right), \quad
	\ if~ 1\leq q\leq \infty$;

${\rm(ii)}~ \left \{ \varphi: \varphi\in L^{\infty}\left(0, T; X\right), \varphi_t\in L^r\left(0, T; Y\right)\right \}\hookrightarrow \hookrightarrow C\left([0, T]; E\right), \quad
	\ if~ 1< r\leq \infty$.
\end{lemma}

In the following, we discuss the incompressible limit of Navier-Stokes/Cahn-Hilliard system and Navier-Stokes/Allen-Cahn system
in Section 2 and 3, respectively.

\section{Navier-Stokes/Cahn-Hilliard system}\label{sec:2}
This section is devoted to incompressible limit of the Navier-Stokes/Cahn-Hilliard system,
i.e. the operator $A_K=A_{CH}$ in $(\ref{q.1.1})$.

\subsection{Statements of the theorems}
When the density is away from vacuum, we rewrite the compressible Navier-Stokes/Cahn-Hilliard system $(\ref{q.1.1})$ with the parameter $\varepsilon$  as follows:
	\begin{align}\label{q.2.1}
	\begin{cases}
		\rho_t^\varepsilon+\nabla\cdot\left(\rho^\varepsilon u^\varepsilon\right)=0,\\
		u^\varepsilon_t+\left(u^\varepsilon\cdot\nabla\right) u^\varepsilon+\dfrac1{\varepsilon^2\rho^\varepsilon}\nabla\left(P\left(\rho^\varepsilon\right)\right)=\dfrac{\nu(\rho^\varepsilon,\phi^\varepsilon)}{\rho^\varepsilon}\Delta u^\varepsilon+\dfrac{\eta(\rho^\varepsilon,\phi^\varepsilon)}{\rho^\varepsilon}\nabla\left(\nabla\cdot u^\varepsilon\right)
		-\dfrac{1}{\rho^\varepsilon}\Delta\phi^\varepsilon\nabla\phi^\varepsilon,\\
		\phi^\varepsilon_t+\left(u^\varepsilon\cdot\nabla\right) \phi^\varepsilon\\
\quad=-\dfrac{1}{\left(\rho^\varepsilon\right)^2}\Delta^2\phi^\varepsilon
-\dfrac2{\rho^\varepsilon}\nabla\left(\dfrac1{\rho^\varepsilon}\right)\cdot\nabla\Delta\phi^\varepsilon
-\dfrac1{\rho^\varepsilon}\Delta\left(\dfrac1{\rho^\varepsilon}\right)\Delta\phi^\varepsilon
+\dfrac1{\rho^\varepsilon}\Delta\left(\left(\phi^\varepsilon\right)^3-\phi^\varepsilon\right).
	\end{cases}	
	\end{align}
Define
\begin{align*}
&E_{s}(\rho,u)(t)
=\displaystyle\sum_{|\alpha|\leq s}\displaystyle\int_{\mathbb{T}^N}\left(\displaystyle\frac{1}{\varepsilon^2}|\nabla^{\alpha} \left(\rho-1\right)|^2+|\nabla^{\alpha}u|^2\right),
\\
&\widetilde{E}_{s}(\rho,u)(t)
=\displaystyle\sum_{|\alpha|\leq s}\displaystyle\int_{\mathbb{T}^N}\left(\displaystyle\frac{P'\left(\rho\right)}{\varepsilon^2\rho}|\nabla^{\alpha} \left(\rho-1\right)|^2+\rho|\nabla^{\alpha}u|^2\right), \nonumber
\end{align*}
and
$$
F_{s}(\nabla\phi)(t)=\displaystyle\sum_{|\alpha|\leq s}\displaystyle\int_{\mathbb{T}^N}|\nabla\nabla^{\alpha}\phi|^2, \qquad
\widetilde{F}_{s}(\nabla\phi)(t)=\displaystyle\sum_{|\alpha|\leq s}\displaystyle\int_{\mathbb{T}^N}\rho|\nabla\nabla^{\alpha}\phi|^2.
$$
Then
\begin{align}\label{q.2.3}
	E_{s}(\rho,u)(t)\backsim \widetilde{E}_{s}(\rho,u)(t) \quad {\rm and}\quad
F_{s}(\nabla\phi)(t)\backsim \widetilde{F}_{s}(\nabla\phi)(t),
\end{align}
provided that $|\rho-1|$ is sufficiently small.

Now, we state the main results of this section.
\begin{theorem}\label{th1}
Consider the compressible Navier-Stokes/Cahn-Hilliard system $(\ref{q.1.1})$ with the following initial data
\begin{align}\label{q.2.4}
\rho^\varepsilon(x, 0)=1+\overline\rho^\varepsilon_0(x), \quad
u^\varepsilon(x, 0)=u_0(x)+\overline u^\varepsilon_0(x), \quad
\phi^\varepsilon(x, 0)=\phi_0(x)+\overline\phi^\varepsilon_0(x),
\end{align}
where $ u_0, \phi_0$ satisfy
%
\begin{align}\label{q.2.5}
u_0\in H^{s+3}\left(\mathbb{T}^N\right), \qquad~
\nabla\cdot u_0=0, \qquad~
\phi_0\in H^{s+4}\left(\mathbb{T}^N\right),
\end{align}
for any $s\ge[\frac N2]+2$, and, for small positive constant $\kappa_0$,
the functions $\overline\rho^\varepsilon_0(x)$, $\overline u^\varepsilon_0(x)$, $\overline\phi^\varepsilon_0(x)$ are assumed to satisfy
\begin{align}\label{q.2.6}
\left\|\overline\rho_0^\varepsilon\right\|_{s+2}\leq \kappa_0\varepsilon^2, \qquad~
\left\|\overline u_0^\varepsilon\right\|_{s+3}\leq \kappa_0\varepsilon, \qquad ~
\left\|\overline\phi_0^\varepsilon\right\|_{s+4}\leq \kappa_0\varepsilon.
\end{align}
Then the following statements hold.

Uniform stability: There exist fixed constants $T_0$ and $C$ independent of $\varepsilon$ such that a unique strong solution
$(\rho^\varepsilon, u^\varepsilon, \phi^\varepsilon)$ of system $(\ref{q.1.1})$ exists
for all small $\varepsilon$ on the time interval $[0, T_0]$ with properties:
\begin{align}\label{q.2.7}
\begin{cases}
E_{s+2}(\rho^\varepsilon,u^\varepsilon)(t)
+F_{s+1}(\nabla\phi^\varepsilon)(t)
+\|\phi^\varepsilon\|^2(t)
\\
\qquad\qquad +\displaystyle\int_0^t\left(\nu_*\left\|\nabla u^\varepsilon\right\|^2_{s+2}+\eta_*\left\|\nabla\cdot u^\varepsilon\right\|^2_{s+2}+\left\|\nabla\phi^\varepsilon\right\|^2_{s+3}\right)\leq C,
\\
E_{s+1}(\rho_t^\varepsilon, u_t^\varepsilon)(t)+F_{s}(\phi^\varepsilon_t)(t)
\\
\qquad\qquad
+\displaystyle\int_0^t\left(\nu_*\left\|\nabla u^\varepsilon_t\right\|^2_{s+1}+\eta_*\left\|\nabla\cdot u^\varepsilon_t\right\|^2_{s+1}+\left\|\phi^\varepsilon_t\right\|^2_{s+2}\right)\leq C,
	\end{cases}
\end{align}
where
$E_{s+1}(\rho_t,u_t)(t)
=\displaystyle\sum_{|\beta|\leq {s+1}}\displaystyle\int_{\mathbb{T}^N}\left(\displaystyle\frac{1}{\varepsilon^2}|\nabla^{\beta} \rho_t|^2+|\nabla^{\beta}u_t|^2\right)$ and
$F_{s}(\phi_t)(t)
=\displaystyle\sum_{|\beta|\leq s}\displaystyle\int_{\mathbb{T}^N}|\nabla^{\beta} \phi_t|^2.$

Local existence of solutions for incompressible Navier-Stokes/Cahn-Hilliard system: There exist functions $u$ and $\phi$ such that
\begin{align}\label{q.2.8}
	\begin{cases}
		\rho^{\varepsilon}\to1~~\text{in}~ L^{\infty}(\left[0,T_0];H^{s+2}\right) \cap  {\rm Lip}(\left[0,T_0\right];H^{s+1}),
\\
		u^{\varepsilon}\stackrel{\mathrm{w}^*}\rightharpoonup u~~\text{in}~L^{\infty}([0,T_0];H^{s+2}) \cap  {\rm Lip}([0,T_0];H^{s+1}),
\\
		u^{\varepsilon}\to u~~\text{in}~ C(\left[0,T_0\right];H^{s'+2}),
\\
		\phi^{\varepsilon}\stackrel{\mathrm{w}^*}\rightharpoonup\phi~~\text{in}~ L^{\infty}(\left[0,T_0\right];H^{s+2})\cap  {\rm Lip}(\left[0,T_0\right];H^{s}),
\\
		\phi^{\varepsilon}\to \phi~~\text{in}~ C(\left[0,T_0\right];H^{{s'}+2})
	\end{cases}
\end{align}
for any $s'\in{[0,s)}$, and the function pair $\left(u,\phi\right)$ is the unique strong solution of
the incompressible Navier-Stokes/Cahn-Hilliard system (\ref{q.1.2}) with the initial data
%
\begin{align}\label{q.2.9}
u(x, 0)=u_0(x),		\ \
\phi(x, 0)=\phi_0(x),
\end{align}
for some $p\in L^\infty\left(\left[0, T_0\right]; H^{s+1}\right)\cap L^2\left(\left[0, T_0\right];H^{s+2}\right)$.
\end{theorem}
\begin{remark}
Here, although the initial data depends on $\varepsilon$, the local existence time $T_{0}$ is independent of $\varepsilon$.
\end{remark}
\begin{theorem}\label{th2}
Consider the strong solutions $(\rho^\varepsilon, u^\varepsilon, \phi^\varepsilon)$ of the Navier-Stokes/Cahn-Hilliard system (\ref{q.1.1}) obtained in Theorem \ref{th1}. Suppose in addition that the initial data satisfy
%
	\begin{align}\label{q.2.10}
	\left\|u_0\right\|_{s+2}^2+\left\|\nabla\phi_0\right\|_{s+1}^2+\left\|\phi_0^2-1\right\|^2\leq\delta,
\end{align}
where $\delta$ is a positive constant. Let $\Theta=\left(\delta+\varepsilon^2\kappa_0^2\right)$. If $\delta$ is sufficiently small, then for any fixed $T>0$, the strong solution $\left(\rho^{\varepsilon},u^{\varepsilon},\phi^{\varepsilon}\right)$ satisfies the estimates
\begin{align}\label{q.2.11}
	&E_{s+2}\left(\rho^\varepsilon,u^\varepsilon\right)(t)
+F_{s+1}(\nabla\phi^\varepsilon)(t)
+\big\|\left(\phi^\varepsilon\right)^2-1\big\|^2(t) \nonumber
\\
&+\int_0^t\left(\nu_*\left\|\nabla u^\varepsilon\right\|^2_{s+2}
+\eta_*\left\|\nabla\cdot u^\varepsilon\right\|^2_{s+2}
+\left\|\nabla\phi^\varepsilon\right\|^2_{s+3}\right)\leq 4\Theta,
\quad t\in[0, T^\varepsilon),
\end{align}
and
\begin{align}\label{q.2.12}
&E_{s+1}(\rho_t^\varepsilon, u_t^\varepsilon)(t)+F_{s}(\phi^\varepsilon_t)(t)   \nonumber
\\
&+\int_{0}^{t}\left(\nu_*\left\|\nabla u^\varepsilon_t\right\|_{s+1}^2+\eta_*\left\|\nabla\cdot u^\varepsilon_t\right\|_{s+1}^2+\left\|\nabla\phi^\varepsilon_t\right\|_{s+1}^2\right)
		\leq C {\rm exp} Ct, \quad t\in[0,T],
\end{align}
where $T^\varepsilon>T$ and $T^\varepsilon\to\infty$ as $\varepsilon\to 0$.

Furthermore, as $\varepsilon\to 0$, $\left(\rho^{\varepsilon}, u^{\varepsilon}, \phi^{\varepsilon}\right)$ converges to the unique global strong solution $\left(1,u,\phi\right)$ of the incompressible Navier-Stokes/Cahn-Hilliard equations $(\ref{q.1.2})$, and
%
\begin{align}\label{q.2.13}
\left\|u\right\|_{s+2}^2(t)
+\left\|\nabla\phi\right\|_{s+1}^2(t)
+\left\|\phi^2-1\right\|^2(t)
+\int_0^t\left(\nu_*\left\|\nabla u\right\|^2_{s+2}+\left\|\nabla\phi\right\|^2_{s+3}\right)\leq C_1\delta,
\end{align}
for any $t>0$, where $C_1$ is a uniform constant independent of $\delta$ and $t$.
\end{theorem}
\begin{theorem}\label{th3} Under the assumptions of Theorem \ref{th1},
we have the following convergence rate estimates as $\varepsilon\to 0$:
\begin{align} \label{q.2.14}
&\left\|u^\varepsilon-u\right\|^2(t)
+\left\|\phi^\varepsilon-\phi\right\|_1^2(t)
+\int_{0}^{t}\left(\left\|u^\varepsilon-u\right\|_1^2+\left\|\phi^\varepsilon-\phi\right\|_3^2\right)\leq C\varepsilon
\end{align} for $t\in[0,T_0]$.
Furthermore, we have
\begin{align}\label{q.2.15}
&\qquad \quad \left\|\rho^\varepsilon-1\right\|_s^2(t)\leq C\varepsilon^2, \quad
\left\|\nabla\left(\rho^\varepsilon-1\right)\right\|_{s}^2(t)\leq C\varepsilon^4,\qquad \forall t\in[0,T_0].
\end{align}
Moreover, the statement also holds for the strong solution given in Theorem \ref{th2} for $t\in \left[0,T^\varepsilon\right)$.
\end{theorem}

In the following subsections, we prove the above theorems one by one.

\subsection{Local existence and uniform stability}\label{subsec:2.2}
In this subsection, we will give the uniform estimates for our results and then prove Theorem \ref{th1}.
Let $U_0=\left(1+\overline\rho_0^\varepsilon, u_0+\overline u_0^\varepsilon, \phi_0+\overline\phi_0^\varepsilon\right)$.
We consider a set of functions
	$B_{T_0}^\varepsilon\left(U_0\right)$ contained in
	$\{(\rho, u, \phi)| (\rho, u)\in L^\infty([0, T_0]; H^{s+2})\cap {\rm Lip}([0, T_0]; H^{s+1}),\ \phi\in L^\infty([0, T_0]; H^{s+2})\cap {\rm Lip}([0, T_0]; H^{s})\}$ with $s\ge 3$
	and defined by
	\begin{align}\label{q.3.1}
		\begin{cases}
			\left|\dfrac{\rho-1}{\varepsilon}\right|+|u-u_0|+\left|\phi-\phi_0\right|<\kappa,
			\\
			E_{s+2}(\rho,u)+F_{s+1}(\nabla\phi)+\left\|\phi\right\|^2
\\
\qquad\qquad+\displaystyle\int_0^t\left(\nu_*\left\|\nabla u\right\|_{s+2}^2+\eta_*\left\|\nabla\cdot u\right\|_{s+2}^2+\left\|\nabla\phi\right\|_{s+3}^2\right)\le K_1,
			\\
			E_{s+1}(\rho_t,u_t)+F_{s}(\phi_t)\\
\qquad\qquad
			+\displaystyle\int_0^t\left(\nu_*\left\|\nabla u_t\right\|_{s+1}^2+\eta_*\left\|\nabla\cdot u_t\right\|_{s+1}^2+\left\|\phi_t\right\|_{s+2}^2\right)\le K_2.
		\end{cases}
	\end{align}

For any $V=\left(\xi^\varepsilon,\ v^\varepsilon,\ \psi^\varepsilon\right)\in B_{T_0}^\varepsilon(U_0)$, define $U=(\rho^\varepsilon,\ u^\varepsilon,\ \phi^\varepsilon)=\Lambda(V)$ as the unique solution of the following linearized problem
	\begin{equation}\label{q.3.2}
	\begin{aligned}
		\begin{cases}
			\rho^\varepsilon_t+(v^\varepsilon\cdot\nabla)\rho^\varepsilon+\xi^\varepsilon\nabla\cdot u^\varepsilon=0,
			\\
			u^\varepsilon_t+(v^\varepsilon\cdot\nabla) u^\varepsilon+\dfrac1{\varepsilon^2}\dfrac{P^{'}(\xi^\varepsilon)}{\xi^\varepsilon}\nabla\rho^\varepsilon
			=\dfrac{\nu(\xi^\varepsilon,\psi^\varepsilon)}{\xi^\varepsilon}\Delta u^\varepsilon+\dfrac{\eta(\xi^\varepsilon,\psi^\varepsilon)}{\xi^\varepsilon}\nabla\left(\nabla\cdot u^\varepsilon\right)
			-\dfrac{1}{\xi^\varepsilon}\Delta\phi^\varepsilon\nabla\phi^\varepsilon,
			\\
			\phi^\varepsilon_t+\dfrac{1}{\left(\xi^\varepsilon\right)^2}\Delta^2\phi^\varepsilon\\
\quad= -\left(v^\varepsilon\cdot\nabla\right) \psi^\varepsilon
-\dfrac2{\xi^\varepsilon}\nabla\left(\dfrac1{\xi^\varepsilon}\right)\cdot\nabla\Delta\psi^\varepsilon
-\dfrac1{\xi^\varepsilon}\Delta\left(\dfrac1{\xi^\varepsilon}\right)\Delta\psi^\varepsilon
+\dfrac1{\xi^\varepsilon}\Delta\left(\left(\psi^\varepsilon\right)^3-\psi^\varepsilon\right).
		\end{cases}
	\end{aligned}
	\end{equation}
for which the existence and uniqueness of the solutions is guaranteed by the standard theory of parabolic equations and Navier-Stokes equations. Now we are to show that for appropriate choices of $T_{0}$, $\kappa$, $K_{1}$, $K_{2}$ independent of $\varepsilon$, $ \Lambda$ maps $B_{T_0}^\varepsilon(U_0)$ into itself and it is a contraction in certain function spaces. We emphasize that the solutions will depend on the value of the parameter $\varepsilon$, but for convenience, the dependence will not always be displayed in this subsection.

\begin{lemma}\label{L3.2}
	Suppose that $B_{T_0}^{\varepsilon}\left(U_0\right)$ is defined by (\ref{q.3.1}) and $\Lambda:V\longrightarrow U$ is defined by the system (\ref{q.3.2}). Then, under the assumptions in Theorem \ref{th1}, there exist constants $T_0,\ K_1$ and $K_2$ independent of $\varepsilon$ such that $\Lambda$ maps $B_{T_0}^{\varepsilon}\left(U_0\right)$ into itself.
	\end{lemma}
\noindent{\it\bfseries Proof.}\quad
Firstly, applying $D^{\alpha_1}$ to the first and second equations of (\ref{q.3.2}) and $D^{\alpha_2}$ to the third one respectively, and then we get
	\begin{equation}\label{q.3.3}
		\begin{aligned}
			\begin{cases}
				\partial_tD^{\alpha_1}\rho+\left(v\cdot\nabla\right)D^{\alpha_1}\rho+\xi\nabla\cdot D^{\alpha_1}u=\Pi_1,\\
				\partial_tD^{\alpha_1}u+\left(v\cdot\nabla\right)D^{\alpha_1}u+\dfrac{1}{\varepsilon^2}\dfrac{P'\left(\xi\right)}{\xi}\nabla D^{\alpha_1}\rho=\Pi_2,\\
				\partial_tD^{\alpha_2}\phi+\dfrac{1}{\xi^2}D^{\alpha_2}\Delta^2\phi=\Pi_3,
			\end{cases}
		\end{aligned}
	\end{equation}
	where
\begin{align*}
 \Pi_1&=-\left[D^{\alpha_1}\left(v\cdot\nabla\rho\right)-\left(v\cdot\nabla\right)D^{\alpha_1}\rho\right]-\left[D^{\alpha_1}\left(\xi\nabla\cdot u\right)-\xi\nabla\cdot D^{\alpha_1}u\right],\\
\Pi_2&=\dfrac{\nu(\xi,\psi)}{\xi}\Delta D^{\alpha_1} u+\dfrac{\eta(\xi,\psi)}{\xi}\nabla D^{\alpha_1}(\nabla\cdot u)-\dfrac{1}{\xi}D^{\alpha_1}\left(\Delta\phi\nabla\phi\right)\\
&\quad-\left[D^{\alpha_1}\left(v\cdot\nabla u\right)-\left(v\cdot\nabla\right)D^{\alpha_1}u\right]-
	\dfrac{1}{\varepsilon^2}\left[D^{\alpha_1}\left(\dfrac{P'(\xi)}{\xi}\nabla \rho\right)-
	\dfrac{P'(\xi)}{\xi}\nabla D^{\alpha_1}\rho\right]\\
&\quad+\left[D^{\alpha_1}\left(\dfrac{\nu(\xi,\psi)}{\xi}\Delta u\right)-
	\dfrac{\nu(\xi,\psi)}{\xi}\Delta D^{\alpha_1} u\right]\\
&\quad
+\left[D^{\alpha_1}\left(\dfrac{\eta(\xi,\psi)}{\xi}\nabla(\nabla\cdot u)\right)-\dfrac{\eta(\xi,\psi)}{\xi}\nabla D^{\alpha_1}(\nabla\cdot u)\right]\\
	&\quad-\left[D^{\alpha_1}\left(\dfrac{1}{\xi}\left(\Delta\phi\nabla\phi\right)\right)
-\dfrac{1}{\xi}D^{\alpha_1}\left(\Delta\phi\nabla\phi\right)\right],\\
	\Pi_3&=-D^{\alpha_2}(v\cdot\nabla \psi)-D^{\alpha_2}\left[\dfrac2{\xi}\nabla\left(\dfrac1{\xi}\right)\cdot\nabla\Delta\psi\right]
-D^{\alpha_2}\left[\dfrac1{\xi}\Delta\left(\dfrac1{\xi}\right)\Delta\psi\right]\\
&\quad-D^{\alpha_2}\left[\dfrac{1}{\xi}\Delta\left(\psi^3-\psi\right)\right]-\left[D^{\alpha_2}\left(\dfrac{1}{\xi^2}\Delta^2\phi\right)-\dfrac{1}{\xi^2}
D^{\alpha_2}\Delta^2\phi\right].
\end{align*}

	Next, we prove that $\Lambda$ maps $B_{T_0}^\varepsilon(U_0)$ into itself by two steps and denote by $C$ the constants independent of $\varepsilon$, $K_{1}$, $K_{2}$ in these two steps. Without loss of generality, we assume that
the positive constants $T_{0}^{-1}$, $\varepsilon^{-1}$, $K_{1}$ and $K_{2}$ are all bigger than $1$.

$\mathbf{Step\ one}$ : Estimates of $\phi$.

On one hand, we rewrite the equation $(\ref{q.3.2})_3$ as
\begin{align}\label{q.3.4}
\xi^2\phi_t+\Delta^2\phi=-\xi^2 v\cdot\nabla\psi-2\xi\nabla\left(\dfrac1{\xi}\right)\cdot\nabla\Delta\psi-\xi\Delta\left(\dfrac1{\xi}\right)\Delta\psi+\xi\Delta(\psi^3-\psi).
\end{align}
For $|\alpha_1|\le s$, applying $\nabla^{\alpha_1}$ to the above equation and multiplying by $\nabla^{\alpha_1}\phi_t$ subsequently,
then integrating the result on $\mathbb{T}^N$ yield
		\begin{align}\label{q.3.5}
			&\dfrac12\dfrac{\rm d}{{\rm d}t}\int_{\mathbb{T}^N}|\nabla^{\alpha_1}\Delta\phi|^2+\int_{\mathbb{T}^N}\xi^2|\nabla^{\alpha_1}\phi_t|^2
\nonumber\\
&=-\int_{\mathbb{T}^N}\left[\nabla^{\alpha_1}(\xi^2\phi_t)-\xi^2\nabla^{\alpha_1}\phi_t\right]
\nabla^{\alpha_1}\phi_t
-\int_{\mathbb{T}^N}\nabla^{\alpha_1}(\xi^2v\cdot\nabla\psi)\nabla^{\alpha_1}\phi_t
\nonumber\\
&\quad-\int_{\mathbb{T}^N}\nabla^{\alpha_1}\left[2\xi\nabla\left(\dfrac1{\xi}\right)\cdot\nabla\Delta\psi\right]\nabla^{\alpha_1}\phi_t
-\int_{\mathbb{T}^N}\nabla^{\alpha_1}\left[\xi\Delta\left(\dfrac1{\xi}\right)\Delta\psi\right]\nabla^{\alpha_1}\phi_t
\nonumber\\
&\quad+\int_{\mathbb{T}^N}\nabla^{\alpha_1}(\xi\Delta(\psi^3-\psi))\nabla^{\alpha_1}\phi_t
=\sum_{i=1}^{5}M_i.
		\end{align}

We will give the estimates of $M_i(i=1,...,5)$. Choose $\kappa$ small enough such that $|\xi-1|\le \frac{1}{2}$.
It follows from Lemma \ref{L3.1}, the Sobolev embedding $H^2(\mathbb{T}^N)$ $\hookrightarrow$ $L^\infty(\mathbb{T}^N)$ for $N =2, 3$, and the Cauchy inequality that
		\begin{align}
			|M_1|&\le C(\|\nabla(\xi^2)\|_{\infty}\|\phi_t\|_{s-1}+\|\nabla(\xi^2)\|_{s-1}\|\phi_t\|_{\infty})\|\nabla^{\alpha_1}\phi_t\|
	\nonumber\\
&\le C\|\nabla(\xi^2)\|_{s-1}\|\phi_t\|_{s-1}\|\nabla^{\alpha_1}\phi_t\|
	\nonumber\\
&\le C\varepsilon(1+\|\xi\|_\infty^{s-1})\left\|\dfrac1\varepsilon\nabla\xi\right\|_{s-1}\|\phi_t\|_{s-1}\|\nabla^{\alpha_1}\phi_t\|\nonumber\\
&	
\leq C(\tau) \varepsilon^2 K_1	\|\phi_t\|_{s-1}^2+\tau\|\nabla^{\alpha_1}\phi_t\|^2
,\label{q.3.6}\\		
	|M_2|&\le C(\|\xi^2v\|_{\infty}\|\nabla\psi\|_s+\|\xi^2v\|_s\|\nabla\psi\|_{\infty})\|\nabla^{\alpha_1}\phi_t\|\nonumber\\
&\le C[\|\xi\|_{\infty}^2\|v\|_{\infty}
+(\|\xi\|_\infty^2\|v\|_s+\|v\|_\infty\|\xi^2\|_s)]\|\nabla\psi\|_{s}\|\nabla^{\alpha_1}\phi_t\|\nonumber\\
&\le C[\|\xi\|_{s}^2\|v\|_{s}
+(\|\xi\|_s^2\|v\|_s+2\|v\|_s\|\xi\|_\infty\|\xi\|_s)]\|\nabla\psi\|_{s}\|\nabla^{\alpha_1}\phi_t\|\nonumber\\
			&  \leq C\|\xi\|_s^2\|v\|_s\|\nabla\psi\|_s\left\|\nabla^{\alpha_1}\phi_t\right\|
\leq CK_1^2\left\|\nabla^{\alpha_1}\phi_t\right\|
\leq C(\tau) K_1^4+\tau\left\|\nabla^{\alpha_1}\phi_t\right\|^2,\label{q.3.7}\\
|M_3|&\le C\left(\left\|\xi\nabla\left(\dfrac1{\xi}\right)\right\|_\infty\|\nabla\Delta\psi\|_s+\|\nabla\Delta\psi\|_\infty\left\|\xi\nabla\left(\dfrac1{\xi}\right)\right\|_s
\right)\|\nabla^{\alpha_1}\phi_t\|\nonumber\\
&\le C\left[\left\|\xi\right\|_s\left\|\nabla\left(\dfrac1{\xi}\right)\right\|_s
+
\left(\|\xi\|_\infty\left\|\nabla\left(\dfrac1{\xi}\right)\right\|_s+\|\xi\|_s\left\|\nabla\left(\dfrac1{\xi}\right)\right\|_\infty\right)
\right]\|\nabla\Delta\psi\|_s\|\nabla^{\alpha_1}\phi_t\|\nonumber\\
&\le C\left\|\xi\right\|_s\left\|\nabla\left(\dfrac1{\xi}\right)\right\|_s\|\nabla\Delta\psi\|_s\|\nabla^{\alpha_1}\phi_t\|
\le C K_1^{\frac12}\left\|\nabla\xi\right\|_{s}\|\nabla\Delta\psi\|_s\|\nabla^{\alpha_1}\phi_t\|
\nonumber\\
&\le C(\tau)\varepsilon^2 K_1^2\|\nabla\Delta\psi\|_s^2+\tau\|\nabla^{\alpha_1}\phi_t\|^2
\le C(\tau)\varepsilon\|\nabla\Delta\psi\|_s^2+\tau\|\nabla^{\alpha_1}\phi_t\|^2,   \label{q.3.6-1}
\\
|M_4|&\le C\left(\left\|\xi\Delta\left(\dfrac1{\xi}\right)\right\|_\infty\|\Delta\psi\|_s+\|\Delta\psi\|_\infty\left\|\xi\Delta\left(\dfrac1{\xi}\right)\right\|_s
\right)\|\nabla^{\alpha_1}\phi_t\|\nonumber\\
&\le C\left[\left\|\xi\right\|_s\left\|\Delta\left(\dfrac1{\xi}\right)\right\|_s+
\left(\left\|\xi\right\|_s\left\|\Delta\left(\dfrac1{\xi}\right)\right\|_\infty+\left\|\xi\right\|_\infty\left\|\Delta\left(\dfrac1{\xi}\right)\right\|_s\right)
\right]\|\Delta\psi\|_s\|\nabla^{\alpha_1}\phi_t\|\nonumber\\
&\le C\left\|\xi\right\|_s\left\|\Delta\left(\dfrac1{\xi}\right)\right\|_s\|\Delta\psi\|_s\|\nabla^{\alpha_1}\phi_t\|
\le CK_1^{\frac12}\left\|\nabla\xi\right\|_{s+1}\|\Delta\psi\|_s\|\nabla^{\alpha_1}\phi_t\|
\nonumber\\
&\le C(\tau)\varepsilon^2  K_1^3+\tau\|\nabla^{\alpha_1}\phi_t\|^2,\label{q.3.7}
\\
			|M_5|&\le C(\|\xi\|_{\infty}\|\Delta(\psi^3-\psi)\|_s+\|\xi\|_s\|\Delta(\psi^3-\psi)\|_{\infty})\|\nabla^{\alpha_1}\phi_t\|\nonumber\\
&\le C\|\xi\|_s\left(\|\Delta\left(\psi^3\right)\|_{s}+\|\Delta\psi\|_s\right)\|\nabla^{\alpha_1}\phi_t\|\nonumber\\
&
\le CK_1^{\frac12}(1+\|\psi\|_\infty^{s+1})\|\nabla\psi\|_{s+1}\|\nabla^{\alpha_1}\phi_t\|\nonumber\\
&\le CK_1\|\nabla^{\alpha_1}\phi_t\|
\leq C(\tau) K_1^2+\tau\left\|\nabla^{\alpha_1}\phi_t\right\|^2,\label{q.3.8}
		\end{align}
where we have used the fact that
\begin{align}\label{q.3.80000}
  \|\xi\|_s\le C\left(1+\varepsilon\left\|\dfrac1{\varepsilon}\nabla\xi\right\|_{s-1}\right)\le CK_1^{\frac12}+C\varepsilon K_1^{\frac12}\le CK_1^{\frac12}.
\end{align}
Substituting the estimates (\ref{q.3.6})-(\ref{q.3.8}) into (\ref{q.3.5}), summing over $\alpha_1$ and then taking $\tau$ and $\varepsilon$ small enough, we obtain
		\begin{align}\label{q.3.9}
			&\dfrac{\rm d}{{\rm d}t}\displaystyle\sum_{|\alpha_1|\le s}\|\nabla^{\alpha_1}\Delta\phi\|^2
+\displaystyle\sum_{|\alpha_1|\le s}\|\nabla^{\alpha_1}\phi_t\|^2
\le C(\tau)\varepsilon\|\nabla\Delta\psi\|_s^2+  C K_1^4.
		\end{align}
Recalling that $\int_{\mathbb{T}^N}\nabla^k\phi=0,(k\geq1)$, the $\rm Poincar\acute{e}$ inequality give $\|\nabla\phi\|_{s+1}(t)\le C\|\Delta\phi\|_{s}(t)$.
Integrating (\ref{q.3.9}) over $[0, t]\subseteq[0, T_0]$ gives
\begin{align*}
\|\nabla\phi\|_{s+1}^2(t)\leq \|\nabla\phi(x,0)\|_{s+1}^2+C(\tau)\varepsilon\int_0^t\|\nabla\Delta\psi\|_s^2+  C T_0K_1^4.
\end{align*}
 Applying (\ref{q.2.4})-(\ref{q.2.6}) and taking $\varepsilon$ small enough, one obtains
	\begin{align}\label{q.3.10}
		\|\nabla\phi\|_{s+1}^2(t)+\int_0^t\|\phi_t\|_{s}^2{\rm d}s\le C,
	\end{align}
	for $t\in[0, T_0]$, provided that $T_0<T_1= K_1^{-4}$.
Then using (\ref{q.3.4}), (\ref{q.3.80000}) and Lemma \ref{L3.1}, taking $\varepsilon$ small enough, one has
\begin{align}\label{qa.3.12}
\int_0^t\|\Delta^2\phi\|_{s}^2
&\le\int_0^t\left\|\xi^2\right\|_{s}^2\left\|\phi_t\right\|_{s}^2
+\int_0^t\left\|\xi^2v\right\|_{s}^2\left\|\nabla\psi\right\|_{s}^2
+\int_0^t\left\|\xi\right\|_s^2\left\|\nabla\left(\dfrac1{\xi}\right)\right\|_s^2\left\|\nabla\Delta\psi\right\|_s^2
\nonumber\\
&\quad+\int_0^t\left\|\xi\right\|_s^2\left\|\Delta\left(\dfrac1{\xi}\right)\right\|_s^2\left\|\Delta\psi\right\|_s^2
+\int_0^t\left\|\xi\right\|_{s}^2\left\|\Delta(\psi^3-\psi)\right\|_{s}^2\nonumber\\
&\le C(1+\varepsilon^2 K_1)\int_0^t\left\|\phi_t\right\|_{s}^2
+CK_1^3 T_0+C\varepsilon^2K_1^2
\int_0^t\left\|\nabla\Delta\psi\right\|_s^2+C\varepsilon K_1^3 T_0
\nonumber\\
&\quad+CT_0(1+\varepsilon^2 K_1)(1+\|\psi\|_\infty^{s+1})\left\|\nabla\psi\right\|_{s+1}^2
\nonumber\\
&\le C(1+\varepsilon^2 K_1)
+C\left(1+\varepsilon+\varepsilon^2\right) K_1^3 T_0+C(1+\varepsilon^2 K_1)K_1T_0\nonumber\\
&\le C(1+\varepsilon^2 K_1)
+C\left(1+\varepsilon+\varepsilon^2\right) K_1^3 K_1^{-4}+C(1+\varepsilon^2 K_1)K_1K_1^{-4}
\leq C,
\end{align}
for $t\in[0, T_0]$.
Hence, it follows from the $\rm Poincar\acute{e}$ inequality, (\ref{q.3.10}) and (\ref{qa.3.12}) that
	\begin{align}\label{q.3.11}
		\|\nabla\phi\|_{s+1}^2(t)+\int_0^t\left(\|\phi_t\|_{s}^2+\|\nabla\phi\|_{s+3}^2\right){\rm d}s\le C.
	\end{align}

On the other hand, applying $\nabla^{\beta}\nabla\partial_t$ to the equation (\ref{q.3.2})$_3$ and multiplying $\nabla^{\beta}\nabla\phi_t$ for $0\le|\beta|\le s-1$, we obtain after integrating the result over $\mathbb{T}^N$ that
	\begin{align}\label{q.3.13}
		&\dfrac12\dfrac{\rm d}{{\rm d}t}\int_{\mathbb{T}^N}|\nabla^{\beta}\nabla\phi_t|^2+\int_{\mathbb{T}^N}\dfrac{1}{\xi^2}|\nabla^{\beta}\nabla\Delta\phi_t|^2
		\nonumber\\ &=-\int_{\mathbb{T}^N}\left[\nabla^{\beta}\nabla\left(\dfrac1{\xi^2}\Delta^2\phi\right)_t-\dfrac1{\xi^2}\nabla^{\beta}\nabla\Delta^2\phi_t\right]\cdot\nabla^{\beta}\nabla\phi_t
		-\int_{\mathbb{T}^N}\nabla^{\beta}\nabla(v\cdot\nabla\psi)_t\cdot\nabla^{\beta}\nabla\phi_t
		\nonumber\\
  &\quad-2\int_{\mathbb{T}^N}\nabla^{\beta}\nabla\left[\dfrac1{\xi}\nabla\left(\dfrac1{\xi}\right)\cdot\nabla\Delta\psi\right]_t\cdot\nabla^{\beta}\nabla\phi_t
  -\int_{\mathbb{T}^N}\nabla^{\beta}\nabla\left[\dfrac1{\xi}\Delta\left(\dfrac1{\xi}\right)\Delta\psi\right]_t\cdot\nabla^{\beta}\nabla\phi_t\nonumber\\
&\quad+\int_{\mathbb{T}^N}\nabla^{\beta}\nabla\left(\dfrac{1}{\xi}\Delta(\psi^3-\psi)\right)_t\cdot\nabla^{\beta}\nabla\phi_t
     -2\int_{\mathbb{T}^N}\nabla\left(\dfrac{1}{\xi^2}\right)\cdot\nabla^{\beta}\nabla\nabla\phi_t\cdot\nabla^{\beta}\nabla\Delta\phi_t
     \nonumber\\
&\quad-\int_{\mathbb{T}^N}\Delta\left(\dfrac{1}{\xi^2}\right)\nabla^{\beta}\nabla\Delta\phi_t\cdot\nabla^{\beta}\nabla\phi_t
\triangleq\sum_{i=1}^7N_i,
	\end{align}
where we have used the following calculation by integration by parts that
\begin{align*}
&\int_{\mathbb{T}^N}\dfrac{1}{\xi^2}\nabla^{\beta}\nabla\Delta^2\phi_t\cdot\nabla^{\beta}\nabla\phi_t
=\int_{\mathbb{T}^N}\dfrac{1}{\xi^2}\nabla^{\beta}\nabla\Delta\nabla_i^2\phi_t\cdot\nabla^{\beta}\nabla\phi_t\nonumber\\
&=-\int_{\mathbb{T}^N}\nabla_i\left(\dfrac{1}{\xi^2}\right)\nabla^{\beta}\nabla\Delta\nabla_i\phi_t\cdot\nabla^{\beta}\nabla\phi_t
-\int_{\mathbb{T}^N}\dfrac{1}{\xi^2}\nabla^{\beta}\nabla\Delta\nabla_i\phi_t\cdot\nabla^{\beta}\nabla\nabla_i\phi_t\nonumber\\
&=\int_{\mathbb{T}^N}\nabla_i^2\left(\dfrac{1}{\xi^2}\right)\nabla^{\beta}\nabla\Delta\phi_t\cdot\nabla^{\beta}\nabla\phi_t
+\int_{\mathbb{T}^N}\nabla_i\left(\dfrac{1}{\xi^2}\right)\nabla^{\beta}\nabla\Delta\phi_t\cdot\nabla^{\beta}\nabla\nabla_i\phi_t\nonumber\\
&\quad+\int_{\mathbb{T}^N}\nabla_i\left(\dfrac{1}{\xi^2}\right)\nabla^{\beta}\nabla\Delta\phi_t\cdot\nabla^{\beta}\nabla\nabla_i\phi_t
+\int_{\mathbb{T}^N}\dfrac{1}{\xi^2}\nabla^{\beta}\nabla\Delta\phi_t\cdot\nabla^{\beta}\nabla\nabla_i^2\phi_t\nonumber\\
&=\int_{\mathbb{T}^N}\Delta\left(\dfrac{1}{\xi^2}\right)\nabla^{\beta}\nabla\Delta\phi_t\cdot\nabla^{\beta}\nabla\phi_t
+2\int_{\mathbb{T}^N}\nabla\left(\dfrac{1}{\xi^2}\right)\cdot\nabla^{\beta}\nabla\nabla\phi_t\cdot\nabla^{\beta}\nabla\Delta\phi_t
+\int_{\mathbb{T}^N}\dfrac{1}{\xi^2}|\nabla^{\beta}\nabla\Delta\phi_t|^2.
\end{align*}

Next, choose $\tau$ small enough, and we calculate $N_i (i=1,...,7)$ as follows:

Firstly, $N_1$ should provide different discussions divided into $s=3$ and $s\geq4$. When $s\geq4$, we obtain after using integration by parts that
		\begin{align}
N_1&=\int_{\mathbb{T}^N}\nabla^{\beta}\left(\left(\dfrac{1}{\xi^2}\right)_t\Delta^2\phi\right)\nabla^{\beta}\Delta\phi_t \underline{+\int_{\mathbb{T}^N}\left[\nabla^{\beta}\left(\dfrac{1}{\xi^2}\Delta^2\phi_t\right)
-\dfrac{1}{\xi^2}\nabla^{\beta}\Delta^2\phi_t\right]\nabla^{\beta}\Delta\phi_t}_{N_{1}^{1}}\nonumber\\
&\quad-\int_{\mathbb{T}^N}\nabla\left(\dfrac{1}{\xi^2}\right)\cdot\nabla\nabla^{\beta}\phi_t\nabla^{\beta}\Delta^2\phi_t\nonumber\\
&=\int_{\mathbb{T}^N}\nabla^{\beta}\left(\left(\dfrac{1}{\xi^2}\right)_t\Delta^2\phi\right)\nabla^{\beta}\Delta\phi_t +\int_{\mathbb{T}^N}\left[\nabla^{\beta}\left(\dfrac{1}{\xi^2}\Delta^2\phi_t\right)
-\dfrac{1}{\xi^2}\nabla^{\beta}\Delta^2\phi_t\right]\nabla^{\beta}\Delta\phi_t\nonumber\\
&\quad+\int_{\mathbb{T}^N}\nabla\nabla^{\beta}\phi_t\cdot\nabla^2\left(\dfrac{1}{\xi^2}\right)\cdot\nabla^{\beta}\nabla\Delta\phi_t
+\int_{\mathbb{T}^N}\nabla\left(\dfrac{1}{\xi^2}\right)\cdot\nabla^{\beta}\nabla^2\phi_t\cdot\nabla^{\beta}\nabla\Delta\phi_t\nonumber\\
		&\le C\left(\left\|\left(\dfrac1{\xi^2}\right)_t\right\|_{\infty}\left\|\Delta^2\phi\right\|_{s-1}
		+\left\|\left(\dfrac1{\xi^2}\right)_t\right\|_{s-1}\left\|\Delta^2\phi\right\|_{\infty}\right)\left\|\nabla^{\beta}\Delta\phi_t\right\|
		\nonumber \\
		&\quad+C\left(\left\|\nabla\left(\dfrac1{\xi^2}\right)\right\|_{\infty}\left\|\Delta^2\phi_t\right\|_{s-2}
		+\left\|\nabla\left(\dfrac1{\xi^2}\right)\right\|_{s-2}\left\|\Delta^2\phi_t\right\|_{\infty}\right)\left\|\nabla^{\beta}\Delta\phi_t\right\|
		\nonumber \\
        &\quad+C\left\|\nabla^2\left(\dfrac{1}{\xi^2}\right)\right\|_2\left\|\nabla\phi_t\right\|_{s-1}\left\|\nabla\Delta\phi_t\right\|_{s-1}
+C\left\|\nabla\left(\dfrac{1}{\xi^2}\right)\right\|_2\left\|\nabla^2\phi_t\right\|_{s-1}\left\|\nabla\Delta\phi_t\right\|_{s-1}
\nonumber\\
&\le C\left\|\nabla^{\beta}\Delta\phi_t\right\|\left(\left\|\left(\dfrac1{\xi^2}\right)_t\right\|_{s-1}\left\|\Delta^2\phi\right\|_{s-1}
		+\left\|\nabla\left(\dfrac1{\xi^2}\right)\right\|_{2}\left\|\Delta^2\phi_t\right\|_{s-2}\right.\nonumber \\
&\quad\left.
		+\left\|\nabla\left(\dfrac1{\xi^2}\right)\right\|_{s-2}\left\|\Delta^2\phi_t\right\|_{2}\right)
+C\left\|\nabla\xi\right\|_{3}\left\|\nabla\phi_t\right\|_{s}\left\|\nabla\Delta\phi_t\right\|_{s-1}
\nonumber\\
&\le C\left((1+\left\|\nabla\xi\right\|_{s-2})\left\|\xi_t\right\|_{s-1}\left\|\Delta^2\phi\right\|_{s-1}
		+\left(\left\|\nabla\xi\right\|_{s-2}+\left\|\nabla\xi\right\|_{3}\right)\left\|\nabla\phi_t\right\|_{s+1}\right)
\left\|\nabla\Delta\phi_t\right\|_{s-1}
		\nonumber \\
&\le C(\tau)\varepsilon^2  K_1K_2\left\|\Delta^2\phi\right\|_{s-1}^2+C(\tau)\varepsilon^2 K_1\left\|\nabla\phi_t\right\|_{s+1}^2+\tau\left\|\nabla\Delta\phi_t\right\|_{s-1}^2,\label{q.3.14}
\end{align}
where we have used the fact by Lemma \ref{L3.1} that
\begin{align}\label{q.3.44}
		\left\|\left(\dfrac{1}{\xi^2}\right)_t\right\|_{s-1}=\left\|-\dfrac{2}{\xi^3}\xi_t\right\|_{s-1}&\leq C\left\|\dfrac{1}{\xi^3}\right\|_\infty\left\|\xi_t\right\|_{s-1}+C\left\|\xi_t\right\|_\infty\left\|\dfrac{1}{\xi^3}\right\|_{s-1} \nonumber\\
& \leq C\left\|\xi_t\right\|_{s-1}(1+\left\|\nabla\xi\right\|_{s-2}).
	\end{align}
When $s=3$, we only estimate the second term $N_{1}^1$ in $N_{1}$ with $|\beta|=2$ since the case $|\beta|=1$ is more easier to estimate with similar way and the case $|\beta|=0$ is in fact $N_{1}^1=0$.
\begin{align}
N_{1}^1(|\beta|=2)&=\int_{\mathbb{T}^N}\left[\nabla_i\nabla_j\left(\dfrac{1}{\xi^2}\Delta^2\phi_t\right)
-\dfrac{1}{\xi^2}\nabla_i\nabla_j\Delta^2\phi_t\right]:\nabla_i\nabla_j\Delta\phi_t\nonumber\\
&=\int_{\mathbb{T}^N}\left[\nabla^2\left(\dfrac{1}{\xi^2}\right)\Delta^2\phi_t
+2\nabla\left(\dfrac{1}{\xi^2}\right)\nabla\Delta^2\phi_t\right]:\nabla^2\Delta\phi_t\nonumber\\
&\leq C\left(\left\|\nabla^2\left(\dfrac{1}{\xi^2}\right)\right\|_\infty\|\Delta^2\phi_t\|
+\left\|\nabla\left(\dfrac{1}{\xi^2}\right)\right\|_\infty\|\nabla\Delta^2\phi_t\|
\right)\left\|\nabla\Delta\phi_t\right\|_1
\nonumber\\
&\leq C\left(C\|\nabla\xi\|_3\|\Delta^2\phi_t\|
+C\|\nabla\xi\|_2\|\nabla\Delta^2\phi_t\|\right)\left\|\nabla\Delta\phi_t\right\|_2
\nonumber\\
&\leq C(\tau)\varepsilon^2 K_1\|\nabla\Delta\phi_t\|_2^2+\tau\left\|\nabla\Delta\phi_t\right\|_2^2.
\end{align}
Next, we deduce the estimate of $N_2-N_7$ by integration by parts that
\begin{align}
     N_2&=-\int_{\mathbb{T}^N}\nabla^{\beta}\nabla\left(v_t\cdot\nabla\psi+v\cdot\nabla\psi_t\right)\cdot\nabla^{\beta}\nabla\phi_t =\int_{\mathbb{T}^N}\nabla^{\beta}\left(v_t\cdot\nabla\psi+v\cdot\nabla\psi_t\right)\nabla^{\beta}\Delta\phi_t
       \nonumber\\
		&\le\left(\|v_t\|_{\infty}\|\nabla\psi\|_{s-1}+\|v_t\|_{s-1}\|\nabla\psi\|_{\infty}
          +\|v\|_{\infty}\|\nabla\psi_t\|_{s-1}+\|v\|_{s-1}\|\nabla\psi_t\|_{\infty}\right)\|\nabla^{\beta}\Delta\phi_t\|
\nonumber\\
&\le CK_1^\frac{1}{2}K_2^\frac{1}{2}\|\Delta\phi_t\|_{s-1}\le \tau\|\Delta\phi_t\|_{s-1}^2+C(\tau)K_1K_2,\label{q.3.15}\\
N_3&=-2\int_{\mathbb{T}^N}\nabla^{\beta}\nabla\left[\dfrac1{\xi}\nabla\left(\dfrac1{\xi}\right)\cdot\nabla\Delta\psi\right]_t\cdot\nabla^{\beta}\nabla\phi_t
=2\int_{\mathbb{T}^N}\nabla^{\beta}\left[\dfrac1{\xi}\nabla\left(\dfrac1{\xi}\right)\cdot\nabla\Delta\psi\right]_t\nabla^{\beta}\Delta\phi_t
\nonumber\\
&=2\int_{\mathbb{T}^N}\nabla^{\beta}\left[\left(\dfrac1{\xi}\right)_t\nabla\left(\dfrac1{\xi}\right)\cdot\nabla\Delta\psi
+\dfrac1{\xi}\nabla\left(\dfrac1{\xi}\right)_t\cdot\nabla\Delta\psi
+\dfrac1{\xi}\nabla\left(\dfrac1{\xi}\right)\cdot\nabla\Delta\psi_t\right]\nabla^{\beta}\Delta\phi_t
\nonumber\\
&\le C\left(\left\|\left(\dfrac1{\xi}\right)_t\right\|_{s-1}\left\|\nabla\left(\dfrac1{\xi}\right)\right\|_{s-1}\left\|\nabla\Delta\psi\right\|_{s-1}
+\left\|\dfrac1{\xi}\right\|_{s-1}\left\|\nabla\left(\dfrac1{\xi}\right)_t\right\|_{s-1}\left\|\nabla\Delta\psi\right\|_{s-1}
\right.\nonumber\\
        &\left.
\quad+\left\|\dfrac1{\xi}\right\|_{s-1}\left\|\nabla\left(\dfrac1{\xi}\right)\right\|_{s-1}\left\|\nabla\Delta\psi_t\right\|_{s-1}\right)
\|\nabla^{\beta}\Delta\phi_t\|
\nonumber\\
&\le C\left[(1+\left\|\nabla\xi\right\|_{s-2}) \left\|\xi_t\right\|_{s-1}\left\|\nabla\xi\right\|_{s-1}\left\|\nabla\Delta\psi\right\|_{s-1}
+(1+\left\|\nabla\xi\right\|_{s-2}) \left\|\nabla\xi\right\|_{s-1}\left\|\nabla\Delta\psi_t\right\|_{s-1}
 \right.\nonumber\\
        &\left.\quad+(1+\left\|\nabla\xi\right\|_{s-2})(1+ \left\|\nabla\xi\right\|_{s-1})\left\|\xi_t\right\|_{s}\left\|\nabla\Delta\psi\right\|_{s-1}
\right]
\|\Delta\phi_t\|_{s-1}
\nonumber\\
&\le C(\tau)\varepsilon^2K_1^3K_2+C(\tau)\varepsilon^2K_1^2\left\|\nabla\Delta\psi_t\right\|_{s-1}^2+\tau\|\Delta\phi_t\|_{s-1}^2+C,
\\
N_4&=-\int_{\mathbb{T}^N}\nabla^{\beta}\nabla\left[\dfrac1{\xi}\Delta\left(\dfrac1{\xi}\right)\Delta\psi\right]_t\cdot\nabla^{\beta}\nabla\phi_t
=\int_{\mathbb{T}^N}\nabla^{\beta}\left[\dfrac1{\xi}\Delta\left(\dfrac1{\xi}\right)\Delta\psi\right]_t\nabla^{\beta}\Delta\phi_t
\nonumber\\
&\le C\left(\left\|\left(\dfrac1{\xi}\right)_t\right\|_{s-1}\left\|\Delta\left(\dfrac1{\xi}\right)\right\|_{s-1}\left\|\Delta\psi\right\|_{s-1}
+\left\|\dfrac1{\xi}\right\|_{s-1}\left\|\Delta\left(\dfrac1{\xi}\right)_t\right\|_{s-1}\left\|\Delta\psi\right\|_{s-1}
\right.\nonumber\\
        &\left.
\quad+\left\|\dfrac1{\xi}\right\|_{s-1}\left\|\Delta\left(\dfrac1{\xi}\right)\right\|_{s-1}\left\|\Delta\psi_t\right\|_{s-1}\right)
\|\nabla^{\beta}\Delta\phi_t\|
\nonumber\\
&\le C\left((1+\left\|\nabla\xi\right\|_{s-2})\left\|\xi_t\right\|_{s-1}\left\|\nabla\xi\right\|_{s}\left\|\Delta\psi\right\|_{s-1}
+(1+\left\|\nabla\xi\right\|_{s-2})(1+\left\|\nabla\xi\right\|_{s})\left\|\xi_t\right\|_{s+1}\left\|\Delta\psi\right\|_{s-1}
 \right.\nonumber\\
        &\left.\quad+(1+\left\|\nabla\xi\right\|_{s-2})\left\|\nabla\xi\right\|_{s}\left\|\Delta\psi_t\right\|_{s-1}
\right)
\|\Delta\phi_t\|_{s-1}
\nonumber\\
&\le C(\tau)\varepsilon^2 K_1^3K_2+C(\tau)\varepsilon^2K_1^2\left\|\Delta\psi_t\right\|_{s-1}^2+\tau\|\Delta\phi_t\|_{s-1}^2,
\\
		N_5&=\int_{\mathbb{T}^N}\nabla^{\beta}\nabla\left(\dfrac{1}{\xi}\Delta(\psi^3-\psi)\right)_t\cdot\nabla^{\beta}\nabla\phi_t \nonumber  \\ &=-\int_{\mathbb{T}^N}\nabla^{\beta}\left[\left(\dfrac{1}{\xi}\right)_t\Delta\left(\psi^3-\psi\right)\right]\nabla^{\beta}\Delta\phi_t
+\int_{\mathbb{T}^N}\nabla^{\beta}\nabla\left[\frac{1}{\xi}\Delta\left(\psi^3-\psi\right)_t\right]\cdot\nabla^{\beta}\nabla\phi_t
		\nonumber \\	
&\le \left\|\left(\dfrac1\xi\right)_t\right\|_{s-1}\left\|\Delta\left(\psi^3-\psi\right)\right\|_{s-1}\|\nabla^{\beta}\Delta\phi_t\|
+\left\|\dfrac1\xi\right\|_{s}\left\|\Delta\left(\psi^3-\psi\right)_t\right\|_{s}\|\nabla^{\beta}\nabla\phi_t\|
\nonumber\\
&\le C(1+\left\|\nabla\xi\right\|_{s-2})\left\|\xi_t\right\|_{s-1}\left\|\nabla\psi\right\|_s\|\nabla^{\beta}\Delta\phi_t\|
+C(1+\left\|\nabla\xi\right\|_{s-1})\left\|\left(3\psi^2-1\right)\psi_t\right\|_{s+2} \|\nabla^{\beta}\nabla\phi_t\|
\nonumber\\
&\le C(\tau)\varepsilon^2K_1^2K_2+\tau\|\Delta\phi_t\|_{s-1}^2
+CK_1^\frac12\left\|3\psi^2-1\right\|_{s+2}\left\|\psi_t\right\|_{s+2}\|\nabla\phi_t\|_{s-1}
\nonumber\\
&\le C(\tau)\varepsilon^2K_1^2K_2+\tau\|\Delta\phi_t\|_{s-1}^2
+CK_1\left\|\psi_t\right\|_{s+2}\|\nabla\phi_t\|_{s-1}
\nonumber\\
&\le C(\tau)\varepsilon^2K_1^2K_2+C(\tau)K_1^2\|\nabla\phi_t\|_{s-1}^2+\tau\|\Delta\phi_t\|_{s-1}^2+\tau\left\|\psi_t\right\|_{s+2}^2,
\\
N_6&=-2\int_{\mathbb{T}^N}\nabla\left(\dfrac{1}{\xi^2}\right)\cdot\nabla^{\beta}\nabla\nabla\phi_t\cdot\nabla^{\beta}\nabla\Delta\phi_t
\nonumber\\
&
\leq C\left\|\nabla\left(\dfrac{1}{\xi^2}\right)\right\|_\infty\left\|\nabla^2\phi_t\right\|_{s-1}\left\|\nabla\Delta\phi_t\right\|_{s-1}
\leq C(\tau)\varepsilon^2 K_1\left\|\nabla\Delta\phi_t\right\|_{s-1}^2+\tau\left\|\nabla\Delta\phi_t\right\|_{s-1}^2,
\\
N_7&=-\int_{\mathbb{T}^N}\Delta\left(\dfrac{1}{\xi^2}\right)\nabla^{\beta}\nabla\Delta\phi_t\cdot\nabla^{\beta}\nabla\phi_t
\nonumber\\
&\leq C\left\|\Delta\left(\dfrac{1}{\xi^2}\right)\right\|_\infty\left\|\nabla\Delta\phi_t\right\|_{s-1}\left\|\nabla\phi_t\right\|_{s-1}
\leq C(\tau)\varepsilon^2 K_1\left\|\nabla\Delta\phi_t\right\|_{s-1}^2+\tau\left\|\nabla\Delta\phi_t\right\|_{s-1}^2.\label{qq.3.23}
\end{align}
	Putting the estimates on $N_i(i=1,...,7)$ into (\ref{q.3.13}), using the Cauchy inequality, and summing over $\beta$, we have after taking $\varepsilon$ and $\tau$ small enough that
	\begin{align}
		\label{cchi-t-2}
		&\dfrac{\rm d}{{\rm d}t}\displaystyle\sum_{|\beta|\le s-1}\left\|\nabla^{\beta}\nabla\phi_t\right\|^2
		+\displaystyle\sum_{|\beta|\le s-1}\|\nabla^{\beta}\nabla\Delta\phi_t\|^2  \nonumber\\
&
		 \quad~\le C(\tau)\varepsilon^2 K_1K_2\left\|\Delta^2\phi\right\|_{s-1}^2+C(\tau)K_1^2\|\nabla\phi_t\|_{s-1}^2+C(\tau)K_1^3K_2\nonumber\\
&\quad+\tau \left\|\psi_t\right\|_{s+2}^2+C(\tau)\varepsilon^2K_1^2\left\|\nabla\psi_t\right\|_{s+1}^2+C.
\end{align}
Then integrating (\ref{cchi-t-2}) over $[0, t]\subseteq[0,  T_0]$,
we find that
\begin{align}\label{cchi-t-22}
		&\left\|\nabla\phi_t\right\|_{s-1}^2(t)+\int_0^t\left\|\nabla\Delta\phi_t\right\|_{s-1}^2{\rm d}s\nonumber\\
&\leq \left\|\nabla\phi_t(x,0)\right\|_{s-1}^2+ C(\tau)\varepsilon^2 K_1K_2\int_0^t\left\|\Delta^2\phi\right\|_{s-1}^2+C(\tau)\int_0^tK_1^2\|\nabla\phi_t\|_{s-1}^2+C(\tau)K_1^3K_2T_0\nonumber\\
&\quad+\tau \int_0^t\left\|\psi_t\right\|_{s+2}^2+C(\tau)\varepsilon^2K_1^2\int_0^t\left\|\nabla\psi_t\right\|_{s+1}^2+CT_0.
\end{align}
From the equation $(\ref{q.3.2})_3$ and the constraints of the initial data, we see that
	\begin{align*}
		&\left\|\nabla\phi_t(x,0)\right\|_{s-1}\le\left\|\phi_t(x,0)\right\|_s\nonumber\\
		&\le C\left(\left\|\dfrac1{\xi^2(x,0)}\Delta^2\phi(x,0)\right\|_s
+\left\|v(x,0)\cdot\nabla\psi(x,0)\right\|_s
        +\left\|\dfrac1{\xi(x,0)}\nabla\left(\dfrac1{\xi(x,0)}\right)\cdot\nabla\Delta\psi(x,0)\right\|_s
         \right.\nonumber\\
        &\left.
        \quad+\left\|\dfrac1{\xi(x,0)}\Delta\left(\dfrac1{\xi(x,0)}\right)\Delta\psi(x,0)\right\|_s
+\left\|\dfrac1{\xi(x,0)}\Delta\left(\psi^3(x,0)-\psi(x,0)\right)\right\|_s\right)
\le C.
\end{align*}
Then applying the Gronwall inequality to (\ref{cchi-t-22}), taking $\tau$, $\varepsilon$ small enough, by (\ref{q.3.11}),  
we deduce that
\begin{align*}
		\left\|\nabla\phi_t\right\|_{s-1}^2(t)&\leq {\rm exp}^{(CK_1^2T_0)}\left(\left\|\nabla\phi_t(x,0)\right\|_{s-1}^2+C(\tau)\varepsilon^2 K_1K_2\int_0^t\left\|\Delta^2\phi\right\|_{s-1}^2+C(\tau)K_1^3K_2T_0\right.\nonumber\\
&\left.\quad+\tau \int_0^t\left\|\psi_t\right\|_{s+2}^2+C(\tau)\varepsilon^2K_1^2\int_0^t\left\|\nabla\psi_t\right\|_{s+1}^2+CT_0\right)\le C,
\end{align*}
 for $t\in[0, T_0]$, provided that $T_0<T_2\triangleq{\rm min}\{T_1,K_1^{-3}K_2^{-1}\}$.
 It follows from (\ref{cchi-t-2}) and $\rm Poincar\acute{e}$ inequality that
\begin{align}\label{chi_t}
		\left\|\nabla\phi_t\right\|_{s-1}^2(t)+\int_0^t\left\|\nabla\phi_t\right\|_{s+1}^2{\rm d}s\le C,
\end{align}
for $t\in[0, T_0]$, provided that $T_0<T_2$.

In addition, operating $\partial_t$ to $(\ref{q.3.2})_3$, multiplying the result equation by $\phi_t$, and then integrating the result over $\mathbb{T}^N$, we have
	\begin{align}\label{qq.3.19}
		&\dfrac12\dfrac{\rm d}{{\rm d}t}\int_{\mathbb{T}^N}|\phi_t|^2+\int_{\mathbb{T}^N}\dfrac{1}{\xi^2}|\Delta\phi_t|^2
		\nonumber\\
		&=-\int_{\mathbb{T}^N}\left(\dfrac{1}{\xi^2}\right)_t\Delta^2\phi\phi_t
-\int_{\mathbb{T}^N}(v\cdot\nabla\psi)_t\phi_t
-2\int_{\mathbb{T}^N}\left(\dfrac{1}{\xi}\nabla\left(\dfrac1{\xi}\right)\cdot\nabla\Delta\psi\right)_t\phi_t
\nonumber\\
		&\quad-\int_{\mathbb{T}^N}\left(\dfrac{1}{\xi}\Delta\left(\dfrac1{\xi}\right)\Delta\psi\right)_t\phi_t
+\int_{\mathbb{T}^N}\left(\dfrac{1}{\xi}\Delta(\psi^3-\psi)\right)_t\phi_t
\nonumber\\
		&\quad-\int_{\mathbb{T}^N}\Delta\left(\dfrac{1}{\xi^2}\right)\Delta\phi_t\phi_t -2\int_{\mathbb{T}^N}\nabla\left(\dfrac{1}{\xi^2}\right)\cdot\nabla\phi_t\Delta\phi_t
\nonumber\\
		&\le C\|\phi_t\|\left[\left\|\left(\dfrac{1}{\xi^2}\right)_t\right\|_\infty\|\Delta^2\phi\|+\|v_t\|\|\nabla\psi\|_\infty+\|v\|_\infty\|\nabla\psi_t\|
+\left\|\left(\frac{1}{\xi}\right)_t\right\|\left\|\nabla\left(\frac{1}{\xi}\right)\right\|_\infty\left\|\nabla\Delta\psi\right\|_\infty
\right.\nonumber\\
        &\left.\quad
        +\left\|\frac{1}{\xi}\right\|_\infty\left\|\nabla\left(\frac{1}{\xi}\right)_t\right\|\left\|\nabla\Delta\psi\right\|_\infty
        +\left\|\frac{1}{\xi}\right\|_\infty\left\|\nabla\left(\frac{1}{\xi}\right)\right\|_\infty\left\|\nabla\Delta\psi_t\right\|
        +\left\|\left(\frac{1}{\xi}\right)_t\right\|\left\|\Delta\left(\frac{1}{\xi}\right)\right\|_\infty\left\|\Delta\psi\right\|_\infty
\right.\nonumber\\
        &\left.\quad
        +\left\|\frac{1}{\xi}\right\|_\infty\left\|\Delta\left(\frac{1}{\xi}\right)_t\right\|\left\|\Delta\psi\right\|_\infty
        +\left\|\frac{1}{\xi}\right\|_\infty\left\|\Delta\left(\frac{1}{\xi}\right)\right\|_\infty\left\|\Delta\psi_t\right\|
+\left\|\left(\frac{1}{\xi}\right)_t\right\|\left\|\Delta\left(\psi^3-\psi\right)\right\|_\infty
\right.\nonumber\\
        &\left.\quad+\left\|\frac{1}{\xi}\right\|_\infty\left\|\Delta\left(\psi^3-\psi\right)_t\right\|
        +\left\|\Delta\left(\frac{1}{\xi^2}\right)\right\|_\infty\|\Delta\phi_t\|
\right]
+C\left\|\nabla\left(\frac{1}{\xi^2}\right)\right\|_\infty\|\nabla\phi_t\|\|\Delta\phi_t\|\nonumber\\
&\le C K_1^3K_2\|\phi_t\|^2+\varepsilon^2K_1\|\Delta\phi_t\|^2+\tau\|\nabla\phi_t\|^2+C
\nonumber\\
&\le C K_1^3K_2\|\phi_t\|^2+\varepsilon^2K_1\|\Delta\phi_t\|^2+\tau\|\nabla^2\phi_t\|^2+C,
	\end{align}
where in the last inequality we have used the $\rm Poincar\acute{e}$ inequality.
Integrating the result over $[0, t]\subseteq[0,  T_0]$, taking $\varepsilon$, $\tau$ small enough, we obtain after using the $\rm Poincar\acute{e}$ inequality and the Gronwall inequality that
	\begin{align}\label{chit}
		\left\|\phi_t\right\|^2(t)+\int_0^t\left\|\nabla\phi_t\right\|_1^2{\rm d}s\le C,
\end{align}
for $t\in[0,\ T_0]$, provided that $T_0<T_2$.

Thus, it follows from (\ref{q.3.11}), (\ref{chi_t}) and (\ref{chit})that
	\begin{align}\label{q.3.19}
		\left\|\phi_t\right\|_{s}^2(t)+\int_0^t\left\|\phi_t\right\|_{s+2}^2{\rm d}s\le C.
\end{align}

 $\mathbf{Step\ two}$: Estimates of $\rho$ and u.

Taking the $L^2$ inner product of the $(\ref{q.3.3})_1$ and $(\ref{q.3.3})_2$ with the $\dfrac{1}{\varepsilon^2}\dfrac{P'(\xi)}{\xi}D^{\alpha_1}(\rho-1)$ and $\xi D^{\alpha_1} u$ respectively, and then integrating by parts, one deduces
	\begin{align}\label{q.3.20}
		&\dfrac{1}{2}\dfrac{\rm d}{{\rm d}t}\int_{\mathbb{T}^N}\left(\dfrac{P'\left(\xi\right)}{\xi}\left|\dfrac{1}{\varepsilon} D^{\alpha_1}\left(\rho-1\right)\right|^2+\xi| D^{\alpha_1} u|^2\right)\nonumber\\
&\quad+\int_{\mathbb{T}^N}\nu(\xi,\psi)| D^{\alpha_1} \nabla u|^2+\int_{\mathbb{T}^N}\eta(\xi,\psi)| D^{\alpha_1} (\nabla\cdot u)|^2=\sum_{i=1}^{10}I_i,
\end{align}
	where
\begin{align}
	&I_1=\cfrac{1}{2}\displaystyle\int_{\mathbb{T}^N}\left[\dfrac{1}{\varepsilon^2}| D^{\alpha_1}\left(\rho-1\right)|^2\partial_t\left(\dfrac{P'\left(\xi\right)}{\xi}\right)+\xi_t| D^{\alpha_1} u|^2\right],\nonumber\\
	&I_2=\cfrac{1}{2}\displaystyle\int_{\mathbb{T}^N}\left[\dfrac{1}{\varepsilon^2}| D^{\alpha_1}\left(\rho-1\right)|^2\nabla\cdot\left(\dfrac{P'(\xi)}{\xi} v\right)+\nabla\cdot\left(\xi v\right)| D^{\alpha_1} u|^2\right],\nonumber\\
	&I_3=\displaystyle\int_{\mathbb{T}^N}\dfrac{1}{\varepsilon^2} P''\left(\xi\right) D^{\alpha_1}\left(\rho-1\right) D^{\alpha_1} u\cdot\nabla\xi,\nonumber\\
	&I_4=-\dfrac{1}{\varepsilon^2}\displaystyle\int_{\mathbb{T}^N}\dfrac{P'(\xi)}{\xi} D^{\alpha_1}(\rho-1)\left\{\left[D^{\alpha_1}(v\cdot\nabla\rho)-(v\cdot\nabla)D^{\alpha_1}\rho\right]
+\left[D^{\alpha_1}(\xi\nabla\cdot u)-\xi\nabla\cdot D^{\alpha_1}u\right]\right\},\nonumber\\
	&I_5=-\displaystyle\int_{\mathbb{T}^N} D^{\alpha_1}\left(\Delta\phi\nabla\phi\right)\cdot D^{\alpha_1} u,\nonumber\\
	&I_6=-\displaystyle\int_{\mathbb{T}^N}\xi \left[D^{\alpha_1}\left(v\cdot\nabla u\right)-\left(v\cdot\nabla\right)D^{\alpha_1}  u\right]\cdot D^{\alpha_1}  u ,\nonumber\\
	&I_7=-\dfrac{1}{\varepsilon^2}\displaystyle\int_{\mathbb{T}^N}\xi \left[D^{\alpha_1}\left(\dfrac{P'\left(\xi\right)}{\xi}\nabla \rho\right)-\dfrac{P'\left(\xi\right)}{\xi}\nabla D^{\alpha_1}\rho\right]\cdot D^{\alpha_1} u,\nonumber\\
	&I_8=\displaystyle\int_{\mathbb{T}^N}\xi \left[D^{\alpha_1}\left(\dfrac{\nu(\xi,\psi)}{\xi}\Delta u\right)-
	\dfrac{\nu(\xi,\psi)}{\xi}\Delta D^{\alpha_1} u\right] \cdot D^{\alpha_1}  u \nonumber\\
&\qquad+\displaystyle\int_{\mathbb{T}^N}\xi \left[D^{\alpha_1}\left(\dfrac{\eta(\xi,\psi)}{\xi}\nabla\left(\nabla\cdot u\right)\right)-\dfrac{\eta(\xi,\psi)}{\xi}\nabla D^{\alpha_1}\left(\nabla\cdot u\right)\right] \cdot D^{\alpha_1}  u ,\nonumber\\
	&I_9=-\displaystyle\int_{\mathbb{T}^N}\xi \left[D^{\alpha_1}\left(\dfrac{1}{\xi}\left(\Delta\phi\nabla\phi\right)\right)-\dfrac{1}{\xi}D^{\alpha_1}\left(\Delta\phi\nabla\phi\right)\right]\cdot D^{\alpha_1}  u ,\nonumber\\
&  I_{10}=-\displaystyle\int_{\mathbb{T}^N}\nabla\nu(\xi,\psi)\cdot \nabla D^{\alpha_1} u\cdot D^{\alpha_1}  u
    -\displaystyle\int_{\mathbb{T}^N}\nabla\eta(\xi,\psi)\cdot D^{\alpha_1}u D^{\alpha_1}( \nabla\cdot u),
\nonumber
\end{align}
including the calculation of singular term that
\begin{align}\label{odd}
&\int_{\mathbb{T}^N}\xi\nabla\cdot D^{\alpha_1}u\left[\dfrac{1}{\varepsilon^2}\dfrac{P'(\xi)}{\xi}D^{\alpha_1}(\rho-1)\right]
+\int_{\mathbb{T}^N}\dfrac{1}{\varepsilon^2}\dfrac{P'\left(\xi\right)}{\xi}\nabla D^{\alpha_1}\rho\cdot\xi D^{\alpha_1}u\nonumber\\
&=\int_{\mathbb{T}^N}\dfrac{1}{\varepsilon^2}P'(\xi)\nabla\cdot D^{\alpha_1}uD^{\alpha_1}(\rho-1)
+\int_{\mathbb{T}^N}\dfrac{1}{\varepsilon^2}P'\left(\xi\right)\nabla D^{\alpha_1}\rho\cdot D^{\alpha_1}u\nonumber\\
&=-\int_{\mathbb{T}^N}\dfrac{1}{\varepsilon^2} P''\left(\xi\right) D^{\alpha_1}\left(\rho-1\right) D^{\alpha_1} u\cdot\nabla\xi
-\int_{\mathbb{T}^N}\dfrac{1}{\varepsilon^2}P'\left(\xi\right) D^{\alpha_1}u\cdot\nabla D^{\alpha_1}(\rho-1)\nonumber\\
&\quad+\int_{\mathbb{T}^N}\dfrac{1}{\varepsilon^2}P'\left(\xi\right)\nabla D^{\alpha_1}\rho\cdot D^{\alpha_1}u\nonumber\\
&=-\int_{\mathbb{T}^N}\dfrac{1}{\varepsilon^2} P''\left(\xi\right) D^{\alpha_1}\left(\rho-1\right) D^{\alpha_1} u\cdot\nabla\xi=-I_3.
\end{align}

Firstly, the estimates of $I_1-I_3$ will be given. Choose $\kappa$ small enough so that $E_{s+2}(\rho,u)\backsim\widetilde{E}_{s+2}(\rho,u)$ hold. Then by the smoothness of the pressure $P(\cdot)$, we get
		\begin{align}
			|I_1|&\leq\left\|\dfrac{\xi P''(\xi)-P'(\xi)}{\xi^2}\right\|_\infty\|\xi_t\|_\infty\left\|\dfrac{1}{\varepsilon}D^{\alpha_1}(\rho-1)\right\|^2+\left\|\xi_t\right\|_\infty\left\|D^{\alpha_1}  u\right\|^2\nonumber \\
			&\leq C\left\|\xi_t\right\|_2\left(\left\|\dfrac{1}{\varepsilon}D^{\alpha_1}\left(\rho-1\right)\right\|^2+\left\|D^{\alpha_1}  u\right\|^2\right)\nonumber \\
			&\leq C\varepsilon K_2^{\frac{1}{2}}\left(\left\|\dfrac{1}{\varepsilon}D^{\alpha_1}(\rho-1)\right\|^2+\left\|D^{\alpha_1}  u\right\|^2\right),\label{q.3.21}\\
			| I_2|&\leq\left(\left\|\dfrac{\xi P''\left(\xi\right)-P'\left(\xi\right)}{\xi^2}\right\|_\infty\left\|\nabla\xi\right\|_\infty\left\|v\right\|_\infty+\left\|\dfrac{P'\left(\xi\right)}{\xi}\right\|_\infty\|\nabla\cdot v\|_\infty\right)\left\|\dfrac{1}{\varepsilon}D^{\alpha_1}\left(\rho-1\right)\right\|^2\nonumber \\
			&\quad+\left(\left\|\xi\right\|_\infty\left\|\nabla\cdot v\right\|_\infty+\left\|v\right\|_\infty\left\|\nabla\xi\right\|_\infty\right)\left\|D^{\alpha_1}  u\right\|^2\nonumber \\		&\leq C\left(\varepsilon\left\|\frac{1}{\varepsilon}\nabla\xi\right\|_2\left\|v\right\|_2+\left\|v\right\|_3\right)\left(\left\|\dfrac{1}{\varepsilon}D^{\alpha_1}
\left(\rho-1\right)\right\|^2+\left\|D^{\alpha_1}  u\right\|^2\right)\nonumber \\
			&\leq C\left(\varepsilon K_1+K_1^\frac{1}{2}\right)\left(\left\|\dfrac{1}{\varepsilon}D^{\alpha_1}(\rho-1)\right\|^2+\left\|D^{\alpha_1}  u\right\|^2\right),\label{q.3.22}\\
			| I_3|&\leq\dfrac{1}{\varepsilon}\left\|P''\left(\xi\right)\right\|_\infty \left\|\nabla\xi\right\|_\infty\left(\left\|\dfrac{1}{\varepsilon}D^{\alpha_1}\left(\rho-1\right)\right\|^2+\left\|D^{\alpha_1}  u\right\|^2\right)\nonumber \\
			&\leq C\left\|\dfrac{1}{\varepsilon}\nabla\xi\right\|_2\left(\left\|\dfrac{1}{\varepsilon}D^{\alpha_1}\left(\rho-1\right)\right\|^2+\left\|D^{\alpha_1}  u\right\|^2\right)\nonumber \\
			&\leq CK_1^\frac{1}{2}\left(\left\|\dfrac{1}{\varepsilon}D^{\alpha_1}\left(\rho-1\right)\right\|^2+\left\|D^{\alpha_1}u\right\|^2\right).\label{q.3.23}
		\end{align}

 Secondly, we show the estimates of the $I_4-I_9$ in the following three cases.

$\mathbf {Case\ 1}$: $|\alpha_1|=0$ or $D^{\alpha_1}=\partial_t$.

If $|\alpha_1|=0$, then the quantities $I_j,\ 4\leq j\leq10$, are all equal to zero expect for $I_5$ and $I_{10}$. Thus it suffices to estimate $I_5$ and $I_{10}$. Choose $\kappa$ small enough so that $|\psi-\psi_0|<\frac12$. Using (\ref{q.3.11}) and Lemma \ref{L3.1}, we have
	\begin{align}
			| I_5|&\leq C\left\|u\right\|\left\|\nabla\phi\right\|_\infty\left\|\Delta\phi\right\|\leq C\|u\|^2+C,\label{q.3.24}\\
    | I_{10}|&\leq C\left\|\nabla\nu(\xi,\psi)\right\|_\infty\|\nabla u\|\|u\|+ C\left\|\nabla\eta(\xi,\psi)\right\|_\infty\| \nabla\cdot u\|\|u\|\nonumber\\
            &\leq C\left\|\nabla\nu(\xi,\psi)\right\|_2\|\nabla u\|\|u\|+ C\left\|\nabla\eta(\xi,\psi)\right\|_2\| \nabla\cdot u\|\|u\|\nonumber\\
            &\leq C(\left\|\nabla_\xi\nu(\xi,\psi)\right\|_2+\left\|\nabla_\psi\nu(\xi,\psi)\right\|_2)\|\nabla u\|\|u\|\nonumber\\
            &\quad+ C(\left\|\nabla_\xi\eta(\xi,\psi)\right\|_2+\left\|\nabla_\psi\eta(\xi,\psi)\right\|_2)\| \nabla\cdot u\|\|u\|\nonumber\\
            &\leq C\left[(1+\|\xi\|_\infty^{2})\|\nabla\xi\|_{2}+(1+\|\psi\|_\infty^{2})\|\nabla\psi\|_{2}\right](\|\nabla u\|+\| \nabla\cdot u\|)\|u\|\nonumber\\
             &\leq CK_1^\frac12(\|\nabla u\|+\| \nabla\cdot u\|)\|u\|\leq C(\tau)K_1\|u\|^2+\tau(\|\nabla u\|^2+\| \nabla\cdot u\|^2).\label{q.3.24000}
             \end{align}

If $D^{\alpha_1}=\partial_t$, then using (\ref{q.3.11}), (\ref{q.3.19}) and the Cauchy inequality lead to
		\begin{align}
			| I_4|&\leq C\left\|\dfrac{1}{\varepsilon}\rho_t\right\|\left(\|v_t\|_\infty\left\|\dfrac{1}{\varepsilon}\nabla\rho\right\|+\left\|\dfrac{1}{\varepsilon}\xi_t\right\|_\infty\|\nabla\cdot u\|\right)\nonumber \\
			&\leq CK_2^{\frac{1}{2}}\left(\left\|\dfrac{1}{\varepsilon}\rho_t\right\|^2+\left\|\dfrac{1}{\varepsilon}\nabla\rho\right\|^2+\left\|\nabla\cdot u\right\|^2\right),\label{q.3.25}\\
			| I_5|&\leq C\left(\left\|\Delta\phi_t\right\|\left\|\nabla\phi\right\|_\infty+\left\|\Delta\phi\right\|_\infty\left\|\nabla\phi_t\right\|\right)\|u_t\|\nonumber \\
			&\leq C\left(\left\|\Delta\phi_t\right\|\left\|\nabla\phi\right\|_2+\left\|\Delta\phi\right\|_2\left\|\nabla\phi_t\right\|\right)\|u_t\|\leq C\left\|u_t\right\|^2+C,\label{q.3.26}\\
			| I_6|&\leq C\left\|v_t\right\|_\infty\left\|\nabla u\right\|\left\|u_t\right\|\leq CK_2^{\frac{1}{2}}\left(\left\|\nabla u\right\|^2+\left\|u_t\right\|^2\right),\label{q.3.27}\\
			| I_7|&\leq C\left\|\dfrac{1}{\varepsilon}\xi_t\right\|_\infty\left\|\dfrac{1}{\varepsilon}\nabla\rho\right\|\left\|u_t\right\|\leq CK_2^{\frac{1}{2}}\left(\left\|\dfrac{1}{\varepsilon}\nabla\rho\right\|^2+\left\|u_t\right\|^2\right),\label{q.3.28}\\
			| I_8|&=\left|\int_{\mathbb{T}^N}\xi \left[\partial_t\left(\dfrac{\nu(\xi,\psi)}{\xi}\Delta u\right)-
	\dfrac{\nu(\xi,\psi)}{\xi}\Delta u_t\right] \cdot  u_t \right.\nonumber\\
&\left.\quad+\int_{\mathbb{T}^N}\xi \left[\partial_t\left(\dfrac{\eta(\xi,\psi)}{\xi}\nabla(\nabla\cdot u)\right)-
	\dfrac{\eta(\xi,\psi)}{\xi}\nabla(\nabla\cdot u)_t\right] \cdot  u_t\right|\nonumber\\
&=\left|\int_{\mathbb{T}^N}\xi \partial_t\left(\dfrac{\nu(\xi,\psi)}{\xi}\right)\Delta u \cdot   u_t
+\int_{\mathbb{T}^N}\xi \partial_t\left(\dfrac{\eta(\xi,\psi)}{\xi}\right)\nabla(\nabla\cdot u) \cdot   u_t\right|\nonumber\\
&=\left|\int_{\mathbb{T}^N} \dfrac{\partial_t\nu(\xi,\psi)\xi-\nu(\xi,\psi)\xi_t}{\xi}\Delta u \cdot   u_t
+\int_{\mathbb{T}^N} \dfrac{\partial_t\eta(\xi,\psi)\xi-\eta(\xi,\psi)\xi_t}{\xi}\nabla(\nabla\cdot u) \cdot   u_t \right|\nonumber\\
&=\left|\int_{\mathbb{T}^N} \left(\partial_t\nu(\xi,\psi)-\dfrac{\nu(\xi,\psi)\xi_t}{\xi}\right)\Delta u \cdot   u_t
+\int_{\mathbb{T}^N} \left(\partial_t\eta(\xi,\psi)-\dfrac{\eta(\xi,\psi)\xi_t}{\xi}\right)\nabla(\nabla\cdot u) \cdot   u_t \right|\nonumber\\
&\le\left|\int_{\mathbb{T}^N} \left(\nu_\xi(\xi,\psi)\xi_t+\nu_\psi(\xi,\psi)\psi_t-\dfrac{\nu(\xi,\psi)\xi_t}{\xi}\right)\Delta u \cdot   u_t \right|\nonumber\\
&\quad+\left|\int_{\mathbb{T}^N} \left(\eta_\xi(\xi,\psi)\xi_t+\eta_\psi(\xi,\psi)\psi_t-\dfrac{\eta(\xi,\psi)\xi_t}{\xi}\right)\nabla(\nabla\cdot u) \cdot   u_t \right| \nonumber\\
&\le\left(\varepsilon\|\nu_\xi(\xi,\psi)\|_\infty\left\|\dfrac{\xi_t}{\varepsilon}\right\|_\infty+\|\nu_\psi(\xi,\psi)\|_\infty\|\psi_t\|_\infty
+\varepsilon\left\|\dfrac{\nu(\xi,\psi)}{\xi}\right\|_\infty\left\|\dfrac{\xi_t}{\varepsilon}\right\|_\infty\right)\|\Delta u\|\| u_t\|\nonumber\\
&\quad+\left(\varepsilon\|\eta_\xi(\xi,\psi)\|_\infty\left\|\dfrac{\xi_t}{\varepsilon}\right\|_\infty+\|\eta_\psi(\xi,\psi)\|_\infty\|\psi_t\|_\infty
+\varepsilon\left\|\dfrac{\eta(\xi,\psi)}{\xi}\right\|_\infty\left\|\dfrac{\xi_t}{\varepsilon}\right\|_\infty\right)\|\nabla(\nabla\cdot u)\|\| u_t\|\nonumber\\
&\le C\left(\varepsilon K_2^{\frac12}+K_2^{\frac12}\right)(\|\Delta u\|+\|\nabla(\nabla\cdot u)\|)\| u_t\|\nonumber\\
&\le C\left(\varepsilon K_2^{\frac12}+K_2^{\frac12}\right)(\| u\|_2^2+\| u_t\|^2),\label{q.3.29}\\
			| I_9|&\leq C\varepsilon\left\|\dfrac{1}{\varepsilon}\xi_t\right\|_\infty\left\|\Delta\phi\right\|\left\|\nabla\phi\right\|_\infty\left\|u_t\right\|\leq C\varepsilon^2K_2\left\|u_t\right\|^2+C,\\
|I_{10}|&=\left|-\displaystyle\int_{\mathbb{T}^N}\nabla\nu(\xi,\psi)\cdot \nabla u_t\cdot  u_t
    -\displaystyle\int_{\mathbb{T}^N} u_t\cdot\nabla\eta(\xi,\psi) ( \nabla\cdot u_t)\right|\nonumber\\
    &\leq C\left\|\nabla\nu(\xi,\psi)\right\|_\infty\|\nabla u_t\|\| u_t\|+ C\left\|\nabla\eta(\xi,\psi)\right\|_\infty\|  \nabla\cdot u_t\|\| u_t\|\nonumber\\
            &\leq C\left\|\nabla\nu(\xi,\psi)\right\|_2\|\nabla u_t\|\|u_t\|+ C\left\|\nabla\eta(\xi,\psi)\right\|_2\| \nabla\cdot u_t\|\|u_t\|\nonumber\\
            &\leq C\left[(1+\|\xi\|_\infty^{2})\|\nabla\xi\|_{2}+(1+\|\psi\|_\infty^{2})\|\nabla\psi\|_{2}\right](\|\nabla u_t\|+\| \nabla\cdot u_t\|)\|u_t\|\nonumber\\
             &\leq CK_1^\frac12(\|\nabla u_t\|+\| \nabla\cdot u_t\|)\|u_t\|\leq C(\tau)K_1\|u_t\|^2+\tau(\|\nabla u_t\|^2+\| \nabla\cdot u_t\|^2)\label{q.3.30},
		\end{align}
where one has used Lemma \ref{L3.1} in the estimate of $I_{10}$.

$\mathbf {Cases\ 2}$: $D^{\alpha_1}=\nabla^{\alpha_1}$ for $1\leq|\alpha_1|\leq s+2$, where $s\geq3$.

It follows from Lemma \ref{L3.1} that
	\begin{align}
		| I_4|&\leq C\left\|\dfrac{1}{\varepsilon}\nabla^{\alpha_1}\left(\rho-1\right)\right\|\left(\left\|\nabla v\right\|_\infty\left\|\dfrac{1}{\varepsilon}\nabla\rho\right\|_{s+1}
		+\left\|\nabla v\right\|_{s+1}\left\|\dfrac{1}{\varepsilon}\nabla\rho\right\|_\infty
		\nonumber \right.\\
&\left.\quad+\left\|\nabla\cdot u\right\|_{s+1}\left\|\dfrac{1}{\varepsilon}\nabla\xi\right\|_\infty
		+\left\|\nabla\cdot u\right\|_\infty\left\|\dfrac{1}{\varepsilon}\nabla\xi\right\|_{s+1}\right)
		\nonumber \\
		&\leq CK_1^{\frac{1}{2}}\left(\left\|\dfrac{1}{\varepsilon}\left(\rho-1\right)\right\|_{s+2}^2+\left\|u\right\|_{s+2}^2\right),\label{q.3.31}\\
		|I_5|&=\left|-\int_{\mathbb{T}^N}\nabla^{\alpha_1}\left(\Delta\phi\nabla\phi\right)\cdot\nabla^{\alpha_1} u\right|
=\left|\int_{\mathbb{T}^N}\nabla^{\alpha_1}\left(\nabla\phi\otimes\nabla\phi-\frac{|\nabla\phi|^2}{2}I_N\right):\nabla\nabla^{\alpha_1} u\right|
  \nonumber\\
		&\leq C\left\|\nabla u\right\|_{s+2}\left\|\nabla\phi\right\|_\infty\left\|\nabla\phi\right\|_{s+2}
		\leq\tau\|\nabla u\|_{s+2}^2+C\left(\tau\right)\left\|\nabla\phi\right\|_{s+2}^2,\label{q.3.32}
\end{align}
where we have used integration by parts and (\ref{q.3.11}) in the estimate of $I_5$.
Similarly as in (\ref{q.3.31}), we get
\begin{align}\label{q.3.33}
		| I_6|&\leq C\left\|u\right\|_{s+2}(\left\|\nabla v\right\|_\infty\left\|\nabla u\right\|_{s+1}+\left\|\nabla v\right\|_{s+1}\left\|\nabla u\right\|_\infty)
		\leq CK_1^{\frac{1}{2}}\left\|u\right\|_{s+2}^2.
\end{align}

On the other hand, Lemma \ref{L3.1} gives
\begin{align}
		| I_7|&\leq C\frac{1}{\varepsilon^{2}}\left( \left\| \frac{\xi P''\left(\xi  \right)-P'\left(\xi \right)}{\xi ^{2} } \nabla\xi \right\|_\infty\left\|\nabla\rho\right\|_{s+1} + \left\|\nabla\rho\right\|_\infty\left\|\nabla\left( \frac{P'(\xi )}{\xi }  \right)\right\|_{s+1} \right) \left\|u\right\|_{s+2}
		\nonumber \\
		&\leq C\left(\left\|\dfrac{1}{\varepsilon}\nabla\xi\right\|_2\left\|\dfrac{1}{\varepsilon}\nabla\rho\right\|_{s+1}+\left\|\dfrac{1}{\varepsilon}\nabla\xi\right\|_{s+1}
\left\|\dfrac{1}{\varepsilon}\nabla\rho\right\|_2\right)\|u\|_{s+2}
		\nonumber \\
		&\leq CK_1^{\frac{1}{2}}\left(\left\|\dfrac{1}{\varepsilon}\left(\rho-1\right)\right\|_{s+2}^2+\left\|u\right\|_{s+2}^2\right),\label{q.3.34}\\
		| I_8|&\le C\left(\left\|\nabla\left(\dfrac{\nu(\xi,\psi)}{\xi}\right)\right\|_{s+1}\|\Delta u\|_\infty+\left\|\nabla\left(\dfrac{\nu(\xi,\psi)}{\xi}\right)\right\|_\infty\|\Delta u\|_{s+1}\right)\|u\|_{s+2}\nonumber\\
&\quad+C\left(\left\|\nabla\left(\dfrac{\eta(\xi,\psi)}{\xi}\right)\right\|_{s+1}\|\nabla(\nabla\cdot u) \|_\infty+\left\|\nabla\left(\dfrac{\eta(\xi,\psi)}{\xi}\right)\right\|_\infty\|\nabla(\nabla\cdot u)\|_{s+1}\right)\|u\|_{s+2}\nonumber\\
&\le C\left\|\nabla\left(\dfrac{\nu(\xi,\psi)}{\xi}\right)\right\|_{s+1}\|\Delta u\|_{s+1}\|u\|_{s+2}+C\left\|\nabla\left(\dfrac{\eta(\xi,\psi)}{\xi}\right)\right\|_{s+1}\|\nabla(\nabla\cdot u)\|_{s+1}\|u\|_{s+2}\nonumber\\
&\le C\left[(1+\|\xi\|_\infty^{s+1})\|\nabla\xi\|_{s+1}+(1+\|\psi\|_\infty^{s+1})\|\nabla\psi\|_{s+1}\right](\|\Delta u\|_{s+1}+\|\nabla(\nabla\cdot u)\|_{s+1})\|u\|_{s+2}\nonumber\\
&\le CK_1^\frac{1}{2}(\|\Delta u\|_{s+1}+\|\nabla(\nabla\cdot u)\|_{s+1})\|u\|_{s+2}\nonumber\\
&\leq C(\tau)K_1\left\|u\right\|_{s+2}^2+\tau\left(\left\|\nabla u\right\|_{s+2}^2+\left\|\nabla\cdot u\right\|_{s+2}^2\right),\label{q.3.35}\\
		| I_9|&\leq C\left\|u\right\|_{s+2}\left(\left\|\nabla \left(\frac{1}{\xi } \right) \right\|_\infty\left\|\Delta\phi\nabla\phi\right\|_{s+1} + \left\|\nabla\left( \frac{1}{\xi  }\right)\right\|_{s+1}\left\|\Delta\phi\nabla\phi\right\|_\infty
		\right)
		\nonumber \\
		&\leq C\varepsilon K_1^{\frac{1}{2}}\left\|u\right\|_{s+2}\left(\left\|\Delta\phi\right\|_2\left\|\nabla\phi\right\|_{s+1}+\left\|\Delta\phi\right\|_{s+1}\left\|\nabla\phi\right\|_2+\left\|\Delta\phi\right\|_2\left\|\nabla\phi\right\|_2\right)
		\nonumber \\
		&\leq C\varepsilon^2K_1\left\|u\right\|_{s+2}^2+C\left\|\nabla\phi\right\|_{s+2}^2,\label{q.3.360000}
	\end{align}
where one has used (\ref{q.3.11}) in the estimate of $I_9$. Moreover, using Lemma \ref{L3.1} and the Cauchy inequality, one obtains
\begin{align}
|I_{10}|&=\left|-\displaystyle\int_{\mathbb{T}^N}\nabla\nu(\xi,\psi)\cdot \nabla D^{\alpha_1} u\cdot D^{\alpha_1}  u
    -\displaystyle\int_{\mathbb{T}^N}\nabla\eta(\xi,\psi)\cdot D^{\alpha_1}u D^{\alpha_1}( \nabla\cdot u)\right|\nonumber\\
    &\leq C\left\|\nabla\nu(\xi,\psi)\right\|_\infty\|\nabla u\|_{s+2}\|  u\|_{s+2}+ C\left\|\nabla\eta(\xi,\psi)\right\|_\infty\|  \nabla\cdot u\|_{s+2}\|  u\|_{s+2}\nonumber\\
            &\leq C\left\|\nabla\nu(\xi,\psi)\right\|_2\|\nabla u\|_{s+2}\|u\|_{s+2}+ C\left\|\nabla\eta(\xi,\psi)\right\|_2\| \nabla\cdot u\|_{s+2}\|u\|_{s+2}\nonumber\\
            &\leq C\left[(1+\|\xi\|_\infty^{2})\|\nabla\xi\|_{2}+(1+\|\psi\|_\infty^{2})\|\nabla\psi\|_{2}\right](\|\nabla u\|_{s+2}+\| \nabla\cdot u\|_{s+2})\|u\|_{s+2}\nonumber\\
             &\leq CK_1^\frac12(\|\nabla u\|_{s+2}+\| \nabla\cdot u\|_{s+2})\|u\|_{s+2}\nonumber\\
             &\leq C(\tau)K_1\|u\|_{s+2}^2+\tau(\|\nabla u\|_{s+2}^2+\| \nabla\cdot u\|_{s+2}^2).\label{q.3.36}
\end{align}
	
Putting the estimates (\ref{q.3.21})-(\ref{q.3.23}) and (\ref{q.3.31})-(\ref{q.3.36}) together, summing over $\alpha_1$, and then choosing $\tau$ small enough, we obtain
	\begin{align}\label{q.3.37}
		&\dfrac{\rm d}{{\rm d}t}\sum_{|\alpha_1|\leq s+2}\int_{\mathbb{T}^N}\left(\dfrac{P'(\xi)}{\xi}\left|\dfrac{1}{\varepsilon} \nabla^{\alpha_1}(\rho-1)\right|^2+\xi|\nabla^{\alpha_1} u|^2\right)\nonumber\\
		&\qquad+\sum_{|\alpha_1|\leq s+2}\nu_*\int_{\mathbb{T}^N}|\nabla^{\alpha_1} \nabla u|^2+\sum_{|\alpha_1|\leq s+2}\eta_*\int_{\mathbb{T}^N} (\nabla\cdot u)|^2	\nonumber\\
		&\leq C\left(K_1+K_2\right)\left(\left\|\dfrac{1}{\varepsilon}\left(\rho-1\right)\right\|_{s+2}^2+\|u\|_{s+2}^2\right)+C\left\|\nabla\phi\right\|_{s+2}^2,
\end{align}
where we has used the fact that $\nu_*\le\nu\le\nu^*,\eta_*\le\eta\le\eta^*$.
Then using (\ref{q.2.3}), (\ref{q.2.5}), (\ref{q.2.6}), (\ref{q.3.11}), (\ref{q.3.37}) and the Gronwall inequality, one deduces
	\begin{align}\label{q.3.39}
		\left\|\dfrac{1}{\varepsilon}(\rho-1)\right\|_{s+2}^2(t)&+\left\|u\right\|_{s+2}^2(t)\nonumber\\
&\leq {\rm exp}\left(C\left(K_1+K_2\right)T_0\right)\left(3\kappa_0^2\varepsilon^2+2\left\|u_0\right\|_{s+2}^2+C\int_0^t\left\|\nabla\phi\right\|_{s+2}^2\right)\leq C
	\end{align}
	for $t\in[0,T_0]$, provided that $T_0<T_2$.

 Furthermore, integrating (\ref{q.3.37}) over $[0, t]\subseteq[0, T_0]$ , we arrive at
	\begin{align}\label{q.3.40}
\nu_*\int_{0}^{t}\|\nabla u\|_{s+2}^2{\rm d}s +\eta_*\int_{0}^{t}\|\nabla\cdot u\|_{s+2}^2{\rm d}s \leq C \quad{\rm for}\ t\in[0,T_0].
	\end{align}

$\mathbf{Cases\ 3}$: $D^{\alpha_1}=\nabla^\beta\partial_t$ for $1\leq|\beta|\leq s+1$, where $s\geq3$.

For the terms $I_4-I_7$ and $I_9-I_{10}$, using integration by parts in $I_5$, we use Lemma \ref{L3.1}, the Cauchy inequality, (\ref{q.3.11}), (\ref{q.3.19}), (\ref{q.3.39}) and a similar calculation as in (\ref{q.3.44}) to get
	\begin{align}
|I_4|&\leq C\left\|\frac{\rho_t}{\varepsilon}\right\|_{s+1}\left(\left\|v_t\right\|_\infty\left\|\frac{\nabla\rho}{\varepsilon}\right\|_{s+1}
+\left\|v_t\right\|_{s+1}\left\|\frac{\nabla\rho}{\varepsilon}\right\|_\infty+\left\|\nabla v\right\|_\infty\left\|\frac{\nabla\rho_t}{\varepsilon}\right\|_{s}+\left\|\nabla v\right\|_{s}\left\|\frac{\nabla\rho_t}{\varepsilon}\right\|_\infty  \nonumber \right.\\
		&\left.\quad+\left\|\frac{\xi_t}{\varepsilon}\right\|_\infty\left\|\nabla\cdot u\right\|_{s+1}+\left\|\frac{\xi_t}{\varepsilon}\right\|_{s+1}\left\|\nabla\cdot u\right\|_\infty
		+\left\|\frac{\nabla\xi}{\varepsilon}\right\|_\infty\left\|\nabla\cdot u_t\right\|_{s}+\left\|\frac{\nabla\xi}{\varepsilon}\right\|_{s}\left\|\nabla\cdot u_t\right\|_\infty\right) \nonumber\\
		&\leq C\left(K_1+K_2\right)\left(\left\|\frac{\rho_t}{\varepsilon}\right\|_{s+1}^2+\left\|u_t\right\|_{s+1}^2\right)+C,\label{qq.3.45}\\
		|I_5|&\leq C\|\nabla u_t\|_{s+1}\left(\left\|\Delta\phi_t\right\|_\infty\left\|\nabla\phi\right\|_{s}+\left\|\Delta\phi_t\right\|_{s}\left\|\nabla\phi\right\|_\infty
+\left\|\Delta\phi\right\|_\infty\left\|\nabla\phi_t\right\|_{s}+\left\|\Delta\phi\right\|_{s}\left\|\nabla\phi_t\right\|_\infty\right)
		\nonumber\\
		&\leq \tau\|\nabla u_t\|_{s+1}^2+C(\tau)\left\|\nabla\phi_t\right\|_{s+1}^2,\label{q.3.41}\\
		|I_6|&\leq C\left\|u_t\right\|_{s+1}\left(\left\|\nabla^\beta\left(v_t\cdot\nabla u\right)+\nabla^\beta\left(v\cdot\nabla u_t\right)-v\nabla^\beta\nabla u_t\right\|\right) \nonumber\\
		&\leq C\|u_t\|_{s+1}\left(\left\|v_t\right\|_\infty\left\|\nabla u\right\|_{s+1}+\left\|v_t\right\|_{s+1}\left\|\nabla u\right\|_\infty+\left\|\nabla v\right\|_\infty\left\|\nabla u_t\right\|_{s}+\left\|\nabla v\right\|_{s}\left\|\nabla u_t\right\|_\infty\right) \nonumber\\
		&\leq\tau\left\|\nabla u_t\right\|_{2}^2+C\left(\tau\right)\left(K_1+K_2\right)\left\|u_t\right\|_{s+1}^2+C,\\\label{q.3.42}
|I_7|&\leq\dfrac{C }{\varepsilon^2} \left\|u_t\right\|_{s+1} \left[\|\xi_t\|_\infty\left\|\nabla\rho\right\|_{s+1}+(1+\|\nabla\xi\|_{s})\|\xi_t\|_{s+1}\left\|\nabla\rho\right\|_\infty
+\left\|\nabla\xi\right\|_\infty\left\|\nabla\rho_t\right\|_{s}
\right.\nonumber\\
		&\quad\left.+\left\|\nabla\xi\right\|_{s}\left\|\nabla\rho_t\right\|_\infty\right]
		 \nonumber\\
		&\leq C\left(K_1^2+K_2^2\right)\left(\left\|\frac{1}{\varepsilon}\rho_t\right\|_{s+1}^2+\left\|u_t\right\|_{s+1}^2\right)+C,\\
		|I_9|&\leq C\left\|\nabla^\beta u_t\right\|\left(\left\|\left(\dfrac{1}{\xi}\right)_t\right\|_\infty\left\|\Delta\phi\nabla\phi\right\|_{s+1}+\left\|\left(\dfrac{1}{\xi}\right)_t\right\|_{s+1}\left\|\Delta\phi\nabla\phi\right\|_\infty\right) \nonumber\\
		&\quad+\left\|\nabla\left(\dfrac{1}{\xi}\right)\right\|_\infty\left\|\left(\Delta\phi\nabla\phi\right)_t\right\|_{s}+\left\|\nabla\left(\dfrac{1}{\xi}\right)\right\|_{s}\left\|(\Delta\phi\nabla\phi)_t\right\|_\infty \nonumber\\
		&\leq C\left\|u_t\right\|_{s+1}[\left\|\xi_t\right\|_\infty\left(\left\|\Delta\phi\right\|_\infty\left\|\nabla\phi\right\|_{s+1}+\left\|\Delta\phi\right\|_{s+1}\left\|\nabla\phi\right\|_\infty\right)  \nonumber\\
		&\quad +\left\|\Delta\phi\nabla\phi\right\|_\infty\left\|\xi_t\right\|_{s+1}\left(1+\left\|\nabla\xi\right\|_{s}\right)  +\left\|\nabla\xi\right\|_\infty\left(\left\|\Delta\phi_t\right\|_\infty\left\|\nabla\phi\right\|_{s}
+\left\|\Delta\phi_t\right\|_{s}\left\|\nabla\phi\right\|_\infty\right.\nonumber\\ &\quad\left.+\left\|\Delta\phi\right\|_\infty\left\|\nabla\phi_t\right\|_{s}
+\left\|\Delta\phi\right\|_{s}\left\|\nabla\phi_t\right\|_\infty\right)+\left\|\nabla\xi\right\|_{s}\left(\left\|\Delta\phi\right\|_\infty\left\|\nabla\phi_t\right\|_\infty
+\left\|\Delta\phi_t\right\|_\infty\left\|\nabla\phi\right\|_\infty\right)]
		 \nonumber\\
&\le C\left(K_1^2+K_2^2\right)\left\|u_t\right\|_{s+1}^2+C\left\|\nabla\phi_t\right\|_{s+1}^2+C\left\|\Delta\phi\right\|_{s+1}^2+C,\\
|I_{10}|
    &\leq C\left\|\nabla\nu(\xi,\psi)\right\|_\infty\|\nabla u_t\|_{s+1}\|  u_t\|_{s+1}+ C\left\|\nabla\eta(\xi,\psi)\right\|_\infty\|  \nabla\cdot u_t\|_{s+1}\|  u_t\|_{s+1}\nonumber\\
    &\leq C\left\|\nabla\nu(\xi,\psi)\right\|_2\|\nabla u_t\|_{s+1}\|u_t\|_{s+1}+ C\left\|\nabla\eta(\xi,\psi)\right\|_2\| \nabla\cdot u_t\|_{s+1}\|u_t\|_{s+1}\nonumber\\
     &\leq C\left[(1+\|\xi\|_\infty^{2})\|\nabla\xi\|_{2}+(1+\|\psi\|_\infty^{2})\|\nabla\psi\|_{2}\right](\|\nabla u_t\|_{s+1}+\| \nabla\cdot u_t\|_{s+1})\|u_t\|_{s+1}\nonumber\\
             &\leq CK_1^\frac12(\|\nabla u_t\|_{s+1}+\| \nabla\cdot u_t\|_{s+1})\|u_t\|_{s+1}\nonumber\\
             &\leq C(\tau)K_1\|u_t\|_{s+1}^2+\tau(\|\nabla u_t\|_{s+1}^2+\| \nabla\cdot u_t\|_{s+1}^2). \label{qq.3.56}
\end{align}
For the term $I_8$, we have
\begin{align}\label{qI8}
I_8
&=\int_{\mathbb{T}^N}\xi \left[\nabla^{\beta}\partial_t\left(\dfrac{\nu(\xi,\psi)}{\xi}\Delta u\right)-
	\dfrac{\nu(\xi,\psi)}{\xi}\Delta \nabla^{\beta} u_t\right] : \nabla^{\beta}  u_t \nonumber\\
&\quad+\int_{\mathbb{T}^N}\xi \left[\nabla^{\beta}\partial_t\left(\dfrac{\eta(\xi,\psi)}{\xi}\nabla(\nabla\cdot u)\right)-
	\dfrac{\eta(\xi,\psi)}{\xi} \nabla^{\beta}\nabla(\nabla\cdot u_t)\right] : \nabla^{\beta}  u_t=I_8^1+I_8^2.
\end{align}
Then using integration by parts and Lemma \ref{L3.1} yield
\begin{align}\label{qI81}
I_8^1
&=\int_{\mathbb{T}^N}\xi \left[\nabla^{\beta}\left(\partial_t\left(\dfrac{\nu(\xi,\psi)}{\xi}\right)\Delta u\right)+\nabla^{\beta}\left(\dfrac{\nu(\xi,\psi)}{\xi}\Delta u_t\right)
-\dfrac{\nu(\xi,\psi)}{\xi}\Delta \nabla^{\beta} u_t\right] : \nabla^{\beta}  u_t \nonumber\\
&=\int_{\mathbb{T}^N}\xi \left\{\nabla^{\beta}\left[\left(\nu_\xi(\xi,\psi)\xi_t+\nu_\psi(\xi,\psi)\psi_t-\dfrac{\nu(\xi,\psi)\xi_t}{\xi^2}\right)\Delta u\right]
\right.\nonumber\\
&\left.\quad
+\nabla^{\beta}\left(\dfrac{\nu(\xi,\psi)}{\xi}\Delta u_t\right)
-\dfrac{\nu(\xi,\psi)}{\xi}\Delta \nabla^{\beta} u_t\right\} : \nabla^{\beta}  u_t \nonumber\\
&=\int_{\mathbb{T}^N}\xi \left\{\nabla^{\beta}\left[\left(\nu_\xi(\xi,\psi)\xi_t-\dfrac{\nu(\xi,\psi)\xi_t}{\xi^2}\right)\Delta u\right]
+\nabla^{\beta}\left(\dfrac{\nu(\xi,\psi)}{\xi}\Delta u_t\right)
-\dfrac{\nu(\xi,\psi)}{\xi}\Delta \nabla^{\beta} u_t\right\} : \nabla^{\beta}  u_t\nonumber\\
&\quad+\int_{\mathbb{T}^N}\xi\nabla^{\beta}\left(\nu_\psi(\xi,\psi)\psi_t\Delta u\right): \nabla^{\beta}  u_t \nonumber\\
%
&=\int_{\mathbb{T}^N}\xi \left\{\nabla^{\beta}\left[\left(\nu_\xi(\xi,\psi)\xi_t-\dfrac{\nu(\xi,\psi)\xi_t}{\xi^2}\right)\Delta u\right]
+\nabla^{\beta}\left(\dfrac{\nu(\xi,\psi)}{\xi}\Delta u_t\right)
-\dfrac{\nu(\xi,\psi)}{\xi}\Delta \nabla^{\beta} u_t\right\} : \nabla^{\beta}  u_t\nonumber\\
&\quad-\int_{\mathbb{T}^N}\xi\nabla^{\beta-1}\left(\nu_\psi(\xi,\psi)\psi_t\Delta u\right) \nabla^{\beta-1} \Delta u_t
-\int_{\mathbb{T}^N}\nabla^{\beta-1}\left(\nu_\psi(\xi,\psi)\psi_t\Delta u\right) \cdot\nabla^{\beta}  u_t\cdot\nabla\xi\nonumber\\
&\le C \left(\left\|\left(\nu_\xi(\xi,\psi)-\dfrac{\nu(\xi,\psi)}{\xi^2}\right)\xi_t\right\|_{s+1}\left\|\Delta u\right\|_{\infty}
+\left\|\left(\nu_\xi(\xi,\psi)-\dfrac{\nu(\xi,\psi)}{\xi^2}\right)\xi_t\right\|_{\infty}\left\|\Delta u\right\|_{s+1} \nonumber\right.\\
		&\left.\quad+\left\|\nabla\left(\dfrac{\nu(\xi,\psi)}{\xi}\right)\right\|_{s}\left\|\Delta u_t\right\|_{\infty}
+\left\|\nabla\left(\dfrac{\nu(\xi,\psi)}{\xi}\right)\right\|_{\infty}\left\|\Delta u_t\right\|_{s}\right)\|u_t\|_{s+1}\nonumber\\
&\quad+ C\left\|\nu_\psi(\xi,\psi)\psi_t\Delta u\right\|_{s}\left\|\Delta u_t\right\|_{s}
+C\left\|\nu_\psi(\xi,\psi)\psi_t\Delta u\right\|_{s}\left\|u_t\right\|_{s+1}\left\|\nabla\xi\right\|_{\infty}\nonumber\\
&\le C \left(\left\|\left(\nu_\xi(\xi,\psi)-\dfrac{\nu(\xi,\psi)}{\xi^2}\right)\xi_t\right\|_{s+1}\left\|\Delta u\right\|_{s+1} +\left\|\nabla\left(\dfrac{\nu(\xi,\psi)}{\xi}\right)\right\|_{s}\left\|\Delta u_t\right\|_{s}\right)\|u_t\|_{s+1}\nonumber\\
&\quad+ C\left\|\nu_\psi(\xi,\psi)\psi_t\Delta u\right\|_{s}\left\|\Delta u_t\right\|_{s}
+C\left\|\nu_\psi(\xi,\psi)\psi_t\Delta u\right\|_{s}\left\|u_t\right\|_{s+1}\left\|\nabla\xi\right\|_{2}\nonumber\\
&\le C \left(\left\|\nu_\xi(\xi,\psi)-\dfrac{\nu(\xi,\psi)}{\xi^2}\right\|_{s+1}\left\|\xi_t\right\|_{s+1}\left\|\Delta u\right\|_{s+1} +\left\|\nabla\left(\dfrac{\nu(\xi,\psi)}{\xi}\right)\right\|_{s}\left\|\Delta u_t\right\|_{s}\right)\|u_t\|_{s+1}\nonumber\\
&\quad+ C\left\|\nu_\psi(\xi,\psi)\right\|_{s}\left\|\psi_t\right\|_{s}\left\|\Delta u\right\|_{s}\left\|\Delta u_t\right\|_{s}
+C\left\|\nu_\psi(\xi,\psi)\right\|_{s}\left\|\psi_t\right\|_{s}\left\|\Delta u\right\|_{s}\left\|u_t\right\|_{s+1}\left\|\nabla\xi\right\|_{2}\nonumber\\
&\le C \left[ \left( \left\|\nu_\xi(\xi,\psi)-\dfrac{\nu(\xi,\psi)}{\xi^2}\right\|+\left\|\nabla\left(\nu_\xi(\xi,\psi)-\dfrac{\nu(\xi,\psi)}{\xi^2}\right)\right\|_{s} \right)  \left\|\xi_t\right\|_{s+1}\left\|\Delta u\right\|_{s+1}\right. \nonumber\\
&\left.\quad+\left\|\nabla\left(\dfrac{\nu(\xi,\psi)}{\xi}\right)\right\|_{s}\left\|\Delta u_t\right\|_{s}\right]\|u_t\|_{s+1}\nonumber\\
&\quad
+C(\left\|\nu_\psi(\xi,\psi)\right\|+\left\|\nabla\nu_\psi(\xi,\psi)\right\|_{s-1})\left\|\psi_t\right\|_{s}\left\|\Delta u\right\|_{s}(\left\|\Delta u_t\right\|_{s}+\left\|u_t\right\|_{s+1}\left\|\nabla\xi\right\|_{2})\nonumber\\
&\le C \left[\left(1+\|\nabla\xi\|_{s}+\|\nabla\psi\|_{s}\right)  \left\|\xi_t\right\|_{s+1}\left\|\Delta u\right\|_{s+1}+\left(\|\nabla\xi\|_{s}+\|\nabla\psi\|_{s}\right)\left\|\Delta u_t\right\|_{s}\right]\|u_t\|_{s+1}\nonumber\\
&\quad+ C\left[1+\|\nabla\xi\|_{s-1}+\|\nabla\psi\|_{s-1}\right]\left\|\psi_t\right\|_{s}\left\|\Delta u\right\|_{s}(\left\|\Delta u_t\right\|_{s}+\left\|u_t\right\|_{s+1}\left\|\nabla\xi\right\|_{2})\nonumber\\
&\le C \left(K_1^\frac12 K_2^\frac12\left\|\nabla u\right\|_{s+2}+K_1^\frac12\left\|\nabla u_t\right\|_{s+1}\right)\|u_t\|_{s+1}
+ CK_1K_2^\frac12(\left\|\nabla u_t\right\|_{s+1}+K_1^\frac12\left\|u_t\right\|_{s+1})\nonumber\\
&\le C(\tau) \left(K_1 K_2+K_1+K_1^3K_2\right)\|u_t\|_{s+1}^2+C\left\|\nabla u\right\|_{s+2}^2+\tau\left\|\nabla u_t\right\|_{s+1}^2+C(\tau)K_1^2K_2.
\end{align}
Similar to $I_8^1$, we also have
\begin{align}\label{qI82}
I_8^2
&\le C(\tau) \left(K_1 K_2+K_1+K_1^3K_2\right)\|u_t\|_{s+1}^2+C\left\|\nabla\cdot u\right\|_{s+2}^2+\tau\left\|\nabla\cdot u_t\right\|_{s+1}^2+C(\tau)K_1^2K_2.
\end{align}
Inserting (\ref{qI81})-(\ref{qI82}) into (\ref{qI8}) gives
\begin{align}\label{qI80}
I_8
&\le C(\tau) \left(K_1^6+ K_2^2\right)\|u_t\|_{s+1}^2+C(\left\|\nabla u\right\|_{s+2}^2+\left\|\nabla\cdot u\right\|_{s+2}^2)
\nonumber\\
&\quad+\tau(\left\|\nabla u_t\right\|_{s+1}^2+\left\|\nabla\cdot u_t\right\|_{s+1}^2)+C(\tau)K_1^2K_2.
\end{align}
 Now, putting the estimates (\ref{q.3.21})-(\ref{q.3.23}), (\ref{q.3.25})-(\ref{q.3.30}), (\ref{qq.3.45})-(\ref{qq.3.56}), (\ref{qI80}) together, summing over $\beta$, and then choosing $\tau$ small enough, one arrival at
	\begin{align*}
		&\dfrac{\rm d}{{\rm d}t}\sum_{|\beta|\leq s+1}\int_{\mathbb{T}^N}\left(\dfrac{P'\left(\xi\right)}{\xi}\left|\dfrac{1}{\varepsilon} \nabla^{\beta}\rho_t\right|^2+\xi|\nabla^{\beta} u_t|^2\right)
		\nonumber\\\nonumber
		&\quad+\sum_{|\beta|\leq s+1}\nu_*\int_{\mathbb{T}^N}|\nabla^{\beta} \nabla u_t|^2+\sum_{|\beta|\leq s+1}\eta_*\int_{\mathbb{T}^N}|\nabla^{\beta} (\nabla\cdot u_t)|^2
		\nonumber\\\nonumber
		&\leq C\left(K_1^6+K_2^2\right)\left(\left\|\dfrac{1}{\varepsilon}\rho_t\right\|_{s+1}^2+\|u_t\|_{s+1}^2\right)\nonumber\\
&\quad+C\left(\|\nabla \phi_t\|_{s+1}^2+\|\Delta\phi\|_{s+1}^2+\|\nabla u\|_{s+2}^2+\|\nabla\cdot u\|_{s+2}^2\right)+C(\tau)K_1^2K_2+C.
	\end{align*}
 Recalling the constraints of the initial data and (\ref{q.2.1}), we get
	\begin{align}\label{q.3.51}
		&\left\|\frac{1}{\varepsilon}\rho_t\left(x,0\right)\right\|_{s+1}^2+\|u_t\left(x,0\right)\|_{s+1}^2
		\nonumber\\
		&\leq C\left(\left\|\frac{1}{\varepsilon}\left(\left(u_0+\overline u_0^\varepsilon\right)\cdot\nabla\right)\overline\rho_0^\varepsilon\right\|_{s+1}^2
		+\left\|\frac{1}{\varepsilon}\left(\overline\rho_0^\varepsilon+1\right)\nabla\cdot\overline u_0^\varepsilon\right\|_{s+1}^2+\left\|\left(\overline\rho_0^\varepsilon+1\right)^{-1}\frac{1}{\varepsilon^2}\nabla\overline\rho_0^\varepsilon\right\|_{s+1}^2
		\right.\nonumber\\
		&\left.\quad+\left\|\left(\left(u_0+\overline u_0^\varepsilon\right)\cdot\nabla\right)\left(u_0+\overline u_0^\varepsilon\right)\right\|_{s+1}^2+\left\|\left(\overline\rho_0^\varepsilon+1\right)^{-1}  \nu(\overline\rho_0^\varepsilon+1, \phi_0+\overline \phi_0^\varepsilon)  \Delta\left(u_0+\overline u_0^\varepsilon\right)\right\|_{s+1}^2
\right.\nonumber\\
		&\left.\quad
+\left\|\left(\overline\rho_0^\varepsilon+1\right)^{-1}  \eta(\overline\rho_0^\varepsilon+1, \phi_0+\overline \phi_0^\varepsilon) \nabla\left(\nabla\cdot\overline u_0^\varepsilon\right)\right\|_{s+1}^2
\right.\nonumber\\
		&\left.\quad+\left\|\left(\overline\rho_0^\varepsilon+1\right)^{-1}\left(\Delta\left(\phi_0+\overline\phi_0^\varepsilon\right)
\nabla\left(\phi_0+\overline\phi_0^\varepsilon\right)\right)\right\|_{s+1}^2\right)\leq C,
	\end{align}
where one has used the fact by Lemma \ref{L3.1} that
\begin{align*}
&\|\nu(\overline\rho_0^\varepsilon+1, \phi_0+\overline \phi_0^\varepsilon)\|_{s+1}
\le \|\nu(\overline\rho_0^\varepsilon+1, \phi_0+\overline \phi_0^\varepsilon)\|
+\|\nabla\nu(\overline\rho_0^\varepsilon+1, \phi_0+\overline \phi_0^\varepsilon)\|_{s}\\
&\le \|\nu(\overline\rho_0^\varepsilon+1, \phi_0+\overline \phi_0^\varepsilon)\|
+(1+\|\overline\rho_0^\varepsilon+1\|_\infty^s)\|\nabla\overline\rho_0^\varepsilon\|_{s}
+(1+\|\phi_0+\overline \phi_0^\varepsilon\|_\infty^s)\|\nabla( \phi_0+\overline \phi_0^\varepsilon)\|_{s}\le C.
\end{align*}
Then using the Gronwall inequality, (\ref{q.3.11}), (\ref{q.3.19}) and (\ref{q.3.40}) yield
	\begin{align}\label{q.3.52}
		\left\|\frac{1}{\varepsilon}\rho_t\right\|_{s+1}^2\left(t\right)+\|u_t\|_{s+1}^2\left(t\right)+\nu_*\int_{0}^{t}\|\nabla u_t\|_{s+1}^2+\eta_*\int_{0}^{t}\|\nabla\cdot u_t\|_{s+1}^2\leq C
	\end{align}
	for $t\in [0, T_0]$, provided that $T_0$ is small enough such that $T_0\leq T_3:={\rm min}\left\{T_2, \left(K_1^6+K_2^2\right)^{-1}\right\}$.

It remains to show the first inequality of (\ref{q.3.1}). It suffices to show  $\left\|\frac{1}{\varepsilon}\left(\rho-1\right)\right\|_s+\|u-u_0\|_s+\|\phi-\phi_0\|_s\leq c_0^{-1}\kappa$ by the Sobolev inequality, where $c_0$ is the Sobolev constant. Let\ $\overline\rho=\rho-1,\ \overline u=u-u_0\ {\rm and}\ \overline\phi=\phi-\phi_0$.

 When $\phi=\overline\phi+\phi_0$ and $|\alpha_1| \leq s$, it follows from (\ref{q.3.4}) that
 \begin{align*}
\bar{\phi}_t+\dfrac1{\xi^2}\Delta^2\bar{\phi}=- v\cdot\nabla\psi-\dfrac2{\xi}\nabla\left(\dfrac1{\xi}\right)\cdot\nabla\Delta\psi
-\dfrac1{\xi}\Delta\left(\dfrac1{\xi}\right)\Delta\psi+\dfrac1{\xi}\Delta(\psi^3-\psi).
\end{align*}
Then operating the above equality by $\nabla^{\alpha_1}$, multiplying the result by $\nabla^{\alpha_1}\overline\phi$, by Lemma \ref{L3.1}, we have
	\begin{align}\label{q.3.53}
			&\dfrac12\dfrac{\rm d}{{\rm d}t}\int_{\mathbb{T}^N}|\nabla^{\alpha_1}\overline\phi|^2+\int_{\mathbb{T}^N}\dfrac1{\xi^2}|\nabla^{\alpha_1}\Delta\overline\phi|^2
\nonumber\\
&=-\int_{\mathbb{T}^N}\left[\nabla^{\alpha_1}\left(\dfrac1{\xi^2}\Delta^2\overline\phi\right)-\dfrac1{\xi^2}\nabla^{\alpha_1}\Delta^2\overline\phi\right]
\nabla^{\alpha_1}\overline\phi
-\int_{\mathbb{T}^N}\nabla^{\alpha_1}(v\cdot\nabla\psi)\nabla^{\alpha_1}\overline\phi
\nonumber\\
&\quad-2\int_{\mathbb{T}^N}\nabla^{\alpha_1}\left(\dfrac1{\xi}\nabla\left(\dfrac1{\xi}\right)\cdot\nabla\Delta\psi\right)\nabla^{\alpha_1}\overline\phi
-\int_{\mathbb{T}^N}\nabla^{\alpha_1}\left(\dfrac1{\xi}\Delta\left(\dfrac1{\xi}\right)\Delta\psi\right)\nabla^{\alpha_1}\overline\phi
\nonumber\\
&\quad+\int_{\mathbb{T}^N}\nabla^{\alpha_1}\left(\dfrac1{\xi}\Delta(\psi^3-\psi)\right)\nabla^{\alpha_1}\overline\phi
-\int_{\mathbb{T}^N}\Delta\left(\dfrac1{\xi^2}\right)\nabla^{\alpha_1}\Delta\overline\phi\nabla^{\alpha_1}\overline\phi
\nonumber\\
&\quad-2\int_{\mathbb{T}^N}\nabla\left(\dfrac1{\xi^2}\right)\cdot\nabla^{\alpha_1}\nabla\overline\phi\nabla^{\alpha_1}\Delta\overline\phi
-\int_{\mathbb{T}^N}\nabla^{\alpha_1}\left(\dfrac1{\xi^2}\Delta^2\phi_0\right)\nabla^{\alpha_1}\overline\phi\nonumber\\
&\leq C\left(\left\|\nabla\left(\dfrac1{\xi^2}\right)\right\|_\infty\|\Delta^2\overline\phi\|_{s-1}
+\|\Delta^2\overline\phi\|_\infty\left\|\nabla\left(\dfrac1{\xi^2}\right)\right\|_{s-1}
+\|v\|_\infty\|\nabla\psi\|_s+\|\nabla\psi\|_\infty\|v\|_s\right.\nonumber\\
&\left.\quad+\left\|\dfrac1{\xi}\right\|_s\left\|\nabla\left(\dfrac1{\xi}\right)\right\|_s\left\|\nabla\Delta\psi\right\|_s
+\left\|\dfrac1{\xi}\right\|_s\left\|\Delta\left(\dfrac1{\xi}\right)\right\|_s\left\|\Delta\psi\right\|_s
\right.\nonumber\\
        &\left.\quad+\left\|\dfrac1{\xi}\right\|_\infty\|\Delta\left(\psi^3-\psi\right)\|_s
        +\|\Delta\left(\psi^3-\psi\right)\|_\infty\left\|\dfrac1{\xi}\right\|_s
        +\left\|\Delta\left(\dfrac1{\xi^2}\right)\right\|_\infty\|\nabla^{\alpha_1}\Delta\overline\phi\|
        \right.\nonumber\\
        &\left.\quad+\left\|\dfrac1{\xi^2}\right\|_\infty\|\Delta^2\phi_0\|_s+\|\Delta^2\phi_0\|_\infty\left\|\dfrac1{\xi^2}\right\|_s\right)
        \|\nabla^{\alpha_1}\overline\phi\|
         +\left\|\nabla\left(\dfrac1{\xi^2}\right)\right\|_\infty\|\nabla^{\alpha_1}\nabla\overline\phi\|\|\nabla^{\alpha_1}\Delta\overline\phi\|
         \nonumber\\
&\leq C\left(\left\|\nabla\xi\right\|_{s-1}\|\Delta^2\overline\phi\|_{s-1}+\left\|v\right\|_s\left\|\nabla\psi\right\|_s
+(1+\left\|\nabla\xi\right\|_{s-1})\left\|\nabla\xi\right\|_s\left\|\nabla\Delta\psi\right\|_{s}
\right.\nonumber\\
&\left.\quad+(1+\left\|\nabla\xi\right\|_{s-1})\left\|\nabla\xi\right\|_{s+1}\left\|\Delta\psi\right\|_{s}
+(1+\left\|\nabla\xi\right\|_{s-1})\left\|\nabla\psi\right\|_{s+1}+\left\|\nabla\xi\right\|_{3}\|\Delta\overline\phi\|_s
\right.\nonumber\\
&\left.\quad
+(1+\left\|\nabla\xi\right\|_{s-1})\|\Delta^2\phi_0\|_s\right)\|\overline\phi\|_s
+C\left\|\nabla\xi\right\|_{2}\|\nabla\overline\phi\|_s\|\Delta\overline\phi\|_s
\nonumber\\
&\leq C\left(\varepsilon K_1^\frac12\|\Delta^2\overline\phi\|_{s-1}
+\varepsilon K_1\left\|\nabla\Delta\psi\right\|_{s}
+ K_1^\frac32
+K_1^\frac12\|\Delta\overline\phi\|_s
+K_1^\frac12\|\Delta^2\phi_0\|_s\right)\|\overline\phi\|_s\nonumber\\
&\quad+C\varepsilon K_1^\frac12\|\nabla\overline\phi\|_s\|\Delta\overline\phi\|_s
\nonumber\\
&\leq  C\varepsilon\|\Delta^2(\phi-\phi_0)\|_{s-1}^2 +CK_1^3\|\overline\phi\|_s^2+C\varepsilon^2\left\|\nabla\Delta\psi\right\|_s^2+\tau\|\Delta\overline\phi\|_s^2+C\varepsilon^2 K_1\|\nabla\overline\phi\|_s^2
+C\nonumber\\
&\leq  C\varepsilon\|\Delta^2\phi\|_{s-1}^2+C\varepsilon\|\Delta^2\phi_0\|_{s-1}^2+CK_1^3\|\overline\phi\|_s^2+\left(C\varepsilon^2 K_1+\tau\right)\|\Delta\overline\phi\|_s^2
+C\varepsilon^2\left\|\nabla\Delta\psi\right\|_s^2+C.
		\end{align}
Taking $\varepsilon, \tau$ and  $T_0(< T_3)$ sufficiently small, and then summing over $\alpha_1$,  by (\ref{q.2.5}), (\ref{q.3.11}) and the Gronwall inequality we deduce
	\begin{align}\label{q.3.55}
		\|\overline\phi\|_s^2
&\leq \exp\left(CK^3T_0\right)\left(\kappa_0^2\varepsilon^2+C\varepsilon\int_0^t\|\Delta^2\phi\|_{s-1}^2+C\varepsilon T_0\|\Delta^2\phi_0\|_{s-1}^2
+C\varepsilon^2\int_0^t\left\|\nabla\Delta\psi\right\|_s^2+CT_0\right)\nonumber\\
		&\leq \exp\left(CK^3T_0\right)\left(\kappa_0^2\varepsilon^2+C\varepsilon+CT_0\right) c_0^{-1}\kappa,
	\end{align}
where we have used the fact by (\ref{q.2.6}) that $\|\overline\phi(x,0)\|_{s}^2\le \kappa_0^2\varepsilon^2$.

	Similarly as in the proof of (\ref{q.3.37}) and (\ref{q.3.53}), we deduce
	\begin{align*}
		&\dfrac{\rm d}{{\rm d}t}\sum_{|\alpha_1|\leq{s}}\int_{\mathbb{T}^N}\left(\dfrac{P'\left(\xi\right)}{\xi}\left|\dfrac{1}{\varepsilon} \nabla^{\alpha_1}\overline\rho\right|^2+\xi| \nabla^{\alpha_1}\overline u|^2\right)\nonumber\\
		&\quad+\sum_{|\alpha_1|\leq{s}}\nu_*\int_{\mathbb{T}^N}| \nabla^{\alpha_1} \nabla \overline u|^2 +\sum_{|\alpha_1|\leq{s}}\eta_*\int_{\mathbb{T}^N}| \nabla^{\alpha_1} (\nabla\cdot \overline u)|^2  \nonumber\\
		&\leq C\left(  K_1 +K_2\right)\left(\left\|\frac{1}{\varepsilon}\overline\rho\right\|_s^2+\left\|\overline u\right\|_s^2\right)+C\left\|\left(v\cdot\nabla\right)u_0\right\|_s^2+C  +\left\|\frac{\nu(\xi,\psi) }{\xi}\Delta u_0\right\|_s^2  \nonumber\\
&\leq C\left(  K_1 +K_2\right)\left(\left\|\frac{1}{\varepsilon}\overline\rho\right\|_s^2+\left\|\overline u\right\|_s^2\right)+C\left\|\left(v\cdot\nabla\right)u_0\right\|_s^2 \nonumber\\
&\quad+C\left(\left\|\frac{\nu(\xi,\psi) }{\xi}\right\|_\infty\left\|\Delta u_0\right\|_s +\left\|\Delta u_0\right\|_\infty \left\|\frac{\nu(\xi,\psi) }{\xi}\right\|_s \right)^2\nonumber\\
		&\leq C\left(K_1+K_2\right)\left(\left\|\frac{1}{\varepsilon}\overline\rho\right\|_s^2+\left\|\overline u\right\|_s^2\right)+CK_1
+C\left\|\frac{\nu(\xi,\psi) }{\xi}\right\|_{s}^2\|\Delta u_0\|_{s}^2\nonumber\\
&\leq C\left(K_1+K_2\right)\left(\left\|\frac{1}{\varepsilon}\overline\rho\right\|_s^2+\left\|\overline u\right\|_s^2\right)+CK_1\nonumber\\
&\quad+C\left[(1+\|\xi\|_\infty^{s-1})\|\nabla\xi\|_{s-1}+(1+\|\psi\|_\infty^{s-1})\|\nabla\psi\|_{s-1}\right]^2\nonumber\\
&\leq C\left(K_1+K_2\right)\left(\left\|\frac{1}{\varepsilon}\overline\rho\right\|_s^2+\left\|\overline u\right\|_s^2\right)+CK_1,
	\end{align*}
where we have used (\ref{q.2.5}).
	By the Gronwall inequality, (\ref{q.2.3}) and (\ref{q.2.6}), we obtain
	\begin{align}\label{q.3.57}
		\left\|\frac{1}{\varepsilon}\overline\rho\right\|_s^2+\|\overline u\|_s^2  \leq C\exp\left(C\left(K_1+K_2\right)T_0\right)\left(\kappa^2\varepsilon^2+K_1T_0\right)
		<c_0^{-1}\kappa,
	\end{align}
	where we have chosen $T_0\left(<T_3\right)$ and $\varepsilon$ small enough such that (\ref{q.3.55}) and (\ref{q.3.57}) hold.
	Then the Sobolev inequality indicates that
		\begin{align}\label{q.3.5700}
	\left|\frac{1}{\varepsilon}\left(\rho-1\right)\right|+\left|u-u_0\right|+\left|\phi-\phi_0\right|<\kappa.
	\end{align}
The proof of Lemma \ref{L3.2} is completed.
\hfill$\Box$

Now we plan to show that $\Lambda$ is a contractive map. In the proof of this part and the following lemma in this subsection, denote by $C$ the constant
depending on the initial data, the domain, $N$, $s$, and the viscosity coefficients $\nu$, $\eta$, but independent of $\varepsilon$.
\begin{lemma} \label{L3.3}
	Under the assumptions in Theorem \ref{th1}, the map $\Lambda:V\longrightarrow U$ is a contraction in the sense that
	%
	\begin{align}\label{aq.3.58}
		&\sup_{0\leq t\leq T_0}	\left(\left\|\dfrac{1}{\varepsilon}\left(\rho-\hat{\rho}\right)\right\|_1^2+\left\|u-\hat{u}\right\|_1^2
+\left\|\phi-\hat{\phi}\right\|_1^2\right)+\int_{0}^{t}\left(\|u-\hat{u}\|_2^2+\|\phi-\hat{\phi}\|_3^2\right)  \nonumber\\
		&\quad\leq\gamma_1\sup_{0\leq t\leq T_0}\left(\left\|\dfrac{1}{\varepsilon}\left(\xi-\hat{\xi}\right)\right\|_1^2+\left\|v-\hat{v}\right\|_1^2+\left\|\psi-\hat{\psi}\right\|_1^2\right)
+\gamma_2\int_0^t\left\|\psi-\hat{\psi}\right\|_3^2
	\end{align}
	for some $0<\gamma_1<1$ and $0<\gamma_2<1$, provided that $T_0$ is small enough.
	\end{lemma}
\noindent{\it\bfseries Proof.}\quad
 Let $U=\Lambda\left(V\right)$ and $\hat{U} =\Lambda\left(\hat{V} \right)$, where $V$, $\hat{V}\in B_{T_0}^{\varepsilon}\left(U_0\right)$. Then by the definition, we arrive at
		\begin{align}\label{q.3.58}
		\begin{cases}
			\left(\rho-\hat{\rho}\right)_t+\left(v\cdot\nabla\right)\left(\rho-\hat{\rho}\right)
+\left(\left(v-\hat{v}\right)\cdot\nabla\right)\hat{\rho}+\xi\nabla\cdot\left(u-\hat{u}\right)+\left(\xi-\hat{\xi}\right)\nabla\cdot\hat{u}=0,\\
			\left(u-\hat{u}\right)_t+\left(v\cdot\nabla\right)\left(u-\hat{u}\right)+\left(\left(v-\hat{v}\right)\cdot\nabla\right)\hat{u}
+\dfrac{1}{\varepsilon^2}\dfrac{P'\left(\xi\right)}{\xi}\nabla\left(\rho-\hat{\rho}\right)\\
\quad+\dfrac{1}{\varepsilon^2}\left(\dfrac{P'\left(\xi\right)}{\xi}-\dfrac{P'\left(\hat{\xi}\right)}{\hat{\xi}}\right)\nabla\hat{\rho}
			=\dfrac{\nu(\xi,\psi)}{\xi}\Delta\left(u-\hat{u}\right)
+\left(\dfrac{\nu(\xi,\psi)}{\xi}-\dfrac{\nu(\hat{\xi},\hat{\psi})}{\hat{\xi}}\right)\Delta\hat{u} \\
\quad
+\dfrac{\eta(\xi,\psi)}{\xi}\nabla\left(\nabla\cdot\left(u-\hat{u}\right)\right)
+\left(\dfrac{\eta(\xi,\psi)}{\xi}-\dfrac{\eta(\hat{\xi},\hat{\psi})}{\hat{\xi}}\right)\nabla\left(\nabla\cdot\hat{u}\right) \\
\quad
			-\dfrac{1}{\xi}\left[\Delta\left(\phi-\hat{\phi}\right)\nabla\phi+\Delta\hat{\phi}\nabla\left(\phi-\hat{\phi}\right)\right]
-\left(\dfrac{1}{\xi}-\dfrac{1}{\hat{\xi}}\right)\Delta\hat\phi\nabla\hat\phi,\\
			\left(\phi-\hat{\phi}\right)_t+\dfrac{1}{\xi^2}\Delta^2\left(\phi-\hat{\phi}\right)+\left(\dfrac{1}{\xi^2}
-\dfrac{1}{\hat{\xi}^2}\right)\Delta^2\hat{\phi}=-\left(v\cdot\nabla\right)\left(\psi-\hat{\psi}\right)
-\left(\left(v-\hat{v}\right)\cdot\nabla\right)\hat\psi\\
\quad-\dfrac{2}{\xi}\nabla\left(\dfrac{1}{\xi}\right)\cdot\nabla\Delta\left(\psi-\hat{\psi}\right)
-\dfrac{2}{\xi}\nabla\left(\dfrac{1}{\xi}-\dfrac{1}{\hat{\xi}}\right)\cdot\nabla\Delta\hat{\psi}
-2\left(\dfrac{1}{\xi}-\dfrac{1}{\hat{\xi}}\right)\nabla\dfrac{1}{\hat{\xi}}\cdot\nabla\Delta\hat{\psi}
\\
\quad-\dfrac{1}{\xi}\Delta\left(\dfrac{1}{\xi}\right)\Delta\left(\psi-\hat{\psi}\right)
-\dfrac{1}{\xi}\Delta\left(\dfrac{1}{\xi}-\dfrac{1}{\hat{\xi}}\right)\Delta\hat{\psi}
-\left(\dfrac{1}{\xi}-\dfrac{1}{\hat{\xi}}\right)\Delta\dfrac{1}{\hat{\xi}}\Delta\hat{\psi}
\\
\quad+\dfrac{1}{\xi}\left[\Delta\left(\psi^3-\psi\right)-\Delta\left(\hat{\psi}^3-\hat{\psi}\right)\right]
+\left(\dfrac{1}{\xi}-\dfrac{1}{\hat{\xi}}\right)\Delta\left(\hat{\psi}^3-\hat\psi\right).
		\end{cases}
		\end{align}

 Firstly, multiplying $(\ref{q.3.58})_3$ by $\left(\phi-\hat{\phi}\right)$ and $\Delta\left(\phi-\hat{\phi}\right)$ respectively and then integrating by parts, we obtain
	\begin{align}\label{q.3.59}
		&\dfrac{1}{2}\dfrac{\rm d}{{\rm d}t}\int_{\mathbb{T}^N}|\phi-\hat{\phi}|^2
+\int_{\mathbb{T}^N}\dfrac{1}{\xi^2}|\Delta\left(\phi-\hat{\phi}\right)|^2
\nonumber\\
		&=-\int_{\mathbb{T}^N}\left(\dfrac{1}{\xi^2}-\dfrac{1}{\hat{\xi}^2}\right)\Delta^2\hat{\phi}\left(\phi-\hat{\phi}\right)
-\int_{\mathbb{T}^N}\left(v\cdot\nabla\right)\left(\psi-\hat{\psi}\right)\left(\phi-\hat{\phi}\right)
\nonumber\\
&\quad-\int_{\mathbb{T}^N}\left(\left(v-\hat{v}\right)\cdot\nabla\right)\hat\psi\left(\phi-\hat{\phi}\right)
-\int_{\mathbb{T}^N}\dfrac{2}{\xi}\nabla\left(\dfrac{1}{\xi}\right)\cdot\nabla\Delta\left(\psi-\hat{\psi}\right)\left(\phi-\hat{\phi}\right)
\nonumber\\
&\quad-\int_{\mathbb{T}^N}\dfrac{2}{\xi}\nabla\left(\dfrac{1}{\xi}-\dfrac{1}{\hat{\xi}}\right)\cdot\nabla\Delta\hat{\psi}\left(\phi-\hat{\phi}\right)
-\int_{\mathbb{T}^N}2\left(\dfrac{1}{\xi}-\dfrac{1}{\hat{\xi}}\right)\nabla\left(\dfrac{1}{\hat{\xi}}\right)\cdot\nabla\Delta\hat{\psi}\left(\phi-\hat{\phi}\right)
\nonumber\\
&\quad-\int_{\mathbb{T}^N}\dfrac{1}{\xi}\Delta\left(\dfrac{1}{\xi}\right)\Delta\left(\psi-\hat{\psi}\right)\left(\phi-\hat{\phi}\right)
-\int_{\mathbb{T}^N}\dfrac{1}{\xi}\Delta\left(\dfrac{1}{\xi}-\dfrac{1}{\hat{\xi}}\right)\Delta\hat{\psi}\left(\phi-\hat{\phi}\right)
\nonumber\\
&\quad-\int_{\mathbb{T}^N}\left(\dfrac{1}{\xi}-\dfrac{1}{\hat{\xi}}\right)\Delta\left(\dfrac{1}{\hat{\xi}}\right)\Delta\hat{\psi}\left(\phi-\hat{\phi}\right)
+\int_{\mathbb{T}^N}\dfrac{1}{\xi}\left[\Delta\left(\psi^3-\psi\right)-\Delta\left(\hat{\psi}^3-\hat{\psi}\right)\right]\left(\phi-\hat{\phi}\right)
\nonumber\\
&\quad+\int_{\mathbb{T}^N}\left(\dfrac{1}{\xi}-\dfrac{1}{\hat{\xi}}\right)\Delta\left(\hat{\psi}^3-\hat\psi\right)\left(\phi-\hat{\phi}\right)
-2\int_{\mathbb{T}^N}\nabla\left(\dfrac{1}{\xi^2}\right)\cdot\nabla\left(\phi-\hat{\phi}\right)\Delta\left(\phi-\hat{\phi}\right)
\nonumber\\
&\quad-\int_{\mathbb{T}^N}\Delta\left(\dfrac{1}{\xi^2}\right)\Delta\left(\phi-\hat{\phi}\right)\left(\phi-\hat{\phi}\right),
	\end{align}
and
\begin{align}\label{q.3.60}
		&\dfrac{1}{2}\dfrac{\rm d}{{\rm d}t}\int_{\mathbb{T}^N}\left|\nabla\left(\phi-\hat{\phi}\right)\right|^2+\int_{\mathbb{T}^N}\dfrac{1}{\xi^2}\left|\nabla\Delta\left(\phi-\hat{\phi}\right)\right|^2 \nonumber\\
		&=\int_{\mathbb{T}^N}\left(\dfrac{1}{\xi^2}-\dfrac{1}{\hat{\xi}^2}\right)\Delta^2\hat{\phi}\Delta\left(\phi-\hat{\phi}\right)
+\int_{\mathbb{T}^N}\left(v\cdot\nabla\right)\left(\psi-\hat{\psi}\right)\Delta\left(\phi-\hat{\phi}\right)
\nonumber\\
&\quad+\int_{\mathbb{T}^N}\left(\left(v-\hat{v}\right)\cdot\nabla\right)\hat\psi\Delta\left(\phi-\hat{\phi}\right)
+\int_{\mathbb{T}^N}\dfrac{2}{\xi}\nabla\left(\dfrac{1}{\xi}\right)\cdot\nabla\Delta\left(\psi-\hat{\psi}\right)\Delta\left(\phi-\hat{\phi}\right)
\nonumber\\
&\quad+\int_{\mathbb{T}^N}\dfrac{2}{\xi}\nabla\left(\dfrac{1}{\xi}-\dfrac{1}{\hat{\xi}}\right)\cdot\nabla\Delta\hat{\psi}\Delta\left(\phi-\hat{\phi}\right)
+\int_{\mathbb{T}^N}2\left(\dfrac{1}{\xi}-\dfrac{1}{\hat{\xi}}\right)\nabla\left(\dfrac{1}{\hat{\xi}}\right)\cdot\nabla\Delta\hat{\psi}\Delta\left(\phi-\hat{\phi}\right)
\nonumber\\
&\quad+\int_{\mathbb{T}^N}\dfrac{1}{\xi}\Delta\left(\dfrac{1}{\xi}\right)\Delta\left(\psi-\hat{\psi}\right)\Delta\left(\phi-\hat{\phi}\right)
+\int_{\mathbb{T}^N}\dfrac{1}{\xi}\Delta\left(\dfrac{1}{\xi}-\dfrac{1}{\hat{\xi}}\right)\Delta\hat{\psi}\Delta\left(\phi-\hat{\phi}\right)
\nonumber\\
&\quad+\int_{\mathbb{T}^N}\left(\dfrac{1}{\xi}-\dfrac{1}{\hat{\xi}}\right)\Delta\left(\dfrac{1}{\hat{\xi}}\right)\Delta\hat{\psi}\Delta\left(\phi-\hat{\phi}\right)
-\int_{\mathbb{T}^N}\dfrac{1}{\xi}\left[\Delta\left(\psi^3-\psi\right)-\Delta\left(\hat{\psi}^3-\hat{\psi}\right)\right]\Delta\left(\phi-\hat{\phi}\right)
\nonumber\\
&\quad-\int_{\mathbb{T}^N}\left(\dfrac{1}{\xi}-\dfrac{1}{\hat{\xi}}\right)\Delta\left(\hat{\psi}^3-\hat\psi\right)\Delta\left(\phi-\hat{\phi}\right)
-\int_{\mathbb{T}^N}\nabla\left(\dfrac{1}{\xi^2}\right)\cdot\nabla\Delta\left(\phi-\hat{\phi}\right)\Delta\left(\phi-\hat{\phi}\right).
	\end{align}
Secondly, applying the operator $\nabla^{\alpha}(|\alpha|\le 1)$ to $(\ref{q.3.58})_1$ and $(\ref{q.3.58})_2$ respectively and then by $\dfrac{1}{\varepsilon^2}\dfrac{P'\left(\xi\right)}{\xi}\nabla^{\alpha}\left(\rho-\hat{\rho}\right)$ and $\xi\nabla^{\alpha}\left(u-\hat{u}\right)$ respectively, one obtains
	\begin{align}\label{q.3.61}
&\dfrac{1}{2}\dfrac{\rm d}{{\rm d}t}\int_{\mathbb{T}^N}\dfrac{P'\left(\xi\right)}{\xi}\left|\dfrac{1}{\varepsilon}\nabla^{\alpha}\left(\rho-\hat{\rho}\right)\right|^2
-\dfrac{1}{2}\int_{\mathbb{T}^N}\left|\dfrac{1}{\varepsilon}\nabla^{\alpha}(\rho-\hat{\rho})\right|^2\left(\dfrac{P'\left(\xi\right)}{\xi}\right)_t
\nonumber\\
&=-\dfrac{1}{\varepsilon^2}\int_{\mathbb{T}^N}\dfrac{P'\left(\xi\right)}{\xi}\nabla^{\alpha}\left(\rho-\hat{\rho}\right)\nabla^{\alpha}\left[v\cdot\nabla\left(\rho-\hat{\rho}\right)\right] 	
-\dfrac{1}{\varepsilon^2}\int_{\mathbb{T}^N}\dfrac{P'\left(\xi\right)}{\xi}\nabla^{\alpha}\left(\rho-\hat{\rho}\right)\nabla^{\alpha}\left[\left(v-\hat{v}\right)
\cdot\nabla\hat{\rho}\right]\nonumber\\
 &\quad-\dfrac{1}{\varepsilon^2}\int_{\mathbb{T}^N}\dfrac{P'\left(\xi\right)}{\xi}\nabla^{\alpha}\left(\rho-\hat{\rho}\right)\nabla^{\alpha}\left[\xi\nabla\cdot\left(u-\hat{u}\right)\right]	-\dfrac{1}{\varepsilon^2}\int_{\mathbb{T}^N}\dfrac{P'\left(\xi\right)}{\xi}\nabla^{\alpha}\left(\rho-\hat{\rho}\right)\nabla^{\alpha}\left[\left(\xi-\hat{\xi}\right)\nabla\cdot\hat{u}\right],
	\end{align}
	and
	\begin{align}\label{q.3.62}
		&\dfrac{1}{2}\dfrac{\rm d}{{\rm d}t}\int_{\mathbb{T}^N}\xi|\nabla^{\alpha}(u-\hat{u})|^2
+\int_{\mathbb{T}^N}\nu(\xi,\psi)|\nabla^{\alpha}\nabla\left(u-\hat{u}\right)|^2+\int_{\mathbb{T}^N}\eta(\xi,\psi)|\nabla^{\alpha}\nabla\cdot\left(u-\hat{u}\right)|^2 \nonumber\\
		&=\int_{\mathbb{T}^N}\xi_t\dfrac{|\nabla^{\alpha}(u-\hat{u})|^2}{2}
-\int_{\mathbb{T}^N}\xi\nabla^{\alpha}\left[v\cdot\nabla\left(u-\hat{u}\right)\right]\cdot\nabla^{\alpha}\left(u-\hat{u}\right)
\nonumber\\
 &\quad-\int_{\mathbb{T}^N}\xi\nabla^{\alpha}\left[\left(v-\hat{v}\right)\cdot\nabla\hat{u}\right]\cdot\nabla^{\alpha}\left(u-\hat{u}\right)	
-\dfrac1{\varepsilon^2}\int_{\mathbb{T}^N}\xi\nabla^{\alpha}
 \left(u-\hat{u}\right)\cdot\nabla^{\alpha}\left(\dfrac{P'\left(\xi\right)}{\xi}\nabla\left(\rho-\hat{\rho}\right)\right)
\nonumber\\
 &\quad
-\dfrac{1}{\varepsilon^2}\int_{\mathbb{T}^N}\xi\nabla^{\alpha}\left(u-\hat{u}\right)\cdot\nabla^{\alpha}\left[\left(\dfrac{P'\left(\xi\right)}{\xi}
-\dfrac{P'\left(\hat{\xi}\right)}{\hat{\xi}}\right)\nabla\hat{\rho}\right]\nonumber\\
&	
\quad+\int_{\mathbb{T}^N}\xi\nabla^{\alpha}\left[\left(\dfrac{\nu(\xi,\psi)}{\xi}-\dfrac{\nu(\hat{\xi},\hat{\psi})}{\hat{\xi}}\right)\Delta\hat{u}\right]\cdot\nabla^{\alpha}\left(u-\hat{u}\right)
\nonumber\\
&\quad+\int_{\mathbb{T}^N}\xi\left\{\nabla^{\alpha}\left[\dfrac{\nu(\xi,\psi)}{\xi}\Delta(u-\hat{u})\right]-\dfrac{\nu(\xi,\psi)}{\xi}\nabla^{\alpha}\Delta(u-\hat{u})\right\}
\cdot\nabla^{\alpha}\left(u-\hat{u}\right)
\nonumber\\
&\quad-\int_{\mathbb{T}^N}\nabla\nu(\xi,\psi)\cdot\nabla^{\alpha}\nabla(u-\hat{u})
\cdot\nabla^{\alpha}\left(u-\hat{u}\right)
\nonumber\\
&\quad+\int_{\mathbb{T}^N}\xi\left\{\nabla^{\alpha}\left[\dfrac{\eta(\xi,\psi)}{\xi}\nabla(\nabla\cdot(u-\hat{u}))\right]
-\dfrac{\eta(\xi,\psi)}{\xi}\nabla^{\alpha}\nabla(\nabla\cdot(u-\hat{u}))\right\}\cdot\nabla^{\alpha}\left(u-\hat{u}\right)
\nonumber\\
&\quad-\int_{\mathbb{T}^N}
\nabla\eta(\xi,\psi)\nabla^{\alpha}(\nabla\cdot(u-\hat{u}))\cdot\nabla^{\alpha}\left(u-\hat{u}\right)
\nonumber\\
&\quad+\int_{\mathbb{T}^N}\xi\nabla^{\alpha}\left[\left(\dfrac{\eta(\xi,\psi)}{\xi}-\dfrac{\eta(\hat{\xi},\hat{\psi})}{\hat{\xi}}\right)\nabla(\nabla\cdot\hat{u})\right]\cdot\nabla^{\alpha}\left(u-\hat{u}\right)
\nonumber\\
&\quad-\int_{\mathbb{T}^N}\xi\nabla^{\alpha}\left[\left(\dfrac{1}{\xi}-\dfrac{1}{\hat{\xi}}\right)\Delta\hat\phi\nabla\hat\phi\right]\cdot\nabla^{\alpha}\left(u-\hat{u}\right)
\nonumber\\
&\quad-\int_{\mathbb{T}^N}\xi\nabla^{\alpha}\left\{\dfrac{1}{\xi}\left[\Delta\left(\phi-\hat{\phi}\right)\nabla\phi
+\Delta\hat{\phi}\nabla\left(\phi-\hat{\phi}\right)\right]\right\}\cdot\nabla^{\alpha}\left(u-\hat{u}\right),
	\end{align}
where we have used the following fact by integration by parts that
\begin{align*}
&\int_{\mathbb{T}^N}\nabla^{\alpha}\left(\dfrac{\nu(\xi,\psi)}{\xi}\Delta\left(u-\hat{u}\right)\right)\cdot\xi\nabla^{\alpha}\left(u-\hat{u}\right)\\
&=\int_{\mathbb{T}^N}\xi\left\{\nabla^{\alpha}\left[\dfrac{\nu(\xi,\psi)}{\xi}\Delta(u-\hat{u})\right]-\dfrac{\nu(\xi,\psi)}{\xi}\nabla^{\alpha}\Delta(u-\hat{u})\right\}
\cdot\nabla^{\alpha}\left(u-\hat{u}\right)\\
&\quad+\int_{\mathbb{T}^N}\nu(\xi,\psi)\nabla^{\alpha}\Delta(u-\hat{u})\cdot\nabla^{\alpha}\left(u-\hat{u}\right)\\
&=\int_{\mathbb{T}^N}\xi\left\{\nabla^{\alpha}\left[\dfrac{\nu(\xi,\psi)}{\xi}\Delta(u-\hat{u})\right]-\dfrac{\nu(\xi,\psi)}{\xi}\nabla^{\alpha}\Delta(u-\hat{u})\right\}
\cdot\nabla^{\alpha}\left(u-\hat{u}\right)\\
&\quad-\int_{\mathbb{T}^N}\nabla\nu(\xi,\psi)\cdot\nabla^{\alpha}\nabla(u-\hat{u})\cdot\nabla^{\alpha}\left(u-\hat{u}\right)
-\int_{\mathbb{T}^N}\nu(\xi,\psi)|\nabla^{\alpha}\nabla\left(u-\hat{u}\right)|^2.
\end{align*}
Then for $(\ref{q.3.61})$ and $(\ref{q.3.62})$, we give the calculations of singular term as follows:\\
For $\alpha=0$, we obtain
\begin{align} \label{C-1}
-\dfrac{1}{\varepsilon^2}\int_{\mathbb{T}^N}\dfrac{P'\left(\xi\right)}{\xi}\nabla^{\alpha}\left(\rho-\hat{\rho}\right)
\nabla^{\alpha}\left[\xi\nabla\cdot\left(u-\hat{u}\right)\right]
=-\dfrac{1}{\varepsilon^2}\int_{\mathbb{T}^N}{P'\left(\xi\right)}\left(\rho-\hat{\rho}\right)
\nabla\cdot\left(u-\hat{u}\right),
\end{align}
and
\begin{align} \label{C-11}
&-\dfrac1{\varepsilon^2}\int_{\mathbb{T}^N}\xi\nabla^{\alpha}
 \left(u-\hat{u}\right)\cdot\nabla^{\alpha}\left(\dfrac{P'\left(\xi\right)}{\xi}\nabla\left(\rho-\hat{\rho}\right)\right)\nonumber\\
 &=-\dfrac1{\varepsilon^2}\int_{\mathbb{T}^N}{P'\left(\xi\right)}
 \left(u-\hat{u}\right)\cdot\nabla\left(\rho-\hat{\rho}\right)\nonumber\\
 &=\dfrac1{\varepsilon^2}\int_{\mathbb{T}^N}\left(u-\hat{u}\right)\cdot\nabla{P'\left(\xi\right)}\left(\rho-\hat{\rho}\right)
 +\dfrac1{\varepsilon^2}\int_{\mathbb{T}^N}{P'\left(\xi\right)}\nabla\cdot\left(u-\hat{u}\right)\left(\rho-\hat{\rho}\right),
 \end{align}
 where we have used integration by parts.
For $\alpha=1$, in a similar way, we have
\begin{align} \label{C-2}
&-\dfrac{1}{\varepsilon^2}\int_{\mathbb{T}^N}\dfrac{P'\left(\xi\right)}{\xi}\nabla^{\alpha}\left(\rho-\hat{\rho}\right)
\nabla^{\alpha}\left[\xi\nabla\cdot\left(u-\hat{u}\right)\right]
=-\dfrac{1}{\varepsilon^2}\int_{\mathbb{T}^N}\dfrac{P'\left(\xi\right)}{\xi}\nabla_i\left(\rho-\hat{\rho}\right)
\nabla_i\left[\xi\nabla_j\left(u-\hat{u}\right)_j\right]\nonumber\\
&=-\dfrac{1}{\varepsilon^2}\int_{\mathbb{T}^N}\dfrac{P'\left(\xi\right)}{\xi}\nabla_i\left(\rho-\hat{\rho}\right)\nabla_i\xi\nabla_j\left(u-\hat{u}\right)_j
-\dfrac{1}{\varepsilon^2}\int_{\mathbb{T}^N}P'\left(\xi\right)\nabla_i\left(\rho-\hat{\rho}\right)
\nabla_i\nabla_j\left(u-\hat{u}\right)_j,
\end{align}
and
\begin{align}\label{qq.3.630}
&-\dfrac1{\varepsilon^2}\int_{\mathbb{T}^N}\xi\nabla^{\alpha}
 \left(u-\hat{u}\right)\cdot\nabla^{\alpha}\left(\dfrac{P'\left(\xi\right)}{\xi}\nabla\left(\rho-\hat{\rho}\right)\right)
 =-\dfrac1{\varepsilon^2}\int_{\mathbb{T}^N}\xi\nabla_i
 \left(u-\hat{u}\right)_j\nabla_i\left(\dfrac{P'\left(\xi\right)}{\xi}\nabla_j\left(\rho-\hat{\rho}\right)\right)\nonumber\\
 &=-\dfrac1{\varepsilon^2}\int_{\mathbb{T}^N}\xi\nabla_i
 \left(u-\hat{u}\right)_j\nabla_i\left(\dfrac{P'\left(\xi\right)}{\xi}\right)\nabla_j\left(\rho-\hat{\rho}\right)
 -\dfrac1{\varepsilon^2}\int_{\mathbb{T}^N}{P'\left(\xi\right)}\nabla_i
 \left(u-\hat{u}\right)_j\nabla_i\nabla_j\left(\rho-\hat{\rho}\right)\nonumber\\
 &=-\dfrac1{\varepsilon^2}\int_{\mathbb{T}^N}\xi\nabla_i
 \left(u-\hat{u}\right)_j\nabla_i\left(\dfrac{P'\left(\xi\right)}{\xi}\right)\nabla_j\left(\rho-\hat{\rho}\right)
 +\dfrac1{\varepsilon^2}\int_{\mathbb{T}^N}\nabla_j{P'\left(\xi\right)}\nabla_i
 \left(u-\hat{u}\right)_j\nabla_i\left(\rho-\hat{\rho}\right)\nonumber\\
 &\quad+\dfrac1{\varepsilon^2}\int_{\mathbb{T}^N}{P'\left(\xi\right)}\nabla_i\nabla_j
 \left(u-\hat{u}\right)_j\nabla_i\left(\rho-\hat{\rho}\right).
\end{align}

 In conclusion, putting (\ref{q.3.59})-(\ref{qq.3.630}) together, processing as before and using the Cauchy inequality, we are led to
	\begin{align}\label{q.3.63}
		&\dfrac{1}{2}\sum_{\alpha\leq 1}\dfrac{\rm d}{{\rm d}t}\int_{\mathbb{T}^N}\left[\dfrac{P'\left(\xi\right)}{\xi}\left|\dfrac{1}{\varepsilon}\nabla^{\alpha}\left(\rho-\hat{\rho}\right)\right|^2
+\xi|\nabla^{\alpha}(u-\hat{u})|^2\right]
+\dfrac{1}{2}\dfrac{\rm d}{{\rm d}t}\int_{\mathbb{T}^N}\left(|\phi-\hat{\phi}|^2+|\nabla\left(\phi-\hat{\phi}\right)|^2\right)
\nonumber\\	
&\quad+\sum_{\alpha\leq 1}\int_{\mathbb{T}^N}\left[\nu_*|\nabla^{\alpha}\nabla\left(u-\hat{u}\right)|^2 + \eta_*|\nabla^{\alpha}\nabla\cdot\left(u-\hat{u}\right)|^2\right]
\nonumber\\	
&\quad+\int_{\mathbb{T}^N}\dfrac{1}{\xi^2}\left[\left|\Delta\left(\phi-\hat{\phi}\right)\right|^2+\left|\nabla\Delta\left(\phi-\hat{\phi}\right)\right|^2\right]  \nonumber\\
		&\leq C\left(\left\|\dfrac{1}{\varepsilon}\left(\rho-\hat{\rho}\right)\right\|_1^2+\left\|u-\hat{u}\|_1^2+\|\phi-\hat{\phi}\right\|_1^2\right)
+C\left(\left\|\dfrac{1}{\varepsilon}\left(\xi-\hat{\xi}\right)\right\|_1^2+\left\|v-\hat{v}\right\|_1^2+\left\|\psi-\hat{\psi}\right\|_1^2\right)\nonumber\\
&\quad+\varepsilon\left\|\Delta(\psi-\hat{\psi})\right\|_1^2+\varepsilon\left\|\Delta\left(\phi-\hat{\phi}\right)\right\|^2
+\tau\left(\left\|\Delta\left(\phi-\hat{\phi}\right)\right\|_1^2
+\|\nabla\left(u-\hat{u}\right)\|_1^2+\|\nabla\cdot\left(u-\hat{u}\right)\|_1^2\right).
	\end{align}
Then noting that$ \left(U-\hat{U}\right)(0)=0$, by the Gronwall inequality and (\ref{q.3.63}), we have
	\begin{align}\label{q.3.66}
		&\sup_{0\le t\le T_0}\left(\left\|\dfrac{1}{\varepsilon}\left(\rho-\hat{\rho}\right)\right\|_1^2+\left\|u-\hat{u}\|_1^2+\|\phi-\hat{\phi}\right\|_1^2\right)\nonumber\\
		&\leq\exp^{CT_0}\left[CT_0\sup_{0\le t\le T_0}\left(\left\|\dfrac{1}{\varepsilon}\left(\xi-\hat{\xi}\right)\right\|_1^2+\left\|v-\hat{v}\|_1^2+\|\psi-\hat{\psi}\right\|_1^2\right)
+\varepsilon\int_0^t\left\|\Delta(\psi-\hat{\psi})\right\|_1^2\right].
	\end{align}
	Moreover, integrating (\ref{q.3.63}) over $[0, t]\subseteq[0, T_0]$, it holds that
	\begin{align}\label{q.3.67}
		&\int_{0}^{t}\left(\left\|u-\hat{u}\right\|_2^2+\left\|\phi-\hat{\phi}\right\|_3^2\right)\nonumber\\
&\leq CT_0\exp^{CT_0}\left[CT_0\sup_{0\le t\le T_0}\left(\left\|\dfrac{1}{\varepsilon}\left(\xi-\hat{\xi}\right)\right\|_1^2+\left\|v-\hat{v}\|_1^2+\|\psi-\hat{\psi}\right\|_1^2\right)
+\varepsilon\int_0^t\left\|\Delta(\psi-\hat{\psi})\right\|_1^2\right]\nonumber\\
&\quad+CT_0\left(\left\|\dfrac{1}{\varepsilon}\left(\xi-\hat{\xi}\right)\right\|_1^2+\left\|v-\hat{v}\right\|_1^2+\left\|\psi-\hat{\psi}\right\|_1^2\right)
+\varepsilon\int_0^t\left\|\Delta(\psi-\hat{\psi})\right\|_1^2.
	\end{align}
 If we take $T_0\left(<T_3\right)$ and $\varepsilon$ small enough such that $C\left(T_0^2\exp^{CT_0}+T_0\right)<1$, $C\varepsilon T_0\exp^{CT_0}+\varepsilon<1$, (\ref{q.3.55}) and (\ref{q.3.57}) are valid, where $C$ is the uniform constants mentioned in (\ref{q.3.66})-(\ref{q.3.67}), then we can prove the contraction.
The proof of Lemma \ref{L3.3} is completed.
\hfill$\Box$

\begin{remark}
The contraction estimates of the map $\Lambda$ for NSCH system in Lemma \ref{L3.3} is one-order higher than NSAC system in Lemma \ref{L3.3-2}.
This is because of the three-order term $\nabla\Delta(\psi-\hat\psi)$ appearing in the forth term in (\ref{q.3.59}).
And it also leads to the last single contraction term in (\ref{aq.3.58}).
\end{remark}

\begin{lemma}\label{L3.4}
	Consider the incompressible system of Navier-Stokes/Cahn-Hilliard (\ref{q.1.2}) with initial condition (\ref{q.2.5}) for $s\geq 2$. Then for $0<T\leq\infty$, there exists at most one strong solution $\left(u, \phi\right)\in\left \{\left(u, \phi\right): u\in L^{\infty}\left(0, T; H^{s+2}\right),\nabla\phi\in L^{\infty}\left(0, T; H^{s+1}\right)\right \} $.
	\end{lemma}
	\noindent{\it\bfseries Proof.}\quad Assume that $\left(u_1, \phi_1\right)$ and $\left(u_2, \phi_2\right)$ are two strong solutions of (\ref{q.1.2}) with the same initial condition (\ref{q.2.5}). Then we obtain
	\begin{align}\label{q.3.69}
		\begin{cases}
			\nabla\cdot(u_1-u_2)=0,\\
			\partial_t\left(u_1-u_2\right)+\left(u_1\cdot\nabla\right)\left(u_1-u_2\right)+\left(\left(u_1-u_2\right)\cdot\nabla\right)u_2+\nabla\left(p_1-p_2\right)\\
			\quad= \nu(\phi_1)\Delta\left(u_1-u_2\right)+(\nu(\phi_1)-\nu(\phi_2))\Delta u_2 \\
\quad\quad-\nabla\cdot\left(\nabla\left(\phi_1-\phi_2\right)\otimes\nabla\phi_1\right)-\nabla\cdot\left(\nabla\phi_2\otimes\nabla\left(\phi_1-\phi_2\right)\right), \\
			\partial_t\left(\phi_1-\phi_2\right)+\left(u_1\cdot\nabla\right)\left(\phi_1-\phi_2\right)+\left(\left(u_1-u_2\right)\cdot\nabla\right)\phi_2\\
			\quad=-\Delta^2\left(\phi_1-\phi_2\right)+\Delta\left(\phi_1^3-\phi_2^3\right)-\Delta\left(\phi_1-\phi_2\right).
		\end{cases}
	\end{align}
Multiplying $(\ref{q.3.69})_2$ by $\left(u_1-u_2\right)$ and multiplying $(\ref{q.3.69})_3$ by $(\phi_1-\phi_2)$, by integration by parts, we have
%
	\begin{align*}
		&\dfrac{1}{2}\dfrac{\rm d}{{\rm d}t}\int_{\mathbb{T}^N}|u_1-u_2|^2+\int_{\mathbb{T}^N}\nu(\phi_1)|\nabla\left(u_1-u_2\right)|^2  \nonumber \\
		&=-\int_{\mathbb{T}^N}\left(u_1\cdot\nabla\right)\left(u_1-u_2\right)\cdot\left(u_1-u_2\right)-\int_{\mathbb{T}^N}\left(\left(u_1-u_2\right)\cdot\nabla\right)u_2\cdot\left(u_1-u_2\right)   \nonumber\\
		&\quad-\int_{\mathbb{T}^N}\nabla\left(p_1-p_2\right)\cdot\left(u_1-u_2\right)-\int_{\mathbb{T}^N}\nabla\cdot\left(\nabla\left(\phi_1-\phi_2\right)\otimes\nabla\phi_1\right)\cdot\left(u_1-u_2\right)   \nonumber\\
		&\quad-\int_{\mathbb{T}^N}\nabla\cdot\left(\nabla\phi_2\otimes\nabla\left(\phi_1-\phi_2\right)\right)\cdot\left(u_1-u_2\right)
-\int_{\mathbb{T}^N}\nabla\nu(\phi_1)\cdot\nabla\left(u_1-u_2\right)\cdot\left(u_1-u_2\right)\nonumber\\
&\quad+\int_{\mathbb{T}^N}(\nu(\phi_1)-\nu(\phi_2))\Delta u_2\cdot\left(u_1-u_2\right),
	\end{align*}
	%
and
	\begin{align*}
		&\dfrac{1}{2}\dfrac{\rm d}{{\rm d}t}\int_{\mathbb{T}^N}|\phi_1-\phi_2|^2+\int_{\mathbb{T}^N}|\Delta\left(\phi_1-\phi_2\right)|^2  \nonumber\\
		&\quad=-\int_{\mathbb{T}^N}\left(u_1\cdot\nabla\right)\left(\phi_1-\phi_2\right)\left(\phi_1-\phi_2\right)-\int_{\mathbb{T}^N}\left(\left(u_1-u_2\right)\cdot\nabla\right)\phi_2\left(\phi_1-\phi_2\right)   \nonumber\\
		&\qquad+\int_{\mathbb{T}^N}\left(\phi_1^3-\phi_2^3\right)\Delta\left(\phi_1-\phi_2\right)-\int_{\mathbb{T}^N}\Delta\left(\phi_1-\phi_2\right)\left(\phi_1-\phi_2\right),
	\end{align*}
	and then using integration by parts, the Mean value theorem, the $\rm Poincar\acute{e}$ inequality and the Cauchy inequality give
	\begin{align*}
		&\dfrac{1}{2}\dfrac{\rm d}{{\rm d}t}\int_{\mathbb{T}^N}\left(|u_1-u_2|^2+|\phi_1-\phi_2|^2\right)+\nu_*\int_{\mathbb{T}^N}|\nabla\left(u_1-u_2\right)|^2    \nonumber\\	
&\quad+\int_{\mathbb{T}^N}|\nabla\left(\phi_1-\phi_2\right)|^2+\int_{\mathbb{T}^N}|\Delta\left(\phi_1-\phi_2\right)|^2   \nonumber\\
		&\leq M\left(\|u_1-u_2\|^2+\|\phi_1-\phi_2\|^2\right)+\dfrac{\nu_*}{2}\|\nabla\left(u_1-u_2\right)\|^2+\dfrac{1}{2}\|\nabla\left(\phi_1-\phi_2\right)\|_1^2,
	\end{align*}
	where $M\geq\sup_{t\in[0, T]}\frac94\left(\|u_i\|_2^4\left(t\right)+\|\phi_i\|_3^4\left(t\right)+\nu_\phi^2(\phi)\right)$ for $i=1, 2$.

Note that $\|u_1\left(x,0\right)-u_2\left(x, 0\right)\|^2+\|\phi_1\left(x,0\right)-\phi_2\left(x,0\right)\|^2=0$. By the Gronwall inequality, it holds that
	\begin{align}\label{q.3.74}
		\|u_1-u_2\|^2+\|\phi_1-\phi_2\|^2=0,
	\end{align}
	which completes the proof of Lemma \ref{L3.4}.
\hfill$\Box$
%
	%
	
 Now we are in a position to prove Theorem \ref{th1}.
 \\
 \\
 \noindent{\it\bfseries Proof of Theorem \ref{th1}.}
 For any fixed $\varepsilon$, the standard procedure produces a sequence $\left\{\left(\rho_i, u_i\right.\right.$, $\left.\left.\phi_i\right)\right\}_{i=0}^{\infty}$ satisfying
	\begin{equation*}
		\begin{aligned}
		\begin{cases}
			\partial_t\rho_{i+1}^{\varepsilon}+\left(u_i^{\varepsilon}\cdot\nabla\right)\rho_{i+1}^{\varepsilon}+\rho_i^{\varepsilon}\nabla\cdot u_{i+1}^{\varepsilon}=0, \nonumber \\
			\partial_t u_{i+1}^{\varepsilon}+\left(u_i^{\varepsilon}\cdot\nabla\right)u_{i+1}^{\varepsilon}+\dfrac{1}{\varepsilon^2}\dfrac{P'\left(\rho_i^{\varepsilon}\right)}{\rho_i^{\varepsilon}}\nabla\rho_{i+1}^{\varepsilon} \nonumber \\
			\quad=\dfrac{\nu(\rho_i^{\varepsilon},\phi_i^\varepsilon)}{\rho_i^{\varepsilon}}\Delta u_{i+1}^{\varepsilon}+\dfrac{\eta(\rho_i^{\varepsilon},\phi_i^\varepsilon)}{\rho_i^{\varepsilon}}\nabla\left(\nabla\cdot u_{i+1}^{\varepsilon}\right)-\dfrac{1}{\rho_i^{\varepsilon}}\left(\Delta\phi_{i+1}^{\varepsilon}\nabla\phi_{i+1}^{\varepsilon}\right),
\nonumber \\
			\partial_t\phi_{i+1}^{\varepsilon}+\left(u_i^{\varepsilon}\cdot\nabla\right)\phi_i^{\varepsilon}\nonumber\\
\quad=-\dfrac{1}{\left(\rho_i^\varepsilon\right)^2}\Delta^2\phi_{i+1}^\varepsilon
-\dfrac2{\rho_i^\varepsilon}\nabla\left(\dfrac1{\rho_i^\varepsilon}\right)\cdot\nabla\Delta\phi_i^\varepsilon
-\dfrac1{\rho_i^\varepsilon}\Delta\left(\dfrac1{\rho_i^\varepsilon}\right)\Delta\phi_i^\varepsilon
+\dfrac1{\rho_i^\varepsilon}\Delta\left(\left(\phi_i^\varepsilon\right)^3-\phi_i^\varepsilon\right)
		\end{cases}
		\end{aligned}
	\end{equation*}
	as well as the uniform estimates
	\begin{align}\label{q.3.75}
	\begin{cases}
		E_{s+2}(\rho_i^\varepsilon,u_i^\varepsilon)+F_{s+1}(\nabla\phi_i^\varepsilon)+\|\phi_i^\varepsilon\|\\
\quad\quad\quad\quad\quad\quad\quad+\displaystyle\int_0^t\left(\nu_*\left\|\nabla u_i^\varepsilon\right\|^2_{s+2}+\eta_*\left\|\nabla\cdot u_i^\varepsilon\right\|^2_{s+2}+\left\|\nabla\phi_i^\varepsilon\right\|^2_{s+3}\right)\leq C,\\
		E_{s+1}(\partial_t(\rho_i^\varepsilon, u_i^\varepsilon))+F_{s}(\partial_t\phi^\varepsilon_i)\\
\quad\quad\quad\quad\quad\quad\quad+\displaystyle\int_0^t\left(\nu_*\left\|\partial_t\nabla u_i^\varepsilon\right\|^2_{s+1}+\eta_*\left\|\partial_t\nabla\cdot u_i^\varepsilon\right\|^2_{s+1}+\left\|\partial_t\phi_i^\varepsilon\right\|^2_{s+2}\right)\leq C.
	\end{cases}
	\end{align}
Let $\hat{\rho}_{i+1}^{\varepsilon}=\rho_{i+1}^{\varepsilon}-\rho_i^{\varepsilon},\ \hat{u}_{i+1}^{\varepsilon}=u_{i+1}^{\varepsilon}-u_i^{\varepsilon}$ and $\hat{\phi}_{i+1}^{\varepsilon}=\phi_{i+1}^{\varepsilon}-\phi_i^{\varepsilon}$. In view of Lemma \ref{L3.3}, we get
	\begin{align}\label{q.3.76}
		&\sum_{i=2}^{\infty}\left\|\hat{\rho}_i^{\varepsilon}\right\|_1<\infty,\qquad\sum_{i=2}^{\infty}\left(\left\|\hat{u}_i^{\varepsilon}\right\|_1
+\int_{0}^{T_0}\left\|\hat{u}_i^{\varepsilon}\right\|_2^2\right)<\infty , \nonumber\\
		&\qquad\sum_{i=2}^{\infty}\left(\left\|\hat{\phi}_i^{\varepsilon}\right\|_1+\int_{0}^{T_0}\left\|\hat{\phi}_i^{\varepsilon}\right\|_3^2\right)<\infty .
	\end{align}
	
 Let $\rho^{\varepsilon}=\rho_1^{\varepsilon}+\sum_{i=2}^{\infty}\hat{\rho}_i^{\varepsilon}$, $u^{\varepsilon}=u_1^{\varepsilon}+\sum_{i=2}^{\infty}\hat{u}_i^{\varepsilon}$ and $\phi^{\varepsilon}=\phi_1^{\varepsilon}+\sum_{i=2}^{\infty}\hat{\phi}_i^{\varepsilon}$,
	and then we rewrite $\left(\rho_i^{\varepsilon}, u_i^{\varepsilon}, \phi_i^{\varepsilon} \right)$ as
\ $\rho_i^{\varepsilon}=\rho_1^{\varepsilon}+\sum_{n=2}^{i}\left(\rho_n^{\varepsilon}-\rho_{n-1}^{\varepsilon}\right),\ u_i^{\varepsilon}=u_1^{\varepsilon}+\sum_{n=2}^{i}\left(u_n^{\varepsilon}-u_{n-1}^{\varepsilon}\right)$ and $\phi_i^{\varepsilon}=\phi_1^{\varepsilon}+\sum_{n=2}^{i}\left(\phi_n^{\varepsilon}-\phi_{n-1}^{\varepsilon}\right)$.
For the density, it follows from (\ref{q.3.76}) that
	\begin{align}\label{q.3.77}
		\left\|\rho_i^{\varepsilon}-\rho^{\varepsilon}\right\|_1
&=\left\|\sum_{n=2}^{i}\left(\rho_n^{\varepsilon}-\rho_{n-1}^{\varepsilon}\right)-\sum_{i=2}^{\infty}\hat{\rho}_i^{\varepsilon}\right\|_1
=\left\|\sum_{n=2}^{i}\left(\rho_n^{\varepsilon}-\rho_{n-1}^{\varepsilon}\right)-\sum_{i=2}^{\infty}\left(\rho_i^{\varepsilon}-\rho_{i-1}^{\varepsilon}\right)\right\|_1 \nonumber\\&\leq\sum_{k=i+1}^{\infty}\left\|\rho_k^{\varepsilon}-\rho_{k-1}^{\varepsilon}\right\|_1
=\sum_{k=i+1}^{\infty}\left\|\hat{\rho}_k^{\varepsilon}\right\|_1\to 0.
	\end{align}
Similar to (\ref{q.3.77}), we can derive
\begin{align}\label{q.3.78}
\left\|u_i^{\varepsilon}-u^{\varepsilon}\right\|_1+\int_{0}^{T_0}\left\|u_i^{\varepsilon}-u^{\varepsilon}\right\|_2^2\to 0,\quad{\rm and}\quad
\left\|\phi_i^{\varepsilon}-\phi^{\varepsilon}\right\|_1+\int_{0}^{T_0}\left\|\phi_i^{\varepsilon}-\phi^{\varepsilon}\right\|_3^2 \to 0.
\end{align}
 Therefore, (\ref{q.3.77}) and (\ref{q.3.78}) indicate
	\begin{align*}
		&\rho_i^{\varepsilon}\to\rho^{\varepsilon},\quad {\rm in}\  L^{\infty}\left([0, T_0]; H^1\right),\\
		&u_i^{\varepsilon}\to u^{\varepsilon},\quad {\rm in}\  L^{\infty}\left([0, T_0]; H^1\right)\cap L^2\left([0, T_0]; H^2\right),\\
		&\phi_i^{\varepsilon}\to\phi^{\varepsilon},\quad {\rm in}\  L^{\infty}\left([0, T_0]; H^1\right)\cap L^2\left([0, T_0]; H^3\right).\nonumber
	\end{align*}
It infers from the estimates (\ref{q.3.75}) and the lower semi-continuity that
	\begin{align*}
	&\rho^{\varepsilon},\ u^{\varepsilon}\in L^{\infty}\left([0, T_0]; H^{s+2}\right)\cap {\rm Lip}\left([0, T_0]; H^{s+1}\right),
\nonumber\\
&\phi\in L^{\infty}\left([0, T_0]; H^{s+2}\right)\cap {\rm Lip}\left([0, T_0]; H^{s}\right).
\end{align*}
 For any $s'+2\in\left [ 0, s+2  \right) $, by the Sobolev interpolation inequalities, we are led to
	\begin{align}
		&\left\|\left(\rho_i^{\varepsilon},u_i^{\varepsilon},\phi_i^{\varepsilon}\right)
-\left(\rho^{\varepsilon},u^{\varepsilon},\phi^{\varepsilon}\right)\right\|_{s'+2}\nonumber\\
		&\leq C\left\|\left(\rho_i^{\varepsilon},u_i^{\varepsilon},\phi_i^{\varepsilon}\right)
-\left(\rho^{\varepsilon},u^{\varepsilon},\phi^{\varepsilon}\right)\right\|^{\theta}
\left(\left\|\left(\rho_i^{\varepsilon},u_i^{\varepsilon},\phi_i^{\varepsilon}\right)\right\|_{s+2}
+\left\|\left(\rho^{\varepsilon},u^{\varepsilon},\phi^{\varepsilon}\right)\right\|_{s+2}\right)^{1-\theta}\to 0,
		\end{align}
	as $i\to\infty$ for some $\theta\in\left(0,1\right)$, where we have used Lemma \ref{L3.3} to get
	\begin{align*}
		\left\|\left(\rho_i^{\varepsilon},u_i^{\varepsilon},\phi_i^{\varepsilon}\right)-\left(\rho^{\varepsilon},u^{\varepsilon},\phi^{\varepsilon}\right)\right\|
&\leq\sum_{j=i+1}^{\infty}\left\|\left(\rho_j^{\varepsilon},u_j^{\varepsilon},\phi_j^{\varepsilon}\right)-\left(\rho_{j-1}^{\varepsilon},u_{j-1}^{\varepsilon},\phi_{j-1}^{\varepsilon}\right)\right\| \\
		&\leq\dfrac{C\gamma_1^i}{1-\gamma_1}.
	\end{align*}
	Hence, $\left(\rho^{\varepsilon}, u^{\varepsilon}, \phi^{\varepsilon}\right)\in C\left([0,T_0];H^{s'+2}\right)$. In addition, by Lemma \ref{L3.5} and (\ref{q.3.75}), one deduces easily that $\left(\rho^{\varepsilon}, u^{\varepsilon}, \phi^{\varepsilon}\right)$ is a strong solution of compressible equation (\ref{q.1.1}).

The uniqueness can be proved by a similar argument as Lemma \ref{L3.3}, which completes the proof of the uniform stability part of Theorem \ref{th1}.

	 Next, we show that $\left(\rho^{\varepsilon}, u^{\varepsilon}, \phi^{\varepsilon}\right)$ converges to the unique strong solution to the corresponding incompressible equation (\ref{q.1.2}) as $\varepsilon\to 0$.
To see this, note that (\ref{q.2.7}) implies that $\rho^{\varepsilon}\to 1$ in $L^{\infty}\left([0, T_0]; H^{s+2}\right)\cap {\rm Lip}\left([0, T_0]; H^{s+1}\right)$, and there exists a subsequence
	$\left \{\left(u^{\varepsilon_j}, \phi^{\varepsilon_j}\right)\right \}_j $ of $\left \{\left(u^{\varepsilon}, \phi^{\varepsilon}\right)\right \}_{\varepsilon} $ with a limit u and $\phi$ such that
	\begin{align}
		\begin{cases}
			u^{\varepsilon_j}\rightharpoonup u\quad {\rm weakly^* \ in}\  L^{\infty}\left([0, T_0]; H^{s+2}\right)\cap {\rm Lip}\left([0, T_0]; H^{s+1}\right),\\
			\phi^{\varepsilon_j}\rightharpoonup\phi\quad {\rm weakly^*\  in} \ L^{\infty}\left([0, T_0]; H^{s+2}\right)\cap {\rm Lip}\left([0, T_0]; H^{s}\right),\\
			u^{\varepsilon_j}\to u\quad {\rm in}\  C\left([0, T_0]; H^{s'+2}\right),\\
			\phi^{\varepsilon_j}\to\phi\quad {\rm in}\  C\left([0, T_0]; H^{s'+2}\right).
		\end{cases}
	\end{align}
	for any $0\leq s'+2<s+2$, where we have the fact the embedding $H^{s+2}\hookrightarrow H^{s'+2}$ is compact and Lemma \ref{L3.5}.

 Now we are to let $j\to\infty\left(\varepsilon_j\to 0\right)$ in (\ref{q.1.1}).

First of all, multiplying $(\ref{q.1.1})_1$ and $(\ref{q.1.1})_3$ by two smooth test functions $\psi_1\left(x,t\right)$ and $\psi_3\left(x,t\right)$ with compact supports in $t\in[0, T_0]$ respectively, we obtain
	\begin{align*}
		\int_{0}^{T_0}\int_{\mathbb{T}^N}\nabla\cdot u^{\varepsilon_j}\psi_1=-\int_{0}^{T_0}\int_{\mathbb{T}^N}\left(\rho_t^{\varepsilon_j}+\left(u^{\varepsilon_j}\cdot\nabla\right)\rho^{\varepsilon_j}
+\left(\rho^{\varepsilon_j}-1\right)
		\nabla\cdot u^{\varepsilon_j}\right)\psi_1,
	\end{align*}
	\begin{align*}
		&\int_{0}^{T_0}\int_{\mathbb{T}^N}\left(\phi_t^{\varepsilon_j}+\left(u^{\varepsilon_j}\cdot\nabla\right)\phi^{\varepsilon_j}
+\dfrac{1}{\left(\rho^{\varepsilon_j}\right)^2}\Delta^2\phi^{\varepsilon_j}
+\dfrac2{\xi^{\varepsilon_j}}\nabla\left(\dfrac1{\xi^{\varepsilon_j}}\right)\cdot\nabla\Delta\psi^{\varepsilon_j}
\right.\nonumber\\
        &\left.\qquad\qquad\qquad\qquad+\dfrac1{\xi^{\varepsilon_j}}\Delta\left(\dfrac1{\xi^{\varepsilon_j}}\right)\Delta\psi^{\varepsilon_j}
-\dfrac{1}{\rho^{\varepsilon_j}}\Delta\left(\left(\phi^{\varepsilon_j}\right)^3-\phi^{\varepsilon_j}\right)\right)\psi_3=0.
	\end{align*}

	Then $(u,\phi)$ satisfies $(\ref{q.1.2})_1$ and $(\ref{q.1.2})_3$ by sending $j\to\infty$.

	 Let $\psi_2\left(x,t\right)$ be a smooth test function of $(\ref{q.1.1})_2$ with compact supports in $t\in[0,T_0]$ and the divergence free condition $\nabla\cdot\psi_2=0$. Then integration by part leads to
	\begin{align*}
		&\int_{0}^{T_0}\int_{\mathbb{T}^N}\left[u_t^{\varepsilon_j}+\left(u^{\varepsilon_j}\cdot\nabla\right)u^{\varepsilon_j}-\dfrac{\nu(\rho^{\varepsilon_j},\phi^{\varepsilon_j})}{\rho^{\varepsilon_j}}\Delta u^{\varepsilon_j}-\dfrac{\eta(\rho^{\varepsilon_j},\phi^{\varepsilon_j})}{\rho^{\varepsilon_j}}\nabla\left(\nabla\cdot u^{\varepsilon_j}\right)\right.
\nonumber\\
 &\left.\qquad\quad~~+\dfrac{1}{\rho^{\varepsilon_j}}\nabla\cdot\left(\nabla\phi^{\varepsilon_j}\otimes\nabla\phi^{\varepsilon_j}\right)
-\dfrac{1}{\rho^{\varepsilon_j}}\nabla\left(\dfrac{|\nabla\phi^{\varepsilon_j}|^2}{2}\right)\right]\cdot\psi_2 \nonumber\\
		&=-\dfrac{1}{\varepsilon_j^2}\int_{0}^{T_0}\int_{\mathbb{T}^N}\dfrac{1}{\rho^{\varepsilon_j}}\psi_2\cdot\nabla\left(P\left(\rho^{\varepsilon_j}\right)\right)
=-\dfrac{1}{\varepsilon^2}\int_{0}^{T_0}\int_{\mathbb{T}^N}\psi_2\cdot\nabla\int_{1}^{\rho^{\varepsilon_j}}\dfrac{P'\left(\xi\right)}{\xi}{\rm d}\xi=0.
	\end{align*}
Then, let $j\to\infty$ and get
		\begin{align*}
		\textbf{P}\left(u_t+\left(u\cdot\nabla\right)u-{\nu(\phi)}\Delta u+\nabla\cdot\left(\nabla\phi\otimes\nabla\phi\right)
-\nabla\left(\dfrac{|\nabla\phi|^2}{2}\right)\right)=0,
	\end{align*}
	where \textbf{P} is the $L^2$-projection on the divergence free vector fields.\\

	If 
\begin{align}\label{q.3.85}
		u_t+\left(u\cdot\nabla\right)u-{\nu(\phi)}\Delta u+\nabla\cdot\left(\nabla\phi\otimes\nabla\phi\right)
-\nabla\left(\dfrac{|\nabla\phi|^2}{2}\right)=-\nabla\widehat{p},
	\end{align}
	for some $\widehat{p}\in L^\infty\left([0, T_0]; H^{s+1}\right)\cap L^2\left([0, T_0]; H^{s+2}\right)$, then we have
	\begin{align}\label{q.3.86}
		\dfrac{1}{\varepsilon_j^2\rho^{\varepsilon_j}}\nabla\left(P\left(\rho^{\varepsilon_j}\right)\right)\to\nabla\widehat{p},\quad {\rm weakly^*\ in}\  L^\infty\left([0, T_0]; H^{s}\right)\cap L^2\left([0, T_0]; H^{s+1}\right).
		\end{align}

Taking $p=\widehat{p}-\dfrac{|\nabla\phi|^2}{2}$, we see that $(\ref{q.1.2})_2$ follows from (\ref{q.3.85}) directly. Actually, Lemma \ref{L3.4} ensures that the convergence is in fact valid for the sequence $\left(\rho^{\varepsilon},u^{\varepsilon},\phi^{\varepsilon}\right)$ themselves. The proof of Theorem \ref{th1} is completed.
\hfill$\Box$

	\begin{remark}\label{R3.1}
	It follows from (\ref{q.2.7}) and (\ref{q.3.86}) that
	\begin{align}\label{q.3.87}
		\left\|\dfrac{1}{\varepsilon^2}\nabla\rho^{\varepsilon}\right\|_{s}+\int_{0}^{t}\left\|\dfrac{1}{\varepsilon^2}\nabla\rho^{\varepsilon}\right\|_{s+1}^2\leq C,\quad t\in[0,T_0].
	\end{align}
\end{remark}
	\begin{remark}\label{R3.2}
It follows from (\ref{q.1.1}) and (\ref{q.2.7}) that
	\begin{align}\label{q.3.88}
		\left\|\dfrac{1}{\varepsilon}\nabla\cdot u^{\varepsilon}\right\|_{s+1}&\leq\left\|\dfrac{1}{\varepsilon}\rho_t^{\varepsilon}\right\|_{s+1}
+\left\|\dfrac{1}{\varepsilon}\left(u^{\varepsilon}\cdot\nabla\right)\rho^{\varepsilon}\right\|_{s+1}
+\left\|\dfrac{1}{\varepsilon}\left(\rho^{\varepsilon}-1\right)\nabla\cdot u^{\varepsilon}\right\|_{s+1}  \nonumber\\
		&\leq C,\quad t\in[0,T_0].
\end{align}
\end{remark}

\subsection{Proof of Theorem \ref{th2}}\label{subsec:2.3}
This subsection is devoted to get a priori estimates (\ref{q.2.11}) and (\ref{q.2.12}) with small initial displacements and small initial data. Then together with Theorem \ref{th1}, Theorem \ref{th2} can be proved by a standard procedure.

 First, we assume that
		\begin{align}\label{q.4.1}
		E_{s+2}\left(\rho,u\right)+F_{s+1}(\nabla\phi)&+\left\|\phi^2-1\right\|^2\nonumber\\
&\quad+\int_0^t\left(\nu_*\left\|\nabla u\right\|^2_{s+2}+\eta_*\left\|\nabla\cdot u\right\|^2_{s+2}+\left\|\nabla\phi\right\|^2_{s+3}\right)\leq 4\Theta
\end{align}
	for $t\in [ 0, T^\varepsilon ] $, and then what we need to do is to prove the following desired estimate
	\begin{align}\label{q.4.2}
		E_{s+2}\left(\rho,u\right)+F_{s+1}(\nabla\phi)&+\left\|\phi^2-1\right\|^2\nonumber\\
&\quad+\int_0^t\left(\nu_*\left\|\nabla u\right\|^2_{s+2}+\eta_*\left\|\nabla\cdot u\right\|^2_{s+2}+\left\|\nabla\phi\right\|^2_{s+3}\right)\leq 3\Theta
	\end{align}
	for $t\in [ 0,T^\varepsilon ] $ and $\Theta=\left(\delta+\varepsilon^2\kappa_0^2\right)$.
Then (\ref{q.2.11}) follows by the standard continuity argument and the fact  $E_{s+2}\left(\rho_0,u_0\right)+F_{s+1}(\nabla\phi_0)+\left\|\phi_0^2-1\right\|^2<4\left(\delta+\varepsilon^2\kappa_0^2\right)$.

	 Firstly, we are to prove (\ref{q.2.12}) under the assumptions of (\ref{q.4.1}).
We get back to (\ref{q.3.14})-(\ref{qq.3.23}), (\ref{qq.3.19}), (\ref{q.3.21})-(\ref{q.3.23}), (\ref{q.3.25})-(\ref{q.3.30}) and (\ref{qq.3.45})-(\ref{qq.3.56}), replace $(\xi,v,\psi)$ by $(\rho,u,\phi)$, make use of (\ref{q.4.1}) and the Cauchy inequality, and then make a similar (just a little different) argument to give the following fact:
\begin{align}\label{q.4.3}
		&\sum_{|\beta|\le s+1,|\alpha_1|\le s} \dfrac{\rm d}{{\rm d}t} \int_{\mathbb{T}^N}\left(\dfrac{P'\left(\rho\right)}{\rho}\left|\dfrac{1}{\varepsilon}\nabla^\beta\rho_t\right|^2+\rho|\nabla^\beta u_t|^2+\left|\nabla^{\alpha_1}\phi_t\right|^2\right)\nonumber\\
		&\quad+\sum_{|\beta|\le s+1}\nu_*\int_{\mathbb{T}^N}|\nabla^\beta \nabla u_t|^2+\sum_{|\beta|\le s+1}\eta_*\int_{\mathbb{T}^N}|\nabla^\beta\left(\nabla\cdot u_t\right)|^2
		+\sum_{|\alpha_1|\le s+1}\int_{\mathbb{T}^N}|\nabla^{\alpha_1}\nabla\phi_t|^2  \nonumber\\
		&\leq C\left(\left\|\nabla u\right\|_{s+2}^2+\left\|\nabla\cdot u\right\|_{s+2}^2+\left\|\nabla\phi\right\|_{s+2}^2+1\right)\left(\left\|\dfrac{1}{\varepsilon}\rho_t\right\|_{s+1}^2+\left\|u_t\right\|_{s+1}^2
+\left\|\phi_t\right\|_{s}^2\right),
\end{align}
	provided that $\delta$ and $\varepsilon$ are both small enough. Here we have used the fact $\left\|\rho_t\right\|_\infty\leq\left\|\nabla\cdot\left(\rho u\right)\right\|_\infty\leq C$ inferred from $(\ref{q.4.1})$ to estimate $I_1$. Using the Gronwall inequality, (\ref{q.2.3}), (\ref{q.3.51}) and (\ref{q.4.1}), we obtain
	\begin{align}\label{q.4.4}
		&\left\|\dfrac{1}{\varepsilon}\rho_t\right\|_{s+1}^2+\left\|u_t\right\|_{s+1}^2+\left\|\phi_t\right\|_{s}^2+\int_{0}^{t}\left(\nu_*\left\|\nabla u_t\right\|_{s+1}^2+\eta_*\left\|\nabla\cdot u_t\right\|_{s+1}^2+\left\|\nabla\phi_t\right\|_{s+1}^2\right) \nonumber\\
		&\quad\leq C\left(1+t\right){\rm exp}Ct\leq C {\rm exp} Ct,
	\end{align}
	for $t\in[0,T^\varepsilon]$, and thus we conclude that Remarks \ref{R3.1} and \ref{R3.2} hold for $t\in [0,T^\varepsilon]$ by the same discussion in subsection 2.2.

On one hand, multiplying $(\ref{q.1.1})_3$ by $\mu$ and using $(\ref{q.1.1})_4$, and then integrating over $\mathbb{T}^N$, one obtains after using integration by parts and (\ref{q.4.1}) that
		\begin{align}\label{qq.4.5}
&\dfrac {\rm d}{{\rm d}t}\int_{\mathbb{T}^N}\left(\dfrac{1}{2}\left|\nabla\phi\right|^2+\dfrac{1}{4}\rho\left|\phi^2-1\right|^2\right)
+\int_{\mathbb{T}^N}\left|\nabla\mu\right|^2
\nonumber\\
&=\int_{\mathbb{T}^N}\left(u\cdot\nabla\right)\phi\Delta\phi\le \|u\|_{\infty}\|\nabla\phi\|\|\Delta\phi\|\le 8\left(\delta+\varepsilon^2\kappa_0^2\right)^{\frac32}.
		\end{align}
Integrating the above result over $[0,t]$ and taking $\delta$ and $\varepsilon$ are both small enough, we have
\begin{align}\label{qq.4.6}
\|\nabla\phi\|^2+\|\phi^2-1\|^2\le 3\left(\delta+\varepsilon^2\kappa_0^2\right),
\end{align}
for $t\in[0,T^\varepsilon]$.

On the other hand, combining $(\ref{q.1.1})_3$ with $(\ref{q.1.1})_4$, we can rewrite the equation as
\begin{align}\label{qq.4.7}
\phi_t+\frac1{\rho}\Delta\left(\frac{\Delta\phi}{\rho}\right)-\frac2{\rho}\Delta\phi
=\frac3{\rho}\left(\phi^2-1\right)\Delta\phi-u\cdot\nabla\phi+\frac{6\phi}{\rho}|\nabla\phi|^2.
\end{align}
Taking $\nabla^\alpha$ to (\ref{qq.4.7}), multiplying the result by $\nabla^\alpha \Delta\phi$ with $|\alpha|\le s+1$, we get after using integration by parts that
\begin{align}\label{qq.4.8}
&\frac12\dfrac {\rm d}{{\rm d}t}\int_{\mathbb{T}^N}|\nabla^\alpha\nabla\phi|^2
+\int_{\mathbb{T}^N}\left|\nabla^\alpha\nabla\left(\frac{\Delta\phi}{\rho}\right)\right|^2
+\int_{\mathbb{T}^N}\frac2{\rho}\left|\nabla^\alpha\Delta\phi\right|^2\nonumber
\\
&=-\int_{\mathbb{T}^N}\nabla^\alpha\Delta\left(\frac{\Delta\phi}{\rho}\right)
\left[\nabla^\alpha\left(\frac1{\rho}\Delta\phi\right)-\frac1{\rho}\nabla^\alpha\Delta\phi\right]
\nonumber\\
&\quad+\int_{\mathbb{T}^N}\left[\nabla^\alpha\left(\frac1{\rho}\Delta\left(\frac{\Delta\phi}{\rho}\right)\right)
-\frac1{\rho}\nabla^\alpha\Delta\left(\frac{\Delta\phi}{\rho}\right)\right]\nabla^\alpha\Delta\phi
\nonumber\\
&\quad-\int_{\mathbb{T}^N}\left[\nabla^\alpha\left(\frac2{\rho}\Delta\phi\right)-\frac2{\rho}\nabla^\alpha\Delta\phi\right]\nabla^\alpha\Delta\phi
-\int_{\mathbb{T}^N}\nabla^\alpha\left[\frac3{\rho}\left(\phi^2-1\right)\Delta\phi\right]\nabla^\alpha\Delta\phi
\nonumber\\
&\quad+\int_{\mathbb{T}^N}\nabla^\alpha\left(u\cdot\nabla\phi\right)\nabla^\alpha\Delta\phi
-\int_{\mathbb{T}^N}\nabla^\alpha\left(\frac{6\phi}{\rho}|\nabla\phi|^2\right)\nabla^\alpha\Delta\phi=\sum_{i=1}^{6}Q_i,
\end{align}
where we have used the following fact from integration by parts that
\begin{align*}
&\int_{\mathbb{T}^N}\nabla^\alpha\left(\frac1{\rho}\Delta\left(\frac{\Delta\phi}{\rho}\right)\right)\nabla^\alpha \Delta\phi\\
&=\int_{\mathbb{T}^N}\left[\nabla^\alpha\left(\frac1{\rho}\Delta\left(\frac{\Delta\phi}{\rho}\right)\right)
-\frac1{\rho}\nabla^\alpha\Delta\left(\frac{\Delta\phi}{\rho}\right)\right]\nabla^\alpha\Delta\phi
+\int_{\mathbb{T}^N}\frac1{\rho}\nabla^\alpha\Delta\left(\frac{\Delta\phi}{\rho}\right)\nabla^\alpha\Delta\phi\\
&=\int_{\mathbb{T}^N}\left[\nabla^\alpha\left(\frac1{\rho}\Delta\left(\frac{\Delta\phi}{\rho}\right)\right)
-\frac1{\rho}\nabla^\alpha\Delta\left(\frac{\Delta\phi}{\rho}\right)\right]\nabla^\alpha\Delta\phi\\
&\quad-\int_{\mathbb{T}^N}\nabla^\alpha\Delta\left(\frac{\Delta\phi}{\rho}\right)
\left[\nabla^\alpha\left(\frac1{\rho}\Delta\phi\right)-\frac1{\rho}\nabla^\alpha\Delta\phi\right]
+\int_{\mathbb{T}^N}\nabla^\alpha\Delta\left(\frac{\Delta\phi}{\rho}\right)
\nabla^\alpha\left(\frac1{\rho}\Delta\phi\right)\\
&=\int_{\mathbb{T}^N}\left[\nabla^\alpha\left(\frac1{\rho}\Delta\left(\frac{\Delta\phi}{\rho}\right)\right)
-\frac1{\rho}\nabla^\alpha\Delta\left(\frac{\Delta\phi}{\rho}\right)\right]\nabla^\alpha\Delta\phi\\
&\quad-\int_{\mathbb{T}^N}\nabla^\alpha\Delta\left(\frac{\Delta\phi}{\rho}\right)
\left[\nabla^\alpha\left(\frac1{\rho}\Delta\phi\right)-\frac1{\rho}\nabla^\alpha\Delta\phi\right]
-\int_{\mathbb{T}^N}\left|\nabla^\alpha\nabla\left(\frac{\Delta\phi}{\rho}\right)\right|^2.
\end{align*}
Next, by the $\rm Poincar\acute{e}$ inequality, Lemma \ref{L3.1} and integration by parts, we give the estimates of $Q_i~(i=1, \cdots, 6)$ as follows:
\begin{align}
Q_1
&=\int_{\mathbb{T}^N}\nabla^\alpha\nabla\left(\frac{\Delta\phi}{\rho}\right)\cdot
\left[\nabla^\alpha\nabla\left(\frac1{\rho}\Delta\phi\right)-\nabla\left(\frac1{\rho}\nabla^\alpha\Delta\phi\right)\right]\nonumber\\
&= \int_{\mathbb{T}^N}\nabla^\alpha\nabla\left(\frac{\Delta\phi}{\rho}\right)\cdot
\left\{\nabla^\alpha\left[\nabla\left(\frac1{\rho}\right)\Delta\phi\right]-\nabla\left(\frac1{\rho}\right)\nabla^\alpha\Delta\phi
+\nabla^\alpha\left(\frac1{\rho}\nabla\Delta\phi\right)-\frac1{\rho}\nabla^\alpha\nabla\Delta\phi\right\}\nonumber\\
&\le \left\|\nabla^\alpha\nabla\left(\frac{\Delta\phi}{\rho}\right)\right\|\left( \left\|\nabla^2\left(\frac1{\rho}\right)\right\|_\infty\|\Delta\phi\|_{s}+\left\|\Delta\phi\right\|_\infty\left\|\nabla^2\left(\frac1{\rho}\right)\right\|_{s}
\right.\nonumber\\
        &\left.
\quad+\left\|\nabla\left(\frac1{\rho}\right)\right\|_\infty\|\nabla\Delta\phi\|_{s}+\|\nabla\Delta\phi\|_\infty\left\|\nabla\left(\frac1{\rho}\right)\right\|_{s}\right)
\nonumber\\
&\le \left\|\nabla\left(\frac{\Delta\phi}{\rho}\right)\right\|_{s+1}\left(\left\|\nabla^2\left(\frac1{\rho}\right)\right\|_{s}\|\Delta\phi\|_{s}
+\left\|\nabla\left(\frac1{\rho}\right)\right\|_{s}\|\nabla\Delta\phi\|_{s}\right)
\nonumber\\
&\le \tau \left\|\nabla\left(\frac{\Delta\phi}{\rho}\right)\right\|_{s+1}^2+ C(\tau)\varepsilon^2 \left\|\dfrac1{\varepsilon}\nabla\rho\right\|_{s+1}^2\|\Delta\phi\|_{s+1}^2\nonumber\\
&\le \tau \left\|\nabla\left(\frac{\Delta\phi}{\rho}\right)\right\|_{s+1}^2+C(\tau) \varepsilon^2 \left(\delta+\varepsilon^2\kappa_0^2\right)\|\Delta\phi\|_{s+1}^2,\label{qq.4.9}\\
%
%
Q_2&\le\left[ \left\|\nabla\left(\frac1{\rho}\right)\right\|_\infty\left\|\Delta\left(\frac{\Delta\phi}{\rho}\right)\right\|_s
+\left\|\Delta\left(\frac{\Delta\phi}{\rho}\right)\right\|_\infty\left\|\nabla\left(\frac1{\rho}\right)\right\|_s\right]\|\Delta\phi\|_{s+1}
\nonumber\\
&\le \tau \left\|\Delta\left(\frac{\Delta\phi}{\rho}\right)\right\|_s^2+ C(\tau)\varepsilon^2 \left\|\dfrac1{\varepsilon}\nabla\rho\right\|_{s}^2\|\Delta\phi\|_{s+1}^2\nonumber\\
&\le \tau \left\|\Delta\left(\frac{\Delta\phi}{\rho}\right)\right\|_s^2+ C(\tau)\varepsilon^2 \left(\delta+\varepsilon^2\kappa_0^2\right)\|\Delta\phi\|_{s+1}^2,\label{qq.4.10}\\
Q_3&\le \left[\left\|\nabla\left(\frac2{\rho}\right)\right\|_\infty\|\Delta\phi\|_s+\|\Delta\phi\|_\infty\left\|\nabla\left(\frac2{\rho}\right)\right\|_s\right]\|\Delta\phi\|_{s+1}
\nonumber\\
&\le \left\|\nabla\left(\frac2{\rho}\right)\right\|_s\|\Delta\phi\|_{s}\|\Delta\phi\|_{s+1}\le \tau \|\Delta\phi\|_{s+1}^2+ C(\tau)\varepsilon^2 \left\|\dfrac1{\varepsilon}\nabla\rho\right\|_{s}^2\|\Delta\phi\|_{s}^2
\nonumber\\
&\le \tau \|\Delta\phi\|_{s+1}^2+ C(\tau)\varepsilon^2 \left(\delta+\varepsilon^2\kappa_0^2\right)\|\Delta\phi\|_{s}^2,\label{qq.4.11}\\
Q_4&\le \left\|\frac3{\rho}\right\|_{s+1}\left\|\phi^2-1\right\|_{s+1}\left\|\Delta\phi\right\|_{s+1}\left\|\Delta\phi\right\|_{s+1}\nonumber\\
&\le \left\|\frac3{\rho}\right\|_{s+1}\left(\left\|\phi^2-1\right\|+\left\|\nabla\left(\phi^2-1\right)\right\|_{s}\right)\left\|\Delta\phi\right\|_{s+1}^2\nonumber\\
&\le C\left(\left\|\phi^2-1\right\|+2\left\|\phi\right\|_{s}\left\|\nabla\phi\right\|_{s}\right)\left\|\Delta\phi\right\|_{s+1}^2\le C\left(\delta+\varepsilon^2\kappa_0^2\right)^\frac12 \left\|\Delta\phi\right\|_{s+1}^2,\label{qq.4.12}\\
Q_5&\le \left(\left\|u\right\|_{s+1}\left\|\nabla\phi\right\|_\infty+\left\|\nabla\phi\right\|_{s+1}\left\|u\right\|_\infty\right)\left\|\Delta\phi\right\|_{s+1}\nonumber\\
&\le C\left\|u\right\|_{s+1}\left\|\nabla\phi\right\|_{s+1}\left\|\Delta\phi\right\|_{s+1}\nonumber\\
&\le C(\tau)\left\|u\right\|_{s+1}^2\left\|\nabla\phi\right\|_{s+1}^2+\tau\left\|\Delta\phi\right\|_{s+1}^2\le C(\tau)\left(\delta+\varepsilon^2\kappa_0^2\right)\left\|\Delta\phi\right\|_{s+1}^2+\tau\left\|\Delta\phi\right\|_{s+1}^2,\label{qq.4.13}\\
Q_6&\le C\left(\left\|\dfrac{\phi}{\rho}\right\|_\infty\left\||\nabla\phi|^2\right\|_{s+1}+\left\||\nabla\phi|^2\right\|_\infty\left\|\dfrac{\phi}{\rho}\right\|_{s+1}\right)
\left\|\Delta\phi\right\|_{s+1}\nonumber\\
&\le C\left[\left\|\dfrac{1}{\rho}\right\|_\infty\left\|\phi\right\|_\infty\left\|\nabla\phi\right\|_\infty\left\|\nabla\phi\right\|_{s+1}
+\left\|\nabla\phi\right\|_\infty^2\left(\left\|\dfrac1{\rho}\right\|_{\infty}\left\|\phi\right\|_{s+1}
+\left\|\phi\right\|_{\infty}\left\|\dfrac1{\rho}\right\|_{s+1}\right)\right]
\left\|\Delta\phi\right\|_{s+1}\nonumber\\
&\le C\left\|\dfrac1\rho\right\|_{s+1}\left\|\phi\right\|_{s+1}\left\|\nabla\phi\right\|_{s+1}^2\left\|\Delta\phi\right\|_{s+1}
\le C\left(\delta+\varepsilon^2\kappa_0^2\right)^\frac12 \left\|\Delta\phi\right\|_{s+1}^2,\label{qq.4.14}
\end{align}
where in the last inequality we have used the fact that
\begin{align*}
\|\phi\|^2=\int_{\mathbb{T}^N}(\phi^2-1+1)\le \|1\|\left\|\phi^2-1\right\|+\left|\mathbb{T}^N\right|\le C.
\end{align*}
Putting (\ref{qq.4.8})-(\ref{qq.4.14}) together, summing over $\alpha$ and taking $\tau$, $\varepsilon$, $\delta$ small enough we get
\begin{align*}
&\sum_{\alpha\le s+1}\dfrac {\rm d}{{\rm d}t}\int_{\mathbb{T}^N}|\nabla^\alpha\nabla\phi|^2
+\sum_{\alpha\le s+1}\int_{\mathbb{T}^N}\left|\nabla^\alpha\nabla\left(\frac{\Delta\phi}{\rho}\right)\right|^2
+\sum_{\alpha\le s+1}\int_{\mathbb{T}^N}\left|\nabla^\alpha\Delta\phi\right|^2\le 0.
\end{align*}
Then integrating the above inequality over $[0,t]$ and by (\ref{q.2.4}), (\ref{q.2.6}) and (\ref{q.2.10}), one obtains
 \begin{align}\label{qq.4.16}
\left\|\nabla\phi\right\|_{s+1}^2
+\int_0^t\left\|\nabla\left(\dfrac{\Delta\phi}{\rho}\right)\right\|_{s+1}^2+\int_0^t\left\|\nabla\phi\right\|_{s+2}^2\le \left\|\nabla\phi_0\right\|_{s+1}^2+\varepsilon^2\kappa_0^2
\le \delta+\varepsilon^2\kappa_0^2.
\end{align}
Note that by (\ref{q.4.1}) we have
\begin{align}\label{qq.4.17}
\left\|\Delta^2\phi\right\|_s=\left\|\Delta\left(\rho\dfrac{\Delta\phi}{\rho}\right)\right\|_s
&\le\left\|\Delta\rho \left(\frac{\Delta\phi}{\rho}\right) \right\|_s
+C\left\|\nabla\rho\cdot\nabla\left(\dfrac{\Delta\phi}{\rho}\right)\right\|_s
+\left\|\rho \Delta\left(\frac{\Delta\phi}{\rho}\right) \right\|_s\nonumber\\
&\le\left\|\Delta\rho \right\|_s\left\| \frac{\Delta\phi}{\rho} \right\|_s
+C\left\|\nabla\rho\right\|_s\left\|\nabla\left(\dfrac{\Delta\phi}{\rho}\right)\right\|_s
+C\left\|\rho\right\|_s\left\| \Delta\left(\frac{\Delta\phi}{\rho}\right) \right\|_s\nonumber\\
&\le C\left(\left\| \frac{\Delta\phi}{\rho} \right\|+\left\| \nabla\left(\frac{\Delta\phi}{\rho}\right) \right\|_{s-1}\right)
+C\left\|\nabla\left(\dfrac{\Delta\phi}{\rho}\right)\right\|_s
+C\left\| \Delta\left(\frac{\Delta\phi}{\rho}\right) \right\|_s\nonumber\\
&\le C\left(\left\|\Delta\phi\right\|+\left\|\nabla\left(\dfrac{\Delta\phi}{\rho}\right)\right\|_{s+1}\right).
\end{align}
In conclusion, collecting (\ref{qq.4.16})-(\ref{qq.4.17}), we derive
\begin{align}\label{qq.4.19}
\left\|\nabla\phi\right\|_{s+1}^2
+\int_0^t\left\|\nabla\phi\right\|_{s+3}^2
\le \delta+\varepsilon^2\kappa_0^2,
\end{align}
for $t\in[0,T^\varepsilon]$.

 In addition, recalling (\ref{q.3.20}) and replacing $(\xi,v)$ by $(\rho,u)$, we derive
	\begin{align}\label{q.4.22}
		&\dfrac{1}{2}\dfrac{\rm d}{{\rm d}t}\int_{\mathbb{T}^N}\left(\dfrac{1}{\varepsilon^2}\dfrac{P'\left(\rho\right)}{\rho}\left|D^{\alpha_1}\left(\rho-1\right)\right|^2
+\rho\left|D^{\alpha_1}u\right|^2\right)\nonumber\\
		&\quad+\int_{\mathbb{T}^N}\nu(\rho,\phi)\left|\nabla D^{\alpha_1}u\right|^2+\int_{\mathbb{T}^N}\eta(\rho,\phi)\left|\nabla\cdot D^{\alpha_1}u\right|^2=\sum_{k=1}^{9}J_k,
	\end{align}
	where there are some slight changes from $I_i$ as follows:
		\begin{align*}
		J_1&=\dfrac{1}{2\varepsilon^2}\int_{\mathbb{T}^N}\dfrac{2P'\left(\rho\right)-P''\left(\rho\right)\rho}{\rho}\nabla\cdot u|D^{\alpha_1}\left(\rho-1\right)|^2,\\
		J_2&=\dfrac{1}{\varepsilon^2}\int_{\mathbb{T}^N}P''\left(\rho\right)D^{\alpha_1}\left(\rho-1\right)D^{\alpha_1}u\cdot\nabla\rho,\\
		J_3&=-\dfrac{1}{\varepsilon^2}\int_{\mathbb{T}^N}\dfrac{P'\left(\rho\right)}{\rho}D^{\alpha_1}\left(\rho-1\right) \left\{\left[D^{\alpha_1}\left(u\cdot\nabla\rho\right)-u\cdot\nabla D^{\alpha_1}\rho \right] \nonumber \right.\\
		&\left.\quad+\left[D^{\alpha_1}\left(\rho\nabla\cdot u\right)-\rho\nabla\cdot D^{\alpha_1}u\right]\right\},\\
		J_4&=-\int_{\mathbb{T}^N}D^{\alpha_1}\left(\Delta\phi\nabla\phi\right)\cdot D^{\alpha_1}u,\\
		J_5&=-\int_{\mathbb{T}^N}\rho\left[D^{\alpha_1}\left(u\cdot\nabla u\right)-u\cdot\nabla D^{\alpha_1}u\right]\cdot D^{\alpha_1}u,\\
		J_6&=-\dfrac{1}{\varepsilon^2}\int_{\mathbb{T}^N}\rho\left[D^{\alpha_1}\left(\dfrac{P'\left(\rho\right)}{\rho}\nabla\rho\right)-\dfrac{P'\left(\rho\right)}{\rho}\nabla D^{\alpha_1}\rho\right]\cdot D^{\alpha_1}u,\\
		J_7&=\int_{\mathbb{T}^N}\rho\left[D^{\alpha_1}\left(\dfrac{\nu(\rho,\phi)}{\rho}\Delta u\right)-\dfrac{\nu(\rho,\phi)}{\rho}\Delta D^{\alpha_1}u\right]\cdot D^{\alpha_1}u\nonumber \\
		&\quad+\int_{\mathbb{T}^N}\rho\left\{D^{\alpha_1}\left[\dfrac{\eta(\rho,\phi)}{\rho}\nabla\left(\nabla\cdot u\right)\right]-\dfrac{\eta(\rho,\phi)}{\rho}\nabla D^{\alpha_1}(\nabla\cdot u)\right\}\cdot D^{\alpha_1}u,\\
		J_8&=-\int_{\mathbb{T}^N}\rho\left\{D^{\alpha_1}\left[\dfrac{1}{\rho}\left(\Delta\phi\nabla\phi\right)\right]
-\dfrac{1}{\rho}D^{\alpha_1}\left(\Delta\phi\nabla\phi\right)\right\}\cdot D^{\alpha_1}u,\\
J_9&=-\int_{\mathbb{T}^N}\nabla\nu(\rho,\phi)\cdot\nabla D^{\alpha_1}u\cdot D^{\alpha_1}u-\int_{\mathbb{T}^N}\nabla\eta(\rho,\phi)\cdot D^{\alpha_1}u\cdot D^{\alpha_1}(\nabla\cdot u),
\end{align*}
	since we have used $(\ref{q.1.1})_1$ to give
	\begin{align*}	&\dfrac{1}{2\varepsilon^2}\int_{\mathbb{T}^N}|D^{\alpha_1}\left(\rho-1\right)|^2\partial_t\left(\dfrac{P'\left(\rho\right)}{\rho}\right)+\dfrac{1}{2\varepsilon^2}\int_{\mathbb{T}^N}|D^{\alpha_1}\left(\rho-1\right)|^2\nabla\cdot\left(\dfrac{P'\left(\rho\right)}{\rho}u\right)  \nonumber\\\nonumber
		&\quad=\dfrac{1}{2\varepsilon^2}\int_{\mathbb{T}^N}\left(\left(\dfrac{P'\left(\rho\right)}{\rho}\right)'\rho_t+\nabla\cdot\left(\dfrac{P'\left(\rho\right)}{\rho}u\right)\right)|D^{\alpha_1}\left(\rho-1\right)|^2 \nonumber\\\nonumber
		&\quad=\dfrac{1}{2\varepsilon^2}\int_{\mathbb{T}^N}\left(\left(\dfrac{P'\left(\rho\right)}{\rho}\right)'\rho_t+\left(\dfrac{P'\left(\rho\right)}{\rho}\right)'\nabla\rho\cdot u+\left(\dfrac{P'\left(\rho\right)}{\rho}\right)'\rho\nabla\cdot u\right)|D^{\alpha_1}\left(\rho-1\right)|^2 \\
		&\qquad+\dfrac{1}{2\varepsilon^2}\int_{\mathbb{T}^N}\left(\dfrac{P'\left(\rho\right)}{\rho}-\left(\dfrac{P'\left(\rho\right)}{\rho}\right)^{'}\rho\right)\nabla\cdot u|D^{\alpha_1}\left(\rho-1\right)|^2\nonumber\\\nonumber
		&\quad=\dfrac{1}{2\varepsilon^2}\int_{\mathbb{T}^N}\left(\dfrac{P'\left(\rho\right)}{\rho}-\left(\dfrac{P'\left(\rho\right)}{\rho}\right)^{'}\rho\right)\nabla\cdot u|D^{\alpha_1}\left(\rho-1\right)|^2\nonumber\\
		&\quad=\dfrac{1}{2\varepsilon^2}\int_{\mathbb{T}^N}\dfrac{2P'\left(\rho\right)-P''\left(\rho\right)\rho}{\rho}\nabla\cdot u|D^{\alpha_1}\left(\rho-1\right)|^2.
	\end{align*}

 Now we give the dispersive estimates about $J_k$ for $k=1,2,...,9$, when $D^{\alpha_1}=\nabla^{\alpha_1}$ with $|\alpha_1|\leq s+2,\ s\geq 3$.\\
First of all, it follows from  Remark \ref{R3.2} that
		\begin{align}\label{q.4.32}	
		|J_1|&\leq C\varepsilon\left\|\dfrac{1}{\varepsilon}\nabla\cdot u\right\|_\infty\left\|\dfrac{1}{\varepsilon}\nabla^{\alpha_1}\left(\rho-1\right)\right\|^2\leq C\varepsilon\left\|\dfrac{1}{\varepsilon}\nabla^{\alpha_1}(\rho-1)\right\|^2.
\end{align}
Next, using the ${\rm H\ddot{o}lder}$ inequality, the Sobolev embedding and Remark \ref{R3.1}, we have
\begin{align}\label{q.4.33}	
		|J_2|&\leq C\dfrac{1}{\varepsilon}\left\|\nabla\rho\right\|_{L^4}\left\|\nabla^{\alpha_1}u\right\|_{L^4}\left\|\dfrac{1}{\varepsilon}\nabla^{\alpha_1}\left(\rho-1\right)\right\|\nonumber \\
		&\leq C\varepsilon\left\|\dfrac{1}{\varepsilon^2}\nabla\rho\right\|_{s}\left(\left\|\nabla^{\alpha_1}u\right\|+\left\|\nabla\nabla^{\alpha_1}u\right\|\right)\left\|\dfrac{1}{\varepsilon}\nabla^{\alpha_1}\left(\rho-1\right)\right\|\nonumber \\
		&\leq\tau\left\|\nabla\nabla^{\alpha_1}u\right\|^2+C\left(\tau\right)\varepsilon\left(\left\|\nabla^{\alpha_1}u\right\|^2+\left\|\dfrac{1}{\varepsilon}\nabla^{\alpha_1}\left(\rho-1\right)\right\|^2\right).
\end{align}
Then by the ${\rm H\ddot{o}lder}$ inequality, the Cauchy inequality, the $\rm Poincar\acute{e}$ inequality and Lemma \ref{L3.1}, we derive
\begin{align}
		|J_3|&\leq C\dfrac{1}{\varepsilon}\left\|\dfrac{1}{\varepsilon}\nabla^{\alpha_1}\left(\rho-1\right)\right\|\left(\left\|\nabla u\right\|_\infty\left\|\nabla\rho\right\|_{s+1}+\left\|\nabla u\right\|_{s+1}\left\|\nabla\rho\right\|_\infty\nonumber \right.\\
		&\left.\quad+\left\|\nabla\cdot u\right\|_\infty\left\|\nabla\rho\right\|_{s+1}+\left\|\nabla\cdot u\right\|_{s+1}\left\|\nabla\rho\right\|_\infty\right)\nonumber \\
		&\leq C\varepsilon\left\|\dfrac{1}{\varepsilon}\nabla^{\alpha_1}\left(\rho-1\right)\right\|\left\|u\right\|_{s+2}\left\|\dfrac{1}{\varepsilon^2}\nabla\rho\right\|_{s+1}\nonumber \\
		&\leq C\varepsilon\left\|\dfrac{1}{\varepsilon}\left(\rho-1\right)\right\|_{s+2}^2+C\varepsilon\left\|\dfrac{1}{\varepsilon^2}\nabla\rho\right\|_{s+1}^2,\label{q.4.34}	\\
		|J_4|&=\left |\int_{\mathbb{T}^N}\nabla^{\alpha_1}\left(\nabla\phi\otimes\nabla\phi-\dfrac{|\nabla\phi|^2}{2}I_N \right):\nabla\nabla^{\alpha_1}u\right |\nonumber \\
		&\leq\tau\left\|\nabla\nabla^{\alpha_1}u\right\|^2+C\left(\tau\right)\left\|\nabla\phi\right\|_\infty^2\left\|\nabla\phi\right\|_{s+2}^2\nonumber \\
		&\leq\tau\left\|\nabla u\right\|_{s+2}^2+C\left(\tau\right)\left(\delta+\varepsilon^2\kappa_0^2\right)\left\|\nabla\phi\right\|_{s+2}^2,\\
		|J_5|&\leq C\left\|\nabla^{\alpha_1}u\right\|\left\|\nabla u\right\|_\infty\left\|\nabla u\right\|_{s+1}\leq C\left(\delta^{\frac{1}{2}}+\varepsilon \kappa_0\right)\left\|\nabla u\right\|_{s+2}^2,\\
		|J_6|&\leq C\dfrac{1}{\varepsilon^2}\left\|\nabla^{\alpha_1}u\right\|\left(\left\|\nabla\rho\right\|_\infty\left\|\nabla\rho\right\|_{s+1}+\left\|\nabla\rho\right\|_{s+1}\left\|\nabla\rho\right\|_\infty\right)\nonumber \\
		&\leq C\dfrac{1}{\varepsilon^2}\left\|\nabla^{\alpha_1}u\right\|\left\|\nabla\rho\right\|_2\left\|\nabla\rho\right\|_{s+1}\nonumber \\
		&\leq C\varepsilon\left\|\nabla^{\alpha_1}u\right\|\left\|\dfrac{1}{\varepsilon^2}\nabla\rho\right\|_{s+1}\left\|\dfrac{1}{\varepsilon}\nabla\rho\right\|_{s+1}\nonumber \\
		&\leq C\varepsilon\left\|u\right\|_{s+2}^2+C\varepsilon\left\|\dfrac{1}{\varepsilon^2}\nabla\rho\right\|_{s+1}^2,\\
		|J_7|&\leq C\left\|\nabla^{\alpha_1}u\right\|\left\{\left\|\nabla\left(\dfrac{\nu(\rho,\phi)}{\rho}\right)\right\|_\infty\left\|\Delta u\right\|_{s+1}
+\left\|\nabla\left(\dfrac{\nu(\rho,\phi)}{\rho}\right)\right\|_{s+1}\left\|\Delta u\right\|_\infty\nonumber \right.\\
		&\left.\quad+\left\|\nabla\left(\dfrac{\eta(\rho,\phi)}{\rho}\right)\right\|_\infty\left\|\nabla\left(\nabla\cdot u\right)\right\|_{s+1}+\left\|\nabla\left(\dfrac{\eta(\rho,\phi)}{\rho}\right)\right\|_{s+1}\left\|\nabla\left(\nabla\cdot u\right)\right\|_\infty\right\}\nonumber \\
&\leq C\left\|\nabla^{\alpha_1}u\right\|\left(\left\|\Delta u\right\|_{s+1} +\left\|\nabla\left(\nabla\cdot u\right)\right\|_{s+1}\right) \left[(1+\|\rho\|_\infty^{s+1})\|\nabla\rho\|_{s+1}+(1+\|\phi\|_\infty^{s+1})\|\nabla\phi\|_{s+1}\right]\nonumber \\
		&\leq C\left\|u\right\|_{s+2}\left\|\nabla u\right\|_{s+2}\left(\varepsilon\left\|\dfrac{1}{\varepsilon}\nabla\rho\right\|_{s+1}+\|\nabla\phi\|_{s+1}\right)\nonumber \\
		&\leq C(\tau,\delta)\varepsilon\left\|u\right\|_{s+2}^2+\tau\left\|\nabla u\right\|_{s+2}^2
+C(\tau)\left(\delta+\varepsilon^2\kappa_0^2\right)\left\|\nabla\phi\right\|_{s+1}^2,\\
		|J_8|&\leq C\left\|\nabla^{\alpha_1}u\right\|\left(\left\|\nabla\rho\right\|_\infty\left\|\Delta\phi\nabla\phi\right\|_{s+1}+\left\|\nabla\rho\right\|_{s+1}\left\|\Delta\phi\right\|_\infty\left\|\nabla\phi\right\|_\infty\right)\nonumber \\
		&\leq C\left\|\nabla^{\alpha_1}u\right\|\left(\left\|\nabla\rho\right\|_2\left(\left\|\Delta\phi\right\|_\infty\left\|\nabla\phi\right\|_{s+1}+\left\|\Delta\phi\right\|_{s+1}\left\|\nabla\phi\right\|_\infty\right)+\left\|\Delta\phi\right\|_2\left\|\nabla\phi\right\|_2\left\|\nabla\rho\right\|_{s+1}\right)\nonumber \\
&\leq C\varepsilon\left\|u\right\|_{s+2}\left\|\dfrac1\varepsilon\nabla\rho\right\|_{s+1}\left\|\nabla\phi\right\|_{s+1}\left\|\Delta\phi\right\|_{s+1}\nonumber\\
		&\leq C\varepsilon\left\|u\right\|_{s+2}^2
+C\left(\delta+\varepsilon^2\kappa_0^2\right)\left\|\nabla\phi\right\|_{s+2}^2,\\
|J_9|&\leq\left[(1+\|\rho\|_\infty^{2})\|\nabla\rho\|_{2}+(1+\|\phi\|_\infty^{2})\|\nabla\phi\|_{2}\right]\left\|u\right\|_{s+2}(\left\|\nabla u\right\|_{s+2}+\left\|\nabla\cdot u\right\|_{s+2})\nonumber\\
&\leq C\varepsilon\left\|\dfrac1\varepsilon\nabla\rho\right\|_{2}\left\|u\right\|_{s+2}(\left\|\nabla u\right\|_{s+2}+\left\|\nabla\cdot u\right\|_{s+2})+C\|\nabla\phi\|_{2}\left\|u\right\|_{s+2}(\left\|\nabla u\right\|_{s+2}+\left\|\nabla\cdot u\right\|_{s+2})\nonumber\\
&\leq C(\tau,\delta)\varepsilon\left\|u\right\|_{s+2}^2+\tau(\left\|\nabla u\right\|_{s+2}^2+\left\|\nabla\cdot u\right\|_{s+2}^2)
+C(\tau)\left(\delta+\varepsilon^2\kappa_0^2\right)\left\|\nabla\phi\right\|_{2}^2.\label{q.4.39}
	\end{align}
 Hence, we put (\ref{q.4.32})-(\ref{q.4.39}) together, sum over $\alpha_1$, and then choose the parameter $\tau,\ \delta$, and $\varepsilon$ small enough to give the following dispersive energy estimate
\begin{align*}
&\dfrac{1}{2}\dfrac{\rm d}{{\rm d}t}\sum_{|\alpha_1|\leq s+2}\int_{\mathbb{T}^N}\left(\dfrac{1}{\varepsilon^2}\dfrac{P'\left(\rho\right)}{\rho}|\nabla^{\alpha_1}\left(\rho-1\right)|^2+\rho|\nabla^{\alpha_1}u|^2\right)\nonumber
\nonumber\\
&\quad+\nu_*\sum_{|\alpha_1|\leq s+2}\int_{\mathbb{T}^N}|\nabla^{\alpha_1}\nabla u|^2+\eta_*\sum_{|\alpha_1|\leq s+2}\int_{\mathbb{T}^N}|\nabla^{\alpha_1}(\nabla\cdot u)|^2 \nonumber\\
		&\leq C\varepsilon\left(\left\|u\right\|_{s+2}^2+\left\|\dfrac{1}{\varepsilon}(\rho-1)\right\|_{s+2}^2\right)
+C\varepsilon\left\|\dfrac{1}{\varepsilon^2}\nabla\rho\right\|_{s+1}^2
+C\left(\delta+\varepsilon^2\kappa_0^2\right)\left\|\nabla\phi\right\|_{s+2}^2.
\end{align*}
Then by the Gronwall inequality, Remark \ref{R3.1} and (\ref{qq.4.19}), one can derive (\ref{q.4.2}) for $t\in [ 0,T^\varepsilon ] $ with $T^\varepsilon={\varepsilon}^{\delta -1}$($\delta<1$ is a small positive constant), provided that $\varepsilon$ and $\delta$ are both small enough.

Finally, by (\ref{q.2.10}) and by (\ref{q.2.12}), we can proceed as in subsection \ref{subsec:2.2} by taking smooth test functions to prove that the limiting function $(1, u,\phi)$ and $\nabla p$ satisfy the incompressible system (\ref{q.1.2}) in the time interval $[0, T]$ with the initial data (\ref{q.2.9}) satisfying the constraints (\ref{q.2.5}) and (\ref{q.2.10}). Since $T$ is an arbitrary positive constant, we have in fact obtained a unique global strong solution
of (\ref{q.1.2}) and proved Theorem \ref{th2}.

\subsection{Proof\ of\ Theorem \ref{th3}}\label{subsec:2.4}
In this subsection, we show the convergence rates of $(\rho^\varepsilon, u^{\varepsilon}, \phi^\varepsilon)$ as ${\varepsilon}\to0$ by the modulated energy method with the help of uniform estimates (\ref{q.2.7}).

We rewrite the compressible Navier-Stokes/Cahn-Hilliard equation with the parameter $\varepsilon$ as follows
\allowdisplaybreaks
	\begin{align}\label{q.5.1}
		\begin{cases}
			\rho_t^\varepsilon+\nabla\cdot\left(\rho^\varepsilon u^\varepsilon\right)=0,\\
			\left(\rho^\varepsilon u^\varepsilon\right)_t+\nabla\cdot\left(\rho^\varepsilon u^\varepsilon\otimes u^\varepsilon\right)+\dfrac1{\varepsilon^2}\nabla P\left(\rho^\varepsilon\right)\\
\quad=\nu(\rho^\varepsilon,\phi^\varepsilon)\Delta u^\varepsilon+\eta(\rho^\varepsilon,\phi^\varepsilon)\nabla\left(\nabla\cdot u^\varepsilon\right)
			-\Delta\phi^\varepsilon\nabla\phi^\varepsilon,\\
\phi^\varepsilon_{t}+ u^\varepsilon\cdot\nabla\phi^\varepsilon=\dfrac1{\rho^\varepsilon}\Delta\mu^\varepsilon,
			\\
			\rho^\varepsilon\mu^\varepsilon
			=-\Delta\phi^\varepsilon+\rho^\varepsilon\left(\left(\phi^\varepsilon\right)^3-\phi^\varepsilon\right).\\	
		\end{cases}	
	\end{align}
 Multiplying $(\ref{q.5.1})_2$ by $u^\varepsilon$, and then integrating over $\mathbb{T}^N$ and using integration by parts, we get that
	\begin{align}\label{q.5.2}
		&\dfrac {\rm d}{{\rm d}t}\int_{\mathbb{T}^N}\left(\dfrac{1}{2}\rho^\varepsilon|u^\varepsilon|^2
+\dfrac{1}{\varepsilon^2}\omega\left(\rho^\varepsilon\right)\right)
+\int_{\mathbb{T}^N}\nu(\rho^\varepsilon,\phi^\varepsilon)|\nabla u^\varepsilon|^2
+\int_{\mathbb{T}^N}\eta(\rho^\varepsilon,\phi^\varepsilon) |\nabla\cdot u^\varepsilon|^2 \nonumber\\
		&=-\int_{\mathbb{T}^N}\left(u^\varepsilon\cdot\nabla\right)\phi^\varepsilon\Delta\phi^\varepsilon
-\int_{\mathbb{T}^N}\nabla\nu(\rho^\varepsilon,\phi^\varepsilon)\cdot\nabla u^\varepsilon\cdot u^\varepsilon
-\int_{\mathbb{T}^N}u^\varepsilon\cdot\nabla\eta(\rho^\varepsilon,\phi^\varepsilon) \left(\nabla\cdot u^\varepsilon\right) ,
	\end{align}
	where $\omega\left(\rho\right)=\rho\int_{1}^{\rho}\frac{P\left(z\right)}{z^2}{\rm d}z$.
 Multiplying $(\ref{q.5.1})_3$ by $\mu^\varepsilon$, and then using $(\ref{q.5.1})_4$ and integrating over  $\mathbb{T}^N$, one obtains after using integration by parts that
		\begin{align}\label{q.5.3}
\dfrac {\rm d}{{\rm d}t}\int_{\mathbb{T}^N}\left(\dfrac{1}{2}\left|\nabla\phi^\varepsilon\right|^2+\dfrac{1}{4}\rho^\varepsilon\left|\left(\phi^\varepsilon\right)^2-1\right|^2\right)
+\int_{\mathbb{T}^N}\left|\nabla\mu^\varepsilon\right|^2=\int_{\mathbb{T}^N}\left(u^\varepsilon\cdot\nabla\right)\phi^\varepsilon\Delta\phi^\varepsilon.
		\end{align}
 Adding (\ref{q.5.2}) and (\ref{q.5.3}) together, it holds
		\begin{align}\label{q.5.4}
		&\dfrac {\rm d}{{\rm d}t}\int_{\mathbb{T}^N}\left(\dfrac{1}{2}\rho^\varepsilon\left|u^\varepsilon\right|^2+\dfrac{1}{\varepsilon^2}\omega\left(\rho^\varepsilon\right)+\dfrac{1}{2}\left|\nabla\phi^\varepsilon \right|^2+\dfrac{1}{4}\rho^\varepsilon\left|\left(\phi^\varepsilon\right)^2-1\right|^2\right) \nonumber\\
		&\qquad+\int_{\mathbb{T}^N}\nu(\rho^\varepsilon,\phi^\varepsilon)\left|\nabla u^\varepsilon \right|^2+\int_{\mathbb{T}^N}\eta(\rho^\varepsilon,\phi^\varepsilon)|\nabla\cdot u^\varepsilon |^2
+\int_{\mathbb{T}^N}\left|\nabla\mu^\varepsilon\right|^2 \nonumber\\
&=-\int_{\mathbb{T}^N}\nabla\nu(\rho^\varepsilon,\phi^\varepsilon)\cdot\nabla u^\varepsilon\cdot u^\varepsilon
-\int_{\mathbb{T}^N}u^\varepsilon\cdot\nabla\eta(\rho^\varepsilon,\phi^\varepsilon) \left(\nabla\cdot u^\varepsilon\right) .
		\end{align}
 Let
		\begin{align*}
\Pi^\varepsilon\left(x,t\right)=\displaystyle{\frac{1}{\varepsilon^2}\left(\omega\left(\rho^\varepsilon\right)-P\left(1\right)\left(\rho^\varepsilon-1\right)\right)},
\end{align*}
it follows from the Taylor series and (\ref{q.2.6}) that
		\begin{align*}
	\Pi^\varepsilon\left(x,t\right)=\displaystyle{\frac{1}{\varepsilon^2}
\left(\frac{P'\left(1\right)}{2}\left(\rho^\varepsilon-1\right)^2\right)+\frac{1}{\varepsilon^2}o\left(\left(\rho^\varepsilon-1\right)^2\right)},
	\end{align*}
	then we get
	$\int_{\mathbb{T}^N}\Pi^\varepsilon\left(x,0\right)\leq C\varepsilon^2$.
Integrating (\ref{q.5.4}) over $[0,t]$ leads to
	\begin{align}\label{q.5.6}
		&\int_{\mathbb{T}^N}\left(\frac{1}{2}\rho^\varepsilon\left|u^\varepsilon \right|^2+\Pi^\varepsilon(x,t)+\frac{1}{2}\left|\nabla\phi^\varepsilon \right|^2+\frac{1}{4}\rho^\varepsilon\left|\left(\phi^\varepsilon\right)^2-1 \right|^2\right)\nonumber\\
		&\quad\quad+\int_{0}^{t}\int_{\mathbb{T}^N}\nu(\rho^\varepsilon,\phi^\varepsilon)\left|\nabla u^\varepsilon \right|^2+\int_{0}^{t}\int_{\mathbb{T}^N}\eta(\rho^\varepsilon,\phi^\varepsilon)\left|\nabla\cdot u^\varepsilon \right|^2+\int_{0}^{t}\int_{\mathbb{T}^N}\left|\nabla\mu^\varepsilon\right|^2\nonumber\\
		&=-\int_0^t\int_{\mathbb{T}^N}\nabla\nu(\rho^\varepsilon,\phi^\varepsilon)\cdot\nabla u^\varepsilon\cdot u^\varepsilon
-\int_0^t\int_{\mathbb{T}^N}u^\varepsilon\cdot\nabla\eta(\rho^\varepsilon,\phi^\varepsilon) \left(\nabla\cdot u^\varepsilon\right) \nonumber\\
&\quad+\int_{\mathbb{T}^N}\left(\frac{1}{2}\rho_0^\varepsilon\left|u_0^\varepsilon \right|^2
+\Pi^\varepsilon(x,0)+\frac{1}{2}\left|\nabla\phi_0^\varepsilon\right|^2
+\frac{1}{4}\rho_0^\varepsilon\left|\left(\phi_0^\varepsilon\right)^2-1\right|^2\right).
	\end{align}
	Similarly, from (\ref{q.1.2}), we get the basic energy law for the incompressible equation as
\begin{align}\label{q.5.7}
		&\int_{\mathbb{T}^N}\left(\frac{1}{2}\left|u\right|^2+\frac{1}{2}\left|\nabla\phi\right|^2
+\frac{1}{4}\left|\phi^2-1\right|^2\right)
+\int_{0}^{t}\int_{\mathbb{T}^N}\nu(\phi)\left|\nabla u\right|^2
+\int_0^t\int_{\mathbb{T}^N}\left|\nabla\mu\right|^2 \nonumber\\
		&\quad=-\int_0^t\int_{\mathbb{T}^N}\nabla\nu(\phi)\cdot\nabla u\cdot u
+\int_{\mathbb{T}^N}\left(\frac{1}{2}\left|u_0\right|^2+\frac{1}{2}\left|\nabla\phi_0\right|^2+\frac{1}{4}\left|\phi_0^2-1\right|^2\right).
	\end{align}
 Then multiplying $(\ref{q.5.1})_2$ by $u$ and integrating over $Q_t$, one derive after using integration by parts that
	%
		\begin{align}\label{q.5.8}
		&\int_{\mathbb{T}^N}\rho^\varepsilon u^\varepsilon\cdot u+\int_{0}^{t}\int_{\mathbb{T}^N}\rho^\varepsilon u^\varepsilon\cdot\left((u\cdot\nabla) u+\nabla p-\nu(\phi)\Delta u+{\rm {\rm div}}\left(\nabla\phi\otimes\nabla\phi\right)\right) \nonumber\\
		&\quad-\int_{0}^{t}\int_{\mathbb{T}^N}\left(\rho^\varepsilon u^\varepsilon\otimes u^\varepsilon\right):\nabla u
		+\int_{0}^{t}\int_{\mathbb{T}^N}\nu(\rho^\varepsilon,\phi^\varepsilon)\nabla u^\varepsilon:\nabla u+\int_{0}^{t}\int_{\mathbb{T}^N}\left(u\cdot\nabla\right)\phi^\varepsilon\Delta\phi^\varepsilon\nonumber\\
		&=-\int_{0}^{t}\int_{\mathbb{T}^N}\nabla\nu(\rho^\varepsilon,\phi^\varepsilon)\cdot\nabla u^\varepsilon\cdot u
-\int_{\mathbb{T}^N}u\cdot\nabla\eta(\rho^\varepsilon,\phi^\varepsilon)\left(\nabla\cdot u^\varepsilon\right)
+\int_{\mathbb{T}^N}\rho_0^\varepsilon u_0^\varepsilon\cdot u_0.
	\end{align}
	%
 Next, multiplying $(\ref{q.5.1})_3$ by $\Delta\phi$, and then integrating over $Q_t$, we obtain after using integration by parts that
	\begin{align}\label{q.5.9}
		&-\int_{\mathbb{T}^N}\nabla\phi^\varepsilon\cdot\nabla\phi
+\int_{0}^{t}\int_{\mathbb{T}^N}\Delta\phi^\varepsilon\left(\left(u\cdot\nabla\right)\phi-\Delta\mu\right)
+\int_{0}^{t}\int_{\mathbb{T}^N}u^\varepsilon\cdot\nabla\phi^\varepsilon\Delta\phi\nonumber\\ &\quad-\int_{0}^{t}\int_{\mathbb{T}^N}\frac{1}{\rho^\varepsilon}\Delta\mu^\varepsilon\Delta\phi =-\int_{\mathbb{T}^N}\nabla\phi_0^\varepsilon\cdot\nabla\phi_0.
	\end{align}
 With (\ref{q.5.6})-(\ref{q.5.9}) at hand, we can derive that
		\begin{align}\label{q.5.10} &\int_{\mathbb{T}^N}\left(\dfrac{1}{2}\left|\sqrt{\rho^\varepsilon}u^\varepsilon-u\right|^2+\dfrac{1}{2}\left|\nabla\phi^\varepsilon-\nabla\phi\right|^2+\dfrac{1}{4}\left|\sqrt{\rho^\varepsilon}\left(\left(\phi^\varepsilon\right)^2-1\right)\right|^2+\dfrac{1}{4}\left|\phi^2-1\right|^2+\Pi^\varepsilon\left(x,t\right) \right)\nonumber\\
		&\quad+\int_{0}^{t}\int_{\mathbb{T}^N}\left(\nu(\rho^\varepsilon,\phi^\varepsilon) \left|\nabla u^\varepsilon-\nabla u\right|^2+\eta(\rho^\varepsilon,\phi^\varepsilon) \left|\nabla\cdot u^\varepsilon\right|^2+\left|\nabla\mu^\varepsilon-\nabla\mu\right|^2\right)\nonumber\\ &=\int_{\mathbb{T}^N}\dfrac{1}{2}\left(\left|\sqrt{\rho_0^\varepsilon}u_0^\varepsilon-u_0\right|^2+\dfrac{1}{2}\left|\nabla\phi_0^\varepsilon-\nabla\phi_0\right|^2
+\dfrac{1}{4}\left|\sqrt{\rho_0^\varepsilon}\left(\left(\phi_0^\varepsilon\right)^2-1\right)\right|^2+\dfrac{1}{4}\left|\phi_0^2-1\right|^2+\Pi^\varepsilon\left(x,0\right)\right)\nonumber\\
		&\quad+\sum_{i=1}^{8}Y_i,
	\end{align}
	and
	\begin{align}
Y_1&=\int_{\mathbb{T}^N}\left(\rho^\varepsilon u^\varepsilon \cdot u-\sqrt{\rho^\varepsilon}u^\varepsilon \cdot u\right)-\int_{\mathbb{T}^N}\left(\rho_0^\varepsilon u_0^\varepsilon \cdot u_0-\sqrt{\rho^\varepsilon}u_0^\varepsilon \cdot u_0\right),\nonumber\\
Y_2&=\int_{0}^{t}\int_{\mathbb{T}^N}\rho^\varepsilon u^\varepsilon\cdot\left(u\cdot\nabla\right)u-\int_{0}^{t}\int_{\mathbb{T}^N}\left(\rho^\varepsilon u^\varepsilon\otimes u^\varepsilon\right):\nabla u ,
\nonumber \\
Y_3&=\int_{0}^{t}\int_{\mathbb{T}^N}\rho^\varepsilon u^\varepsilon\cdot\nabla p,
	\nonumber \\
Y_4&=-\int_{0}^{t}\int_{\mathbb{T}^N}\nu(\rho^\varepsilon,\phi^\varepsilon)\nabla u^\varepsilon:\nabla u
+\int_{0}^{t}\int_{\mathbb{T}^N}\left[\nu(\rho^\varepsilon,\phi^\varepsilon)-\nu(\phi)\right]|\nabla u|^2\nonumber\\
&\quad-\int_0^t\int_{\mathbb{T}^N}\nabla\nu(\phi)\cdot\nabla u\cdot u-\int_{0}^{t}\int_{\mathbb{T}^N}\nu(\phi)\rho^\varepsilon u^\varepsilon\cdot\Delta u\nonumber\\
&\quad-\int_0^t\int_{\mathbb{T}^N}\nabla\nu(\rho^\varepsilon,\phi^\varepsilon)\cdot\nabla u^\varepsilon\cdot u^\varepsilon
+\int_{0}^{t}\int_{\mathbb{T}^N}\nabla\nu(\rho^\varepsilon,\phi^\varepsilon)\cdot\nabla u^\varepsilon\cdot u
	\nonumber \\
&=-\int_{0}^{t}\int_{\mathbb{T}^N}\nu(\rho^\varepsilon,\phi^\varepsilon)\nabla \left(u^\varepsilon-u\right):\nabla u
-\int_0^t\int_{\mathbb{T}^N}\nabla(\nu(\phi)u):\nabla u\nonumber\\
&\quad-\int_{0}^{t}\int_{\mathbb{T}^N}\nu(\phi)\left(\rho^\varepsilon-1\right) u^\varepsilon\cdot\Delta u
-\int_{0}^{t}\int_{\mathbb{T}^N}\nu(\phi) u^\varepsilon\cdot\Delta u\nonumber\\
&\quad-\int_0^t\int_{\mathbb{T}^N}\nabla\nu(\rho^\varepsilon,\phi^\varepsilon)\cdot\nabla u^\varepsilon\cdot \left(u^\varepsilon-u\right)
\nonumber\\
&=-\int_{0}^{t}\int_{\mathbb{T}^N}\nu(\rho^\varepsilon,\phi^\varepsilon)\nabla \left(u^\varepsilon-u\right):\nabla u
+\int_0^t\int_{\mathbb{T}^N}\nu(\phi) u\cdot\Delta u
\nonumber\\
&\quad-\int_{0}^{t}\int_{\mathbb{T}^N}\nu(\phi)\left(\rho^\varepsilon-1\right) u^\varepsilon\cdot\Delta u
-\int_{0}^{t}\int_{\mathbb{T}^N}\nu(\phi) u^\varepsilon\cdot\Delta u
\nonumber\\
&\quad+\int_0^t\int_{\mathbb{T}^N}\nu(\rho^\varepsilon,\phi^\varepsilon)\nabla u^\varepsilon:\nabla \left(u^\varepsilon-u\right)
+\int_0^t\int_{\mathbb{T}^N}\nu(\rho^\varepsilon,\phi^\varepsilon)\Delta u^\varepsilon\cdot \left(u^\varepsilon-u\right)
\nonumber\\
&=\int_{0}^{t}\int_{\mathbb{T}^N}\nu(\rho^\varepsilon,\phi^\varepsilon)\left|\nabla\left( u^\varepsilon-u\right)\right|^2
+\int_{0}^{t}\int_{\mathbb{T}^N}\nu(\phi) \Delta\left( u^\varepsilon- u\right) \cdot\left(u^\varepsilon-u\right)
\nonumber\\
&\quad-\int_{0}^{t}\int_{\mathbb{T}^N}\nu(\phi)\left(\rho^\varepsilon-1\right) u^\varepsilon\cdot\Delta u+\int_0^t\int_{\mathbb{T}^N}\left(\nu(\rho^\varepsilon,\phi^\varepsilon)-\nu(\phi)\right)\Delta u^\varepsilon\cdot \left(u^\varepsilon-u\right)
\nonumber\\
&=\int_{0}^{t}\int_{\mathbb{T}^N}\nu(\rho^\varepsilon,\phi^\varepsilon)\left|\nabla\left( u^\varepsilon-u\right)\right|^2
-\int_{0}^{t}\int_{\mathbb{T}^N}\nu(\phi)\left| \nabla\left( u^\varepsilon- u\right) \right|^2\nonumber\\
&\quad
-\int_{0}^{t}\int_{\mathbb{T}^N}\nabla\nu(\phi)\cdot \nabla\left( u^\varepsilon- u\right) \cdot\left(u^\varepsilon-u\right)-\int_{0}^{t}\int_{\mathbb{T}^N}\nu(\phi)\left(\rho^\varepsilon-1\right) u^\varepsilon\cdot\Delta u\nonumber\\
&\quad+\int_0^t\int_{\mathbb{T}^N}\left(\nu(\rho^\varepsilon,\phi^\varepsilon)-\nu(\phi)\right)\Delta u^\varepsilon\cdot \left(u^\varepsilon-u\right)
\nonumber\\
&=\int_{0}^{t}\int_{\mathbb{T}^N}\nu(\rho^\varepsilon,\phi^\varepsilon)\left|\nabla\left( u^\varepsilon-u\right)\right|^2-\int_{0}^{t}\int_{\mathbb{T}^N}\nu(\phi)\left| \nabla\left( u^\varepsilon- u\right) \right|^2
\nonumber\\
&\quad
-\int_{0}^{t}\int_{\mathbb{T}^N}\nu(\phi)\left(\rho^\varepsilon-1\right) u^\varepsilon\cdot\Delta u
-\int_{0}^{t}\int_{\mathbb{T}^N}\nabla\nu(\phi)\cdot \nabla\left( u^\varepsilon- u\right) \cdot\left(u^\varepsilon-\sqrt{\rho^\varepsilon}u^\varepsilon\right)\nonumber\\
&\quad-\int_{0}^{t}\int_{\mathbb{T}^N}\nabla\nu(\phi)\cdot \nabla\left( u^\varepsilon- u\right) \cdot\left(\sqrt{\rho^\varepsilon}u^\varepsilon-u\right)\nonumber\\
&\quad+\int_0^t\int_{\mathbb{T}^N}\left(\nu(\rho^\varepsilon,\phi^\varepsilon)-\nu(\phi)\right)\Delta u^\varepsilon\cdot \left(u^\varepsilon-\sqrt{\rho^\varepsilon}u^\varepsilon\right)\nonumber\\
&\quad+\int_0^t\int_{\mathbb{T}^N}\left(\nu(\rho^\varepsilon,\phi^\varepsilon)-\nu(\phi)\right)\Delta u^\varepsilon\cdot \left(\sqrt{\rho^\varepsilon}u^\varepsilon-u\right),\nonumber\\
Y_5&=-\int_0^t\int_{\mathbb{T}^N}\left(u^\varepsilon-u\right)\cdot\nabla\eta(\rho^\varepsilon,\phi^\varepsilon) \left(\nabla\cdot u^\varepsilon\right)\nonumber\\
&=-\int_0^t\int_{\mathbb{T}^N}\left(u^\varepsilon-\sqrt{\rho^\varepsilon}u^\varepsilon\right)\cdot\nabla\eta(\rho^\varepsilon,\phi^\varepsilon) \left(\nabla\cdot u^\varepsilon\right)
-\int_0^t\int_{\mathbb{T}^N}\left(\sqrt{\rho^\varepsilon}u^\varepsilon-u\right)\cdot\nabla\eta(\rho^\varepsilon,\phi^\varepsilon) \left(\nabla\cdot u^\varepsilon\right),\nonumber\\
Y_6&=\int_{0}^{t}\int_{\mathbb{T}^N}\rho^\varepsilon u^\varepsilon
\cdot\left(\nabla\cdot\left(\nabla\phi\otimes\nabla\phi\right)\right)+\int_{0}^{t}\int_{\mathbb{T}^N}\left(u\cdot\nabla\right)\phi^\varepsilon\Delta\phi^\varepsilon\nonumber\\
		&\quad-\int_{0}^{t}\int_{\mathbb{T}^N}\left(u\cdot\nabla\right)\phi\Delta\phi^\varepsilon-\int_{0}^{t}\int_{\mathbb{T}^N}\left(u^\varepsilon\cdot\nabla\right)\phi^\varepsilon\Delta\phi,
\nonumber\\
&=\int_{0}^{t}\int_{\mathbb{T}^N}\left(\rho^\varepsilon-1\right)\left(u^\varepsilon\cdot\nabla\right)\phi\Delta\phi
+\int_{0}^{t}\int_{\mathbb{T}^N}\rho_t^\varepsilon\dfrac{|\nabla\phi|^2}{2}
\nonumber\\
&\quad-\int_{0}^{t}\int_{\mathbb{T}^N}\left(u^\varepsilon\cdot\nabla\right)\left(\phi^\varepsilon-\phi\right)\Delta\phi
+\int_{0}^{t}\int_{\mathbb{T}^N}\left(u\cdot\nabla\right)\left(\phi^\varepsilon-\phi\right)\Delta\phi^\varepsilon
\nonumber \\
&=\int_{0}^{t}\int_{\mathbb{T}^N}\left(\rho^\varepsilon-1\right)\left(u^\varepsilon\cdot\nabla\right)\phi\Delta\phi+\int_{0}^{t}\int_{\mathbb{T}^N}\rho_t^\varepsilon\dfrac{|\nabla\phi|^2}{2} \nonumber\\
&\quad-\int_{0}^{t}\int_{\mathbb{T}^N}\left(\left(u^\varepsilon-u\right)\cdot\nabla\right)\left(\phi^\varepsilon-\phi\right)\Delta\phi\nonumber \\
		&\quad+\int_{0}^{t}\int_{\mathbb{T}^N}\left(u\cdot\nabla\right)\left(\phi^\varepsilon-\phi\right)\left(\Delta\phi^\varepsilon-\Delta\phi\right)\nonumber \\
		&=\int_{0}^{t}\int_{\mathbb{T}^N}\left(\rho^\varepsilon-1\right)\left(u^\varepsilon\cdot\nabla\right)\phi\Delta\phi+\int_{0}^{t}\int_{\mathbb{T}^N}\rho_t^\varepsilon\dfrac{|\nabla\phi|^2}{2} \nonumber \\
		&\quad-\int_{0}^{t}\int_{\mathbb{T}^N}\left(1-\sqrt{\rho^\varepsilon}\right)\left(u^\varepsilon\cdot\nabla\right)\left(\phi^\varepsilon-\phi\right)\Delta\phi -\int_{0}^{t}\int_{\mathbb{T}^N}\left(\left(\sqrt{\rho^\varepsilon}u^\varepsilon-u\right)\cdot\nabla\right)\left(\phi^\varepsilon-\phi\right)\Delta\phi \nonumber \\
		&\quad-\int_{0}^{t}\int_{\mathbb{T}^N}\nabla\left( u\cdot\nabla\right)\left(\phi^\varepsilon-\phi\right)\cdot\left(\nabla\phi^\varepsilon-\nabla\phi\right),\nonumber\\
Y_7&=-\int_{0}^{t}\int_{\mathbb{T}^N}\left\{\left(\dfrac{1}{\rho^\varepsilon}-1\right)\Delta\phi^{\varepsilon}\Delta\mu
+\left(1-\dfrac{1}{\rho^\varepsilon}\right)\Delta\mu^\varepsilon\Delta\phi\right\}
,\nonumber\\
Y_8&=-2\int_{0}^{t}\int_{\mathbb{T}^N}\nabla\mu^\varepsilon\cdot\nabla\mu
+\int_{0}^{t}\int_{\mathbb{T}^N}\dfrac{1}{\rho^\varepsilon}\Delta\phi^{\varepsilon}\Delta\mu
+\int_{0}^{t}\int_{\mathbb{T}^N}\Delta\mu^{\varepsilon}\Delta\phi\nonumber\\
&=\int_{0}^{t}\int_{\mathbb{T}^N}\mu^\varepsilon\Delta\mu
+\int_{0}^{t}\int_{\mathbb{T}^N}\Delta\mu^\varepsilon\mu
+\int_{0}^{t}\int_{\mathbb{T}^N}\dfrac{1}{\rho^\varepsilon}\Delta\phi^{\varepsilon}\Delta\mu
+\int_{0}^{t}\int_{\mathbb{T}^N}\Delta\mu^{\varepsilon}\Delta\phi\nonumber\\
&=\int_{0}^{t}\int_{\mathbb{T}^N}\left(\mu^\varepsilon+\dfrac{1}{\rho^\varepsilon}\Delta\phi^{\varepsilon}\right)\Delta\mu
+\int_{0}^{t}\int_{\mathbb{T}^N}\Delta\mu^\varepsilon\left(\mu+\Delta\phi\right)\nonumber\\
&=\int_{0}^{t}\int_{\mathbb{T}^N}\left\{\left((\phi^\varepsilon)^3-\phi^\varepsilon\right)\Delta\mu
+\Delta\mu^\varepsilon(\phi^3-\phi)\right\}\nonumber
\nonumber\\
&=\int_{0}^{t}\int_{\mathbb{T}^N}\left[\left((\phi^\varepsilon)^3-\phi^\varepsilon\right)-(\phi^3-\phi)\right]\Delta\mu
+\int_{0}^{t}\int_{\mathbb{T}^N}\Delta\mu^\varepsilon\left((\phi^3-\phi)-\left((\phi^\varepsilon)^3-\phi^\varepsilon\right)\right)
\nonumber\\
&\quad+\int_{0}^{t}\int_{\mathbb{T}^N}\left(\phi^3-\phi\right)\Delta\mu+\int_{0}^{t}\int_{\mathbb{T}^N}\Delta\mu^\varepsilon\left((\phi^\varepsilon)^3-\phi^\varepsilon\right)
\nonumber\\
&=\int_{0}^{t}\int_{\mathbb{T}^N}\nabla\left[\left((\phi^\varepsilon)^3-\phi^\varepsilon\right)-(\phi^3-\phi)\right]\cdot
\nabla\left(\mu^\varepsilon-\mu\right)
\nonumber\\
&\quad+\int_{0}^{t}\int_{\mathbb{T}^N}\left(\phi^3-\phi\right)(\phi_t+u\cdot\nabla \phi)
+\int_{0}^{t}\int_{\mathbb{T}^N}\left(\rho^\varepsilon\phi_t^\varepsilon+\rho^\varepsilon u^\varepsilon\cdot\nabla \phi^\varepsilon\right)\left((\phi^\varepsilon)^3-\phi^\varepsilon\right)
\nonumber\\
&=\int_{0}^{t}\int_{\mathbb{T}^N}\nabla\left[\left((\phi^\varepsilon)^3-\phi^\varepsilon\right)-(\phi^3-\phi)\right]\cdot
\nabla\left(\mu^\varepsilon-\mu\right)
\nonumber\\
&\quad+\int_{0}^{t}\int_{\mathbb{T}^N}\left[\left(\dfrac{\phi^4}4-\dfrac{\phi^2}2\right)_t+u\cdot\nabla \left(\dfrac{\phi^4}4-\dfrac{\phi^2}2\right)\right]+\int_{0}^{t}\int_{\mathbb{T}^N}\left[\rho^\varepsilon\left(\dfrac{(\phi^\varepsilon)^4}4-\dfrac{(\phi^\varepsilon)^2}2\right)\right]_t\nonumber\\
&\quad-
\int_{0}^{t}\int_{\mathbb{T}^N}\rho_t^\varepsilon\left(\dfrac{(\phi^\varepsilon)^4}4-\dfrac{(\phi^\varepsilon)^2}2\right)
+\int_{0}^{t}\int_{\mathbb{T}^N}\rho^\varepsilon u^\varepsilon\cdot\nabla \left(\dfrac{(\phi^\varepsilon)^4}4-\dfrac{(\phi^\varepsilon)^2}2\right)
\nonumber\\
&=\int_{0}^{t}\int_{\mathbb{T}^N}\nabla\left[\left((\phi^\varepsilon)^3-\phi^\varepsilon\right)-(\phi^3-\phi)\right]\cdot
\nabla\left(\mu^\varepsilon-\mu\right)
\nonumber\\
&\quad+\int_{\mathbb{T}^N}\dfrac14\left(\phi^2-1\right)^2-\int_{\mathbb{T}^N}\dfrac14\left(\phi_0^2-1\right)^2
+\int_{\mathbb{T}^N}\dfrac14\rho^\varepsilon\left(\left(\phi^\varepsilon\right)^2-1\right)^2-\int_{\mathbb{T}^N}\dfrac14\rho_0^\varepsilon\left(\left(\phi_0^\varepsilon\right)^2-1\right)^2,
\label{Y8s}
\end{align}
where we have used integration by parts, (\ref{q.1.2}) and (\ref{q.5.1}).

 Now we estimate each $Y_i (i=1,...,8)$.
Firstly, it follows from Lemma \ref{L3.1}, the $\rm H\ddot{o}lder$ inequality, (\ref{q.2.6}), (\ref{q.2.7}) that
	\begin{align} |Y_1|&\leq\|u\|_{\infty}\left(\int_{\mathbb{T}^N}\rho^\varepsilon{\left|u^\varepsilon\right|}^2\right)^{\frac{1}{2}}\left(\int_{\mathbb{T}^N}\left|\sqrt{\rho^\varepsilon}-1 \right|^2\right)^{\frac{1}{2}}+\|u_0\|_\infty\left(\int_{\mathbb{T}^N}\rho_0^\varepsilon\left|u_0^\varepsilon\right|^2\right)^{\frac{1}{2}}\left(\int_{\mathbb{T}^N}\left|\sqrt{\rho_0^\varepsilon}-1 \right|^2\right)^{\frac{1}{2}}\nonumber \\	
		&\leq C\left \|\dfrac{1}{\sqrt{\rho^\varepsilon}}+1\right \|_\infty\left(\int_{\mathbb{T}^N}\left|\rho^\varepsilon-1 \right|^2\right)^{\frac{1}{2}}+C\left\|\dfrac{1}{\sqrt{\rho_0^\varepsilon}}+1\right\|_\infty\left(\int_{\mathbb{T}^N}\left|\rho_0^\varepsilon-1 \right|^2\right)^{\frac{1}{2}}\nonumber \\	
		&\leq C\varepsilon,	\label{q.5.11}\\	
		|Y_2|&= \left|-\int_{0}^{t}\int_{\mathbb{T}^N}\left(\left(\sqrt{\rho^\varepsilon}u^\varepsilon-u\right)\otimes\left(\sqrt{\rho^\varepsilon}u^\varepsilon-u\right)\right):\nabla u\nonumber \right.\\	
		&\left.\quad+\int_{0}^{t}\int_{\mathbb{T}^N}\left(\rho^\varepsilon-\sqrt{\rho^\varepsilon}\right)u^\varepsilon\cdot\left(\left(u\cdot\nabla\right)u\right)
-\int_{0}^{t}\int_{\mathbb{T}^N}\left(\sqrt{\rho^\varepsilon}u^\varepsilon-u\right)\cdot\nabla\left(\dfrac{|u|^2}{2}\right) \right|\nonumber \\
	&\leq C\varepsilon+C\int_{0}^{t}\int_{\mathbb{T}^N}\left|\sqrt{\rho^\varepsilon}u^\varepsilon-u\right|^2
+\int_{0}^{t}\int_{\mathbb{T}^N}\left|\sqrt{\rho^\varepsilon}\left(1-\sqrt{\rho^\varepsilon}\right)u^\varepsilon\cdot\nabla\left(\dfrac{|u|^2}2\right)
+\rho_t^\varepsilon\cdot\dfrac{|u|^2}2\right|
\nonumber\\
		&\leq C\varepsilon+C\int_{0}^{t}\int_{\mathbb{T}^N}\left|\sqrt{\rho^\varepsilon}u^\varepsilon-u\right|^2,\label{q.5.13}
\end{align}
where we have used in $(\ref{q.5.13})$ the following fact from the incompressibility $\nabla\cdot u=0$ that
\begin{align*}
&\left|\int_{0}^{t}\int_{\mathbb{T}^N}\left(\sqrt{\rho^\varepsilon}u^\varepsilon-u\right)\cdot\nabla\left(\dfrac{|u|^2}{2}\right) \right|
=\left|\int_{0}^{t}\int_{\mathbb{T}^N}\sqrt{\rho^\varepsilon}u^\varepsilon\cdot\nabla\left(\dfrac{|u|^2}{2}\right) \right|
\nonumber\\
&=\left|-\int_{0}^{t}\int_{\mathbb{T}^N}\sqrt{\rho^\varepsilon}(\sqrt{\rho^\varepsilon}-1)u^\varepsilon\cdot\nabla\left(\dfrac{|u|^2}{2}\right)
+\int_{0}^{t}\int_{\mathbb{T}^N}\rho^\varepsilon u^\varepsilon\cdot\nabla\left(\dfrac{|u|^2}{2}\right) \right|\nonumber\\
&=\left|-\int_{0}^{t}\int_{\mathbb{T}^N}\sqrt{\rho^\varepsilon}(\sqrt{\rho^\varepsilon}-1)u^\varepsilon\cdot\nabla\left(\dfrac{|u|^2}{2}\right)
+\int_{0}^{t}\int_{\mathbb{T}^N}\rho_t^\varepsilon\dfrac{|u|^2}{2} \right|\nonumber\\
&\le C\left(\int_0^t\int_{\mathbb{T}^N}\left|\rho^\varepsilon-1 \right|^2\right)^{\frac{1}{2}}
+C\left(\int_0^t\int_{\mathbb{T}^N}\left|\rho_t^\varepsilon \right|^2\right)^{\frac{1}{2}}\le C\varepsilon.
\end{align*}
In a similar way, we get
\begin{align}
		|Y_3|&=\left|\int_{0}^{t}\int_{\mathbb{T}^N}\rho_t^\varepsilon p\right|\leq C\left(\int_{0}^{t}\int_{\mathbb{T}^N}\left|\rho_t^\varepsilon\right|^2\right)^{\frac{1}{2}} \left(\int_{0}^{t}\int_{\mathbb{T}^N}p^2\right)^{\frac{1}{2}}\leq C \varepsilon.\label{q.5.14}
\end{align}
Now we turn to estimate $Y_4$ as follows:
\begin{align}
Y_4
&\leq\int_{0}^{t}\int_{\mathbb{T}^N}\nu(\rho^\varepsilon,\phi^\varepsilon)\left|\nabla\left( u^\varepsilon-u\right)\right|^2-\int_{0}^{t}\int_{\mathbb{T}^N}\nu(\phi)\left| \nabla\left( u^\varepsilon- u\right) \right|^2
\nonumber\\
&\quad+C\left(\int_{0}^{t}\int_{\mathbb{T}^N}\left|1-\rho^\varepsilon\right|^2\right)^{\frac{1}{2}}\left(\int_{0}^{t}\int_{\mathbb{T}^N}\left|u^\varepsilon\right|^2\left|\Delta u\right|^2\right)^{\frac{1}{2}}
\nonumber\\
&\quad+\dfrac{C}{\nu_*}\|\nabla\nu(\phi)u^\varepsilon\|_\infty\left(\int_{0}^{t}\int_{\mathbb{T}^N}\nu(\phi)|\nabla\left(u^\varepsilon-u\right)|^2\right)^{\frac{1}{2}}
\left(\int_{0}^{t}\int_{\mathbb{T}^N}\left|1-\sqrt{\rho^\varepsilon}\right|^2\right)^\frac12\nonumber\\
&\quad+\dfrac{C}{\nu_*}\|\nabla\nu(\phi)\|_\infty\left(\int_{0}^{t}\int_{\mathbb{T}^N}\nu(\phi)|\nabla\left(u^\varepsilon-u\right)|^2\right)^{\frac{1}{2}}
\left(\int_{0}^{t}\int_{\mathbb{T}^N}|\sqrt{\rho^\varepsilon}u^\varepsilon-u|^2\right)^{\frac{1}{2}}
\nonumber  \\
&\quad+ C\|\nu(\rho^\varepsilon,\phi^\varepsilon)-\nu(\phi)\|_\infty\left(\int_{0}^{t}\int_{\mathbb{T}^N}|u^\varepsilon\Delta u^\varepsilon|^2\right)^\frac12\left(\int_{0}^{t}\int_{\mathbb{T}^N}\left|1-\sqrt{\rho^\varepsilon}\right|^2\right)^\frac12\nonumber\\
&\quad+C\left(\int_{0}^{t}\int_{\mathbb{T}^N}|\rho^\varepsilon-1|^2\right)^{\frac{1}{2}}
\left(\int_{0}^{t}\int_{\mathbb{T}^N}|\sqrt{\rho^\varepsilon}u^\varepsilon-u|^2\right)^{\frac{1}{2}}
\nonumber  \\
&\quad+C\left(\int_{0}^{t}\int_{\mathbb{T}^N}|\phi^\varepsilon-\phi|^2\right)^{\frac{1}{2}}
\left(\int_{0}^{t}\int_{\mathbb{T}^N}|\sqrt{\rho^\varepsilon}u^\varepsilon-u|^2\right)^{\frac{1}{2}}
\nonumber  \\
&\leq\int_{0}^{t}\int_{\mathbb{T}^N}\nu(\rho^\varepsilon,\phi^\varepsilon)\left|\nabla\left( u^\varepsilon-u\right)\right|^2-\int_{0}^{t}\int_{\mathbb{T}^N}\nu(\phi)\left| \nabla\left( u^\varepsilon- u\right) \right|^2
+ C\varepsilon\nonumber\\
&\quad
+C\int_{0}^{t}\int_{\mathbb{T}^N}\left(|\sqrt{\rho^\varepsilon}u^\varepsilon-u|^2+|\phi^\varepsilon-\phi|^2
\right)+\int_{0}^{t}\int_{\mathbb{T}^N}\dfrac{\nu(\phi)}2\left|\nabla\left( u^\varepsilon-u\right)\right|^2\nonumber\\
&\leq C\varepsilon^2
+C\int_{0}^{t}\int_{\mathbb{T}^N}\left(|\sqrt{\rho^\varepsilon}u^\varepsilon-u|^2+|\phi^\varepsilon-\phi|^2
\right)
\nonumber\\
&\quad+\int_{0}^{t}\int_{\mathbb{T}^N}\nu(\rho^\varepsilon,\phi^\varepsilon)\left|\nabla\left( u^\varepsilon-u\right)\right|^2 -\int_{0}^{t}\int_{\mathbb{T}^N}\dfrac{\nu(\phi)}2\left|\nabla\left( u^\varepsilon-u\right)\right|^2,\label{Y4}
\end{align}
where we have used the Mean value theorem and the Cauchy inequality.
Similarly, we derive that
\begin{align}
		|Y_5|
&\leq \dfrac{C}{\eta_*}\left\|\nabla\eta(\rho^\varepsilon,\phi^\varepsilon)u^\varepsilon\right\|_\infty
\left(\int_{0}^{t}\int_{\mathbb{T}^N}\left|1-\sqrt{\rho^\varepsilon}\right|^2\right)^\frac12
\left(\int_{0}^{t}\int_{\mathbb{T}^N}\eta(\rho^\varepsilon,\phi^\varepsilon)\left|\nabla\cdot u^\varepsilon\right|^2\right)^\frac12\nonumber\\
&\quad+ \dfrac{C}{\eta_*}\left\|\nabla\eta(\rho^\varepsilon,\phi^\varepsilon)\right\|_\infty
\left(\int_{0}^{t}\int_{\mathbb{T}^N}\left|\sqrt{\rho^\varepsilon}u^\varepsilon-u\right|^2\right)^\frac12
\left(\int_{0}^{t}\int_{\mathbb{T}^N}\eta(\rho^\varepsilon,\phi^\varepsilon)\left|\nabla\cdot u^\varepsilon\right|^2\right)^\frac12\nonumber\\
&\leq C\varepsilon^2+C\int_{0}^{t}\int_{\mathbb{T}^N}\left|\sqrt{\rho^\varepsilon}u^\varepsilon-u\right|^2
+\dfrac12\int_{0}^{t}\int_{\mathbb{T}^N}\eta(\rho^\varepsilon,\phi^\varepsilon)\left|\nabla\cdot u^\varepsilon\right|^2,\\
		|Y_6|&\leq C\left(\int_{0}^{t}\int_{\mathbb{T}^N}|\rho^\varepsilon-1|^2\right)^{\frac{1}{2}}\left(\int_{0}^{t}\int_{\mathbb{T}^N}|u^\varepsilon|^2|\nabla\phi|^2|\Delta\phi|^2\right)^{\frac{1}{2}}
\nonumber\\
&
\quad+C\left(\int_{0}^{t}\int_{\mathbb{T}^N}|\rho_t^\varepsilon|^2\right)^{\frac{1}{2}}\left(\int_{0}^{t}\int_{\mathbb{T}^N}\dfrac{|\nabla\phi|^4}{4}\right)^{\frac{1}{2}} \nonumber \\		&\quad+C\|u^\varepsilon\|_\infty\|\Delta\phi\|_\infty\left(\int_{0}^{t}\int_{\mathbb{T}^N}|1-\sqrt{\rho^\varepsilon}|^2\right)^{\frac{1}{2}}\left(\int_{0}^{t}\int_{\mathbb{T}^N}|\nabla\left(\phi^\varepsilon-\phi\right)|^2\right)^{\frac{1}{2}}\nonumber  \\		&\quad+C\|\Delta\phi\|_\infty\left(\int_{0}^{t}\int_{\mathbb{T}^N}|\sqrt{\rho^\varepsilon}u^\varepsilon-u|^2\right)^{\frac{1}{2}}\left(\int_{0}^{t}\int_{\mathbb{T}^N}|\nabla\left(\phi^\varepsilon-\phi\right)|^2\right)^{\frac{1}{2}}\nonumber  \\
		&\quad+C\|\nabla u\|_\infty\int_{0}^{t}\int_{\mathbb{T}^N}|\nabla\left(\phi^\varepsilon-\phi\right)|^2\nonumber \\
		&\leq C\left(\varepsilon+\varepsilon^2\right)+C\int_{0}^{t}\int_{\mathbb{T}^N}\left(|\sqrt{\rho^\varepsilon}u^\varepsilon-u|^2+|\nabla\left(\phi^\varepsilon-\phi\right)|^2\right),\\
		|Y_7|&\leq C\left(\int_0^t\int_{\mathbb{T}^N}\left|\rho^\varepsilon-1 \right|^2\right)^{\frac{1}{2}}\leq C\varepsilon,
		\\
		Y_8&\leq \dfrac12\int_{0}^{t}\int_{\mathbb{T}^N}\left|\nabla\left(\mu^\varepsilon-\mu\right)\right|^2
+C\int_{0}^{t}\int_{\mathbb{T}^N}\left[\left|\nabla(\phi^\varepsilon-\phi)\right|^2+\left|\phi^\varepsilon-\phi\right|^2\right]\nonumber\\
&\quad+\int_{\mathbb{T}^N}\dfrac14\left(\phi^2-1\right)^2-\int_{\mathbb{T}^N}\dfrac14\left(\phi_0^2-1\right)^2
+\int_{\mathbb{T}^N}\dfrac14\rho^\varepsilon\left(\left(\phi^\varepsilon\right)^2-1\right)^2-\int_{\mathbb{T}^N}\dfrac14\rho_0^\varepsilon\left(\left(\phi_0^\varepsilon\right)^2-1\right)^2.\label{q.5.17}
\end{align}
 In conclusion, adding up (\ref{q.5.10}) to (\ref{q.5.17}) leads to
		\begin{align}\label{q.5.19}
		&\int_{\mathbb{T}^N}\left(\dfrac{1}{2}\left|\sqrt{\rho^\varepsilon}u^\varepsilon-u\right|^2+\Pi^\varepsilon\left(x,t\right)
+\dfrac{1}{2}\left|\nabla\phi^\varepsilon-\nabla\phi\right|^2\right)+\int_{0}^{t}\int_{\mathbb{T}^N}\dfrac{\nu(\phi)}{2}\left|\nabla u^\varepsilon-\nabla u\right|^2 \nonumber\\
		&\quad+\int_{0}^{t}\int_{\mathbb{T}^N}\left(\dfrac{\eta(\rho^\varepsilon,\phi^\varepsilon)}2\left|\nabla\cdot u^\varepsilon\right|^2
+\left|\nabla\mu^\varepsilon-
\nabla\mu\right|^2\right) \nonumber\\
		&\leq\int_{\mathbb{T}^N}\left(\dfrac{1}{2}\left|\sqrt{\rho_0^\varepsilon}u_0^\varepsilon-u_0\right|^2+\Pi^\varepsilon\left(x,0\right)
+\dfrac{1}{2}\left|\nabla\phi_0^\varepsilon-\nabla\phi_0\right|^2\right)+C\varepsilon \nonumber\\
		&\quad+C\int_{0}^{t}\int_{\mathbb{T}^N}\left(\left|\sqrt{\rho^\varepsilon}u^\varepsilon-u\right|^2
+\left|\nabla\left(\phi^\varepsilon-\phi\right)\right|^2+\left|\phi^\varepsilon-\phi\right|^2\right) .
		\end{align}
In addition, subtracting $(\ref{q.2.1})_3$ from $(\ref{q.1.2})_3$, using $(\ref{q.2.1})_4$ and multiplying the result equality by $\phi^\varepsilon-\phi$, we obtain after integrating over $Q_t$ that
\begin{align}\label{q.5.18}
&\frac12\int_{\mathbb{T}^N}\left|\phi^\varepsilon-\phi\right|^2+\int_0^t\int_{\mathbb{T}^N}\left|\Delta\left(\phi^\varepsilon-\phi\right)\right|^2
\nonumber\\
&=
-\int_0^t\int_{\mathbb{T}^N}\frac{\left(1-\rho^\varepsilon\right)\left(1+\rho^\varepsilon\right)}{\left(\rho^\varepsilon\right)^2}\Delta^2\phi^\varepsilon\left(\phi^\varepsilon-\phi\right)
-\int_0^t\int_{\mathbb{T}^N}\left(u^\varepsilon-u\right)\cdot\nabla\phi^\varepsilon\left(\phi^\varepsilon-\phi\right)
\nonumber\\
&\quad-\int_0^t\int_{\mathbb{T}^N}u\cdot\nabla\left(\phi^\varepsilon-\phi\right)\left(\phi^\varepsilon-\phi\right)
-\int_0^t\int_{\mathbb{T}^N}\frac{2}{\rho^\varepsilon}\nabla\left(\frac{1}{\rho^\varepsilon}\right)\nabla\Delta\phi^\varepsilon\left(\phi^\varepsilon-\phi\right)
\nonumber\\
&\quad-\int_0^t\int_{\mathbb{T}^N}\frac{1}{\rho^\varepsilon}\Delta\left(\frac{1}{\rho^\varepsilon}\right)\Delta\phi^\varepsilon\left(\phi^\varepsilon-\phi\right)
+\int_0^t\int_{\mathbb{T}^N}\frac{1-\rho^\varepsilon}{\rho^\varepsilon}\Delta\left(\left(\phi^\varepsilon\right)^3-\phi^\varepsilon\right)
\left(\phi^\varepsilon-\phi\right)
\nonumber\\
&\quad+\int_0^t\int_{\mathbb{T}^N}
       \left[\Delta\left(\left(\phi^\varepsilon\right)^3-\phi^\varepsilon\right)-\Delta\left(\phi^3-\phi\right)\right]\left(\phi^\varepsilon-\phi\right)
+\int_{\mathbb{T}^N}\frac12\left|\phi_0^\varepsilon-\phi_0\right|^2
\nonumber\\
&\leq C\varepsilon+C\int_0^t\int_{\mathbb{T}^N}\left|u^\varepsilon-u\right|^2
+C\int_0^t\int_{\mathbb{T}^N}\left(\left|\phi^\varepsilon-\phi\right|^2+\left|\nabla\left(\phi^\varepsilon-\phi\right)\right|^2\right)\nonumber\\
&\quad
+\int_0^t\int_{\mathbb{T}^N}\frac14\left|\Delta\left(\phi^\varepsilon-\phi\right)\right|^2
+\int_{\mathbb{T}^N}\frac12\left|\phi_0^\varepsilon-\phi_0\right|^2\nonumber\\
&\leq C\varepsilon+C\int_0^t\int_{\mathbb{T}^N}\left|\sqrt{\rho^\varepsilon}u^\varepsilon-u\right|^2
+C\int_0^t\int_{\mathbb{T}^N}\left(\left|\phi^\varepsilon-\phi\right|^2+\left|\nabla\left(\phi^\varepsilon-\phi\right)\right|^2\right)\nonumber\\
&\quad
+\frac14\int_0^t\int_{\mathbb{T}^N}\left|\Delta\left(\phi^\varepsilon-\phi\right)\right|^2
+\frac12\int_{\mathbb{T}^N}\left|\phi_0^\varepsilon-\phi_0\right|^2.
\end{align}
 By the following facts indicated from (\ref{q.2.6}) and (\ref{q.2.7}) that
		\begin{align}\label{q.5.20}
		\int_{\mathbb{T}^N}\left|\sqrt{\rho_0^\varepsilon}u_0^\varepsilon-u_0\right|^2=\int_{\mathbb{T}^N}\left|\sqrt{\rho_0^\varepsilon}u_0^\varepsilon-u_0^\varepsilon\right|^2+\int_{\mathbb{T}^N}\left|u_0^\varepsilon-u_0\right|^2\leq C\left(\varepsilon^4+\varepsilon^2\right),
	\end{align}
	and
	\begin{align}\label{q.5.21}
		\int_{\mathbb{T}^N}\left|u^\varepsilon-u\right|^2\leq\int_{\mathbb{T}^N}\left|\sqrt{\rho^\varepsilon}u^\varepsilon-u^\varepsilon\right|^2+\int_{\mathbb{T}^N}\left|\sqrt{\rho^\varepsilon}u^\varepsilon-u\right|^2\leq C\varepsilon^2+\int_{\mathbb{T}^N}\left|\sqrt{\rho^\varepsilon}u^\varepsilon-u\right|^2.
		\end{align}
 By (\ref{q.5.19})-(\ref{q.5.21}), the $\rm Poincar\acute{e}$ inequality and the Gronwall inequality, we have
	\begin{align}\label{q.5.22}
		\left\|u^\varepsilon-u\right\|^2+\left\|\phi^\varepsilon-\phi\right\|_1^2\leq C\varepsilon.
	\end{align}	
Then inserting (\ref{q.5.22}) into (\ref{q.5.19}), we get
	\begin{align}\label{q.5.23}
		\int_{0}^{t}\left(\nu(\phi)\left\|\nabla\left(u^\varepsilon-u\right)\right\|^2+\left\|\nabla\mu^\varepsilon-\nabla\mu\right\|^2
+\left\|\nabla\phi^\varepsilon-\nabla\phi\right\|_1^2\right)\leq C\varepsilon.
		\end{align}

	Finally, subtracting $(\ref{q.1.1})_4$ from $(\ref{q.1.2})_4$, operating the result equality by $\nabla$, we have
	\begin{align*}
		&\nabla\Delta\phi^\varepsilon-\nabla\Delta\phi\\
&=\nabla[\left(1-\rho^\varepsilon\right)\mu^\varepsilon]
-\nabla\left(\mu^\varepsilon-\mu\right)+\nabla\left[\left(\rho^\varepsilon -1\right)\left(\left(\phi^\varepsilon\right)^3-\phi^3\right)\right]
+\nabla\left(\left(\phi^\varepsilon\right)^3-\phi^3\right)-\nabla\left(\phi^\varepsilon-\phi\right).
		\end{align*}
	Then the $\rm Poincar\acute{e}$ inequality implies  
\begin{align}\label{q.5.25}
		\left\|\nabla\phi^\varepsilon-\nabla\phi\right\|_2^2
&\leq C\left\|\nabla\rho^\varepsilon\right\|^2\left\|\mu^\varepsilon\right\|_\infty^2
+C\left\|1-\rho^\varepsilon\right\|^2\left\|\nabla\mu^\varepsilon\right\|_\infty^2
+C\left\|\nabla\mu^\varepsilon-\nabla\mu\right\|^2\nonumber\\
&\quad+ C\left\|\nabla\rho^\varepsilon\right\|^2\left\|\left(\phi^\varepsilon\right)^3-\phi^3\right\|_\infty^2
+C\left\|1-\rho^\varepsilon\right\|^2\left\|\nabla\left(\left(\phi^\varepsilon\right)^3-\phi^3\right)\right\|_\infty^2\nonumber\\
&\quad+C\left\|\phi^\varepsilon-\phi\right\|^2+C\left\|\nabla\phi^\varepsilon-\nabla\phi\right\|^2,
		\end{align}
	for $t\in[0,T^0]$.
	And then (\ref{q.5.25}), together with (\ref{q.5.23}) and (\ref{q.5.22}), implies that
\begin{align*}
		\int_{0}^{t}\left\|\nabla\left(\phi^\varepsilon-\phi\right)\right\|_2^2\leq C\varepsilon.
	\end{align*}
Furthermore, it follows from (\ref{q.3.5700}) and (\ref{q.3.87}) that (\ref{q.2.15}).
This completes the proof Theorem \ref{th3}.
\hfill$\Box$

\section{Navier-Stokes/Allen-Cahn system}\label{sec:3}
This section is devoted to incompressible limit of the Navier-Stokes/Allen-Cahn system,
i.e. the operator $A_K=A_{AC}$ in $(\ref{q.1.1})$.

\subsection{Statements of the theorems.}
When the density is away from vacuum, we rewrite the compressible Navier-Stokes/Allen-Cahn system $(\ref{q.1.1})$ as follows
	\begin{align}\label{qqq.2.1}
	\begin{cases}
		\rho_t^\varepsilon+\nabla\cdot\left(\rho^\varepsilon u^\varepsilon\right)=0,\\
		u^\varepsilon_t+\left(u^\varepsilon\cdot\nabla\right) u^\varepsilon+\dfrac1{\varepsilon^2\rho^\varepsilon}\nabla\left(P\left(\rho^\varepsilon\right)\right)=\dfrac{\nu(\rho^\varepsilon,\phi^\varepsilon)}{\rho^\varepsilon}\Delta u^\varepsilon+\dfrac{\eta(\rho^\varepsilon,\phi^\varepsilon)}{\rho^\varepsilon}\nabla\left(\nabla\cdot u^\varepsilon\right)
		-\dfrac{1}{\rho^\varepsilon}\Delta\phi^\varepsilon\nabla\phi^\varepsilon,\\
		\phi^\varepsilon_t+\left(u^\varepsilon\cdot\nabla\right) \phi^\varepsilon=\dfrac{1}{\left(\rho^\varepsilon\right)^2}\Delta\phi^\varepsilon-\dfrac1{\rho^\varepsilon}\left(\left(\phi^\varepsilon\right)^3-\phi^\varepsilon\right).\\
	\end{cases}	
	\end{align}
Define
\begin{align*}
	\begin{cases}
		E_s\left(U\left(t\right)\right)=\dfrac{1}{2}\displaystyle\sum_{|\alpha|\leq s}\displaystyle\int_{\mathbb{T}^N}\left(\displaystyle\frac{1}{\varepsilon^2}|\nabla^{\alpha} \left(\rho-1\right)|^2+|\nabla^{\alpha}u|^2+|\nabla\nabla^{\alpha}\phi|^2\right),\\
		\widetilde{E}_{s}\left(U\left(t\right)\right)=\dfrac{1}{2}\displaystyle\sum_{|\alpha|\leq s}\displaystyle\int_{\mathbb{T}^N}\left(\displaystyle\frac{1}{\varepsilon^2}\frac{P'\left(\rho\right)}{\rho}|\nabla^{\alpha} \left(\rho-1\right)|^2+\rho|\nabla^{\alpha}u|^2+\rho^2|\nabla\nabla^{\alpha}\phi|^2\right),\\
	\end{cases}
\end{align*}
where $U=\left(\rho,u,\phi\right)$. Then
\begin{align}\label{qqq.2.3}
	E_s\left(U\left(t\right)\right)\backsim\widetilde{E}_{s}\left(U\left(t\right)\right),
\end{align}
provided that $|\rho-1|$ is sufficiently small.

Now, we state the main results of this section.
\begin{theorem}\label{th1-2}
Consider the compressible Navier-Stokes/Allen-Cahn system $(\ref{q.1.1})$ with the following initial data
	\begin{align}\label{qqq.2.4}
	\rho^\varepsilon\left(x,0\right)=1+\overline\rho^\varepsilon_0(x), \quad~
	u^\varepsilon\left(x,0\right)=u_0(x)+\overline u^\varepsilon_0(x), \quad~
	\phi^\varepsilon\left(x,0\right)=\phi_0(x)+\overline\phi^\varepsilon_0(x),
\end{align}
where $u_0$ and $\phi_0$ satisfy
	\begin{align}\label{qqq.2.5}
	u_0\in H^{s+1}\left(\mathbb{T}^N\right),\quad\quad~ \nabla\cdot u_0=0,\quad\quad~\phi_0\in H^{s+2}\left(\mathbb{T}^N\right),
\end{align}
for any $s\ge[\frac N2]+2$. Moreover, for small positive constant $\kappa_0$, the functions
$\overline\rho^\varepsilon_0(x)$, $\overline u^\varepsilon_0(x)$, $\overline\phi^\varepsilon_0(x)$ are assumed to satisfy
	\begin{align}\label{qqq.2.6}
	\left\|\overline\rho_0^\varepsilon(x)\right\|_{s}\leq \kappa_0\varepsilon^2,
	\quad\quad~
	\left\|\overline u_0^\varepsilon(x)\right\|_{s+1}\leq \kappa_0\varepsilon,
\quad\quad~
	\left\|\overline\phi_0^\varepsilon(x)\right\|_{s+2}\leq \kappa_0\varepsilon.
\end{align}
Then the following statements hold.

Uniform stability: There exist constants $T_0$ and $C$ independent of $\varepsilon$ such that a unique strong solution
$(\rho^\varepsilon, u^\varepsilon, \phi^\varepsilon)$ of system $(\ref{q.1.1})$ exists for all small $\varepsilon$ on the time interval $[0, T_0]$ with properties:
\begin{align}\label{qqq.2.7}
	\begin{cases}
		E_s\left(U^\varepsilon\left(t\right)\right)+\left\|\phi^\varepsilon\right\|^2+\displaystyle\int_0^t\left(\nu_*\left\|\nabla u^\varepsilon\right\|^2_s+\eta_*\left\|\nabla\cdot u^\varepsilon\right\|^2_s+\left\|\nabla\phi^\varepsilon\right\|^2_{s+1}\right)\leq C,
\\
		E_{s-1}\left(\partial_t U^\varepsilon\left(t\right)\right)+\left\|\phi^\varepsilon_t\right\|^2+\displaystyle\int_0^t\left(\nu_*\left\|\nabla \partial_tu^\varepsilon\right\|^2_{s-1}+\eta_*\left\|\nabla\cdot \partial_t u^\varepsilon\right\|^2_{s-1}+\left\|\nabla\partial_t\phi^\varepsilon\right\|^2_s\right)\leq C,
	\end{cases}
\end{align}
where $E_{s-1}\left(\partial_tU\left(t\right)\right)=\dfrac{1}{2}\displaystyle\sum_{|\beta|\leq {s-1}}\displaystyle\int_{\mathbb{T}^N}\left(\displaystyle\frac{1}{\varepsilon^2}|\nabla^{\beta} \partial_t\rho|^2+|\nabla^{\beta} \partial_t u|^2+|\nabla^{\beta} \partial_t\nabla\phi|^2\right)$.

Local existence of solutions for incompressible Navier-Stokes/Allen-Cahn system: There exist functions $u$ and $\phi$ such that
\begin{align}\label{qqq.2.8}
\begin{cases}
\rho^{\varepsilon}\to1~~\text{in}~ L^{\infty}(\left[0,T_0\right];H^{s}) \cap  {\rm Lip}(\left[0,T_0\right];H^{s-1}),
\\
u^{\varepsilon}\stackrel{\mathrm{w}^*}\rightharpoonup u~~\text{in}~L^{\infty}(\left[0,T_0\right];H^{s})
\cap  {\rm Lip}(\left[0,T_0\right];H^{s-1}),
\\
u^{\varepsilon}\to u~~\text{in}~ C(\left[0,T_0\right];H^{s'}),
\\
\phi^{\varepsilon}\stackrel{\mathrm{w}^*}\rightharpoonup\phi~~\text{in}~ L^{\infty}(\left[0,T_0\right];H^{s+1})
\cap  {\rm Lip}(\left[0,T_0\right];H^{s}),
\\
\phi^{\varepsilon}\to \phi~~\text{in}~ C(\left[0,T_0\right];H^{{s'}+1})
\end{cases}
\end{align}
for any $s'\in{[0,s)}$, and the function pair $\left(u,\phi\right)$ is the unique strong solution of
the incompressible Navier-Stokes/Allen-Cahn system (\ref{q.1.2}) with the initial data
\begin{equation}\label{qqq.2.9}
u(x, 0)=u_0(x),		\quad\quad
\phi(x, 0)=\phi_0(x),
\end{equation}
for some $p\in L^\infty\left(\left[0, T_0\right];\, H^{s-1}\right)\cap L^2\left(\left[0, T_0\right];\, H^s\right)$.
\end{theorem}

\begin{theorem}\label{th2-2} Consider the strong solutions $(\rho^\varepsilon, u^\varepsilon, \phi^\varepsilon)$
of the Navier-Stokes/Allen-Cahn system (\ref{q.1.1}) obtained in Theorem \ref{th1-2}. Suppose in addition that the initial data satisfies
%
	\begin{align}\label{qqq.2.10}
	\left\|u_0\right\|_s^2+\left\|\phi_0^2-1\right\|_s^2+\left\|\nabla\phi_0\right\|_s^2\leq\delta,
\end{align}
where $\delta$ is a positive constant. If $\delta$ is sufficiently small
and $\bar{M}\triangleq\left(\delta+\varepsilon^2\kappa_0^2\right)$, then for any fixed $T>0$, the strong solution $\left(\rho^{\varepsilon},u^{\varepsilon},\phi^{\varepsilon}\right)$ satisfies the estimates:
\begin{align}\label{qqq.2.11}
	&E_s\left(U^\varepsilon\left(t\right)\right)
+\left\|\left({\phi^\varepsilon}\right)^2-1\right\|^2_s(t)
\nonumber\\
&+\int_0^t\left(\nu_*\left\|\nabla u^\varepsilon\right\|^2_s(t)
+\eta_*\left\|\nabla\cdot u^\varepsilon\right\|^2_s
+\left\|\left({\phi^\varepsilon}\right)^2-1\right\|^2_{s+1}+\left\|\nabla\phi^\varepsilon\right\|^2_{s+1}\right)\leq 4\bar{M},
\end{align}
for $t\in[0, T^\varepsilon)$, and
\begin{align}\label{qqq.2.12}
&E_{s-1}\left(\partial_t U^\varepsilon\left(t\right)\right)
+\left\|\phi^\varepsilon_t\right\|^2(t)
\nonumber\\
&\quad+\int_0^t\left(\nu_*\left\|\nabla \partial_t u^\varepsilon\right\|^2_{s-1}+\eta_*\left\|\nabla\cdot \partial_tu^\varepsilon\right\|^2_{s-1}+\left\|\nabla\partial_t\phi^\varepsilon\right\|^2_s\right)\leq C {\rm exp}Ct,
\end{align}
for $ t\in[0,T]$, where $T^\varepsilon>T$ and $T^\varepsilon\to\infty$ as $\varepsilon\to 0$.

Furthermore, as $\varepsilon\to 0$, $\left(\rho^{\varepsilon},\ u^{\varepsilon},\ \phi^{\varepsilon}\right)$
converges to the unique global strong solution $\left(1,u,\phi\right)$ of the incompressible Navier-Stokes/Allen-Cahn system $(\ref{q.1.2})$, and
	\begin{align}\label{qqq.2.13}
	\left\|u\right\|_s^2(t)
+\left\|\phi^2-1\right\|_s^2(t)
+\left\|\nabla\phi\right\|_s^2(t)
+\int_0^t\left(\nu_*\left\|\nabla u\right\|^2_s+\left\|\phi^2-1\right\|_{s+1}^2+\left\|\nabla\phi\right\|^2_{s+1}\right)\leq C_1\delta
\end{align}
for any $t>0$, where $C_1$ is a uniform constant independent of $\delta$ and $t$.
\end{theorem}

\begin{theorem}\label{th3-2} Under the assumptions of Theorem \ref{th1-2}, the convergence rate of $\rho^\varepsilon$, $u^\varepsilon$ and $\phi^\varepsilon$ $(\varepsilon\to 0)$ are deduced as
\begin{align}\label{qqq.2.14}
\left\|u^\varepsilon-u\right\|^2(t)+
\left\|\phi^\varepsilon-\phi\right\|_2^2(t)
+\int_{0}^{t}\left(\left\|u^\varepsilon-u\right\|_1^2+\left\|\phi^\varepsilon-\phi\right\|_3^2\right)
\leq C\varepsilon,
\end{align}
for $t\in[0,T_0]$. Furthermore, we have
\begin{align}\label{qqq.2.15}
\left\|\rho^\varepsilon-1\right\|_s^2(t)\leq C\varepsilon, \quad\quad
\left\|\nabla\left(\rho^\varepsilon-1\right)\right\|_{s-2}^2\leq C\varepsilon^4, \qquad \forall t\in[0,T_0].
\end{align}
 The statement also holds for the strong solution given in Theorem \ref{th2-2} for $t\in \left[0,T^\varepsilon\right)$.
\end{theorem}

In the following proof, we only state the different part and omit the analogous part.

\subsection{Local existence and uniform stability.}\label{subsec:3.2}
\quad~~In this subsection, we will give the uniform estimates for our results and then prove Theorem \ref{th1-2}. Let $U_0=\left(1+\overline\rho_0^\varepsilon,\ u_0+\overline u_0^\varepsilon,\ \phi_0+\overline\phi_0^\varepsilon\right)$. We consider a set of functions
	$B_{T_0}^\varepsilon\left(U_0\right)$ contained in
	$\{(\rho, u, \phi)~|~ (\rho, u, \nabla\phi)\in L^\infty([0, T_0]; H^s)\cap {\rm Lip}([0, T_0]; H^{s-1})\}$ with $s\ge 3$
	and defined by
	\begin{align}\label{qqq.3.1}
		\begin{cases}
			\left|\dfrac{\rho-1}{\varepsilon}\right|+|u-u_0|+\left|\phi-\phi_0\right|<\kappa,
			\\
			E_s\left(U\left(t\right)\right)+\left\|\phi\right\|^2+\displaystyle\int_0^t\left(\nu_*\left\|\nabla u\right\|_s^2+\eta_*\left\|\nabla\cdot u\right\|_s^2+\left\|\nabla\phi\right\|_{s+1}^2\right)\le \tilde K_1,
			\\
			E_{s-1}\left(\partial_tU\left(t\right)\right)+\left\|\phi_t\right\|^2
			+\displaystyle\int_0^t\left(\nu_*\left\|\nabla u_t\right\|_{s-1}^2+\eta_*\left\|\nabla\cdot u_t\right\|_{s-1}^2+\left\|\nabla\phi_t\right\|_s^2\right)\le \tilde K_2.
		\end{cases}
	\end{align}

For any $V=\left(\xi^\varepsilon, v^\varepsilon, \psi^\varepsilon\right)\in B_{T_0}^\varepsilon(U_0)$, define $U=(\rho^\varepsilon, u^\varepsilon, \phi^\varepsilon)=\Lambda(V)$ as the unique solution of the following linearized problem
	\begin{align}\label{qqq.3.2}
		\begin{cases}
			\rho^\varepsilon_t+(v^\varepsilon\cdot\nabla)\rho^\varepsilon+\xi^\varepsilon\nabla\cdot u^\varepsilon=0,
			\\
			u^\varepsilon_t+(v^\varepsilon\cdot\nabla) u^\varepsilon+\dfrac1{\varepsilon^2}\dfrac{P^{'}(\xi^\varepsilon)}{\xi^\varepsilon}\nabla\rho^\varepsilon
			=\dfrac{\nu(\xi^\varepsilon,\psi^\varepsilon)}{\xi^\varepsilon}\Delta u^\varepsilon+\dfrac{\eta(\xi^\varepsilon,\psi^\varepsilon)}{\xi^\varepsilon}\nabla\left(\nabla\cdot u^\varepsilon\right)
			-\dfrac{1}{\xi^\varepsilon}\Delta\phi^\varepsilon\nabla\phi^\varepsilon,
			\\
			\phi^\varepsilon_t-\dfrac{1}{\left(\xi^\varepsilon\right)^2}\Delta\phi^\varepsilon
			=-(v^\varepsilon\cdot\nabla)\psi^\varepsilon-\dfrac{1}{\xi^\varepsilon}\left[\left(\psi^\varepsilon\right)^3-\psi^\varepsilon\right],
		\end{cases}
	\end{align}
for which the existence and uniqueness of the solutions is guaranteed by the standard theory of parabolic equations and Navier-Stokes equations. Now we are to show that for appropriate choices of $T_{0}$, $\kappa$, $\tilde K_{1}$, $\tilde K_{2}$ independent of $\varepsilon$, $ \Lambda$ maps $B_{T_0}^\varepsilon(U_0)$ into itself and it is a contraction in certain function spaces. We emphasize that the solutions will depend on the value of the parameter $\varepsilon$, but for convenience, the dependence will not always be displayed in this subsection.

\begin{lemma}\label{L3.2-2}
	Suppose that $B_{T_0}^{\varepsilon}\left(U_0\right)$ is defined by (\ref{qqq.3.1}) and $\Lambda:V\longrightarrow U$ is defined by the system(\ref{qqq.3.2}). Then, under the assumptions in Theorem \ref{th1-2}, there exist constants $T_0$, $\kappa$, $\tilde K_1$ and $\tilde K_2$ independent of $\varepsilon$ such that $\Lambda$ maps $B_{T_0}^{\varepsilon}\left(U_0\right)$ into itself.
	\end{lemma}

\noindent{\it\bfseries Proof.}\quad Firstly, we apply $D^{\alpha_1}$ to the first and second equations of $(\ref{qqq.3.2})$ and $D^{\alpha_2}$ to the third one respectively, and then we get
%
		\begin{align}\label{qqq.3.3}
			\begin{cases}
				\partial_tD^{\alpha_1}\rho+\left(v\cdot\nabla\right)D^{\alpha_1}\rho+\xi\nabla\cdot D^{\alpha_1}u=\Pi_1,\\
				\partial_tD^{\alpha_1}u+\left(v\cdot\nabla\right)D^{\alpha_1}u+\dfrac{1}{\varepsilon^2}\dfrac{P'\left(\xi\right)}{\xi}\nabla D^{\alpha_1}\rho=\Pi_2,\\
				\partial_tD^{\alpha_2}\phi-\dfrac{1}{\xi^2}D^{\alpha_2}\Delta\phi=\Pi_3,
			\end{cases}
		\end{align}
	where
\begin{align*}
 \Pi_1&=-\left[D^{\alpha_1}\left(v\cdot\nabla\rho\right)-\left(v\cdot\nabla\right)D^{\alpha_1}\rho\right]-\left[D^{\alpha_1}\left(\xi\nabla\cdot u\right)-\xi\nabla\cdot D^{\alpha_1}u\right],\\
\Pi_2&=\dfrac{\nu(\xi,\psi)}{\xi}\Delta D^{\alpha_1} u+\dfrac{\eta(\xi,\psi)}{\xi}\nabla D^{\alpha_1}(\nabla\cdot u)-\dfrac{1}{\xi}D^{\alpha_1}\left(\Delta\phi\nabla\phi\right)\\
&\quad-\left[D^{\alpha_1}\left(v\cdot\nabla u\right)-\left(v\cdot\nabla\right)D^{\alpha_1}u\right]-
	\dfrac{1}{\varepsilon^2}\left[D^{\alpha_1}\left(\dfrac{P'(\xi)}{\xi}\nabla \rho\right)-
	\dfrac{P'(\xi)}{\xi}\nabla D^{\alpha_1}\rho\right]\\
&\quad+\left[D^{\alpha_1}\left(\dfrac{\nu(\xi,\psi)}{\xi}\Delta u\right)-
	\dfrac{\nu(\xi,\psi)}{\xi}\Delta D^{\alpha_1} u\right]
\nonumber\\
&\quad+\left[D^{\alpha_1}\left(\dfrac{\eta(\xi,\psi)}{\xi}\nabla(\nabla\cdot u)\right)-\dfrac{\eta(\xi,\psi)}{\xi}\nabla D^{\alpha_1}(\nabla\cdot u)\right]\\
	&\quad-\left[D^{\alpha_1}\left(\dfrac{1}{\xi}\left(\Delta\phi\nabla\phi\right)\right)
-\dfrac{1}{\xi}D^{\alpha_1}\left(\Delta\phi\nabla\phi\right)\right],\\
	\Pi_3&=-D^{\alpha_2}(v\cdot\nabla \psi)-D^{\alpha_2}\left[\dfrac{1}{\xi}\left(\psi^3-\psi\right)\right]+\left[D^{\alpha_2}\left(\dfrac{1}{\xi^2}\Delta\phi\right)-\dfrac{1}{\xi^2}
D^{\alpha_2}\Delta\phi\right].
\end{align*}

	We will prove that $\Lambda$ maps $B_{T_0}^\varepsilon(U_0)$ into itself by two steps and denote by $C$ the constants independent of $\varepsilon$, $\tilde K_{1}$, $\tilde K_{2}$ in these two steps. Without loss of generality, we assume that $T_{0}^{-1}$, $\varepsilon$, $\tilde K_{1}$ and $\tilde K_{2}$ are all bigger than $1$.

 $\mathbf{Step\ one}$ : Estimates of $\phi$.

On one hand, we rewrite the equation $(\ref{qqq.3.2})_3$ as
\begin{align}\label{qqq.3.4}
\xi^2\phi_t-\Delta\phi=-\xi^2 v\cdot\nabla\psi-\xi(\psi^3-\psi).
\end{align}
For $|\alpha_1|\le s$, applying $\nabla^{\alpha_1}$ to the above equation and multiplying by $\nabla^{\alpha_1}\phi_t$,
then integrating the result on $\mathbb{T}^N$ yield
		\begin{align}\label{qqq.3.5}
			&\dfrac12\dfrac{\rm d}{{\rm d}t}\int_{\mathbb{T}^N}|\nabla^{\alpha_1}\nabla\phi|^2+\int_{\mathbb{T}^N}\xi^2|\nabla^{\alpha_1}\phi_t|^2
\nonumber\\
&=-\int_{\mathbb{T}^N}\left[\nabla^{\alpha_1}(\xi^2\phi_t)-\xi^2\nabla^{\alpha_1}\phi_t\right]
\nabla^{\alpha_1}\phi_t
-\int_{\mathbb{T}^N}\nabla^{\alpha_1}(\xi^2v\cdot\nabla\psi)\nabla^{\alpha_1}\phi_t
\nonumber\\
			&\quad-\int_{\mathbb{T}^N}\nabla^{\alpha_1}[\xi(\psi^3-\psi)]\nabla^{\alpha_1}\phi_t
=\sum_{i=1}^{3}\tilde M_i.
		\end{align}

Next, we will give the estimates of $\tilde M_1$, $\tilde M_2$, $\tilde M_3$. Choose $\kappa$ small enough so that $|\xi-1|\le \frac{1}{2}$.

It follows from Lemma \ref{L3.1}, the Sobolev embedding $H^2(\mathbb{T}^N)$ $\hookrightarrow$ $L^\infty(\mathbb{T}^N)$ for $N =2, 3$, and the Cauchy inequality that
		\begin{align}
			|\tilde M_1|&\le C(\|\nabla(\xi^2)\|_{\infty}\|\phi_t\|_{s-1}+\|\nabla(\xi^2)\|_{s-1}\|\phi_t\|_{\infty})\|\nabla^{\alpha_1}\phi_t\|
	\nonumber\\
&\leq C \varepsilon \tilde K_1^\frac12	\|\phi_t\|_{s-1}\|\nabla^{\alpha_1}\phi_t\|	
\leq C(\tau)	\tilde K_1\|\phi_t\|_{s-1}^2+\tau\|\nabla^{\alpha_1}\phi_t\|^2	
,\label{qqq.3.6}\\		
	|\tilde M_2|&\le C(\|\xi^2v\|_{\infty}\|\nabla\psi\|_s+\|\xi^2v\|_s\|\nabla\psi\|_{\infty})\|\nabla^{\alpha_1}\phi_t\|\nonumber\\
			&\leq C \tilde K_1^2\left\|\nabla^{\alpha_1}\phi_t\right\|
\leq C(\tau)\tilde K_1^4+\tau\|\nabla^{\alpha_1}\phi_t\|^2,\label{qqq.3.7}\\
			|\tilde M_3|&\le C(\|\xi\|_{\infty}\|\psi^3-\psi\|_s+\|\xi\|_s\|\psi^3-\psi\|_{\infty})\|\nabla^{\alpha_1}\phi_t\|\nonumber\\
			&\leq C \tilde K_1(\tilde K_1+1)\left\|\nabla^{\alpha_1}\phi_t\right\|
\leq C(\tau)\tilde K_1^4+\tau\|\nabla^{\alpha_1}\phi_t\|^2.\label{qqq.3.8}
		\end{align}
where we have used
\begin{align*}
\|\xi\|_s\le C\left(1+\varepsilon\left\|\dfrac1{\varepsilon}\nabla\xi\right\|_{s-1}\right)\le C\tilde K_1^{\frac12}+C\varepsilon \tilde K_1^{\frac12}\le C\tilde K_1^{\frac12}.
\end{align*}
By using the equation (\ref{qqq.3.4}), we have
\begin{align}\label{aaa.3.9}
\|\phi_t\|_{s-1}
\le&\left\|\dfrac1{\xi^2}\Delta\phi\right\|_{s-1}+\|v\cdot\nabla\psi\|_{s-1}+\left\|\dfrac1\xi(\psi^3-\psi)\right\|_{s-1}
\nonumber\\
\le&\left\|\dfrac1{\xi^2}\right\|_{\infty}\|\Delta\phi\|_{s-1}+\left\|\dfrac{1}{\xi^2}\right\|_{s-1}\|\Delta\phi\|_{\infty}
+\|v\|_{\infty}\|\nabla\psi\|_{s-1}+\|v\|_{s-1}\|\nabla\psi\|_{\infty}
\nonumber\\
&+\left\|\dfrac1\xi\right\|_{\infty}\|\psi^3-\psi\|_{s-1}+\left\|\dfrac1\xi\right\|_{s-1}\|\psi^3-\psi\|_{\infty}
\nonumber\\
\le&C(\tilde K_1+1)\|\nabla\phi\|_{s}+C \tilde K_1( \tilde K_1+1).
\end{align}
Substituting the estimates (\ref{qqq.3.6})-(\ref{aaa.3.9}) into (\ref{qqq.3.5}), summing over $\alpha_1$, and then taking $\varepsilon$ small enough we obtain
		\begin{align*}
			&\dfrac{\rm d}{{\rm d}t}\displaystyle\sum_{|\alpha_1|\le s}\|\nabla^{\alpha_1}\nabla\phi\|^2
+\displaystyle\sum_{|\alpha_1|\le s}\|\nabla^{\alpha_1}\phi_t\|^2
\le C \tilde K_1^3\|\nabla\phi\|_{s}^2+C \tilde K_1^5.
		\end{align*}
Then integrating the above inequality over $[0, t]\subseteq[0, T_0]$, and applying (\ref{qqq.2.4})-(\ref{qqq.2.6}) and the Gronwall inequality yield
	%
	\begin{align}\label{qqq.3.10}
		\|\nabla\phi\|_{s}^2(t)+\int_0^t\|\phi_t\|_{s}^2{\rm d}s\le C,
	\end{align}
	for $t\in[0, T_0]$, provided that $T_0<T_1^*= \tilde K_1^{-3}$.
By using the equation (\ref{qqq.3.4}), one has
\begin{align*}
\|\Delta\phi\|_{s}
\le&\left\|\xi^2\right\|_{s}\left\|\phi_t\right\|_{\infty}+\left\|\xi^2\right\|_{\infty}\left\|\phi_t\right\|_{s}
+\left\|\xi^2v\right\|_{\infty}\left\|\nabla\psi\right\|_{s}+\left\|\xi^2v\right\|_{s}\left\|\nabla\psi\right\|_{\infty}
\nonumber\\
&+\left\|\xi\right\|_{\infty}\left\|\psi^3-\psi\right\|_{s}+\left\|\xi\right\|_{s}\left\|\psi^3-\psi\right\|_{\infty}
\nonumber\\
\le&C(1+\varepsilon^2 \tilde K_1)\|\phi_t\|_{s}+C \tilde K_1\left(\tilde K_1+1\right).
\end{align*}
Integrating the above equation in $[0, t]\subseteq[0, T_0]$, by (\ref{qqq.3.10}) we are led to
\begin{align}
  \int_{0}^{t}\|\Delta\phi\|_{s}^2{\rm d}s &\leq \int_{0}^{t}\left[C(1+\varepsilon^2 \tilde K_1)^2\|\phi_t\|_{s}^2 +C \tilde K_1^4\right]{\rm d}s
  \leq C,  \nonumber
\end{align}
for $t\in[0, T_0]$, provided that $T_0<T_1^*$.
And then we get
\begin{align}\label{qqq.3.11}
 \int_{0}^{t}\|\nabla\phi\|_{s+1}^2 {\rm d}s&\leq C\int_{0}^{t}\|\Delta\phi\|_{s}^2{\rm d}s\leq C,
\end{align}
where we have used the fact that
$\|\nabla\phi\|_{s+1}\leq C\|\Delta\phi\|_{s}$.
Thus putting (\ref{qqq.3.10}) and (\ref{qqq.3.11}) together, we arrival at
\begin{align}\label{qqq.3.12}
  \|\nabla\phi\|_{s}^2(t)+\int_0^t\|\phi_t\|_{s}^2{\rm d}s+\int_{0}^{t}\|\nabla\phi\|_{s+1}^2 {\rm d}s\leq C.
\end{align}

On the other hand, applying $\nabla^{\alpha_1}\partial_t$ to the equation (\ref{qqq.3.2})$_3$, multiplying $\nabla^{\alpha_1}\phi_t$ for $1\le|\alpha_1|\le s$,
and then integrating the result over $\mathbb{T}^N$, we get
	\begin{align}\label{qqq.3.13}
		&\dfrac12\dfrac{\rm d}{{\rm d}t}\int_{\mathbb{T}^N}|\nabla^{\alpha_1}\phi_t|^2+\int_{\mathbb{T}^N}\dfrac{1}{\xi^2}|\nabla^{\alpha_1}\nabla\phi_t|^2
		\nonumber\\
		&=-\int_{\mathbb{T}^N}\nabla\left(\dfrac{1}{\xi^2}\right)\cdot\nabla^{\alpha_1}\nabla\phi_t\cdot\nabla^{\alpha_1}\phi_t
		-\int_{\mathbb{T}^N}\nabla^{\alpha_1}(v\cdot\nabla\psi)_t\cdot\nabla^{\alpha_1}\phi_t
		\nonumber\\
&\quad-\int_{\mathbb{T}^N}\nabla^{\alpha_1}\left(\dfrac{1}{\xi}(\psi^3-\psi)\right)_t\cdot\nabla^{\alpha_1}\phi_t
+\int_{\mathbb{T}^N}\left[\nabla^{\alpha_1}\left(\dfrac1{\xi^2}\Delta\phi\right)_t-\dfrac1{\xi^2}\nabla^{\alpha_1}\Delta\phi_t\right]\cdot\nabla^{\alpha_1}\phi_t
	=\displaystyle\sum_{i=1}^4\tilde N_i.
	\end{align}
Next, choose $\tau$ small enough, and we calculate $\tilde N_i (i=1,...,4)$ as follows:
		\begin{align}
		\tilde N_1&\le C\left\|\nabla\left(\dfrac1{\xi^2}\right)\right\|_{\infty}\left\|\nabla^{\alpha_1}\nabla\phi_t\right\| \left\|\nabla^{\alpha_1}\phi_t\right\|
		\le C\varepsilon \tilde K_1^\frac{1}{2}\left\|\nabla^2\phi_t\right\|_{s-1}\left\|\nabla^{\alpha_1}\phi_t\right\| \nonumber\\
		&\leq\tau\left\|\nabla^2\phi_t\right\|_{s-1}^2+C(\tau)\varepsilon^2 \tilde K_1\left\|\nabla^{\alpha_1}\phi_t\right\|^2,\label{qqq.3.14}\\
    \tilde N_2&=-\int_{\mathbb{T}^N}\nabla^{\alpha_1}\left(v_t\cdot\nabla\psi+v\cdot\nabla\psi_t\right)\cdot\nabla^{\alpha_1}\phi_t =\int_{\mathbb{T}^N}\nabla^{\alpha_1-1}\left(v_t\cdot\nabla\psi+v\cdot\nabla\psi_t\right)\nabla^{\alpha_1-1}\Delta\phi_t
       \nonumber\\
		&\le C\left(\|v_t\|_{\infty}\|\nabla\psi\|_{s-1}+\|v_t\|_{s-1}\|\nabla\psi\|_{\infty}
          +\|v\|_{\infty}\|\nabla\psi_t\|_{s-1}+\|v\|_{s-1}\|\nabla\psi_t\|_{\infty}\right)\|\nabla^{\alpha_1-1}\Delta\phi_t\|
\nonumber\\
&\le C\tilde K_1^\frac{1}{2}\tilde K_2^\frac{1}{2}\|\Delta\phi_t\|_{s-1}\le \tau\|\Delta\phi_t\|_{s-1}^2+C(\tau)\tilde K_1\tilde K_2,\label{qqq.3.15}\\
		\tilde N_3&=-\int_{\mathbb{T}^N}\nabla^{\alpha_1}\left(\dfrac{1}{\xi}(\psi^3-\psi)\right)_t\cdot\nabla^{\alpha_1}\phi_t \nonumber  \\ &=\int_{\mathbb{T}^N}\nabla^{{\alpha_1}-1}\left[\left(\dfrac{1}{\xi}\right)_t\left(\psi^3-\psi\right)+\frac{1}{\xi}\left(3\psi^2-1\right)\psi_t\right]
\nabla^{\alpha_1-1}\Delta\phi_t
		\nonumber \\
		&\le\left(\left\|\left(\dfrac1\xi\right)_t\right\|_\infty\|\psi^3-\psi\|_{s-1}+\left\|\left(\dfrac1\xi\right)_t\right\|_{s-1}\|\psi^3-\psi\|_{\infty}
\nonumber\right.\\
        &\left.\quad  +\left\|\dfrac1\xi(3\psi^2-1)\right\|_{\infty}\|\psi_t\|_{s-1}+\left\|\dfrac1\xi(3\psi^2-1)\right\|_{s-1}\|\psi_t\|_{\infty}\right)
     \|\nabla^{\alpha_1-1}\Delta\phi_t\|
\nonumber\\
&\le C \tilde K_1^2\tilde K_2^\frac{1}{2}\|\Delta\phi_t\|_{s-1}
\le\tau\|\Delta\phi_t\|_{s-1}^2+ C(\tau) \tilde K_1^4\tilde K_2,\label{qqq.3.16}\\
		\tilde N_4&=-\int_{\mathbb{T}^N}\nabla^{\alpha_1-1}\left(\left(\dfrac{1}{\xi^2}\right)_t\Delta\phi\right)\cdot\nabla^{\alpha_1-1}\Delta\phi_t
		+\int_{\mathbb{T}^N}\left[\nabla^{\alpha_1}\left(\dfrac{1}{\xi^2}\Delta\phi_t\right)-\dfrac{1}{\xi^2}\nabla^{\alpha_1}\Delta\phi_t\right]\cdot\nabla^{\alpha_1}\phi_t\nonumber\\
		&\le C\left(\left\|\left(\dfrac1{\xi^2}\right)_t\right\|_{\infty}\left\|\Delta\phi\right\|_{s-1}
		+\left\|\left(\dfrac1{\xi^2}\right)_t\right\|_{s-1}\left\|\Delta\phi\right\|_{\infty}\right)\left\|\nabla^{\alpha_1-1}\Delta\phi_t\right\|
		\nonumber \\
		&\quad+C\left(\left\|\nabla\left(\dfrac1{\xi^2}\right)\right\|_{\infty}\left\|\Delta\phi_t\right\|_{s-1}
		+\left\|\nabla\left(\dfrac1{\xi^2}\right)\right\|_{s-1}\left\|\Delta\phi_t\right\|_{\infty}\right)\left\|\nabla^{\alpha_1}\phi_t\right\|
		\nonumber \\
&\le C \tilde K_1^\frac12\tilde K_2^\frac{1}{2}\left\|\Delta\phi_t\right\|_{s-1}+C\tilde K_1^\frac12\left\|\Delta\phi_t\right\|_{s-1}\left\|\nabla^{\alpha_1}\phi_t\right\|
		\nonumber \\
&\le C(\tau) \tilde K_1\left\|\nabla^{\alpha_1}\phi_t\right\|^2+\tau\left\|\Delta\phi_t\right\|_{s-1}^2+C(\tau)\tilde K_1\tilde K_2,\label{qqq.3.17}
\end{align}
	where we have used (\ref{qqq.3.12}).
Putting the estimates $(\ref{qqq.3.14})$-(\ref{qqq.3.17}) into (\ref{qqq.3.13}) and summing over $\alpha_1$, we obtain
	\begin{align}
		\label{chi-t-2}
		&\dfrac{\rm d}{{\rm d}t}\displaystyle\sum_{1\le|\alpha_1|\le s}\left\|\nabla^{\alpha_1}\phi_t\right\|^2
		+\displaystyle\sum_{1\le|\alpha_1|\le s}\|\nabla^{\alpha_1}\nabla\phi_t\|^2
		 \le C\tilde K_1\left\|\nabla^{\alpha_1}\phi_t\right\|^2+C\tilde K_1^4\tilde K_2.
\end{align}
From the equation $(\ref{qqq.3.2})_3$ and the constraints of the initial data, we see that
	\begin{align*}
		&\left\|\nabla\phi_t(x,0)\right\|_{s-1}\le\left\|\phi_t(x,0)\right\|_s\nonumber\\
		&\le C\left(\left\|\dfrac1{\xi^2(x,0)}\Delta\phi(x,0)\right\|_s
		+\left\|v(x,0)\cdot\nabla\psi(x,0)\right\|_s+\left\|\dfrac1{\xi(x,0)}\left(\psi^3(x,0)-\psi(x,0)\right)\right\|_s\right)
\le C.
\end{align*}
Then integrating (\ref{chi-t-2}) over $[0, t]\subseteq[0,  T_0]$ and using the Gronwall inequality yield
%
	\begin{align}\label{qqq.3.1900000}
		\left\|\nabla\phi_t\right\|_{s-1}^2(t)+\int_0^t\left\|\nabla^2\phi_t\right\|_{s-1}^2{\rm d}s\le C,
\end{align}
 for $t\in[0, T_0]$, provided that $T_0<T_2^*:={\rm min}\{T_1^*,\tilde K_1^{-4}\tilde K_2^{-1}\}$.

In the case $|\alpha_1|=0$, (\ref{qqq.3.13}) is simplified as follows:
	\begin{align}\label{qqq.3.19}
		&\dfrac12\dfrac{\rm d}{{\rm d}t}\int_{\mathbb{T}^N}|\phi_t|^2+\int_{\mathbb{T}^N}\dfrac{1}{\xi^2}|\nabla\phi_t|^2
		\nonumber\\
		&=-\int_{\mathbb{T}^N}\nabla\left(\dfrac{1}{\xi^2}\right)\cdot\nabla\phi_t\phi_t
		-\int_{\mathbb{T}^N}(v\cdot\nabla\psi)_t\phi_t
		-\int_{\mathbb{T}^N}\left(\dfrac{1}{\xi}(\psi^3-\psi)\right)_t\phi_t
		+\int_{\mathbb{T}^N}\left(\dfrac1{\xi^2}\right)_t\Delta\phi\phi_t \nonumber\\
		&\le C\|\phi_t\|\left(\|\nabla\xi\|_\infty\|\nabla\phi_t\|+\|v_t\|\|\nabla\psi\|_\infty+\|v\|_\infty\|\nabla\psi_t\|+\|\xi_t\|\|\psi^3-\psi\|_\infty
\right.\nonumber\\
        &\left.\quad+\left\|\frac{1}{\xi}\right\|_\infty\left\|\left(\psi^3-\psi\right)_t\right\|+\|\xi_t\|_\infty\|\Delta\phi\|
\right)\nonumber\\
&\le C\|\phi_t\|\left(\varepsilon \tilde K_1^\frac{1}{2}\|\nabla\phi_t\|+ \tilde K_1^\frac{1}{2}\tilde K_2^\frac{1}{2}\right)\nonumber\\
&\le C \tilde K_1\tilde K_2\|\phi_t\|^2+\tau\|\nabla\phi_t\|^2+C.
	\end{align}
Integrating the result over $[0, t]\subseteq[0,  T_0]$,
taking $\tau$ small enough and using the Gronwall inequality, we obtain
	\begin{align*}
		\left\|\phi_t\right\|^2(t)+\int_0^t\left\|\nabla\phi_t\right\|^2{\rm d}s\le C,
\end{align*}
for $t\in[0, T_0]$, provided that $T_0<T_2^*$.
Combining the above inequality with (\ref{qqq.3.1900000}), we get
%
	\begin{align}\label{qqq.3.19}
		\left\|\phi_t\right\|_{s}^2(t)+\int_0^t\left\|\nabla\phi_t\right\|_s^2{\rm d}s\le C.
\end{align}

 $\mathbf{Step\ two}$: Estimates of $\rho$ and u.

 $\mathbf{Cases\ 1}$-$\mathbf{Cases\ 2}$:
 It is similar to subsection \ref{subsec:2.2}, so we will not give details for the part. Obviously, we have the following inequality:
	\begin{align}\label{qqq.3.37}
		&\dfrac{\rm d}{{\rm d}t}\sum_{|\alpha_1|\leq s}\int_{\mathbb{T}^N}\left(\dfrac{P'(\xi)}{\xi}\left|\dfrac{1}{\varepsilon} \nabla^{\alpha_1}(\rho-1)\right|^2+\xi|\nabla^{\alpha_1} u|^2\right)\nonumber\\
		&\qquad+\sum_{|\alpha_1|\leq s}\nu_*\int_{\mathbb{T}^N}|\nabla^{\alpha_1} \nabla u|^2+\sum_{|\alpha_1|\leq s}\eta_*\int_{\mathbb{T}^N}|\nabla^{\alpha_1} (\nabla\cdot u)|^2	\nonumber\\
		&\leq C\left(\tilde K_1+\tilde K_2\right)\left(\left\|\dfrac{1}{\varepsilon}\left(\rho-1\right)\right\|_s^2+\|u\|_s^2\right)+C.
\end{align}
Then by the constraints of the initial data and the Gronwall inequality, we obtain
	\begin{align}\label{qqq.3.39}
		\left\|\dfrac{1}{\varepsilon}(\rho-1)\right\|_s^2(t)+\left\|u\right\|_s^2(t)\leq  {\rm exp}\left(C\left(\tilde K_1+\tilde K_2\right)T_0\right)\left(3\kappa_0^2\varepsilon^2+2\left\|u_0\right\|_s^2+CT_0\right)\leq C
	\end{align}
	for $t\in[0, T_0]$, provided that $T_0<T_2^*$.
  Furthermore, integrating corresponding inequalities over $[0, t]\subseteq[0, T_0]$ , and then we arrive at
	\begin{align}\label{qqq.3.40}
		\nu_*\int_{0}^{t}\|\nabla u\|_s^2{\rm d}s +\eta_*\int_{0}^{t}\|\nabla\cdot u\|_s^2{\rm d}s \leq C \quad{\rm for}\ t\in[0, T_0].
	\end{align}

$\mathbf{Cases\ 3}$: $D^{\alpha_1}=\nabla^\beta\partial_t$ for $1\leq|\beta|\leq s-1$, where $s\geq3$.

Similarly, we only state the terms $\tilde I_4$, $\tilde I_7$ and $\tilde I_8$, some extra discussions have to be shown since the methods to estimate these three terms are different between the cases s=3 and $s\geq4$. Indeed, if $s\geq4$, then it follows from Lemma \ref{L3.1}, the Sobolev embedding, the Cauchy inequality, (\ref{qqq.3.12}), (\ref{qqq.3.19}), (\ref{qqq.3.39}) and a similar calculation as in (\ref{q.3.44}) that
	\begin{align*}
		|\tilde I_4|&\leq C\left\|\frac{\rho_t}{\varepsilon}\right\|_{s-1}\left(\left\|v_t\right\|_\infty\left\|\frac{\nabla\rho}{\varepsilon}\right\|_{s-1}
+\left\|v_t\right\|_{s-1}\left\|\frac{\nabla\rho}{\varepsilon}\right\|_\infty
+\left\|\nabla v\right\|_\infty\left\|\frac{\nabla\rho_t}{\varepsilon}\right\|_{s-2}
\nonumber \right.\\
		&\left.\quad+\left\|\nabla v\right\|_{s-2}\left\|\frac{\nabla\rho_t}{\varepsilon}\right\|_\infty
+\left\|\frac{\xi_t}{\varepsilon}\right\|_\infty\left\|\nabla\cdot u\right\|_{s-1}+\left\|\frac{\xi_t}{\varepsilon}\right\|_{s-1}\left\|\nabla\cdot u\right\|_\infty
		\nonumber \right.\\
		&\left.\quad+\left\|\frac{\nabla\xi}{\varepsilon}\right\|_\infty\left\|\nabla\cdot u_t\right\|_{s-2}+\left\|\frac{\nabla\xi}{\varepsilon}\right\|_{s-2}\left\|\nabla\cdot u_t\right\|_\infty\right) \nonumber\\
		&\leq C\left(\tilde K_1+\tilde K_2\right)\left(\left\|\frac{\rho_t}{\varepsilon}\right\|_{s-1}^2+\left\|u_t\right\|_{s-1}^2\right)+C,\\
		|\tilde I_7|&\leq C\dfrac{1}{\varepsilon^2}\left\|u_t\right\|_{s-1} \left(\|\xi_t\|_\infty\left\|\nabla\rho\right\|_{s-1}+\|\xi_t\|_{s-1}(1+\|\nabla\xi\|_{s-2})\left\|\nabla\rho\right\|_\infty\right.\nonumber\\
        &\left.\quad+\left\|\nabla\xi\right\|_\infty\left\|\nabla\rho_t\right\|_{s-2}+\left\|\nabla\xi\right\|_{s-2}\left\|\nabla\rho_t\right\|_\infty\right)
		 \nonumber\\
		&\leq C\left(\tilde K_1^2+\tilde K_2^2\right)\left(\left\|\frac{\rho_t}{\varepsilon}\right\|_{s-1}^2+\left\|u_t\right\|_{s-1}^2\right)+C.
	\end{align*}
For the term $\tilde I_8$, it can be clearly written as
\begin{align*}
\tilde I_8&=\int_{\mathbb{T}^N}\xi \left[\nabla^{\beta}\partial_t\left(\dfrac{\nu(\xi,\psi)}{\xi}\Delta u\right)-
	\dfrac{\nu(\xi,\psi)}{\xi}\Delta \nabla^{\beta} u_t\right] : \nabla^{\beta}  u_t \nonumber\\
&\quad+\int_{\mathbb{T}^N}\xi \left[\nabla^{\beta}\partial_t\left(\dfrac{\eta(\xi,\psi)}{\xi}\nabla(\nabla\cdot u)\right)-
	\dfrac{\eta(\xi,\psi)}{\xi} \nabla^{\beta}\nabla(\nabla\cdot u_t)\right] : \nabla^{\beta}  u_t=I_8^1+I_8^2.
\end{align*}
By Lemma \ref{L3.1}, we deduce
\begin{align*}
\tilde I_8^1
&=\int_{\mathbb{T}^N}\xi \left\{\nabla^{\beta}\left[\left(\nu_\xi(\xi,\psi)\xi_t+\nu_\psi(\xi,\psi)\psi_t-\dfrac{\nu(\xi,\psi)\xi_t}{\xi^2}\right)\Delta u\right]
\right.\nonumber\\
&\left.\quad
+\nabla^{\beta}\left(\dfrac{\nu(\xi,\psi)}{\xi}\Delta u_t\right)
-\dfrac{\nu(\xi,\psi)}{\xi}\Delta \nabla^{\beta} u_t\right\} : \nabla^{\beta}  u_t \nonumber\\
&\le C\left[\left(\left\|\nu_\xi(\xi,\psi)-\dfrac{\nu(\xi,\psi)}{\xi^2}\right\|_{s-1}\|\xi_t\|_{s-1}
+\left\|\nu_\psi(\xi,\psi)\right\|_{s-1}\|\psi_t\|_{s-1}\right)\|\Delta u\|_{s-1}\right.\nonumber\\
&\left.\quad+\left\|\nabla\left(\dfrac{\nu(\xi,\psi)}{\xi}\right)\right\|_{s-2}\|\Delta u_t\|_{s-2}\right]\|u_t\|_{s-1}\nonumber\\
&\le C\left[(1+\|\nabla\xi\|_{s-2}+\left\|\nabla\psi\right\|_{s-2})(\|\xi_t\|_{s-1}+\|\psi_t\|_{s-1})\|\nabla u\|_{s}\right.\nonumber\\
&\left.\quad+C(\|\nabla\xi\|_{s-2}+\left\|\nabla\psi\right\|_{s-2})\|\nabla u_t\|_{s-1}\right]\|u_t\|_{s-1}\nonumber\\
&\le C\left(\tilde K_1^{\frac12}\tilde K_2^{\frac12}\|\nabla u\|_{s}+ \tilde K_1^{\frac12}\|\nabla u_t\|_{s-1}\right)\|u_t\|_{s-1}\nonumber\\
&\le C(\tau)\tilde K_1\tilde K_2\|u_t\|_{s-1}^2+C\|\nabla u\|_{s}^2+\tau\|\nabla u_t\|_{s-1}^2,
\end{align*}
and the term $\tilde I_8^2$ also has similar result. Then we get
\begin{align}
\tilde I_8
&\le C(\tau)\tilde K_1\tilde K_2\|u_t\|_{s-1}^2+C\left(\|\nabla u\|_{s}^2+\|\nabla\cdot u\|_{s}^2\right)+\tau\left(\|\nabla u_t\|_{s-1}^2+\|\nabla\cdot u_t\|_{s-1}^2\right).
\end{align}

When $s=3$, for the terms $\tilde I_4$, $\tilde I_7$, and $\tilde I_8$, one merely tackles terms with $|\beta|=s-1=2$ $ \left(D^{\alpha_1}=\nabla_i\nabla_j\partial_t\right)$ since it is more easier to draw for the case $|\beta|=1$ in a similar method, or one can get the estimate of $\int_{\mathbb{T}^N}(|\varepsilon^{-1}\nabla^\beta\rho_t|^2+|\nabla^\beta u_t|^2)$ by the interpolation since we have tackled the cases $D^{\alpha_1}=\partial_t$ and
$D^{\alpha_1}=\nabla_i\nabla_j\partial_t$. It follows from Lemma \ref{L3.1}, the Sobolev embedding $H^1\left(\mathbb{T}^N\right)\hookrightarrow L^4\left(\mathbb{T}^N\right)$ and $H^2\left(\mathbb{T}^N\right)\hookrightarrow L^\infty\left(\mathbb{T}^N\right)$, (\ref{qqq.3.39}) and the Cauchy inequality that
	\begin{align*} |\tilde I_4|&\leq\dfrac{1}{\varepsilon^2}\|\nabla_i\nabla_j\rho_t\|\left(\|\nabla_i\nabla_j\left(v\cdot\nabla\rho_t\right)-\left(v\cdot\nabla\right)\nabla_i\nabla_j\rho_t\|+\|\nabla_i\nabla_j\left(v_t\cdot\nabla\rho\right)\|
		\right.\nonumber\\
		&\left.\quad+\|\nabla_i\nabla_j\left(\xi\nabla\cdot u_t\right)-\xi\nabla\cdot\nabla_i\nabla_ju_t\|+\|\nabla_i\nabla_j(\xi_t\nabla\cdot u)\|\right)
		\nonumber\\
		&\leq\dfrac{1}{\varepsilon^2}\|\nabla^2\rho_t\|\left(\|\nabla_i\nabla_jv\cdot\nabla\rho_t+\nabla_j\nabla\rho_t\cdot\nabla_iv+\nabla_jv\cdot\nabla_i\nabla\rho_t\|+\|\nabla^2\left(v_t\cdot\nabla\rho\right)\|
		\right.\nonumber\\
		&\left.\quad+\|\nabla_i\xi\nabla_j\left(\nabla\cdot u_t\right)+\nabla_i\nabla_j\xi\left(\nabla\cdot u_t\right)+\nabla_j\xi\nabla_i\left(\nabla\cdot u_t\right)\|+\|\nabla^2\left(\xi_t\cdot\nabla u\right)\|\right)
		\nonumber\\
		&\leq C\dfrac{1}{\varepsilon}\left\|\dfrac{1}{\varepsilon}\rho_t\right\|_2
		\left(\|\nabla^2v\|_{L^4}\|\nabla\rho_t\|_{L^4}+\|\nabla v\|_\infty\|\nabla^2\rho_t\|+\|v_t\|_\infty\|\nabla\rho\|_2+\|v_t\|_2\|\nabla\rho\|_\infty
		\right.\nonumber\\
		&\left.\quad+\|\nabla^2\xi\|_{L^4}\|\nabla\cdot u_t\|_{L^4}+\|\nabla\xi\|_\infty\|\nabla\left(\nabla\cdot u_t\right)\|+\|\xi_t\|_\infty\|\nabla u\|_2+\|\xi_t\|_2\|\nabla u\|_\infty\right)
		\nonumber\\
		&\leq C\dfrac{1}{\varepsilon}\left\|\dfrac{1}{\varepsilon}\rho_t\right\|_2
		\left(\|\nabla^2v\|_1\|\nabla\rho_t\|_1+\|\nabla v\|_2\|\nabla^2\rho_t\|+\|v_t\|_2\|\nabla\rho\|_2+\|\nabla\rho\|_2 \|v_t\|_2
		\right.\nonumber\\
		&\left.\quad+\|\nabla^2\xi\|_1\|\nabla\cdot u_t\|_1+\|\nabla\xi\|_2\|\nabla\left(\nabla\cdot u_t\right)\|+\|\xi_t\|_2\|\nabla u\|_2+\|\nabla u\|_2\|\xi_t\|_2\right)
		\nonumber\\
		&\leq C\left(\tilde K_1+\tilde K_2\right)\left(\left\|\dfrac{1}{\varepsilon}\rho_t\right\|_2^2+\|u_t\|_2^2\right)+C.
	\end{align*}
A simple calculation indicates that
$$
|\tilde I_7|=-\dfrac{1}{\varepsilon^2}\int_{\mathbb{T}^N}\xi\left\{\left[\nabla_i\nabla_j\left(\dfrac{P'\left(\xi\right)}{\xi}\nabla\rho_t\right)-\dfrac{P'\left(\xi\right)}{\xi}\nabla\nabla_i\nabla_j\rho_t\right]
+\nabla_i\nabla_j\left[\left(\dfrac{P'\left(\xi\right)}{\xi}\right)_t\nabla\rho\right]\right\}\cdot\nabla_i\nabla_ju_t.
$$
 Based on the calculation as follows:
	\begin{align*}	&\nabla_i\nabla_j\left(\dfrac{P'\left(\xi\right)}{\xi}\nabla\rho_t\right)-\dfrac{P'\left(\xi\right)}{\xi}\nabla\nabla_i\nabla_j\rho_t \nonumber\\\nonumber
		&=\left(\dfrac{P'''\left(\xi\right)}{\xi}-\dfrac{2P''\left(\xi\right)}{\xi^2}+\dfrac{2 P'\left(\xi\right)}{\xi^3}\right)\nabla_i\xi\nabla_j\xi\nabla\rho_t+\dfrac{\xi P''\left(\xi\right)-P'\left(\xi\right)}{\xi^2}\nabla_i\nabla_j\xi\nabla\rho_t\nonumber\\\nonumber
		&\quad+\dfrac{\xi P''\left(\xi\right)-P'\left(\xi\right)}{\xi^2}\nabla_j\xi\nabla\nabla_i\rho_t
		+\dfrac{\xi P''\left(\xi\right)-P'\left(\xi\right)}{\xi^2}\nabla_i\xi\nabla\nabla_j\rho_t,\nonumber
		\end{align*}
one deduces
	\begin{align*}
		|\tilde I_7|&\leq C\dfrac{1}{\varepsilon^2}\|u_t\|_2\left(\|\nabla\xi\|_\infty^2\|\nabla\rho_t\|+\|\nabla^2\xi\|_{L^4}\|\nabla\rho_t\|_{L^4} +\|\nabla\xi\|_\infty\|\nabla^2\rho_t\|
		\right.\nonumber\\
		&\left.\quad+\|\xi_t\|_\infty\|\nabla\rho\|_2+\|\nabla\rho\|_\infty\|\xi_t\|_2\left(1+\|\nabla\xi\|_1\right)\right)
		\nonumber\\
		&\leq C\left(\tilde K_1^2+\tilde K_2^2\right)\left(\left\|\dfrac{1}{\varepsilon}\rho_t\right\|_2^2+\|u_t\|_2^2\right)+C.
		\end{align*}
Finally, for the term $\tilde I_8$, we clearly write as
\begin{align*}
		\tilde I_8&=\int_{\mathbb{T}^N}\xi \left[\nabla_i\nabla_j\partial_t\left(\dfrac{\nu(\xi,\psi)}{\xi}\Delta u\right)-
	\dfrac{\nu(\xi,\psi)}{\xi}\Delta \nabla_i\nabla_j u_t\right] \vdots \nabla_i\nabla_j  u_t \nonumber\\
&\quad+\int_{\mathbb{T}^N}\xi \left[\nabla_i\nabla_j\partial_t\left(\dfrac{\eta(\xi,\psi)}{\xi}\nabla(\nabla\cdot u)\right)-
	\dfrac{\eta(\xi,\psi)}{\xi} \nabla_i\nabla_j\nabla(\nabla\cdot u_t)\right] \vdots \nabla_i\nabla_j  u_t=I_8^1+I_8^2,
	\end{align*}
	and then it follows from Lemma \ref{L3.1} that
	\begin{align}\label{IIIq.I-81}
		I_8^1&=\int_{\mathbb{T}^N}\xi \left[\nabla_i\nabla_j\left(\partial_t\left(\dfrac{\nu(\xi,\psi)}{\xi}\right)\Delta u\right)+\nabla_i\nabla_j\left(\dfrac{\nu(\xi,\psi)}{\xi}\Delta u_t\right)-
	\dfrac{\nu(\xi,\psi)}{\xi}\Delta \nabla_i\nabla_j u_t\right] \vdots \nabla_i\nabla_j  u_t \nonumber\\
&=\int_{\mathbb{T}^N}\xi \left\{\nabla_i\nabla_j\left[\left(\nu_\xi(\xi,\psi)\xi_t+\nu_\psi(\xi,\psi)\psi_t-\dfrac{\nu(\xi,\psi)\xi_t}{\xi^2}\right)\Delta u\right]
\right.\nonumber\\
&\left.\quad
+\nabla_i\nabla_j\left(\dfrac{\nu(\xi,\psi)}{\xi}\right)\Delta u_t
+\nabla_j\left(\dfrac{\nu(\xi,\psi)}{\xi}\right)\nabla_i\Delta u_t
+\nabla_i\left(\dfrac{\nu(\xi,\psi)}{\xi}\right)\nabla_j\Delta u_t
\right\} \vdots \nabla_i\nabla_j  u_t \nonumber\\
&\le C\left[\left(\left\|\nu_\xi(\xi,\psi)-\dfrac{\nu(\xi,\psi)}{\xi^2}\right\|_{2}\|\xi_t\|_{2}
+\left\|\nu_\psi(\xi,\psi)\right\|_{2}\|\psi_t\|_{2}\right)\|\Delta u\|_{2}\right.\nonumber\\
&\left.\quad+\left\|\nabla^2\left(\dfrac{\nu(\xi,\psi)}{\xi}\right)\right\|_{L^4}\|\Delta u_t\|_{L^4}
+\left\|\nabla\left(\dfrac{\nu(\xi,\psi)}{\xi}\right)\right\|_{\infty}\|\nabla\Delta u_t\|\right]\|u_t\|_{2}\nonumber\\
&\le C(1+\|\nabla\xi\|_{1}+\left\|\nabla\psi\right\|_{1})(\|\xi_t\|_{2}+\|\psi_t\|_{2})\|\nabla u\|_{3}\|u_t\|_{2}\nonumber\\
&\quad+\left\|\nabla^2\left(\dfrac{\nu(\xi,\psi)}{\xi}\right)\right\|_{1}\|\Delta u_t\|_{1}
+\left(\left\|\nabla\left(\dfrac{\nu(\xi,\psi)}{\xi}\right)\right\|_{2}\|\nabla\Delta u_t\|\right)\|u_t\|_{2}\nonumber\\
&\le C(1+\|\nabla\xi\|_{1}+\left\|\nabla\psi\right\|_{1})(\|\xi_t\|_{2}+\|\psi_t\|_{2})\|\nabla u\|_{3}\|u_t\|_{2}\nonumber\\
&\quad+C(\|\nabla\xi\|_{2}+\left\|\nabla\psi\right\|_{2})\|\nabla u_t\|_{2}\|u_t\|_{2}\nonumber\\
&\le C(\tau)\tilde K_1\tilde K_2\|u_t\|_{2}^2+C\|\nabla u\|_{3}^2+\tau\|\nabla u_t\|_{2}^2.
	\end{align}
After applying a similar way as (\ref{IIIq.I-81}) to $\tilde I_8^2$, we deduce the estimate of $\tilde I_8$ that
\begin{align}
\tilde I_8
&\le C(\tau)\tilde K_1\tilde K_2\|u_t\|_{2}^2+C\left(\|\nabla u\|_{3}^2+\|\nabla\cdot u\|_{3}^2\right)+\tau\left(\|\nabla u_t\|_{2}^2+\|\nabla\cdot u_t\|_{2}^2\right).
\end{align}
The constraints of the initial data and (\ref{q.2.1}) give
	\begin{align}\label{qqq.3.51}
		&\left\|\frac{1}{\varepsilon}\rho_t\left(x,0\right)\right\|_{s-1}^2+\|u_t\left(x,0\right)\|_{s-1}^2
		\nonumber\\
		&\leq C\left(\left\|\frac{1}{\varepsilon}\left(\left(u_0+\overline u_0^\varepsilon\right)\cdot\nabla\right)\overline\rho_0^\varepsilon\right\|_{s-1}^2
		+\left\|\frac{1}{\varepsilon}\left(\overline\rho_0^\varepsilon+1\right)\nabla\cdot\overline u_0^\varepsilon\right\|_{s-1}^2+\left\|\left(\overline\rho_0^\varepsilon+1\right)^{-1}\frac{1}{\varepsilon^2}\nabla\overline\rho_0^\varepsilon\right\|_{s-1}^2
		\right.\nonumber\\
		&\left.\quad+\|\left(\left(u_0+\overline u_0^\varepsilon\right)\cdot\nabla\right)\left(u_0+\overline u_0^\varepsilon\right)\|_{s-1}^2+\|\left(\overline\rho_0^\varepsilon+1\right)^{-1}\nu(\overline\rho_0^\varepsilon+1,\phi_0+\overline\phi_0^\varepsilon)\Delta\left(u_0+\overline u_0^\varepsilon\right)\|_{s-1}^2
\right.\nonumber\\
&\left.\quad+\|\left(\overline\rho_0^\varepsilon+1\right)^{-1}\eta(\overline\rho_0^\varepsilon+1,\phi_0+\overline\phi_0^\varepsilon)\nabla\left(\nabla\cdot\overline u_0^\varepsilon\right)\|_{s-1}^2
		\right.\nonumber\\
		&\left.\quad+\|\left(\overline\rho_0^\varepsilon+1\right)^{-1}\Delta\left(\phi_0+\overline\phi_0^\varepsilon\right)
\nabla\left(\phi_0+\overline\phi_0^\varepsilon\right)\|_{s-1}^2\right)\leq C.
	\end{align}
Finally, summing over the above $\tilde I_1$-$\tilde I_{10}$ and combining with the Gronwall inequality, (\ref{qqq.3.19}) and (\ref{qqq.3.40}) yield
	\begin{align}\label{qqq.3.52}
		\left\|\frac{1}{\varepsilon}\rho_t\right\|_{s-1}^2\left(t\right)+\|u_t\|_{s-1}^2\left(t\right)+\nu_*\int_{0}^{t}\|\nabla u_t\|_{s-1}^2+\eta_*\int_{0}^{t}\|\nabla\cdot u_t\|_{s-1}^2\leq C,
	\end{align}
	for $t\in [0, T_0]$, provided that $T_0$ is small enough such that $T_0\leq T_3^*:={\rm min}\left\{T_2^*, \left(\tilde K_1^2+\tilde K_2^2\right)^{-1}\right\}$.

It remains to show the first inequality of (\ref{qqq.3.1}). It suffices to show  $\left\|\frac{1}{\varepsilon}\left(\rho-1\right)\right\|_s+\|u-u_0\|_s+\|\phi-\phi_0\|_s\leq c_0^{-1}\kappa$ by the Sobolev inequality, where $c_0$ is the Sobolev constant. Let\ $\overline\rho=\rho-1,\ \overline u=u-u_0\ {\rm and}\ \overline\phi=\phi-\phi_0$.

 When $\phi=\overline\phi+\phi_0$ and $|\alpha_1| \leq s$, it follows from (\ref{qqq.3.4})
 \begin{align}\label{qqq.3.4-2}
\phi_t-\dfrac1{\xi^2}\Delta(\phi+\phi_0)=- v\cdot\nabla\psi-\dfrac1{\xi}(\psi^3-\psi).
\end{align}
Operating $\nabla^{\alpha_1}$ to (\ref{qqq.3.4-2}) and multiplying the result by $\nabla^{\alpha_1}\overline\phi$, we obtain after using integration by parts that
	\begin{align}\label{qqq.3.53-1}
			&\dfrac12\dfrac{\rm d}{{\rm d}t}\int_{\mathbb{T}^N}|\nabla^{\alpha_1}\overline\phi|^2+\int_{\mathbb{T}^N}\dfrac1{\xi^2}|\nabla^{\alpha_1}\nabla\overline\phi|^2
\nonumber\\
&=\int_{\mathbb{T}^N}\left[\nabla^{\alpha_1}\left(\dfrac1{\xi^2}\Delta\overline\phi\right)-\dfrac1{\xi^2}\nabla^{\alpha_1}\Delta\overline\phi\right]
\nabla^{\alpha_1}\overline\phi
-\int_{\mathbb{T}^N}\nabla^{\alpha_1}(v\cdot\nabla\psi)\nabla^{\alpha_1}\overline\phi
\nonumber\\
&\quad-\int_{\mathbb{T}^N}\nabla^{\alpha_1}\left(\dfrac1{\xi}(\psi^3-\psi)\right)\nabla^{\alpha_1}\overline\phi
-\int_{\mathbb{T}^N}\nabla\left(\dfrac1{\xi^2}\right)\cdot\nabla^{\alpha_1}\nabla\overline\phi\nabla^{\alpha_1}\overline\phi
+\int_{\mathbb{T}^N}\nabla^{\alpha_1}\left(\dfrac1{\xi^2}\Delta\phi_0\right)\nabla^{\alpha_1}\overline\phi\nonumber\\
&\leq C\left(\left\|\nabla\left(\dfrac1{\xi^2}\right)\right\|_\infty\|\Delta\overline\phi\|_{s-1}+\|\Delta\overline\phi\|_\infty\left\|\nabla\left(\dfrac1{\xi^2}\right)\right\|_{s-1}
+\|v\|_\infty\|\nabla\psi\|_s+\|\nabla\psi\|_\infty\|v\|_s\right.\nonumber\\
        &\left.\quad+\left\|\dfrac1{\xi}\right\|_\infty\|\psi^3-\psi\|_s+\|\psi^3-\psi\|_\infty\left\|\dfrac1{\xi}\right\|_s+\left\|\nabla\left(\dfrac1{\xi^2}\right)\right\|_\infty
        \|\nabla^{\alpha_1}\nabla\overline\phi\|\right.\nonumber\\
        &\left.\quad+\left\|\dfrac1{\xi^2}\right\|_\infty\|\Delta\phi_0\|_s+\|\Delta\phi_0\|_\infty\left\|\dfrac1{\xi^2}\right\|_s\right)
        \|\nabla^{\alpha_1}\overline\phi\|\nonumber\\
&\leq C\left(\tilde K_1^{\frac1{2}}\|\nabla\overline\phi\|_s+\tilde K_1^2\right)\|\nabla^{\alpha_1}\overline\phi\| \leq C\tilde K_1^4\|\nabla^{\alpha_1}\overline\phi\|^2+\tau\|\nabla\overline\phi\|_s+C.
		\end{align}
Summing over $\alpha_1$ and taking $\tau$ small enough, we get
		\begin{align*}
			&\frac{\rm d}{{\rm d}t}\sum_{|\alpha_1|\leq{s}}\|\nabla^{\alpha_1}\overline\phi\|^2+\sum_{|\alpha_1|\leq{s}}\|\nabla^{\alpha_1}\nabla\overline\phi\|^2
			\leq C\tilde K_1^4\|\overline\phi\|_s^2+C.
		\end{align*}
Since (\ref{qqq.2.6}) implies that $\|\overline\phi(x,0)\|_{s}^2\le \kappa_0^2\varepsilon^2$, we conclude from the Gronwall inequality that
	\begin{align}\label{qqq.3.55}
		\|\overline\phi\|_s^2
		&\leq \exp\left(C\tilde K_1^4T_0\right)(\kappa_0^2\varepsilon^2+CT_0)\leq c_0^{-1}\kappa,
	\end{align}
where we have chosen $\varepsilon$ and  $T_0(< T_3^*)$ sufficiently small such that (\ref{qqq.3.55}).

	Similarly as in the proof of (\ref{qqq.3.37}) and (\ref{qqq.3.53-1}), we get
	\begin{align}\label{qqq.3.56}
		&\dfrac{\rm d}{{\rm d}t}\sum_{|\alpha_1|\leq{s}}\int_{\mathbb{T}^N}\left(\dfrac{P'\left(\xi\right)}{\xi}\left|\dfrac{1}{\varepsilon} \nabla^{\alpha_1}\overline\rho\right|^2+\xi| \nabla^{\alpha_1}\overline u|^2\right)\nonumber\\
		&\quad+\sum_{|\alpha_1|\leq{s}}\nu_*\int_{\mathbb{T}^N}| \nabla^{\alpha_1} \nabla \overline u|^2 +\sum_{|\alpha_1|\leq{s}}\eta_*\int_{\mathbb{T}^N}| \nabla^{\alpha_1} (\nabla\cdot \overline u)|^2  \nonumber\\
		&\leq C\left(  \tilde K_1 +\tilde K_2\right)\left(\left\|\frac{1}{\varepsilon}\overline\rho\right\|_s^2+\left\|\overline u\right\|_s^2\right)+C\left\|\left(v\cdot\nabla\right)u_0\right\|_s^2+C\left\|\nabla u_0\right\|_s^2+C  \nonumber\\\nonumber
		&\quad+C\sum_{|\alpha_1|\leq{s}}\left\|\nabla^{\alpha_1}\left(\frac{\nu(\xi,\psi) }{\xi}\Delta u_0\right)-\frac{  \nu(\xi,\psi) }{\xi}\nabla^{\alpha_1}\Delta u_0\right\|^2  \nonumber\\
		&\leq C\left(\tilde K_1+\tilde K_2\right)\left(\left\|\frac{1}{\varepsilon}\overline\rho\right\|_s^2+\left\|\overline u\right\|_s^2\right)
+C\left\|\nabla\left(\frac{\nu(\xi,\psi) }{\xi}\right)\right\|_{s-1}^2\|\Delta u_0\|_{s-1}^2+C\tilde K_1\nonumber\\
&\leq C\left(\tilde K_1+\tilde K_2\right)\left(\left\|\frac{1}{\varepsilon}\overline\rho\right\|_s^2+\left\|\overline u\right\|_s^2\right)+C\left(\|\nabla\xi\|_{s-1}+\|\nabla\psi\|_{s-1}\right)^2+C\tilde K_1\nonumber\\
&\leq C\left(\tilde K_1+\tilde K_2\right)\left(\left\|\frac{1}{\varepsilon}\overline\rho\right\|_s^2+\left\|\overline u\right\|_s^2\right)+C\tilde K_1.
	\end{align}
	Then the Gronwall inequality, (\ref{qqq.2.3}) and (\ref{qqq.2.6}) give
	\begin{align}\label{qqq.3.57}
		\left\|\frac{1}{\varepsilon}\overline\rho\right\|_s^2+\|\overline u\|_s^2  \leq \exp\left(C\left(\tilde K_1+\tilde K_2\right)T_0\right)\left(\kappa^2\varepsilon^2+\tilde K_1T_0\right)
		<c_0^{-1}\kappa,
	\end{align}
	where we have chosen $T_0\left(<T_3^*\right)$ and $\varepsilon$ small enough (\ref{qqq.3.55}) and (\ref{qqq.3.57}) hold.
The proof of Lemma $\ref{L3.2-2}$ is completed.
\hfill$\Box$

The following Lemma indicate that $\Lambda$ is a contractive map. Its proof is similar to the Navier-Stokes/Cahn-Hilliard system, even rather simpler than it  and we will not elaborate further.
\begin{lemma} \label{L3.3-2}
	Under the assumptions in Theorem \ref{th1-2}, the map $\Lambda:V\longrightarrow U$ is a contraction in the sense that
	%
	\begin{align}\label{qqq.3.53}
		&\sup_{0\leq t\leq T_0}	\left(\left\|\dfrac{1}{\varepsilon}\left(\rho-\hat{\rho}\right)\right\|^2+\left\|u-\hat{u}\right\|^2+\left\|\phi-\hat{\phi}\right\|_1^2\right)+\int_{0}^{t}\left(\|u-\hat{u}\|_1^2+\|\phi-\hat{\phi}\|_2^2\right)  \nonumber\\
		&\leq\gamma_3\sup_{0\leq t\leq T_0}\left(\left\|\dfrac{1}{\varepsilon}\left(\xi-\hat{\xi}\right)\right\|^2+\left\|v-\hat{v}\right\|^2+\left\|\psi-\hat{\psi}\right\|_1^2\right).
	\end{align}
	for some $0<\gamma_3<1$, provided that $T_0$ is small enough.
	\end{lemma}

The following Lemma is necessary for the proof of Theorem \ref{th1-2}.
\begin{lemma}\label{L3.4-2}
	Consider the incompressible system of Navier-Stokes/Allen-Cahn (\ref{q.1.2}) with initial condition (\ref{qqq.2.5}) for $s\geq 2$. Then for $0<T\leq\infty$, there exists at most one strong solution $\left(u, \phi\right)\in\left \{\left(u, \phi\right): \left(u, \nabla\phi\right)\in L^{\infty}\left(0, T; H^s\right)\right \} $.
	\end{lemma}

By combining Lemma \ref{L3.2-2}-\ref{L3.4-2} and following the proof steps of Theorem \ref{th1}, Theorem \ref{th1-2} can be proven.

	\begin{remark}\label{R3.1-2}
	It follows
that
	\begin{align*}
		\left\|\dfrac{1}{\varepsilon^2}\nabla\rho^{\varepsilon}\right\|_{s-2}+\int_{0}^{t}\left\|\dfrac{1}{\varepsilon^2}\nabla\rho^{\varepsilon}\right\|_{s-1}^2\leq C,\quad t\in[0,T_0].
	\end{align*}
\end{remark}
	\begin{remark}\label{R3.2-2}
It follows
that
	\begin{align*}
		\left\|\dfrac{1}{\varepsilon}\nabla\cdot u^{\varepsilon}\right\|_{s-1}&\leq\left\|\dfrac{1}{\varepsilon}\rho_t^{\varepsilon}\right\|_{s-1}+\left\|\dfrac{1}{\varepsilon}\left(u^{\varepsilon}\cdot\nabla\right)\rho^{\varepsilon}\right\|_{s-1}+\left\|\dfrac{1}{\varepsilon}\left(\rho^{\varepsilon}-1\right)\nabla\cdot u^{\varepsilon}\right\|_{s-1}  \nonumber\\
		&\leq C,\quad t\in[0,T_0].
\end{align*}
\end{remark}
\allowdisplaybreaks
\subsection{Proof of Theorem \ref{th2-2}.}\label{subsec:3.3}
	\quad~~In this subsection, we devote ourselves to getting a priori energy estimates (\ref{qqq.2.11}) and (\ref{qqq.2.12}) with small initial displacements and small initial data. Then together with Theorem \ref{th1-2}, Theorem \ref{th2-2} can be proved by a standard procedure. Since the estimation of $\rho$ and
$u$ is not significantly different from subsection \ref{subsec:2.3}, we only give the details for estimate of $\phi$.

 First of all, we assume that
		\begin{align}\label{qqq.4.1}
		E_s\left(U\left(t\right)\right)&+\left\|\phi^2-1\right\|^2_s\nonumber\\
&+\int_0^t\left(\nu_*\left\|\nabla u\right\|^2_s+\eta_*\left\|\nabla\cdot u\right\|^2_s+\left\|\phi^2-1\right\|^2_{s+1}+\left\|\nabla\phi\right\|^2_{s+1}\right)\leq 4\left(\delta+\varepsilon^2\kappa_0^2\right)
\end{align}
	for $t\in [ 0,\ T^\varepsilon ] $, and then what we need to do is to prove the following desired estimates
	\begin{align}\label{qqq.4.2}
		E_s\left(U\left(t\right)\right)&+\left\|\phi^2-1\right\|^2_s
\nonumber\\
&+\int_0^t\left(\nu_*\left\|\nabla u\right\|^2_s+\eta_*\left\|\nabla\cdot u\right\|^2_s+\left\|\phi^2-1\right\|^2_{s+1}+\left\|\nabla\phi\right\|^2_{s+1}\right)\leq 3\left(\delta+\varepsilon^2\kappa_0^2\right)
	\end{align}
	for $t\in [ 0,T^\varepsilon ] $.
The (\ref{q.2.11}) follows by the standard continuity argument and the fact  $E_s\left(U(0)\right)<4\left(\delta+\varepsilon^2\kappa_0^2\right)$.

	 Firstly,
recalling Section \ref{sec:2} and making a similar (just a little different) argument to give the following fact:
	\begin{align}\label{qqq.4.4}
		&\left\|\dfrac{1}{\varepsilon}\rho_t\right\|_{s-1}^2+\left\|u_t\right\|_{s-1}^2+\left\|\phi_t\right\|_{s}^2+\int_{0}^{t}\left(\nu_*\left\|\nabla u_t\right\|_{s-1}^2+\eta_*\left\|\nabla\cdot u_t\right\|_{s-1}^2+\left\|\nabla\phi_t\right\|_{s}^2\right) \nonumber\\
		&\quad\leq C\left(1+t\right){\rm exp}Ct\leq C {\rm exp} Ct,
	\end{align}
	for $t\in[0,T^\varepsilon]$, and thus we conclude that Remarks \ref{R3.1-2} and \ref{R3.2-2} hold for $t\in [0,T^\varepsilon]$ by the same discussion in subsection \ref{subsec:3.2}.

Next, we begin to prove (\ref{qqq.4.2}).
On one hand, multiplying $(\ref{qqq.2.1})_3$ by $2\phi$, we obtain
\begin{align*}
(\phi^2-1)_t-\dfrac1{\rho^2}\Delta\left(\phi^2-1\right)+\dfrac2{\rho}\left(\phi^2-1\right)
=-u\cdot\nabla\left(\phi^2-1\right)-\dfrac2{\rho}\left(\phi^2-1\right)^2
-\dfrac2{\rho^2}|\nabla\phi|^2.
\end{align*}
Then operating the above equality by $\nabla^{\alpha_1}$ with $|\alpha_1|\le s$, we deduce
\begin{align*}
&\nabla^{\alpha_1}(\phi^2-1)_t-\dfrac1{\rho^2}\nabla^{\alpha_1}\Delta\left(\phi^2-1\right)+\dfrac2{\rho}\nabla^{\alpha_1}\left(\phi^2-1\right)
\nonumber\\
&=\left[\nabla^{\alpha_1}\left(\dfrac1{\rho^2}\Delta\left(\phi^2-1\right)\right)-\dfrac1{\rho^2}\nabla^{\alpha_1}\Delta\left(\phi^2-1\right)\right]
-\left[\nabla^{\alpha_1}\left(\dfrac2{\rho}\left(\phi^2-1\right)\right)-\dfrac2{\rho}\nabla^{\alpha_1}\left(\phi^2-1\right)\right]
\nonumber\\
&\quad-\nabla^{\alpha_1}\left(u\cdot\nabla\left(\phi^2-1\right)\right)-\nabla^{\alpha_1}\left(\dfrac2{\rho}\left(\phi^2-1\right)^2\right)
-\nabla^{\alpha_1}\left(\dfrac2{\rho^2}|\nabla\phi|^2\right).
\end{align*}
Multiplying the above equality by $\nabla^{\alpha_1}\left(\phi^2-1\right)$, we obtain after integrating over $\mathbb{T}^N$ that
\begin{align}\label{qqq.4.7}
&\dfrac{1}{2}\dfrac{\rm d}{{\rm d}t}\int_{\mathbb{T}^N}\left|\nabla^{\alpha_1}(\phi^2-1)\right|^2 +\int_{\mathbb{T}^N}\dfrac1{\rho^2}\nabla^{\alpha_1}\left|\nabla\left(\phi^2-1\right)\right|^2     +\int_{\mathbb{T}^N}\dfrac2{\rho}\left|\nabla^{\alpha_1}\left(\phi^2-1\right)\right|^2
\nonumber\\
&=\int_{\mathbb{T}^N}\left[\nabla^{\alpha_1}\left(\dfrac1{\rho^2}\Delta\left(\phi^2-1\right)\right)-\dfrac1{\rho^2}\nabla^{\alpha_1}\Delta\left(\phi^2-1\right)\right]\nabla^{\alpha_1}\left(\phi^2-1\right)
\nonumber\\
&\quad-\int_{\mathbb{T}^N}\left[\nabla^{\alpha_1}\left(\dfrac2{\rho}\left(\phi^2-1\right)\right)-\dfrac2{\rho}\nabla^{\alpha_1}\left(\phi^2-1\right)\right]\nabla^{\alpha_1}\left(\phi^2-1\right)
\nonumber\\
&\quad-\int_{\mathbb{T}^N}\nabla^{\alpha_1}\left(u\cdot\nabla\left(\phi^2-1\right)\right)\nabla^{\alpha_1}\left(\phi^2-1\right)
-\int_{\mathbb{T}^N}\nabla^{\alpha_1}\left(\dfrac2{\rho}\left(\phi^2-1\right)^2\right)\nabla^{\alpha_1}\left(\phi^2-1\right)
\nonumber\\
&\quad-\int_{\mathbb{T}^N}\nabla^{\alpha_1}\left(\dfrac2{\rho^2}|\nabla\phi|^2\right)\nabla^{\alpha_1}\left(\phi^2-1\right)
-\int_{\mathbb{T}^N}\nabla\left(\dfrac1{\rho^2}\right)\cdot\nabla^{\alpha_1}\nabla\left(\phi^2-1\right)\nabla^{\alpha_1}\left(\phi^2-1\right)=\sum_{i}^{6}A_i.
\end{align}
%
Next, by (\ref{qqq.4.1}), we estimate each $A_i(i=1,2,...6)$ as follows:
\begin{align}
A_1&\leq C\left[\left\|\nabla\left(\dfrac1{\rho^2}\right)\right\|_\infty\left\|\Delta\left(\phi^2-1\right)\right\|_{s-1}
+\left\|\Delta\left(\phi^2-1\right)\right\|_\infty\left\|\nabla\left(\dfrac1{\rho^2}\right)\right\|_{s-1}
\right] \left\|\nabla^{\alpha_1}\left(\phi^2-1\right)\right\|
\nonumber\\
&\leq C\varepsilon\left\|\dfrac1{\varepsilon}\nabla\rho\right\|_{s-1}\left\|\nabla\left(\phi^2-1\right)\right\|_{s}\left\|\phi^2-1\right\|_{s}
\nonumber\\
&\leq C(\tau)\varepsilon^2(\delta+\varepsilon^2\kappa_0)\left\|\nabla\left(\phi^2-1\right)\right\|_{s}^2+\tau\left\|\phi^2-1\right\|_{s}^2,
\\
A_2&\leq C\left[\left\|\nabla\left(\dfrac2{\rho}\right)\right\|_\infty\left\|\phi^2-1\right\|_{s-1}
+\left\|\phi^2-1\right\|_\infty\left\|\nabla\left(\dfrac2{\rho}\right)\right\|_{s-1}\right] \left\|\nabla^{\alpha_1}\left(\phi^2-1\right)\right\|
\nonumber\\
&\leq C(\tau)\varepsilon^2(\delta+\varepsilon^2\kappa_0)\left\|\phi^2-1\right\|_{s-1}^2+\tau\left\|\phi^2-1\right\|_{s}^2,
\\
A_3&\leq C\left[\|u\|_s\left\|\nabla\left(\phi^2-1\right)\right\|_\infty+\|u\|_\infty\left\|\nabla\left(\phi^2-1\right)\right\|_s
\right] \left\|\nabla^{\alpha_1}\left(\phi^2-1\right)\right\|
\nonumber\\
&\leq C(\tau)(\delta+\varepsilon^2\kappa_0)\left\|\nabla\left(\phi^2-1\right)\right\|_{s}^2+\tau\left\|\phi^2-1\right\|_{s}^2,
\\
A_4&\leq C\left[\left\|\dfrac2{\rho}\right\|_\infty\left\|\phi^2-1\right\|_s^2+\left\|\phi^2-1\right\|_\infty^2\left\|\dfrac2{\rho}\right\|_s
\right] \left\|\nabla^{\alpha_1}\left(\phi^2-1\right)\right\|
\nonumber\\
&\leq C\left(1+\varepsilon\left\|\dfrac1{\varepsilon}\nabla\rho\right\|_{s-1}\right)\left\|\phi^2-1\right\|_{s}^3
\nonumber\\
&\leq C(\delta+\varepsilon^2\kappa_0)^\frac12\left\|\phi^2-1\right\|_{s}^2,
\\
A_5&\leq C\left(\left\|\dfrac2{\rho^2}\right\|_\infty\|\nabla\phi\|_s^2+\|\nabla\phi\|_\infty^2\left\|\dfrac2{\rho^2}\right\|_s
\right) \left\|\nabla^{\alpha_1}\left(\phi^2-1\right)\right\|
\nonumber\\
&\leq C(\tau)\varepsilon^2(\delta+\varepsilon^2\kappa_0)\left\|\nabla\phi\right\|_{s}^2+\tau\left\|\phi^2-1\right\|_{s}^2,
\\
A_6&\leq C\left\|\nabla\left(\dfrac1{\rho^2}\right)\right\|_\infty\|\nabla\left(\phi^2-1\right)\|_s\left\|\nabla^{\alpha_1}\left(\phi^2-1\right)\right\|
\nonumber\\
&\leq C(\tau)\varepsilon(\delta+\varepsilon^2\kappa_0)\left\|\nabla\left(\phi^2-1\right)\right\|_{s}^2+\tau\left\|\phi^2-1\right\|_{s}^2.
\end{align}
 On the other hand, multiplying $(\ref{qqq.2.1})$ by $\rho$, operating $\nabla_i\nabla^{\alpha_1}(|\alpha_1|\le s)$ to the result, we obtain
\begin{align*}
\rho\nabla_i\nabla^{\alpha_1}\phi_t-\dfrac{1}{\rho}\nabla_i\nabla^{\alpha_1}\Delta\phi
&=-\left[\nabla_i\nabla^{\alpha_1}(\rho\phi_t)-\rho\nabla_i\nabla^{\alpha_1}\phi_t\right]
+\left[\nabla_i\nabla^{\alpha_1}\left(\dfrac{1}{\rho}\Delta\phi\right)-\dfrac{1}{\rho}\nabla_i\nabla^{\alpha_1}\Delta\phi\right]
\nonumber\\
&\quad-\nabla_i\nabla^{\alpha_1}\left(\rho u\cdot\nabla\phi\right)
-\nabla_i\nabla^{\alpha_1}\left(\phi^3-3\phi\right)
-2\nabla_i\nabla^{\alpha_1}\phi.
\end{align*}
Taking the $L^2$ inner product of the above equation 
with $\nabla_i\nabla^{\alpha_1}\phi$
and using integration by parts and then we get
\allowdisplaybreaks
		\begin{align}\label{qqq.4.14}
		&\dfrac{1}{2}\dfrac{\rm d}{{\rm d}t}\int_{\mathbb{T}^N}\rho|\nabla_i\nabla^{\alpha_1}\phi|^2
+\int_{\mathbb{T}^N}\dfrac{1}{\rho}|\nabla^{\alpha_1}\Delta\phi|^2
+2\int_{\mathbb{T}^N}|\nabla_i\nabla^{\alpha_1}\phi|^2\nonumber\\
		&=-\int_{\mathbb{T}^N}\left[\nabla_i\nabla^{\alpha_1}(\rho\phi_t)-\rho\nabla_i\nabla^{\alpha_1}\phi_t\right]\cdot\nabla_i\nabla^{\alpha_1}\phi
+\int_{\mathbb{T}^N}\left[\nabla_i\nabla^{\alpha_1}\left(\dfrac{\Delta\phi}{\rho}\right)
-\dfrac{\nabla_i\nabla^{\alpha_1}\Delta\phi}{\rho}\right]\cdot\nabla_i\nabla^{\alpha_1}\phi
\nonumber\\
&\quad-\int_{\mathbb{T}^N}\nabla_i\nabla^{\alpha_1}\left(\rho u\cdot\nabla\phi\right)\cdot\nabla_i\nabla^{\alpha_1}\phi
-\int_{\mathbb{T}^N}\nabla_i\nabla^{\alpha_1}\left(\phi^3-3\phi\right)\cdot\nabla_i\nabla^{\alpha_1}\phi
\nonumber\\
&\quad-\int_{\mathbb{T}^N}\nabla_i\left(\dfrac{1}{\rho}\right)\cdot\nabla_i\nabla^{\alpha_1}\phi\nabla^{\alpha_1}\Delta\phi
+\dfrac{1}{2}\int_{\mathbb{T}^N}\rho_t|\nabla_i\nabla^{\alpha_1}\phi|^2
=\sum_{i=1}^{6}F_i.
		\end{align}
	 Next, by (\ref{qqq.4.1}) and using integration by parts, we give the estimates of $F_i(i=1,...,6)$ as follows:
	\begin{align}
		|F_1|&=\int_{\mathbb{T}^N}\left[\nabla^{\alpha_1}(\rho\phi_t)-\rho\nabla^{\alpha_1}\phi_t\right]\nabla^{\alpha_1}\Delta\phi
-\int_{\mathbb{T}^N}\nabla\rho\cdot\nabla\nabla^{\alpha_1}\phi\nabla^{\alpha_1}\phi_t\nonumber\\
&\leq C\left\|\Delta\phi\right\|_s\left(\left\|\nabla\rho\right\|_\infty\left\|\phi_t\right\|_{s-1}+\left\|\nabla\rho\right\|_{s-1}\left\|\phi_t\right\|_\infty\right)
+C\left\|\nabla\rho\right\|_\infty\left\|\nabla\phi\right\|_s\left\|\phi_t\right\|_s\nonumber  \\
 &\leq\tau\left\|\nabla\phi\right\|_{s+1}^2+C(\tau)\varepsilon^2\left\|\dfrac1{\varepsilon}\nabla\rho\right\|_{s-1}^2\left\|\phi_t\right\|_s^2
 \nonumber\\
 &\leq\tau\left\|\nabla\phi\right\|_{s+1}^2+C(\tau)\varepsilon^2\left(\left\|u\cdot\nabla\phi\right\|_s^2+\left\|\dfrac{1}{\rho^2}\Delta\phi\right\|_s^2
 +\left\|\dfrac{1}{\rho}\left(\phi^2-1\right)\phi\right\|_s^2\right)
 \nonumber\\
 &\leq\tau\left\|\nabla\phi\right\|_{s+1}^2+C(\tau)\varepsilon^2\left(\left\|\nabla\phi\right\|_s^2+\left\|\Delta\phi\right\|_s^2
 +\left\|\phi^2-1\right\|_s^2\right),\\
|F_2|
		&\leq C\left\|\nabla\phi\right\|_{s}\left(\left\|\nabla\left(\dfrac1{\rho}\right)\right\|_\infty\left\|\Delta\phi\right\|_{s}
+\left\|\nabla\left(\dfrac1{\rho}\right)\right\|_{s}\left\|\Delta\phi\right\|_\infty\right)\nonumber  \\
&\leq C\varepsilon\left\|\dfrac1{\varepsilon}\nabla\rho\right\|_{s}\left\|\nabla\phi\right\|_{s}\left\|\Delta\phi\right\|_{s}
\nonumber  \\
&\leq \tau\left\|\Delta\phi\right\|_s^2+C(\tau)\varepsilon^2\left(\delta+\varepsilon^2\kappa_0^2\right)\left\|\nabla\phi\right\|_{s}^2, \\
		|F_3|&\leq C\left\|\rho u\cdot \nabla\phi\right\|_s\left\|\Delta\phi\right\|_s\leq C\left\|\rho \right\|_s\left\|u\right\|_s\left\| \nabla\phi\right\|_s\left\|\Delta\phi\right\|_s\nonumber \\
&\leq \tau\left\|\Delta\phi\right\|_s^2+C(\tau)\left(\delta+\varepsilon^2\kappa_0^2\right)\left\|\nabla\phi\right\|_{s}^2, \\
|F_4|&=-3\int_{\mathbb{T}^N}\nabla^{\alpha_1}\left[\left(\phi^2-1\right)\nabla\phi\right]\cdot\nabla\nabla^{\alpha_1}\phi
\nonumber \\
&\leq C\left\|\phi^2-1\right\|_s\left\|\nabla\phi\right\|_s\left\|\nabla\phi\right\|_s\nonumber\\
&\leq \tau\left\|\nabla\phi\right\|_s^2+C(\tau)\left(\delta+\varepsilon^2\kappa_0^2\right)\left\|\nabla\phi\right\|_{s}^2,\\
|F_5|&\leq C\varepsilon\left\|\dfrac1{\varepsilon}\nabla\rho\right\|_{s}\left\|\nabla\phi\right\|_s\left\|\Delta\phi\right\|_s
\nonumber\\
&\leq \tau\left\|\Delta\phi\right\|_s^2+C(\tau)\varepsilon^2\left(\delta+\varepsilon^2\kappa_0^2\right)\left\|\nabla\phi\right\|_{s}^2,\\
|F_6|&\leq C\left\|\rho_t\right\|_\infty\left\|\nabla\phi\right\|_s^2\leq C\left\|{\rm div}(\rho u)\right\|_2\left\|\nabla\phi\right\|_s^2\leq C\varepsilon^2\left(\delta+\varepsilon^2\kappa_0^2\right)\left\|\nabla\phi\right\|_s^2.
\label{qqq.4.19}
	\end{align}
 Putting (\ref{qqq.4.7})-(\ref{qqq.4.19}) together, summing over $\alpha_1$, we get
		\begin{align}\label{qqq.4.20}
		&\dfrac{1}{2}\sum_{|\alpha_1|\le s}\dfrac{\rm d}{{\rm d}t}\left(\int_{\mathbb{T}^N}\rho|\nabla_i\nabla^{\alpha_1}\phi|^2+\int_{\mathbb{T}^N}\left|\nabla^{\alpha_1}(\phi^2-1)\right|^2 \right) +\sum_{|\alpha_1|\le s}\int_{\mathbb{T}^N}\dfrac1{\rho}\nabla^{\alpha_1}\left|\nabla\left(\phi^2-1\right)\right|^2     \nonumber\\
&\quad+\sum_{|\alpha_1|\le s}\int_{\mathbb{T}^N}\dfrac2{\rho}\left|\nabla^{\alpha_1}\left(\phi^2-1\right)\right|^2
+\sum_{|\alpha_1|\le s}\int_{\mathbb{T}^N}\dfrac{1}{\rho}|\nabla^{\alpha_1}\Delta\phi|^2
+2\sum_{|\alpha_1|\le s}\int_{\mathbb{T}^N}|\nabla^{\alpha_1}\nabla\phi|^2
\nonumber\\
&\leq\tau\left(\left\|\phi^2-1\right\|_{s}^2+\left\|\nabla\phi\right\|_{s+1}^2\right)
+C(\tau)(\delta+\varepsilon^2\kappa_0)^\frac12\left(\left\|\phi^2-1\right\|_{s+1}^2
+\left\|\nabla\phi\right\|_s^2\right)\nonumber\\
&\quad+C\varepsilon^2\left[\left\|\nabla\phi\right\|_s^2+\left\|\Delta\phi\right\|_s^2
 +\left\|\phi^2-1\right\|_s^2\right].
	\end{align}
	Then taking $\varepsilon$ and $\tau$ small enough, we obtain after integrating $(\ref{qqq.4.20})$ over $[0,t]$ that
		\begin{align*}
		&\left\|\phi^2-1\right\|_{s}^2(t)+\left\|\nabla\phi\right\|_s^2\left(t\right)
+\int_{0}^{t}\left(\left\|\phi^2-1\right\|_{s+1}^2
+\left\|\Delta\phi\right\|_s^2+\left\|\nabla\phi\right\|_{s+1}^2\right){\rm d}s
\nonumber\\
&\leq\left\|\left(\phi^2-1\right)(t,0)\right\|_{s}^2+\left\|\nabla\phi\left(x,0\right)\right\|_s^2\leq 2\left(\delta+\varepsilon^2\kappa_0^2\right),
	\end{align*}
	for $t\in[0,T^\varepsilon]$.

Furthermore, for $\rho$ and $u$, we derive
\begin{align*}
&\dfrac{1}{2}\dfrac{\rm d}{{\rm d}t}\sum_{|\alpha_1|\leq s}\int_{\mathbb{T}^N}\left(\dfrac{1}{\varepsilon^2}\dfrac{P'\left(\rho\right)}{\rho}|\nabla^{\alpha_1}\left(\rho-1\right)|^2+\rho|\nabla^{\alpha_1}u|^2\right)\nonumber
\nonumber\\
&\quad+\nu_*\sum_{|\alpha_1|\leq s}\int_{\mathbb{T}^N}|\nabla^{\alpha_1}\nabla u|^2+\eta_*\sum_{|\alpha_1|\leq s}\int_{\mathbb{T}^N}|\nabla^{\alpha_1}(\nabla\cdot u)|^2 \nonumber\\
		&\leq C\varepsilon\left(\left\|u\right\|_s^2+\left\|\dfrac{1}{\varepsilon}(\rho-1)\right\|_{s}^2\right)+C\varepsilon\left\|\dfrac{1}{\varepsilon^2}\nabla\rho\right\|_{s-1}^2
+C\left(\varepsilon_0+\varepsilon^2\kappa_0^2\right)\left\|\nabla^2\phi\right\|_s^2.
\end{align*}
Then it is easy to complete the proof of Theorem \ref{th2-2}, similar to subsection \ref{subsec:2.3}.
\allowdisplaybreaks
\subsection{Proof of Theorem \ref{th3-2}}
Proceeding the similar steps in subsection \ref{subsec:2.4}, we can make an argument to give the following fact
	\begin{align}\label{qqq.5.22}
		\left\|u^\varepsilon-u\right\|^2+\left\|\phi^\varepsilon-\phi\right\|_1^2
+\int_{0}^{t}\left(\nu(\phi)\left\|\nabla\left(u^\varepsilon-u\right)\right\|^2
+\left\|\nabla\left(\phi^\varepsilon-\phi\right)\right\|^2+\left\|\mu^\varepsilon-\mu\right\|^2\right) \leq C\varepsilon.
	\end{align}
Furthermore, subtracting $(\ref{q.1.1})_3$ from $(\ref{q.1.2})_3$, we obtain after using $(\ref{q.1.1})_4$ and $(\ref{q.1.2})_4$,that
	\begin{align}\label{qqq.5.24}
		\left(\phi^\varepsilon-\phi\right)_t-\Delta\left(\phi^\varepsilon-\phi\right)&=-\left(u^\varepsilon\cdot\nabla\right)\left(\phi^\varepsilon-\phi\right)
-\left(u^\varepsilon-u\right)\cdot\nabla\phi+\left(\dfrac1{\left(\rho^\varepsilon\right)^2}-1\right)\Delta\phi^\varepsilon\nonumber\\
&\quad-\left(\dfrac1{\rho^\varepsilon}-1\right)\left(\left(\phi^\varepsilon\right)^3-\phi^\varepsilon\right)-\left(\left(\phi^\varepsilon\right)^3-\phi^3\right)+\left(\phi^\varepsilon-\phi\right).
		\end{align}
	Then the parabolic theory implies that %
\begin{align}\label{qqq.5.25}
		\left\|\phi^\varepsilon-\phi\right\|_2^2&\leq C\left\|u^\varepsilon\right\|_\infty^2\left\|\nabla\left(\phi^\varepsilon-\phi\right)\right\|^2+C\left\|\nabla\phi\right\|_\infty^2\left\|u^\varepsilon-u\right\|^2
\nonumber\\
		&\quad+C\left\|\dfrac1{\left(\rho^\varepsilon\right)^2}-1\right\|^2\left\||\Delta\phi\right\|_\infty^2
+\left\|\dfrac1{\rho^\varepsilon}-1\right\|^2\left\|\left(\phi^\varepsilon\right)^3-\phi^\varepsilon\right\|_\infty^2
\nonumber\\
		&\quad+C\left\|\phi^\varepsilon-\phi\right\|^2\left\|\left(\phi^\varepsilon\right)^2+\phi^\varepsilon\phi
+\phi^2\right\|_\infty^2+C\left\|\phi^\varepsilon-\phi\right\|^2\leq C\varepsilon,
		\end{align}
	for $t\in[0,T^0]$.
 Finally, applying $\nabla$ to (\ref{qqq.5.24}), we have
	\begin{align*}
&\nabla\left(\phi^\varepsilon-\phi\right)_t-\nabla\Delta\left(\phi^\varepsilon-\phi\right)  \nonumber
\\
&=-\nabla\left(u^\varepsilon\cdot\nabla\right)\left(\phi^\varepsilon-\phi\right)
-\left(u^\varepsilon\cdot\nabla\right)\nabla\left(\phi^\varepsilon-\phi\right)
-\nabla\left(u^\varepsilon-u\right)\cdot\nabla\phi
-\left(u^\varepsilon-u\right)\cdot\nabla^2\phi    \nonumber
\\
&\quad +\nabla\left(\dfrac1{\left(\rho^\varepsilon\right)^2}-1\right)\Delta\phi^\varepsilon
+\left(\dfrac1{\left(\rho^\varepsilon\right)^2}-1\right)\nabla\Delta\phi^\varepsilon
-\nabla\left(\dfrac1{\rho^\varepsilon}-1\right)\left(\left(\phi^\varepsilon\right)^3-\phi^\varepsilon\right) \nonumber
\\
&\quad
-\left(\dfrac1{\rho^\varepsilon}-1\right)\nabla\left(\left(\phi^\varepsilon\right)^3-\phi^\varepsilon\right)\nonumber
-\nabla\left(\left(\phi^\varepsilon\right)^3-\phi^3\right)+\nabla\left(\phi^\varepsilon-\phi\right).
\end{align*}
	And then (\ref{qqq.5.22}), together with (\ref{qqq.5.25}), implies that
\begin{align*}
		\int_{0}^{t}\left\|\nabla\left(\phi^\varepsilon-\phi\right)\right\|_2^2\leq C\varepsilon.
	\end{align*}
This completes the proof Theorem \ref{th3-2}.
\hfill$\Box$
\section*{Acknowledgments}
Li's work is supported by the National Natural Science Foundation of China (No.12371205).





\begin{thebibliography}{12}
\bibitem{A-F}
H. Abels, E. Feireisl,
On a diffuse interface model for a two-phase flow of compressible viscous fluids,
{\it Indiana Univ. Math. J.}, {\bf57} (2008), no. 2, 659–698.

\bibitem{A-09}
H. Abels,
On a diffuse interface model for two-phase flows of viscous incompressible fluids with matched densities,
{\it Arch. Ration. Mech. Anal.}, {\bf 194} (2009), no.~2, 463--506.

\bibitem{A-L-N}
H. Abels, Y. Liu, $\rm\breve{S}$ Ne$\rm\breve{c}$asov$\rm\acute{a}$,
Low Mach number limit of a diffuse interface model for two‐phase flows of compressible viscous fluids,
{\it GAMM-Mitteilungen}, (2024), DOI: 10.1002/gamm.202470008

\bibitem{B}
T. Blesgen,
A generalizaion of the Navier-Stokes equations to two-phase flow,
{\it J. Phys. D Appl. Phys.}, {\bf32} (1999), 1119--1123.

\bibitem{B-99}
F. Boyer,
Mathematical study of multi-phase flow under shear through order parameter formulation,
{\it Asymptot. Anal.}, {\bf 20} (1999), no.~2, 175--212.



\bibitem{C-G}
M. Chen, X. Guo,
Global large solutions for a coupled compressible Navier-Stokes/Allen-Cahn system with initial vacuum,
{\it Nonlinear Anal. Real World Appl.}, {\bf37} (2017), 350–373.



\bibitem{C-H-M-S}
Y. Chen, Q. He, M. Mei, X. Shi,
Asymptotic stability of solutions for 1-D compressible Navier-Stokes-Cahn-Hilliard system,
{\it J. Math. Anal. Appl.}, {\bf467} (2018), no. 1, 185–206.

\bibitem{C-H-S}
Y. Chen, H. Hong, X. Shi,
The asymptotic stability of phase separation states for compressible immiscible two-phase flow in 3D,
{\it Acta Math. Sci. Ser. B (Engl. Ed.)}, {\bf43} (2023), no. 5, 2133–2158.



\bibitem{C-Z}
S. Chen, C. Zhu,
Blow-up criterion and the global existence of strong/classical solutions to Navier-Stokes/Allen-Cahn system,
{\it Z. Angew. Math. Phys.}, {\bf72} (2021), no. 1, Paper No. 14, 24 pp.

\bibitem{D-H-W-Z}
S. Ding, J. Huang, H. Wen, R. Zi,
Incompressible limit of the compressible nematic liquid crystal flow,
{\it J. Funct. Anal.}, {\bf264} (2013), no. 7, 1711–1756.

\bibitem{D-L}
S. Ding, Y. Li,
Global solutions of a diffuse interface model for the two-phase flow of compressible viscous fluids in 1D,
{\it Commun. Math. Sci.}, {\bf18} (2020), no. 4, 1055–1086.

\bibitem{D-L-L}
S. Ding, Y. Li, W. Luo,
Global solutions for a coupled compressible Navier-Stokes/Allen-Cahn system in 1D,
{\it J. Math. Fluid Mech.}, {\bf15} (2013), no. 2, 335–360.


\bibitem{Feireisl}
E. Feireisl,
Incompressible limits and propagation of acoustic waves in large domains with boundaries,
{\it Comm. Math. Phys.}, {\bf 294} (2010), no. 1, 73--95.

\bibitem{F-P-R-S}
E. Feireisl, H. Petzeltov$\rm\acute{a}$, E.~Rocca and G.~Schimperna,
Analysis of a phase-field model for two-phase compressible fluids,
{\it Math. Models Methods Appl. Sci.}, {\bf 20} (2010), no.~7, 1129--1160.


\bibitem{F-N}
E. Feireisl, A. Novotný,
The low Mach number limit for the full Navier-Stokes-Fourier system,
{\it Arch. Ration. Mech. Anal.}, {\bf 186} (2007), no. 1, 77--107.

\bibitem{F-N-P}
E. Feireisl, A. Novotný, H. Petzeltová,
On the incompressible limit for the Navier-Stokes-Fourier system in domains with wavy bottoms,
{\it Math. Models Methods Appl. Sci.}, {\bf18} (2008), no. 2, 291--324.

\bibitem{F-J-A}
E. Feireisl, J. Málek, A. Novotný,
Navier's slip and incompressible limits in domains with variable bottoms,
{\it Discrete Contin. Dyn. Syst. Ser. S}, {\bf1} (2008), no. 3, 427--460.


\bibitem{F-P-P}
E. Feireisl, M. Petcu and D. Pra\v{z}\'{a}k,
Relative energy approach to a diffuse interface model of a compressible two-phase flow,
Math. Methods Appl. Sci. {\bf 42} (2019), no.~5, 1465--1479.


\bibitem{F-P-R-S}
E. Feireisl, H. Petzeltová, E. Rocca, G. Schimperna,
Analysis of a phase-field model for two-phase compressible fluids,
{\it Math. Models Methods Appl. Sci.}, {\bf20} (2010), no. 7, 1129–1160.

\bibitem{G-G-10}
C.~G. Gal and M. Grasselli,
Asymptotic behavior of a Cahn-Hilliard-Navier-Stokes system in 2D,
{\it Ann. Inst. H. Poincar\'e C Anal. Non Lin\'eaire}, {\bf 27} (2010), no.~1, 401--436.

\bibitem{G-G-11}
C.~G. Gal and M. Grasselli,
Instability of two-phase flows: a lower bound on the dimension of the global attractor of the Cahn-Hilliard--Navier-Stokes system,
{\it Phys. D}, {\bf 240} (2011), no.~7, 629--635.

\bibitem{G-M-T}
A. Giorgini, A. Miranville and R.~M. Temam,
Uniqueness and regularity for the Navier-Stokes-Cahn-Hilliard system,
{\it SIAM J. Math. Anal.}, {\bf 51} (2019), no.~3, 2535--2574.

\bibitem{H-L-95}
T.~M. Hagstrom and J. Lorenz,
All-time existence of smooth solutions to PDEs of mixed type and the invariant subspace of uniform states,
{\it Adv. in Appl. Math.}, {\bf 16} (1995), no.~2, 219--257.

\bibitem{H-H-77}
P. C. Hohenberg, B. I. Halperin, Theory of dynamic critical phenomena,
{\it Rev. Mod. Phys.}, {\bf49}(1977), 435--479.

\bibitem{H-H-W}
B. Huang, J. Huang, H. Wen,
Low Mach number limit of the compressible Navier-Stokes-Smoluchowski equations in multi-dimensions,
{\it J. Math. Phys.}, {\bf60} (2019), no. 6, 061501, 20 pp.

\bibitem{H-L}
Y. Hao, X. Liu,
Incompressible limit of a compressible liquid crystals system,
{\it Acta Math. Sci. Ser. B (Engl. Ed.)}, {\bf33} (2013), no. 3, 781–796.


\bibitem{H-W}
X. Hu, D. Wang,
Low Mach number limit of viscous compressible magnetohydrodynamic flows,
{\it SIAM J. Math. Anal.}, {\bf 41} (2009), no.~3, 1272--1294.



\bibitem{L-T}
J. Lowengrub, L. Truskinovsky,
Quasi-incompressible Cahn-Hilliard fluids and topological transitions,
{\it R. Soc. Lond. Proc. Ser. A Math. Phys. Eng. Sci.}, {\bf 454} (1998), no.~1978, 2617--2654.

\bibitem{J-J-L}
S. Jiang, Q. Ju, F. Li,
Incompressible limit of the compressible magnetohydrodynamic equations with periodic boundary conditions,
{\it Comm. Math. Phys.}, {\bf 297} (2010), no. 2, 371–400.

\bibitem{J-J-L2}
S. Jiang, Q. Ju, F. Li,
Incompressible limit of the compressible magnetohydrodynamic equations with vanishing viscosity coefficients,
{\it SIAM J. Math. Anal.}, {\bf42} (2010), no. 6, 2539–2553.

\bibitem{J-J-L-X}
S. Jiang, Q. Ju, F. Li, Z. Xin,
Low Mach number limit for the full compressible magnetohydrodynamic equations with general initial data,
{\it Adv. Math.}, {\bf259} (2014), 384–420.

\bibitem{Kotschote1}
M. Kotschote,
Strong solutions of the Navier-Stokes equations for a compressible fluid of Allen-Cahn type,
{\it Arch. Ration. Mech. Anal.}, {\bf206} (2012), no. 2, 489–514.


\bibitem{Kotschote}
M. Kotschote,
Mixing rules and the Navier-Stokes-Cahn-Hilliard equations for compressible heat-conductive fluids,
{\it Bull. Braz. Math. Soc. (N.S.)}, {\bf47} (2016), no. 2, 457–471.

\bibitem{K-M}
S. Klainerman, A. Majda,
Singular limits of quasilinear hyperbolic systems with large parameters and the incompressible limit of compressible fluids,
{\it Comm. Pure Appl. Math.}, {\bf34} (1981), no. 4, 481–524.

\bibitem{K-M-2}
S. Klainerman, A. Majda,
Compressible and incompressible fluids,
{\it Comm. Pure Appl. Math.}, {\bf 35} (1982), no. 5, 629–651.


\bibitem{Lei}
Z. Lei,
Global existence of classical solutions for some Oldroyd-B model via the incompressible limit,
{\it Chinese Ann. Math. Ser. B}, {\bf27} (2006), no. 5, 565–580.


\bibitem{L-Z}
Z. Lei, Y. Zhou,
Global existence of classical solutions for the two-dimensional Oldroyd model via the incompressible limit,
{\it SIAM J. Math. Anal.}, {\bf37} (2005), no. 3, 797–814.

\bibitem{L-L-95}
F. Lin and C. Liu,
Nonparabolic dissipative systems modeling the flow of liquid crystals,
{\it Comm. Pure Appl. Math.}, {\bf 48} (1995), no.~5, 501--537.

\bibitem{L-L-96}
F. Lin and C. Liu,
Partial regularity of the dynamic system modeling the flow of liquid crystals,
{\it Discrete Contin. Dynam. Systems}, {\bf 2} (1996), no.~1, 1--22.

\bibitem{SVN}
V.~N. Starovo\u itov, On the motion of a two-component fluid in the presence of capillary forces,
{\it Math. Notes}, {\bf 62} (1997), no.~1-2, 244--254; translated from  {\it Mat. Zametki}, {\bf 62} (1997), no.~2, 293--305.


\bibitem{Y-Z-Z}
L. Yao, C. Zhu and R. Zi,
Incompressible limit of viscous liquid-gas two-phase flow model,
{\it SIAM J. Math. Anal.}, {\bf 44} (2012), no.~5, 3324--3345.





\end{thebibliography}
\end{document}